\documentclass[12pt]{amsart}
\usepackage[margin=1.2in]{geometry}
\usepackage{parskip}
\usepackage{amssymb,mathtools}
\usepackage{mathrsfs,stmaryrd}
\usepackage{mathbbol}
\usepackage{graphicx}

\PassOptionsToPackage{no-math}{fontspec}
\usepackage[lining]{ebgaramond}

\ifluatex
\makeatletter
\let\mathbfcal\@undefined
\makeatother
\usepackage[cal=cm]{mathalpha}
\fi
\usepackage[final]{microtype}

\DeclareSymbolFontAlphabet{\mathbb}{AMSb}
\DeclareSymbolFontAlphabet{\mathbbl}{bbold}
\usepackage{tikz-cd}
\usepackage{xcolor}

\theoremstyle{definition}
\newtheorem{definition}{Definition}[section]
\newtheorem{example}[definition]{Example}
\newtheorem{assumption}[definition]{Assumption}

\theoremstyle{plain}
\newtheorem{theorem}[definition]{Theorem}
\newtheorem{corollary}[definition]{Corollary}
\newtheorem{proposition}[definition]{Proposition}
\newtheorem{lemma}[definition]{Lemma}

\newtheorem{innercustomthm}{Theorem}
\newenvironment{customthm}[1]
  {\renewcommand\theinnercustomthm{#1}\innercustomthm}
  {\endinnercustomthm}

\theoremstyle{remark}
\newtheorem{remark}[definition]{Remark}

\usepackage{etoolbox}

\makeatletter
\appto\th@remark{
  \addtolength{\@totalleftmargin}{0.5em}
  \addtolength{\linewidth}{-1em}
  \parshape 1 \@totalleftmargin \linewidth
}
\makeatother
\definecolor{bodytext}{gray}{0}
\AtBeginDocument{\color{bodytext}}
\newcommand{\keepstatementblack}[1]{\AtBeginEnvironment{#1}{\color{black}}}
\forcsvlist{\keepstatementblack}{definition,example,assumption,theorem,corollary,proposition,lemma,innercustomthm,customthm,remark}

\newcommand{\define}[4]{\expandafter#1\csname#3#4\endcsname{#2{#4}}}

\forcsvlist{\define{\DeclareMathOperator}{}{}}{coker,AnSpec,AnSpf,Spec,Spf,colim,Spa,an,Fil,cons,Fitt}

\forcsvlist{\define{\newcommand}{\mathsf}{}}
  {Mod,Vect,Fun,Rep,Coh,QCoh,Perf,Loc,Ind}

\makeatletter
\let\D\@undefined
\makeatother
\forcsvlist{\define{\newcommand}{\mathrm}{}}
{R,D,M,I,E,B,X,Y,V,W,U}
\forcsvlist{\define{\newcommand}{\mathrm}{}}{Hom,id,GL,End,Map,op,T,cris,perf,gr,RHom,Tot,fib,cofib,fgu}
\newcommand{\RcalHom}{\underline{\mathrm{RHom}}}
\forcsvlist{\define{\newcommand}{\mathbf}{}}{N,Z,Q,F,A,C}

\newcommand{\rA}{\mathrm{A}}

\newcommand{\rC}{\mathrm{C}}
\newcommand{\rF}{\mathrm{F}}
\newcommand{\rK}{\mathrm{K}}

\newcommand{\rN}{\mathrm{N}}
\newcommand{\rP}{\mathrm{P}}
\newcommand{\rQ}{\mathrm{Q}}
\newcommand{\rS}{\mathrm{S}}
\newcommand{\rZ}{\mathrm{Z}}
\newcommand{\rG}{\mathrm{G}}
\newcommand{\rJ}{\mathrm{J}}

\newcommand{\m}{\mathfrak{m}}
\newcommand{\cyc}{\mathrm{cyc}}
\newcommand{\et}{\ensuremath{\mathrm{\acute{e}t}}}

\newcommand{\proet}{\ensuremath{\mathrm{pro\acute{e}t}}}
\newcommand{\proket}{\ensuremath{\mathrm{prok\acute{e}t}}}

\newcommand{\Prism}{{\mathbbl{\Delta}}}
\newcommand{\Ainf}{\mathbb{A}_{\inf}}
\newcommand{\RGamma}{\mathrm{R}\Gamma}
\newcommand{\Log}{\mathrm{Log}}
\newcommand{\Laur}{\O_\Prism[1/\I_\Prism]_p^{\wedge}}
\newcommand{\heart}{\ensuremath\heartsuit}
\newcommand{\sq}{\mathord{\scaledmathsymbol{\square}}}
\newcommand{\laur}[1]{\mathopen{(\mkern-3mu(}#1\mathclose{)\mkern-3mu)}}
\newcommand{\dquot}{\mathord{/\mkern-6mu/}}

\newenvironment{display}{\begin{center}\begin{tikzcd}}
  {\end{tikzcd}\end{center}}

\usepackage[
  hypertexnames=false,
  colorlinks=true,
  linkcolor={red!50!black},
  citecolor={blue!50!black},
  urlcolor={blue!80!black}
]{hyperref}
\makeatletter
\renewcommand{\l@subsection}{\@tocline{2}{0pt}{3pc}{}{}}
\makeatother

\renewcommand{\O}{\ensuremath{\mathcal{O}}}
\renewcommand{\H}{\mathrm{H}}
\renewcommand{\L}{\mathrm{L}}
\renewcommand{\P}{\mathbf{P}}

\newcommand{\scaledmathsymbol}[1]{\mathchoice{\scalebox{0.75}{$\displaystyle#1$}}{\scalebox{0.75}{$\textstyle#1$}}{\scalebox{0.75}{$\scriptstyle#1$}}{\scalebox{0.75}{$\scriptscriptstyle#1$}}}
\renewcommand{\Box}{\mathord{\scaledmathsymbol{\blacksquare}}}

\title[]{A Prismatic Riemann-Hilbert correspondence}

\author{Dylan Pentland}
\email{dpentland@g.harvard.edu}

\begin{document}

\begingroup
\setlength{\parskip}{2pt}
\begin{abstract}
We study the category of constructible $p$-Kummer-\'{e}tale sheaves on the log adic generic fiber $(\mathrm{X},\mathcal{M}_{\mathrm{D}})_\eta$ for the stratification induced by a horizontal normal crossings divisor $\mathrm{D}$, showing that it is antiequivalent to a category of finitely generated unit Laurent crystals on the log prismatic site. We give a local description of finitely generated unit crystals on a framed log affine in terms of $(\varphi,\Gamma)$-modules, and show pushforward compatibility for strict log smooth proper morphisms.
\end{abstract}

\maketitle

\vspace{-1em}
\tableofcontents
\endgroup

\newpage

\section{Introduction}

Let $k$ be a perfect field of characteristic $p$. For smooth $k$-schemes, Emerton and Kisin gave a contravariant Riemann-Hilbert correspondence describing constructible $\F_p$ sheaves using finitely generated unit Frobenius modules in \cite{EmertonKisin2}. Their result is a contravariant generalization of Katz's correspondence in \cite{KatzPAdicProperties}, which describes the full subcategory of \'{e}tale local systems.

In mixed characteristic, Katz's theorem plays an important role in describing local systems on the adic generic fiber of a $p$-adic formal scheme using $(\varphi,\Gamma)$-modules. This idea is originally due to Fontaine in \cite{Fontaine1990}, and was later generalized by Andreatta-Iovita in \cite{AndreattaIovita2008}. Using modern prismatic methods, globalizations were given in \cite{MinWangRelativePhiGamma} and \cite{BhattScholzePrismaticFCrystals} using Laurent $\rF$-crystals.

In this article we give a Riemann-Hilbert correspondence in mixed characteristic extending these results to describe constructible sheaves on a fixed normal crossings stratification. At a high level, our construction replaces the role played by Katz's correspondence in the theory of $(\varphi,\Gamma)$-modules by the Emerton-Kisin correspondence. 

In our setup, $p$ is odd and $\X/\Spf \W(k)$ is a smooth quasicompact $p$-adic formal scheme, with log structure $\mathcal{M} = \mathcal{M}_\D$ induced by a horizontal normal crossings divisor $\D$. We write $\D^{(b)}_{\cons,\D}((\X,\mathcal{M}_\D)_{\eta},\Z_p)$ for the bounded derived category of derived $p$-complete Kummer-\'{e}tale $\Z_p$-sheaves constructible along the stratification induced by $\D$. As in the correspondences for Kummer-\'{e}tale local systems of \cite{KoshikawaYaoCorrigendum} and \cite[Theorem 7.5.4]{KedlayaLiu2016}, objects $\mathscr{F}$ are $p$-Kummer when the prime-to-$p$ log inertia $\I_{\overline{x}}^{\log,(p)'}$ acts trivially on $\mathscr{F}_{\overline{x}}$ at every log geometric point $\overline{x}$ (Definition \ref{def:p-Kum}).

We let $(\X,\mathcal{M})_\Prism$ denote the absolute log prismatic site of $(\X,\mathcal{M})$, and let $\Laur$ be the sheaf assigning to a log prism $(\rA,\I,\M_{\Spf \rA})^a$ the ring $\rA[1/\I]_p^{\wedge}$. In Definition \ref{def:fgu} we introduce the category $\D_{\fgu}((\X,\mathcal{M})_\Prism, \Laur)$ of finitely generated unit Laurent crystals, which are certain crystals $\mathcal{E}$ over $\Laur$ equipped with an isomorphism $\varphi^*\mathcal{E}\simeq \mathcal{E}$ and a finiteness condition similar to \cite{EmertonKisin2}. The main result of this article is then the following.

\begin{customthm}{A}[Theorem \ref{thm:RH_main}, Corollary \ref{cor:other_bases}]
Consider $(\X,\mathcal{M}_\D)$ with $\X$ smooth and quasicompact over $\Spf \W(k)$ as above. There is a symmetric monoidal equivalence
\[\T_\et: \D_{\fgu}((\X,\mathcal{M}_\D)_\Prism, \Laur)^{\op} \longrightarrow \D^{(b),\mathrm{pKum}}_{\cons,\D}((\X,\mathcal{M}_\D)_{\eta},\Z_p)\]
where the target consists of the $p$-Kummer complexes in  $\D^{(b)}_{\cons,\D}((\X,\mathcal{M}_{\D})_\eta,\Z_p)$. If $\X_\eta$ is pure of dimension $d$, then $\T_\et[d]$ identifies the standard $t$-structure on the source with the $p^+$-perverse $t$-structure on the target.
\label{thm:main_introduction}
\end{customthm}

One does not need to check on all prisms to verify the finitely generated unit condition. It in fact suffices to check on a covering family of logarithmic $q$-de Rham prisms, as in Remark \ref{rem:qdR_test_sufficient}. The theorem also holds over $\W(k)[\zeta_{p^\ell}]$ for $\ell\ge 0$, as explained in Corollary \ref{cor:other_bases}.

We can also study ordinary constructible sheaves within this correspondence. Indeed there is a morphism of sites
\[\epsilon:(\X,\mathcal{M})_{\eta,\proket}\longrightarrow\X_{\eta,\proet}\]
whose pullback $\epsilon^*$ is fully faithful on constructible sheaves and lands in the $p$-Kummer sheaves in the heart of $\D^{(b)}_{\cons,\D}((\X,\mathcal{M})_\eta,\Z_p)$. Thus ordinary constructible sheaves (for the stratification induced by $\D$) form a full subcategory at the level of hearts. See Remark \ref{rem:eff} for the corresponding condition on the crystal side of the correspondence.

This observation that one may extract ordinary constructible sheaves from the correspondence also gives a potential approach to describing all Zariski constructible sheaves via de Jong alterations. If one extended Theorem \ref{thm:main_introduction} to cover the case where the log formal model is a strict semistable pair, then using a strict alteration hypercover (effectively $h$-sheafifying the finitely generated unit side of the correspondence) we could reduce a general stratification to this case, as it seems reasonable that the abelian category $\mathrm{Cons}_{\D}(\X_\eta,\Z_p)$ satisfies similar descent properties. Then passing to a colimit over stratifications, one can access all constructible sheaves in the heart. Over $\O_{\C_p}$ the extension to semistable pairs can already be done via the same arguments, as explained in Remark \ref{rem:perfectoid_base}. The remaining difficulty lies entirely in descending the correspondence to a finite extension of $\Q_p$.

Finitely generated unit Laurent crystals also admit an explicit local description in terms of $(\varphi,\Gamma)$-modules. We state the variant over $\W(k)$, although we can handle $\W(k)[\zeta_{p^\ell}]$ as well. In what follows a \emph{framed log affine} is a $\W(k)$-algebra $\R$ that is $p$-completely \'{e}tale over $\W(k)\langle \T_1,\ldots,\T_r,\T_{r+1}^{\pm 1},\ldots,\T_n^{\pm 1}\rangle$ with log structure $\mathcal{M}$ induced by $\D=\V(\T_1 \cdots \T_r)$.

\begin{customthm}{B}[Variant of Theorem \ref{thm:qdR_crys_Gamma_mod}]
\label{thm:local_introduction}
Let $(\Spf \R,\mathcal{M})$ be a framed log affine, with framing denoted by $\sq$. Evaluation on the log $q$-de Rham prism $(\R^{\sq}\llbracket q-1 \rrbracket, [p]_q, \N^r)^a$ over $\R$ induces an equivalence
\[\D_{\fgu}((\Spf \R,\mathcal{M})_\Prism, \Laur) \simeq \D_{\fgu}(\R^{\sq}\laur{q-1}^{\wedge}_p)^{h\Gamma},\]
where $(a,b_1,\ldots, b_n)\in \Gamma=\Z_p^\times \ltimes \Z_p(1)^n$ acts by $q\mapsto q^a$ and $\T_i\mapsto q^{pb_i}\T_i$, and ``fgu'' or finitely generated unit means that the complex has a unit Frobenius and is perfect over $\R^{\sq}\laur{q-1}^{\wedge}_p\langle \rF\rangle$.
\end{customthm}

We also use this to make $\T_\et$ completely explicit locally in Lemma \ref{lem:factorization} (although in principle the definition in Definition \ref{def:Tet} is already computable). In this formalism it is relatively straightforward to write down some simple examples. In the case of $\Spf \R=\Spf \W(k)\langle \T \rangle$ with $\Gamma = \Z_p^\times \ltimes \Z_p(1)$ and log structure $\mathcal{M}_{\V(\T)}$, we describe the objects corresponding to $j_! \underline{\Z_p}$ and $i_* \underline{\Z_p}$ under Theorem \ref{thm:qdR_crys_Gamma_mod}. Here $j$ is the inclusion of the punctured disk $\mathbb{D}^\times \to (\X,\mathcal{M}_{\V(\T)})_\eta$, and $i$ the complementary inclusion of the origin (with its pullback log structure).

In the notation of Theorem \ref{thm:local_introduction}, $\M=\R^{\sq}\laur{q-1}_p^{\wedge}[1/\T]_p^{\wedge}$ corresponds to $j_! \underline{\Z_p}$. We equip it with an action by $\Gamma = \Z_p^\times \ltimes \Z_p(1)$ where $q \mapsto q^x$ for $x\in \Z_p^\times$ and $\T \mapsto q^{px} \T$ for $x\in \Z_p(1)$. Moreover it has a unit Frobenius, sending $\T \mapsto \T^p$ and $q \mapsto q^p$. The $(\varphi,\Gamma)$-module $\R^{\sq}\laur{q-1}_p^{\wedge}[1/\T]_p^{\wedge}/\R^{\sq}\laur{q-1}_p^{\wedge}$ corresponds to $i_* \underline{\Z_p}[-1]$ (for the unshifted version of $\T_\et$). In Example \ref{ex:nontrivial_cons}, we give a more detailed calculation corresponding to a constructible sheaf (which has a nontrivial gluing datum between strata).

Our third main result is pushforward compatibility.

\begin{customthm}{C}[Theorem \ref{thm:fgu_pushforward_compat}]
Let $f:(\X,\mathcal{M})\to(\Y,\mathcal{N})$ be a strict log smooth proper morphism of pure relative dimension $d$ between smooth quasicompact $p$-adic formal schemes over $\W(k)[\zeta_{p^\ell}]$ for some $\ell\ge0$, and let $\D_\X$ and $\D_\Y$ be horizontal normal crossings divisors inducing $\mathcal{M}$ and $\mathcal{N}$. Suppose the adic generic fibers have pure dimensions $d_\X$ and $d_\Y$. Then there is a commutative diagram
\begin{display}
\D_{\fgu}((\X,\mathcal{M})_\Prism, \Laur)^\op \ar{r}{\T_{\et,\X}[d_\X]} \ar{d}{\R f_* \{d\} [d]} & \D^{(b)}_{\cons,\D_\X}((\X,\mathcal{M})_\eta,\Z_p) \ar{d}{\R f_{\eta,*}} \\
\D_{\fgu}((\Y,\mathcal{N})_\Prism, \Laur)^\op \ar{r}{\T_{\et,\Y}[d_\Y]} & \D^{(b)}_{\cons,\D_\Y}((\Y,\mathcal{N})_\eta,\Z_p).
\end{display}
\end{customthm}

The shift and twist arise from a logarithmic Poincar\'{e} duality result for perfect Laurent crystals proved in Corollary \ref{cor:laurent_prismatic_pd}. For the Laurent $\rF$-crystal equivalences in either \cite{BhattScholzePrismaticFCrystals} or \cite{KoshikawaYaoCorrigendum} the analogous pushforward compatibility theorem does not appear in the literature, the problem being that it is not obvious that $\R f_*$ preserves $\Perf(\X_{\Prism},\Laur)$. For this reason existing pushforward compatibility results (such as the result stated in \cite{GuoLiFrobeniusHeight}) assume that the Laurent $\rF$-crystal arises from a prismatic $\rF$-crystal, as in that situation perfectness can be deduced from the Hodge-Tate comparison. As Theorem \ref{thm:fgu_pushforward_compat} is compatible with a contravariant version of \cite{KoshikawaYaoCorrigendum}, we also obtain pushforward compatibility for this correspondence (as well as the one in \cite[Corollary 3.7]{BhattScholzePrismaticFCrystals} as a degenerate case).

There are also other compatibilities we prove through more straightforward arguments. Pullback compatibility is implied by Corollary \ref{cor:pullback}, and Theorem \ref{thm:open_stratum} gives compatibility with $j_!$ where $j$ is the inclusion of $\X_\eta \setminus \D_\eta$. Locally, through Kashiwara's lemma (Lemma \ref{lem:kashiwara}) we also give compatibility with $i_*$ and $i^!$. Finally, due to the cohomological smoothness results in \cite{Kipp2} we also have compatibility with $f^!$ for strict log smooth proper maps of pure relative dimension $d$ by combining pushforward compatibility with $f^! \Laur \simeq \Laur \{d\}[2d]$ from Lemma \ref{lem:prismatic_stack_cohom_smooth_proper}.

\textbf{Related work}. While not directly providing Riemann-Hilbert correspondences, some of the motivation for this work came from comparison theorems in \cite{Faltings1989} and \cite{yamashita}. The article \cite{Faltings1989} in fact originally motivated Emerton and Kisin to produce \cite{EmertonKisin2} in the hope of obtaining a result of the type in Theorem \ref{thm:main_introduction}, and our correspondence shares a number of similarities with what Faltings does in this article. The comparison theorems in \cite{Faltings1989} allow various support conditions related to a normal crossings divisor $\D$, and point towards it being fruitful to work in a logarithmic setting with the stratification induced by $\D$ just as we do in this article (indeed Faltings uses a somewhat similar logarithmic setup). Various other similarities arise, for example an analogue of the $p$-Kummer condition in \cite[II(i)]{Faltings1989} and some constructions analogous to Remark \ref{rem:eff} in \cite[III(d)]{Faltings1989}. It seems plausible there could be a more direct relationship, which we speculate further about in Remark \ref{rem:comparison}.

There are also several related works which are more recent. In \cite{BhattPadicHodgeNotes}, Bhatt and Lurie construct a fully faithful covariant Riemann-Hilbert functor
\[\mathrm{RH}_\Prism:\D_{p\text{-comp}}((\X_\eta)_\et,\Z_p)^{\mathrm{oc}}\longrightarrow\D(\X^{\Prism,\mathrm{pfd}})^{a,\varphi=1}\]
under the assumption that $\X/\Spf \O_K$ is locally of finite type and $K/\Q_p$ is a nonarchimedean field of finite $p$-cohomological dimension. The source of this functor is the category of derived $p$-complete overconvergent sheaves, as indicated by the superscript `$\mathrm{oc}$'. On the target, `$a$' denotes passage to almost coefficients. Their functor is moreover compatible with their earlier correspondence in positive characteristic in \cite{BhattLurieModpRH}.

Another related Riemann-Hilbert functor was constructed by Mann in \cite{MannPadicSixFunctors}. In \cite[Theorem 3.9.23]{MannPadicSixFunctors} Mann constructs a fully faithful functor from the category of overconvergent \'{e}tale $\F_p$-sheaves on a small $v$-stack $\X$ to $\D_{\Box}(\X, \O^+/\pi)^{\varphi=1}$ with $\pi \mid p$ a pseudouniformizer. Here the target category uses solid almost modules over $\O^+/\pi$.

We point out that our correspondence has two practical advantages when compared to these very general Riemann-Hilbert functors. The first is that we can locally formulate the crystal side of the correspondence as a generalization of $(\varphi,\Gamma)$-modules, which is then amenable to computation as we work over a deperfected base (see Example \ref{ex:nontrivial_cons}). The second is that the functor is an \emph{equivalence} onto the full subcategory of $p$-Kummer constructible sheaves, so there is a clear and easily verifiable condition on the crystal side corresponding to constructibility. While we do not pursue it here, we expect our correspondence (up to a duality) to be realized as a deperfection of a logarithmic variant of the Bhatt-Lurie correspondence describing overconvergent $p$-Kummer sheaves.

The fact that the correspondence is deperfected makes it in principle possible to access $\rC_{\mathrm{st}}$-type comparison theorems for constructible sheaves, which we hope to pursue in future work. This sort of application does not seem possible with other Riemann-Hilbert correspondences due to taking perfections. More is said about this in Remark \ref{rem:comparison}.

\textbf{Overview of the proofs}. All arguments work first over $\W(k)[\zeta_p]$, and then we extend to other bases once Theorem \ref{thm:main_introduction} is proved.

In Definition \ref{def:fgu} we introduce the category of finitely generated unit Laurent crystals (which we denote by $\D_{\fgu}((\X,\mathcal{M})_\Prism,\Laur)$), and later construct a functor
\[\T_\et: \D_{\fgu}((\X,\mathcal{M})_\Prism,\Laur)^\op \longrightarrow \D((\X,\mathcal{M})_{\eta,\proket}, \Z_p).\]
Theorem \ref{thm:main_introduction} requires us to show this is fully faithful and determine the essential image as the $p$-Kummer constructible sheaves along $\D$.

We work in the logarithmic setting, so in place of \'{e}tale descent we use strict \'{e}tale descent. The main idea behind proving Theorem \ref{thm:main_introduction} is to prove strict \'{e}tale descent of both sides and then locally use the explicit $(\varphi,\Gamma)$-module description in Theorem \ref{thm:local_introduction}.

We defer explaining both of these results, and instead focus on how one proves the main theorem using the $(\varphi,\Gamma)$-module description in Theorem \ref{thm:local_introduction}. In Lemma \ref{lem:factorization}, we show that $\T_\et$ has a local description similar to how one constructs the functor for ordinary $(\varphi,\Gamma)$-modules -- from this one can see that the functor lands in the $p$-Kummer subcategory. The task is then to explain where constructibility along $\D$ comes from.

Using \cite{BlickleGamma}, we deduce that modulo $p$ any finitely generated unit $(\varphi,\Gamma)$-module has a presentation $\M \simeq \colim_n \varphi^{*n} \M_0$ where $\M_0$ is a coherent $\Gamma$-module with a $\Gamma$-equivariant injection or ``root map'' $\beta: \M_0 \to \varphi^* \M_0$. In Lemma \ref{lem:toric_stable_ideal} we classify possible Fitting ideals for $\M_0,\varphi^* \M_0$ and $\coker \beta$, allowing one to conclude that on each stratum induced by $\D$ the root map is an isomorphism. Recalling Katz's equivalence one deduces the result is lisse on every log stratum.

The argument for full faithfulness and essential surjectivity also proceeds by examining these strata, the most important being the open stratum $\U_\eta := \X_\eta \setminus \D_\eta$. This part we also do globally in Theorem \ref{thm:open_stratum}. We begin by constructing a crystal $\Laur(\ast \D)$ which is sent to $j_! \underline{\Z_p}$ by $\T_\et$. Deducing essential surjectivity uses the results of \cite{KoshikawaYaoCorrigendum}, as one can use L\"{u}tkebohmert's theorem as recast in \cite[Corollary 6.3.4]{DiaoLanLiuZhuLogAdic} to extend a local system on $\U_\eta$ to a Kummer-\'{e}tale local system. After a finite prime-to-$p$ Kummer cover $(\X_{\rN},\mathcal{M}_{\rN})$ with Galois group $\mu_{\rN}^r$ on the open stratum, this extended local system becomes $p$-Kummer and has a corresponding Laurent $\mathrm{F}$-crystal $\mathcal{E}$. We tensor with $\Laur(\ast \D)$, and keeping track of $\mu_{\rN}^r$ descent data, descend back to a Laurent crystal on $(\X,\mathcal{M})_\Prism$. Full faithfulness on the open stratum uses a decompletion statement that base change to the perfection is fully faithful. We must also understand full faithfulness for analytification on perfectoid spaces, where we use \cite{BhattPadicHodgeNotes}.

To work towards the local version of Theorem \ref{thm:main_introduction} one proves a variant of Kashiwara's lemma in Lemma \ref{lem:kashiwara}. This allows us to handle closed strata with their hollow pullback log structures, and we study a variant of the open stratum result with hollow log structures locally. Full faithfulness follows by identifying thick generators corresponding to local systems supported on locally closed strata and checking full faithfulness for maps between these: the remaining case is roughly to understand maps between local systems supported on closed and open strata. This is handled in Lemma \ref{lem:orthog_red_ket} and Lemma \ref{lem:orthog_red} on each side of the correspondence, which allow us to reduce this assertion to an $\RHom$ on a single stratum. The general outline of the argument is summarized in the toy model in the beginning of \S 5.2.

Now we return to the two ingredients needed to reduce the global theorem to a theorem about $(\varphi,\Gamma)$-modules. The first claim is Theorem \ref{thm:local_introduction} expressing $\D_{\fgu}((\X,\mathcal{M})_\Prism,\Laur)$ locally via $(\varphi,\Gamma)$-modules. This is handled in \S 4.1. Proposition \ref{prop:Gamma_mod_ff} shows that evaluation on the log $q$-de Rham prism yields a fully faithful functor into semilinear $\Gamma$-modules. For essential surjectivity in Theorem \ref{thm:qdR_crys_Gamma_mod} one uses roots to present finitely generated unit $(\varphi,\Gamma)$-modules. Modulo $p$, \cite{BlickleGamma} gives a unique minimal coherent root in the sense of Definition \ref{def:roots}. Minimality forces this to be $\Gamma$-equivariant, and we use Theorem \ref{thm:laur_root_crys} to carry a $\Gamma$-equivariant root to a coherent Laurent crystal. This shows essential surjectivity modulo $p$, which is sufficient for the argument. This idea of upgrading distinguished minimal roots is central to proving both strict \'{e}tale descent and pushforward compatibility.

We now explain strict \'{e}tale descent of $\D_{\fgu}((\X,\mathcal{M})_\Prism,\Laur)$. We work with a solid variant $\D_{\Box}((\X,\mathcal{M})_\Prism,\Laur)$ of Laurent crystals as the ambient category finitely generated unit Laurent crystals reside in, as explained in Definition \ref{def:laur_category}. Roughly using that \'{e}tale covers are descendable, it is not difficult to show $\D_{\Box}((\X,\mathcal{M})_\Prism,\Laur)$ is a stack for the strict \'{e}tale topology. One must then show that the finitely generated unit condition itself is strict \'{e}tale local. We do this by again using roots: by this point we have shown that one can promote these minimal roots \emph{as modules} from log $q$-de Rham evaluations to minimal roots in \emph{coherent Laurent crystals} on framed log affines. Once this is done, their uniqueness as subcrystals makes them glue and proves descent (since coherent Laurent crystals form a stack, using \cite{drinfeldmathew}).

Finally, for pushforward compatibility in Theorem \ref{thm:fgu_pushforward_compat}, we use results from \cite{Kipp2} to construct a stacky approach to solid Laurent crystals in the log setting and equip them with a six functor formalism. The idea of the argument is to show that the pushforward $\R f_* \{d\} [d]$ preserves the finitely generated unit condition, which formally yields the desired compatibility after using Poincar\'{e} duality. To deduce this preservation one again uses roots: it can be checked modulo $p$, and $\R f_* \mathcal{E} \simeq \R f_! (\colim_n \varphi^{*n} \mathcal{E}_0) \simeq \colim_n \R f_* \varphi^{*n} \mathcal{E}_0$. Finally we show that $\R f_* \varphi^{*n} \mathcal{E}_0 \simeq \varphi^{*n}\R f_* \mathcal{E}_0$ using a variant of Cartier descent, yielding $\R f_* \mathcal{E} \simeq \colim_n \varphi^{*n} \R f_* \mathcal{E}_0$ which shows the pushforward is finitely generated unit. The result that $\R f_* \varphi^{*n} \mathcal{E}_0 \simeq \varphi^{*n}\R f_* \mathcal{E}_0$ can also be thought of as related to the Frobenius isogeny property for prismatic cohomology, but in the setting of Laurent crystals.

\textbf{Notation}. All categories are $\infty$-categories, unless otherwise stated. We will generally assume $\X$ is a smooth quasicompact $p$-adic formal scheme over $\Spf\W(k)[\zeta_{p^\ell}]$, $\ell\ge0$, with log structure $\mathcal{M}=\mathcal{M}_\D$ induced by a horizontal normal crossings divisor $\D$. 

Our convention is that $(\X,\mathcal{M})_\Prism$ denotes the strict absolute log prismatic site as developed in \cite{KoshikawaYao2025} (there are other variants, e.g. the saturated site, but see \ref{rem:crystal_indep} for why this does not matter). We use the term `Laurentization' for the passage $(\rA,\I)\mapsto\rA[1/\I]_p^{\wedge}$ for a prism $(\rA,\I)$ or log prism. As in \cite{KoshikawaYao2025} we also use $(\R,\M)^a$ to denote the log ring associated to a prelog ring, and similarly for a prism. We write $\D_{\Box}$ and $\D_{\fgu}$ for solid Laurent crystals and finitely generated unit Laurent crystals, respectively (Definitions \ref{def:laur_category} and \ref{def:fgu}). We use $\T_\et$ for the \'{e}tale realization in Theorem \ref{thm:main_introduction} (see Definition \ref{def:Tet}).

We also use several stacks: $(\X,\mathcal{M})^\Prism$ denotes the log prismatic stack (see \ref{app:log_prismatic_stacks}), and we also use $(\X,\mathcal{M})^{\Prism,\Box}$ and $(\X,\mathcal{M})^{\mathcal{L},\Box}$ for log prismatic and Laurent analytic stacks, following the constructions in \cite{Kipp} (and analytic stacks as defined in \cite{clausenscholze}). The Laurent stack is constructed to capture $\D_{\Box}((\X,\mathcal{M})_\Prism,\Laur)$ as $\D_{\Box}((\X,\mathcal{M})^{\mathcal{L},\Box})$. For a prism $(\rA^+,\I)$ with $\rA=(\rA^+[1/\I])_p^{\wedge}$, we write 
\[\AnSpf \rA := \colim_{n\ge 1} \Spa((\rA^+\dquot p^n)[1/\I],\rA^+\dquot p^n),\]  
where $\Spa$ denotes the analytic stack associated to a Huber pair.

We will frequently reduce to the situation of a framed log affine, as in Definition \ref{def:fr_log_aff}. This has an associated log $q$-de Rham prism we will use, with associated Laurentization $\rA=\R_0^{\sq}\laur{q-1}_p^{\wedge}$. We generally denote the ring of the actual prism by $\rA^+$. We use $\sq$ to denote the framing, which determines the Frobenius and action of $\Gamma = (1+p\Z_p) \ltimes \Z_p(1)^n$ (or other variants depending on the base). We often write $\overline{\rA}=\rA/p$, and use $(-)^{h\Gamma}$ for the totalization of the continuous cobar construction.

We write the internal derived Hom as $\RcalHom$, and $\RHom$ for the usual derived Hom. We also write $\M\dquot f:=\cofib(f:\M\to\M)$, with several elements denoting iterated cofibers; single slashes denote underived quotients.

\textbf{Acknowledgements}. The author would like to thank Mark Kisin for suggesting the existence of a prismatic Riemann-Hilbert correspondence compatible with \cite{EmertonKisin2} to the author and for his guidance throughout the project. The author would also like to thank Oakley Edens, Daishi Kiyohara, Niklas Kipp, and Kush Singhal for helpful conversations.

\textbf{Statement on the use of AI}. AI tools (GPT 5.6 and GPT 6 Astra) were used to check for mathematical errors during and after writing this article. They were also used to help proofread and make minor revisions.

AI tools were used for literature search, which resulted in several useful suggestions. Many of these are simply locating references for standard statements, but some had more impact. The ones which had the most impact were locating \cite{Wigner1973}, which removed a delicate point-set topology argument, and the papers \cite{HeyerMannSixFunctors} and \cite{MikamiDuality} when asked about six functor formalism results in the literature to check properness.

Beyond reference suggestions or corrections, the following arguments in this paper were produced using AI tools: the calculation in Lemma \ref{lem:Cech_cohomology} and its helper Lemma \ref{lem:Cq_criterion}, the helper lemmas for Proposition \ref{prop:open_decompletion_hom} computing group cohomology, and simplifying an earlier example in Remark \ref{rem:perf_invar_counter}. In the helper lemmas for Proposition \ref{prop:open_decompletion_hom}, the idea of using techniques from \cite[\S 5.6]{KedlayaLiu2016} for decompletion is due to the author. AI tools worked out the explicit contracting homotopies needed to actually carry out this approach, and the author carefully checked this and simplified the argument. The role in Lemma \ref{lem:Cech_cohomology} and its helper lemma \ref{lem:Cq_criterion} is similar. These two lemmas, particularly Lemma \ref{lem:Cq_criterion}, have variants elsewhere in the article (for example in Proposition \ref{prop:laurent_cartier_descent}).

\section{The Emerton-Kisin correspondence}

In this section, we will review (and slightly generalize) some results about the Emerton-Kisin correspondence in \cite{EmertonKisin2}.

As a special case of the Emerton-Kisin correspondence for lisse sheaves, we begin with Katz's correspondence which was introduced in \cite[Proposition 4.1.1]{KatzPAdicProperties}. There the result was proved for a regular $\F_p$-algebra, but actually this was not an essential ingredient in the proof as observed in \cite[Proposition 3.4]{BhattScholzePrismaticFCrystals}. Indeed, the result holds for a general $\F_p$-algebra. 

However the version in \cite[Proposition 3.4]{BhattScholzePrismaticFCrystals} only describes $\D^{(b)}_{\mathrm{lisse}}(\Spec \R, \F_p)$, and we will need a variant for $\Z_p$ coefficients in the case $\R$ is a perfect $\F_p$-algebra. We write $\D^{(b)}_{\mathrm{lisse}}(\Spec \R, \Z_p)$ for the full subcategory of $\D(\Spec \R_{\proet},\Z_p)$ consisting of locally bounded derived $p$-complete complexes of sheaves $\mathscr{F}$ such that $\mathscr{F}\dquot p$ has all lisse cohomology sheaves.

\begin{theorem}[Katz] \label{thm:katz}
Let $\R$ be a perfect $\F_p$-algebra. Then there is an equivalence of categories
\[\Perf(\W(\R))^{\varphi=1} \simeq \D^{(b)}_{\mathrm{lisse}}(\Spec \R, \Z_p),\]
induced by the solution functor
\[\M \longmapsto \left(\M\otimes_{\W(\R)}\W(\O_{\Spec\R,\proet})\right)^{\varphi=1}.\]
\end{theorem}

\begin{proof}
We apply Lemma \ref{lem:categorical_nakayama}, which is a variant of derived Nakayama for categories that we will use frequently (and prove shortly).

The hypotheses of this lemma ask that both source and target can be written as $\N$-indexed limits of idempotent complete $\Z_p$-linear categories. For the source we use $\Perf(\W_n(\R))^{\varphi=1}$ and for the target we use $\D^{(b)}_{\mathrm{lisse}}(\Spec \R, \Z/p^n)$ (we impose a Tor boundedness condition now). It is easy to check idempotent completeness here since both sit in ambient idempotent complete categories and their defining conditions are closed under retracts. 

Finally, the reduction modulo $p$ becomes \cite[Proposition 3.4]{BhattScholzePrismaticFCrystals}. Thus by Lemma \ref{lem:categorical_nakayama} the desired claim follows.
\end{proof}

We note \cite[Proposition 3.4]{BhattScholzePrismaticFCrystals} gives this result without any hypothesis on the $\F_p$-algebra $\R$, with $\F_p$ coefficients.

\begin{remark} \label{rem:katz_contra}
The Emerton-Kisin correspondence is contravariant. Since $\Perf$ is self-dual, the previous correspondence also admits a dual or contravariant version. The solution functor is $\M \mapsto \RcalHom_{\Perf^{\varphi}(\W(\R)_\proet)}(\M_\proet, \W(\R)_\proet)$ (regarded as a pro\'{e}tale sheaf). Similarly other correspondences reducing to Katz's theorem, such as \cite[Corollary 3.7]{BhattScholzePrismaticFCrystals} or \cite{KoshikawaYaoCorrigendum} admit contravariant versions.
\end{remark}

Before we go ahead with proving Lemma \ref{lem:categorical_nakayama}, we will need some basic setup for the lemma to work. We thank Kush Singhal for suggesting trying Verdier localization, which allowed the author to simplify their first (much messier) proof.

\textbf{Setup}. Let $\mathcal{C}$ be an idempotent complete $\Perf(\Z_p)$-linear stable $\infty$-category. Write 
\[\mathcal{C}_n:=\mathcal{C}\otimes_{\Perf(\Z_p)}\Perf(\Z/p^n).\]
We say that $\mathcal{C}$ is derived $p$-complete if the natural functor $\mathcal{C}\to \lim_n\mathcal{C}_n$ is an equivalence.

We ask that functors between derived $p$-complete categories $\rF: \mathcal{C} \to \mathcal{D}$ in this setup are induced by a system of compatible functors $\rF_n: \mathcal{C}_n \to \mathcal{D}_n$. Here $\rF$ is $\Z_p$-linear and the $\rF_n$ are induced by scalar extension.

\begin{lemma} \label{lem:categorical_nakayama}
Let $\mathcal{C} = \lim \mathcal{C}_n$ and $\mathcal{D} = \lim \mathcal{D}_n$ be derived $p$-complete idempotent complete categories as in the previous setup. Consider an exact functor
\[\rF:\mathcal{C}\to\mathcal{D}\]
induced by a compatible system $\{\rF_n\}$ of functors $\mathcal{C}_n \to \mathcal{D}_n$.

If $\rF_1$ is an equivalence, then $\rF$ is an equivalence.
\end{lemma}

\begin{proof}
Full faithfulness reduces to the ordinary derived Nakayama lemma. Indeed,
\[\RHom_{\mathcal{C}_1}(\X_1,\Y_1)\simeq \RHom_{\mathcal{C}}(\X,\Y\dquot p) \simeq \RHom_{\mathcal{C}}(\X,\Y) \otimes_{\Z_p}^{\L} \F_p.\]
Note that by our setup $\RHom_{\mathcal{C}}(\X,\Y)$ is derived $p$-complete. Similarly one also gets full faithfulness of $\rF_n$.

For essential surjectivity, we observe that $\rF$ is exact and fully faithful. The essential image will be a thick subcategory of the target $\mathcal{D}$ using idempotent completeness of $\mathcal{C}$; in fact the same is true for each $\rF_n$ (this point is why we impose idempotent completeness).

Thus we can consider the Verdier quotient
\[\mathcal{Q} = \cofib(\rF)\]
where cofibers are taken in $\mathsf{Cat}^{\mathrm{perf}}_{\Z_p}$. Set also $\mathcal{Q}_n = \cofib(\rF_n)$, although a priori it need not be the case $\mathcal{Q}=\lim_n \mathcal{Q}_n$ (and we will not need this).

It will suffice to show each $\mathcal{Q}_n = \mathcal{Q}\otimes_{\Perf(\Z_p)} \Perf(\Z/p^n)$ is zero, using the fact that $\mathcal{Q}_1$ is zero. Because we showed that the essential image in $\mathcal{D}_n$ is thick, this is enough to show essential surjectivity for each $n$ and hence for $\rF$. Otherwise, we would only deduce that $\rF_n(\mathcal{C}_n)$ generates $\mathcal{D}_n$.

Considering some $\M \in \mathcal{Q}$, observe that
\[\cofib(p \cdot \id_{\M}) = \M \otimes^{\L}_{\Z_p} \Z/p\Z.\]
But there is an adjunction $\mathrm{red}_1 \leftrightarrows \mathrm{for}_1$ for the reduction functor $\mathrm{red}_1: \mathcal{Q} \to \mathcal{Q}_1$ (the right adjoint being forgetful), which has the property that $\mathrm{for}_1 \mathrm{red}_1$ recovers $- \otimes_{\Z_p} \Z/p\Z$. Since in this case $\mathrm{red}_1$ is zero we conclude that $\cofib(p \cdot \id_{\M})$ is zero. 

Thus $p$ acts invertibly on $\mathcal{Q}$ and therefore each $\mathcal{Q}_n$, which can only be the case if each $\mathcal{Q}_n$ is the zero category. The claim then follows.
\end{proof}

Katz's theorem is crucial in proving the equivalence
\[\Perf(\X_{\Prism},\Laur)^{\varphi=1} \simeq \D^{(b)}_{\mathrm{lisse}}(\X_{\eta}^\diamond,\Z_p)\]
for bounded $p$-adic formal schemes $\X$ in \cite{BhattScholzePrismaticFCrystals}. Precisely, this means that when $\X=\Spf \R$ is perfectoid, the result reduces to Katz's equivalence; their theorem is deduced by a reduction to this case. Our result for finitely generated unit Laurent crystals will be compatible in a similar way with the Emerton-Kisin correspondence for constructible sheaves.

We now extend the Emerton-Kisin correspondence of \cite{EmertonKisin2} to $\Z_p$ coefficients, deduced from the mod $p$ theorem by Lemma \ref{lem:categorical_nakayama}. If $\mathscr{X}/\Spf \W(k)$ is a smooth formal scheme, we let $\D_{\mathscr{X}}$ denote the sheaf of infinitesimal differential operators. Locally, if $\X \to \A^n_k$ is \'{e}tale with formal lift $\mathscr{X}$, this can be pulled back from the sheaf associated to $\D_{\Spf \W(k) \langle \T_1,\ldots,\T_n \rangle} = \W(k)\langle \T_1, \ldots, \T_n, \partial_1^{[i_1]}, \ldots, \partial_n^{[i_n]}\rangle$, where all $\partial_i$ are given full divided powers. We put
\[\D_{\rF,\mathscr{X}} := \widehat{\bigoplus}_{i\ge 0} \rF^{*i}\D_{\mathscr{X}}.\]
Modules over $\D_{\rF,\mathscr{X}}$ are infinitesimal crystals $\mathcal{E}$ with semilinear Frobenius $\varphi_{\mathcal{E}}$. There is also a useful Taylor formula comparing Frobenius lifts which makes $\D_{\rF,\mathscr{X}}$ canonical and global when constructed this way; see \cite[\S 6.1-6.2]{EmertonKisin1}.

There is also a perverse $t$-structure that the Emerton-Kisin correspondence is compatible with. We define this as follows, using \cite{GabberTStructures} which produces a $t$-structure $p$ on $\D^{(b)}_{\cons}(\X,\F_p)$.

\begin{definition} \label{def:EK_pplus}
We may then construct the $p^+$-perverse $t$-structure, specified by
\[{}^{p^+}\D^{\ge 0}=\{\mathscr{F} \mid \mathscr{F}\dquot p\in{}^p\D^{\ge 0}(\X_\proet,\F_p)\}.\]
We define the $\le 0$ half by orthogonality.
\end{definition}

We use the notation $p^+$ due to the similarity with \cite[Construction 4.3 and Proposition 4.6]{BhattHansenCons}.

This definition does not yet show the $t$-structure is well-defined, but this will follow once we extend the Emerton-Kisin correspondence to $\Z_p$ coefficients and check compatibility on the $\ge 0$ part, as then the $t$-structure is transported across an antiequivalence from the standard $t$-structure. We use the same convention as \cite{HemoRicharzScholbach2023} for constructible sheaves, so that in our situation for $\Z/p^n$ they become the same as bounded constructible sheaves with finite Tor amplitude over $\Z/p^n$. We will not concern ourselves with boundedness versus local boundedness in what follows, as $\X$ is assumed quasicompact.

In this setting one defines $\D^{b}_{\fgu}(\D_{\rF,\mathscr{X}})$ as the full subcategory of $\D_{\mathrm{unit}}(\D_{\rF,\mathscr{X}})$ (i.e. we insist $\rF^* \M \simeq \M$) consisting of objects which are bounded and have cohomology sheaves finitely generated over $\D_{\rF,\mathscr{X}}$.

We assume for convenience in the following theorem that $\mathscr{X}/\W(k)$ is smooth and quasicompact, although this assumption can surely be weakened at least to regularity.

\begin{theorem} \label{thm:EK_formal_lift}
There is an equivalence
\[\mathrm{RSol}: \D^{b}_{\fgu}(\D_{\rF,\mathscr{X}})^\op \simeq \D^{b}_{\cons}(\X, \Z_p),\]
where $\X=\mathscr{X}_{p=0}$ and 
\[\mathrm{RSol}_{\X}(\M) := \RcalHom_{\D_{\rF,\mathscr{X}_{\proet}}}(\M_\proet, \O_{\mathscr{X},\proet}) \in \D(\X_{\proet},\Z_p).\] 
In the case that $\X$ is pure of dimension $d$, then using the standard $t$-structure on the source we see $\mathrm{RSol}_{\X}[d]$ is $p^+$-perverse $t$-exact.
\end{theorem}

\begin{proof}
We first note that one has
\[\D^{b}_{\fgu}(\D_{\rF,\mathscr{X}}) \simeq \lim_n\D^{b}_{\fgu,\mathrm{tf}}(\D_{\rF,\mathscr{X}_{p^n=0}})\]
and similarly
\[\D^{b}_{\cons}(\X,\Z_p) \simeq \lim_n\D^{b}_{\cons}(\X_\proet,\Z/p^n) = \lim_n \D^{b}_{\mathrm{ctf}}(\X_\proet,\Z/p^n).\]
Here $\mathrm{tf}$ means Tor finite. The second equivalence follows from \cite[Theorem 1.6(2)]{HemoRicharzScholbach2023}, and uses their Remark 7.2 to identify with $\D^{b}_{\mathrm{ctf}}(\X_\proet,\Z/p^n)$.

Because each side here sits in an obvious ambient category that is idempotent complete and the imposed conditions are stable under retracts, we get the idempotent completeness requirement for Lemma \ref{lem:categorical_nakayama}. It is straightforward to show that the categories we have written here in each limit are obtained by applying $\otimes_{\Perf(\Z_p)} \Perf(\Z/p^n)$.

Moreover, \cite[Proposition 15.4.1 and Equation (16.1.3)]{EmertonKisin2} give
\[\mathrm{RSol}_{\X}(\M)\dquot p^n \simeq \mathrm{RSol}_{\X,\Z/p^n}(\M \dquot p^n).\]
Thus the solution functor is naturally obtained by a compatible system of solution functors modulo $p^n$. Upon taking mod $p$ reductions we have
\[\D^{b}_{\cons}(\X_\proet, \F_p) \simeq \D^{b}_{\cons}(\X_\et,\F_p)\]
by \cite[Theorem 1.6(1)]{HemoRicharzScholbach2023}. Thus by applying the mod $p$ Emerton-Kisin correspondence (\cite[Corollary 16.2.8]{EmertonKisin2}) we can prove the equivalence via Lemma \ref{lem:categorical_nakayama}.

The perverse $t$-exactness claim remains. In light of Definition \ref{def:EK_pplus}, our strategy will be to reduce to checking it modulo $p$ where all the work is done in \cite{EmertonKisin2}. In \cite{EmertonKisin2}, perverse $t$-exactness is actually proved for a variant with $\O_{\rF,\X}$-modules. This is not an issue, however: one notes that modulo $p$ \cite[Proposition 15.4.3]{EmertonKisin2} identifies $\D^b_{\mathrm{lfgu}}(\D_{\rF,\mathscr{X}})\simeq \D^b_{\mathrm{lfgu}}(\O_{\rF,\X})$, and \cite[16.1.5]{EmertonKisin2} shows the corresponding solution functors are compatible under this identification. In particular we may work with $\O_{\rF,\X}$-modules, even for perverse $t$-exactness, because the identification leaves the underlying $\O_{\rF,\X}$-module unchanged (that is, the identification is $t$-exact).

Assume now that $\X$ has pure dimension $d$, and consider $\mathrm{RSol}_{\X}(\M)[d]$. By derived Nakayama, $\M\in\D^{\le0}$ if and only if $\M\dquot p\in\D^{\le0}$. This compatibility with $-\dquot p$ and perverse $t$-exactness modulo $p$ from \cite[Theorem 11.5.4]{EmertonKisin2} (which applies to $\O_{\rF,\X}$-modules) identifies the condition $\M\in\D^{\le0}$ with $(\mathrm{RSol}_{\X}(\M)[d])\dquot p\in{}^p\D^{\ge0}$, which holds if and only if $\mathrm{RSol}_{\X}(\M)[d] \in {}^{p^+}\D^{\ge0}$ by Definition \ref{def:EK_pplus}. Thus $\mathrm{RSol}[d]$ identifies these halves of the two $t$-structures. Orthogonality identifies the other halves, and thus proves the claim.
\end{proof}

There is also a variant for perfect rings. In this setting we use the definition from \cite{BhattLurieModpRH} to make sense of the finitely generated unit condition: one asks for modules to be in $\mathrm{Thick}(\W(\R)\langle \rF \rangle_p)$ and to be unit.

\begin{theorem}
\label{thm:EK_perfect}
Let $\R$ be a perfect $\F_p$-algebra. Then there is an equivalence
\[\mathrm{RSol}: \D_{\fgu}(\W(\R))^\op \simeq \D^{b}_{\cons}(\Spec \R, \Z_p)\]
given by $\RcalHom_{\D^{\varphi}(\W(\R)_\proet)}(\M_\proet, \W(\R)_\proet) \in \D((\Spec \R)_{\proet},\Z_p)$. The equivalence is moreover symmetric monoidal.
\end{theorem}

\begin{proof}
This is essentially already proven in \cite{BhattLurieModpRH}, but they work largely with $\F_p$ coefficients. The extension to $\Z_p$ coefficients is similar to the previous case. Both sides are stable subcategories of their ambient derived categories closed under retracts, hence are idempotent complete and $\Perf(\Z_p)$-linear. Their reductions modulo $p^n$ (in the sense of the setting of Lemma \ref{lem:categorical_nakayama}) are the finitely generated unit categories over $\W_n(\R)$ and the constructible complexes whose stalks have finite Tor-dimension over $\Z/p^n$. Applying Lemma \ref{lem:categorical_nakayama} now reduces the claim to \cite[Theorem 11.4.4]{BhattLurieModpRH}, together with \cite[Theorems 12.5.4 and 12.6.1]{BhattLurieModpRH}.

Thus all that remains is to check the equivalence is symmetric monoidal. One may essentially copy the argument in \cite[Theorem 8.4.1]{BhattLurieModpRH}. The formula for $\mathrm{RSol}$ gives a lax symmetric monoidal structure:
\[\mathrm{RSol}(\M) \widehat{\otimes}_{\Z_p}^{\L} \mathrm{RSol}(\rN) \longrightarrow \mathrm{RSol}(\M \widehat{\otimes}_{\W(\R)}^{\L} \rN).\] 
We want to check this is in fact symmetric monoidal, which by derived Nakayama may be checked modulo $p$. By using pullback compatibility of the correspondence we may further reduce to checking symmetric monoidality on geometric stalks. Over an algebraically closed perfect field $\kappa$, we then use the dual or contravariant form of Katz's theorem described in the remark following Theorem \ref{thm:katz}. The fact that this is symmetric monoidal can be checked in the same way as \cite[Theorem 8.4.1]{BhattLurieModpRH}: the category is equivalent to $\Perf(\F_p)$, and reducing to the heart (which is fine since it thickly generates) one can decompose objects as finite direct sums of the unit where the claim is obvious. 

Finally, we must also check closure under tensor products. In the usual Emerton-Kisin correspondence with $\O_{\mathrm{F},\X}$-modules, one can check this on the heart using roots. We carry out a similar argument in the present setting, using Remark \ref{rem:fgu_mod_p} to reduce us to working with $\D_{\fgu}(\R)$. There, using \cite[Theorem 12.5.4]{BhattLurieModpRH} we can express every finitely generated unit complex as $\D(\rP)$ for a Frobenius module $\rP\in \Perf^{\varphi}(\R)$, where $\D$ denotes a \emph{weak dual}. By \cite[Construction 11.2.7 and the proof of Proposition 12.5.1]{BhattLurieModpRH}, $\D(\rP)\simeq\colim_n \varphi^{*n}(\rP^\vee)$. Therefore we deduce 
\[\D(\rP) \otimes_{\R}^{\L} \D(\rQ) \simeq \colim_{\N \times \N} \varphi^{*n} \rP^\vee \otimes_\R^{\L} \varphi^{*m} \rQ^\vee \simeq \colim_n \varphi^{*n}(\rP^\vee \otimes_\R^{\L} \rQ^{\vee}) \simeq \D(\rP\otimes_{\R}^{\L}\rQ).\]
Since $\rP \otimes_{\R}^{\L} \rQ$ is again perfect over $\R$, its weak dual is finitely generated unit by \cite[Proposition 12.3.11]{BhattLurieModpRH}.
\end{proof}

\section{Some results on Laurent crystals and \texorpdfstring{$\Gamma$}{Gamma}-modules}

Let $(\X,\mathcal{M})$ be a bounded fs log $p$-adic formal scheme. Associated to $(\X,\mathcal{M})$ is the absolute log prismatic site $(\X,\mathcal{M})_\Prism$, the opposite category of log prisms over $(\X,\mathcal{M})$ (see \cite{KoshikawaLogPrismaticI} and \cite{KoshikawaYao2025}). 

On this site, there is a sheaf of rings $\Laur: (\rA,\I,\M_{\Spf \rA})^a \mapsto \rA[1/\I]_p^{\wedge}$. We call prismatic crystals over $\Laur$ \emph{Laurent crystals}. The goal of this section is to study these Laurent crystals in the setting where $(\X,\mathcal{M})$ is a \emph{framed log affine}, which will be the standard local situation we reduce to in the article.

\begin{assumption} \label{ass:horizontal_nc}
We will frequently work under the hypothesis that $\X/\O_K$ be smooth and quasicompact, and set $\D$ to be a horizontal normal crossings divisor. By this, we mean that $\D\cup\X_s$ is a normal crossings divisor and no irreducible component of $\D$ is contained in the special fiber. Associated to $\D$, there is a log structure $\mathcal{M}=\mathcal{M}_\D$ we put on $\X$.
\end{assumption}

We fix throughout a compatible system $(\zeta_{p^m})_{m\ge 1}$ extending $\zeta_p$, and can use it to identify $\Z_p(1)$ with $\Z_p$ if we choose. In this section we assume $p$ is odd and work over $\W(k)[\zeta_p]$. At the end of this section in Corollary \ref{cor:coh_laur_other_bases}, we explain how to make modifications to obtain the result over $\W(k)$ (and also with more roots of unity). The assumption that $p$ is odd can surely be removed, but we do not address it here as the paper is long enough.

The first well-known result we show is that locally log formal schemes $(\X,\mathcal{M})$ satisfying Assumption \ref{ass:horizontal_nc} look like framed log affines as defined below.

\begin{definition} \label{def:fr_log_aff}
Let $K/\W(k)[1/p]$ be a finite totally ramified extension, where $k$ is a perfect field of characteristic $p$. A framing $\sq$ of $(\Spf \R,\mathcal{M})$ over $\O_K$ is a $p$-completely \'{e}tale map 
\[\sq: \Spf \R \longrightarrow \Spf \O_K\langle \T_1,\ldots,\T_r,\T_{r+1}^{\pm1},\ldots,\T_n^{\pm1}\rangle\] 
under which $\mathcal{M}$ is the pullback of the log structure on the target induced by the divisor $\V(\T_1\cdots \T_r)$. We write 
\[\mathbb{T}_r^n:=(\Spf \O_K\langle \T_1,\ldots,\T_r,\T_{r+1}^{\pm1},\ldots,\T_n^{\pm1}\rangle, \mathcal{M}_{\V(\T_1 \cdots \T_r)})\]
and call $(\Spf \R,\mathcal{M})$ a \emph{framed log affine}.
\end{definition}

\begin{remark}
\label{rem:fr_log_aff_W(k)_model}
For a framed log affine, we use $n$ and $r$ as the corresponding integers in Definition \ref{def:fr_log_aff}. We sometimes call the \'{e}tale lifts of $\T_1,\ldots,\T_r$ ``boundary coordinates''.

If $K/\W(k)[1/p]$ is a finite extension with residue field $k$, the framing can actually be done over $\W(k)$. By this we mean that there is a $p$-completely \'{e}tale $\W(k)\langle \T_1,\ldots,\T_r,\T_{r+1}^{\pm1},\ldots,\T_n^{\pm1}\rangle$-algebra $\R_0$ with $\R\simeq\R_0\widehat{\otimes}_{\W(k)}\O_K$, compatibly with the framing.
\end{remark}

\begin{lemma} \label{lem:strict_etale_local_framing}
Let $(\X,\mathcal{M})$ satisfy Assumption \ref{ass:horizontal_nc} over $\O_K$. Then $(\X,\mathcal{M})$ admits a strict \'{e}tale cover by finitely many framed log affines over $\O_K$.
\end{lemma}

\begin{proof}
This is surely well-known, but we include a proof for completeness. The non-logarithmic case is sketched in \cite[Lemma 4.9]{bhattframing}, and we follow a similar approach. 

By working strict \'{e}tale locally, we can separate the local branches of $\D$ and reduce to the situation where $\D$ is a simple normal crossings divisor. Once we are in this situation, let us zoom in around a closed point $x\in \X$ (which we equip with its pullback log structure) and denote the irreducible components of $\D$ passing through $x$ by $\D_1, \ldots, \D_r$. We can pick a small affine Zariski neighborhood $\Spf \R$ of $x$, so no other components besides these $\D_i$ meet this neighborhood.

In this situation we must now produce a framing for $\Spf \R$. Since $\D$ is a simple normal crossings divisor, after possibly shrinking further there are elements $t_1,\ldots,t_r,s_{r+1},\ldots,s_n\in\R$ such that $\D_i=\V(t_i)$ and
\[\mathrm{d}t_1,\ldots,\mathrm{d}t_r,\mathrm{d}s_{r+1},\ldots,\mathrm{d}s_n\]
is a basis of $\widehat{\Omega}^1_{\R/\O_K}$. We may choose constants $c_j\in\O_K$ such that $u_j=s_j+c_j$ has nonzero image in $\kappa(x)$; thus after Zariski localizing further, each $u_j$ is a unit. We can even choose $c_j$ so they are zero or one as in \cite[Lemma 4.9]{bhattframing}. Using the functions $t_i$ and $u_j$ for the images of the coordinates in $\mathbb{T}_r^n$ we obtain a morphism
\[\Spf \R \xlongrightarrow{\sq} \mathbb{T}_r^n\]
induced by the ring map sending $\T_i$ to $t_i$ for $i\le r$ and $\T_j$ to $u_j$ for $j>r$. The pullback of $\V(\T_1\cdots\T_r)$ is $\V(t_1\cdots t_r)=\D\cap\Spf\R$, hence $\sq$ is strict as a morphism of log formal schemes. Since on stalks by construction the pullback induced by $\sq$ on $\Omega^1$ is an isomorphism this shows that the underlying morphism is $p$-completely \'{e}tale.

Quasicompactness gives a finite cover.
\end{proof}

To relate Laurent crystals on a framed log affine to modules with a semilinear $\Gamma$-action, we evaluate them on the logarithmic $q$-de Rham prism attached to the framing.

Let $(\Spf \R,\mathcal{M})$ be a framed log affine with framing $\sq$ over $\W(k)[\zeta_p]$ and let $\R_0$ be the $\W(k)$-algebra explained in Remark \ref{rem:fr_log_aff_W(k)_model}. Put $\rA^+=\R_0^{\sq}\llbracket q-1\rrbracket := \R_0\llbracket q-1\rrbracket$ and $\rA=\R_0^{\sq}\laur{q-1}^{\wedge}_p:=(\rA^+[1/[p]_q])^{\wedge}_p$. Here, the decoration `$\sq$' indicates the data of the framing, which is needed to give this a log prism structure. Indeed for $\mathbb{T}_r^n$ we can define
\[\W(k)\langle \T_1,\ldots,\T_r,\T_{r+1}^{\pm1},\ldots,\T_n^{\pm1}\rangle \llbracket q-1 \rrbracket\]
with log structure given by the chart $\N^r$ with basis elements $e_i \mapsto \T_i$ and $\delta_{\log}$-structure induced by the Frobenius sending $q\mapsto q^p$ and $\T_i \mapsto \T_i^p$. Then to generalize to a framed log affine, we lift both of these along the framing $\sq$ to induce a log structure and $\delta_{\log}$ structure. Note $\rA^+/[p]_q \simeq \R$, not $\R_0$.

Put $\Gamma=(1+p\Z_p)\ltimes \Z_p(1)^n$. For $\gamma=(a,b_1,\ldots,b_n)\in\Gamma$, regard $b_i$ as an element of $\Z_p$ through the identification $\Z_p(1)\simeq \Z_p$. In these coordinates, the action of $\Gamma$ sends $q\mapsto q^a$ and $\T_i\mapsto q^{pb_i}\T_i$ for $i\le n$. Both of these actions reduce to the identity modulo $[p]_q$, therefore lifting uniquely along the framing $\sq$ to log prism automorphisms over $(\Spf\R,\mathcal{M})$. Thus the crystal property gives an action of $\Gamma$ on every evaluation in $\Coh(\rA)$ for a coherent Laurent crystal.

More generally, in Definition \ref{def:laur_category} we construct a category $\D_{\Box}((\X,\mathcal{M})_\Prism, \Laur)$ of solid Laurent crystals. We will need to work with this larger category of crystals, and our first result will be showing in Proposition \ref{prop:Gamma_mod_ff} that evaluation $\rho_{\rA}^{\Gamma,*}$ on the log $q$-de Rham prism induces a fully faithful functor
\[\rho_{\rA}^{\Gamma,*}: \D_{\Box}((\Spf\R,\mathcal{M})_\Prism,\Laur) \longrightarrow \D_{\Box}(\rA)^{h\Gamma}.\]
By the target, we will always mean the totalization $\lim_{[n]\in \Delta} \D_{\Box}(\rC(\Gamma^n, \rA))$.

We also establish a more down to earth result about roots. By a root, we mean a finite semilinear $\Gamma$-module or coherent Laurent crystal with an injective map $\mathcal{E} \hookrightarrow \varphi^* \mathcal{E}$. This will be crucial in giving a workable local description of finitely generated unit crystals, as the finitely generated unit condition will have a characterization in terms of roots. In particular we prove the following theorem:

\begin{theorem} \label{thm:laur_root_crys}
Let $(\Spf \R,\mathcal{M})$ be a framed log affine and put $\overline{\rA}=\rA/p$. Let $\M$ be a finitely generated $\overline{\rA}$-module equipped with a continuous semilinear $\Gamma$-action as above. Suppose there is a $\Gamma$-equivariant injection 
\[\M\hookrightarrow\varphi^*\M.\] 
Then there is a coherent Laurent crystal $\mathcal{E}$ killed by $p$, with $\rho_{\rA}^\ast(\mathcal{E})=\M$, and an injection $\mathcal{E}\hookrightarrow\varphi^*\mathcal{E}$ whose evaluation is the given map.
\end{theorem}

\subsection{Laurent crystals and \texorpdfstring{$\Gamma$}{Gamma}-modules}

For this subsection, we will fix a framed log affine $(\Spf \R,\mathcal{M})$ over $\Spf \W(k)[\zeta_p]$ throughout. We let $\rA^+ = \R_0^{\sq} \llbracket q-1 \rrbracket$ denote the $q$-de Rham prism attached to the framing $\sq$. The main task will be to construct a functor
\[\rho_{\rA}^{\Gamma,*}: \D_{\Box}((\Spf\R,\mathcal{M})_\Prism,\Laur) \longrightarrow \D_{\Box}(\rA)^{h\Gamma}\]
and subsequently verify that it is fully faithful. As before, $(-)^{h\Gamma}$ denotes the totalization of the continuous cobar construction.

We explain some of the content of Appendix \ref{app:solid_laur} here, as we will be working in this formalism. Throughout we work with analytic rings and stacks as defined by Clausen and Scholze in \cite{clausenscholze} -- see also \cite{RodriguezCamargoSolidGeometry}. We will frequently cite \cite{Kipp} for more background material, as largely our use of analytic stacks comes from the results in this paper. If $(\B,(d))$ is a prism, we give $\B[1/d]_p^{\wedge}$ an associated analytic stack $\AnSpf \B[1/d]_p^{\wedge}$ in Definition \ref{def:anspf}, and thus can make sense of $\D_{\Box}(\AnSpf \B[1/d]_p^{\wedge})$ which we abbreviate to $\D_{\Box}(\B[1/d]_p^{\wedge})$. For dealing with semilinear $\Gamma$-modules we will need a similar construction for $\rC(\Gamma,-)$ as well (for the mod $p^n$ reductions, as a topological ring we use the compact-open topology). The precise definition of $\D_{\Box}((\Spf\R,\mathcal{M})_\Prism,\Laur)$ is given in Definition \ref{def:laur_category} using these constructions, as a limit over all log prisms $\rA$ in $(\Spf \R,\mathcal{M})_\Prism^{\mathrm{op}}$.

We will show that actually there is a natural stacky interpretation of this functor, using the results in Appendix \ref{app:log_prismatic_stacks} and Appendix \ref{app:solid_laur}. For the reader not wishing to wade through the technical details needed to properly set up this theory, the essential point is that $(\Spf \R,\mathcal{M})^{\mathcal{L},\Box}$ is an analytic stack such that $\D_{\Box}((\Spf \R,\mathcal{M})^{\mathcal{L},\Box}) \simeq \D_{\Box}((\Spf\R,\mathcal{M})_\Prism,\Laur)$. We postpone explaining details about this Laurent stack until \S 6.1 as we will only use relatively formal properties of it here.

It is possible to carry the full faithfulness argument out entirely site-theoretically (indeed this was the first proof the author obtained) but the strategy is far easier to understand at a high level from a stacky perspective. We will realize $\rho_{\rA}^{\Gamma,*}$ as a pullback functor along
\[\rho^{\Gamma}_{\rA}: [\AnSpf \rA/\Gamma] \longrightarrow (\Spf \R,\mathcal{M})^{\mathcal{L},\Box}.\]
The stacky interpretation then gives a natural strategy to attack the problem.  In such a situation to verify full faithfulness what one needs to do is to check that the morphism of stacks is cohomologically proper, and moreover that the canonical map $\O \to \R \rho^{\Gamma}_{\rA,*} \O$ is an isomorphism. Once this is established full faithfulness is a purely formal consequence of the projection formula.

The remaining content of the argument is then a calculation of $\R \rho^{\Gamma}_{\rA,*} \O$, which we reinterpret as a cohomology calculation in Lemma \ref{lem:Cech_cohomology}. For this we will need an explicit description of the \v{C}ech nerve of $\rA^+$ and how its Laurentization sits inside the algebra of continuous functions $\rC(\Gamma^n,\rA)$ (see Remark \ref{rem:la_Cq} for a more conceptual description). The following lemma computes $\rA^{+,(1)}$ as a log prismatic envelope, which we use to obtain this description.

\begin{lemma} \label{lem:qdr_actual_cech_envelope}
Let $\rA^+=\R_0^{\sq}\llbracket q-1\rrbracket$, and let $(\rA^+)^{(1)}$ be the self-coproduct of $(\rA^+,[p]_q,\mathcal{M})^a$ in the category of log prisms $(\Spf \R,\mathcal{M})_\Prism^{\mathrm{op}}$. Put $u_0=q_1/q_0$, and let $u_i=\T_{i,1}/\T_{i,0}$ for $1\le i\le n$. Then
\[(\rA^+)^{(1)} \simeq \R_0^{\sq}\llbracket q_0-1\rrbracket[u_0^{\pm1},u_1^{\pm1},\ldots,u_n^{\pm1}]\left\{\frac{u_0-1}{[p]_{q_0}},\ldots,\frac{u_n-1}{[p]_{q_0}}\right\}_{\delta,(p,[p]_{q_0})}^{\wedge}.\]
Here $\delta(u_j)=0$ for $0\le j\le n$.
\end{lemma}

\begin{proof}
This can be checked directly by matching the universal properties of the prismatic envelope. We omit this here, since the argument is the same as Lemma \ref{lem:prism_coproduct}.
\end{proof}

Our next task will be to understand the higher self-coproducts $\rA^{(s)} := \rA^{+,(s)}[1/[p]_q]_p^{\wedge}$ explicitly. We will also enhance the evaluation functor $\rho_{\rA}^\ast$ for the log $q$-de Rham prism to keep track of the $\Gamma$ action coming from a Laurent crystal structure. That is, we will show evaluation on the log $q$-de Rham prism induces a functor
\[\rho_{\rA}^{\Gamma,*}: \D_{\Box}((\Spf\R,\mathcal{M})_\Prism,\Laur) \longrightarrow \D_{\Box}(\rA)^{h\Gamma}.\]
The crystal property does give an abstract action of $\Gamma$ since each $\gamma\in \Gamma$ gives an automorphism of $(\R_0^{\sq} \llbracket q-1 \rrbracket, [p]_q, \M)^a$ as a log prism, but does not obviously ensure continuity! 

We will therefore first need to check whether $\rho_{\rA}^{\Gamma,*}$ is even well-defined. Lemma \ref{lem:qdr_laur_cech} writes the source $\D_{\Box}((\Spf\R,\mathcal{M})_\Prism,\Laur)$ as the limit
\[\lim\left(
\begin{tikzcd}[shift left/.default=+0.8ex]
  \D_{\Box}(\rA)
  \ar[shift left]{r}
  \ar[shift right]{r}
& \D_{\Box}(\rA^{(1)})
  \ar{l}
  \ar[shift right, shift right]{r}
  \ar{r}
  \ar[shift left, shift left]{r}
& \D_{\Box}(\rA^{(2)})
  \ar[shift right]{l}
  \ar[shift left]{l}
  \ar[shift left, shift left, shift left]{r}
  \ar[shift left]{r}
  \ar[shift right]{r}
  \ar[shift right, shift right, shift right]{r}
& \cdots
  \ar[shift left, shift left]{l}
  \ar{l}
  \ar[shift right, shift right]{l}
\end{tikzcd}
\right).
\]
The target of our desired functor $\rho_{\rA}^{\Gamma,*}$ can be thought of as solid $\rA$-modules with continuous semilinear $\Gamma$-action, although Lemma \ref{lem:nuclear_continuous_action} tells us this is only reasonable under a nuclearity hypothesis (which will be satisfied in all applications). We will see the $\Gamma$-action on $\rA^+=\R_0^{\sq}\llbracket q-1 \rrbracket$, after $p$-completely inverting $[p]_q$ induces maps
\[\theta_s:\rA^{(s)} \longrightarrow \rC(\Gamma^s,\rA).\]
A priori these land only in all functions, but we will show they land in continuous functions by studying $\mathrm{im}(\theta_s)$ (and showing it is injective as well). One then obtains a functor
\[\D_{\Box}((\Spf\R,\mathcal{M})_\Prism,\Laur) \xlongrightarrow{\rho_{\rA}^{\Gamma,\ast}} \lim_{[s]\in\Delta}\D_{\Box}(\rC(\Gamma^s,\rA)) =: \D_{\Box}(\rA)^{h\Gamma}.\]

We now describe a specific subring $\rC_q(\Gamma^s,\rA) \subset \rC(\Gamma^s,\rA)$, which will later be identified with $\mathrm{im}(\theta_s)$.

\begin{definition}
We use slightly different coordinates for $\Gamma$, where the coordinate $x_0$ for $a\in 1+p\Z_p$ records $x_0=\frac{a-1}{p}$ and for $b_i \in \Z_p(1)\simeq \Z_p$ we use the coordinate $x_i$ for $1\le i\le n$. Thus we obtain a system of coordinates $x=(x_0,\ldots, x_n)$ for $\Gamma$.

Set $\eta_{0,m}(x_0)=q^{m+p\binom{m}{2}}(q-1)^m\binom{x_0}{m}_{q^p}$ and $\eta_{i,m}(x_i)=q^{p\binom{m}{2}}(q-1)^m\binom{x_i}{m}_{q^p}$ for $i>0$. Write $\eta_\alpha$ for the products of these functions in the $(n+1)$ variables and define
\[\rC_q(\Gamma,\rA) = \left(\bigoplus_{\alpha\in\N^{n+1}}\rA^+\eta_\alpha\right)^{\wedge}_{(p,q-1)}[1/[p]_q]^{\wedge}_p.\]
We define $\rC_q(\Gamma^s,\rA)$ similarly, instead using indices $\alpha\in \N^{s(n+1)}$.\footnote{There is also an alternative formula which will be of use. Considering the divided powers $[m]_{q^p}!\gamma_{m,q^p}(x)=\prod_{j=0}^{m-1}((q^p)^x-q^{pj})$ where $[m]_{q^p}!=\prod_{j=1}^m[j]_{q^p}$, we have $\eta_{0,m}(x_0)=q^m[p]_q^{-m}\gamma_{m,q^p}(x_0)$ and $\eta_{i,m}(x_i)=[p]_q^{-m}\gamma_{m,q^p}(x_i)$ for $i>0$.}
\label{def:Cq}
\end{definition}

In the previous definition, the formula directly shows $\eta_\alpha\in (q-1)^{|\alpha|}\rC(\Gamma,\rA^+)$ since $\binom{x_i}{m}_{q^p}$ is in $\rC(\Gamma,\rA^+)$. Thus any expansion in $\eta_\alpha$ whose coefficients tend to zero $(p,q-1)$-adically in $|\alpha|$ converges uniformly, and we obtain a continuous map $\rC_q(\Gamma^s,\rA) \to \rC(\Gamma^s,\rA)$ with the target having the usual compact-open topology. Here we have given $\rC_q(\Gamma^s,\rA)$ the natural topology induced by the completion in Definition \ref{def:Cq}, followed by Laurentization and $p$-adic completion. Lemma \ref{lem:Cech_function_calc} shows that this map is injective, but we warn the topology on the source is not the subspace topology.

In what follows we use $\D_{\Box}(\rA^{(s)})$ as in Definition \ref{def:laur_category} to mean 
\[\D_{\Box}(\AnSpf \rA^{(s)}):= \D_{\Box}(\colim_n \Spa((\rA^{+(s)}/p^n)[1/[p]_q], (\rA^{+(s)}/p^n)))\] 
for the log prism $\rA^{+,(s)}$ in the logarithmic \v{C}ech nerve. Here $\AnSpf$, as used below, is defined in Definition \ref{def:anspf}.

For $\rC_q(\Gamma^s,\rA)$ we make a similar definition using the completed direct sum presentation. The following lemma is more or less a rephrasing of the bases produced in \cite{GrosLeStumQuiros2023} and \cite{GrosLeStumQuirosNegative} for prismatic envelopes.

\begin{lemma} \label{lem:Cech_function_calc}
The map $\theta_s$ identifies $\rA^{(s)}$ with the subring $\rC_q(\Gamma^s,\rA)\subset\rC(\Gamma^s,\rA)$:
\[\rA^{(s)}\xlongrightarrow{ \sim }\rC_q(\Gamma^s,\rA).\]
This also induces $\AnSpf \rA^{(s)} \simeq \AnSpf \rC_q(\Gamma^s,\rA)$.
\end{lemma}

\begin{proof}
As above, we use $\rA=\R_0^{\sq}\laur{q-1}_p^{\wedge}$. We first treat $s=1$, and we will then see this case covers all higher $s$ formally. It suffices to show there is a natural isomorphism of adic rings
\[\theta_1: \rA^{+,(1)} \xlongrightarrow{ \sim } \left(\bigoplus_{\alpha\in\N^{n+1}}\rA^+\eta_\alpha\right)^{\wedge}_{(p,q-1)}.\]
After $p$-completely inverting $[p]_q$ or equivalently $q-1$, this isomorphism of ordinary adic rings implies the desired isomorphism of analytic stacks.

We do this by using the explicit log prismatic envelope presentation of $\rA^{+,(1)}$. Recall $(\rA^+)^{(1)}$ from Lemma \ref{lem:qdr_actual_cech_envelope} is the $(p,q-1)$-adic (equivalently $(p,[p]_q)$) completion of
\[\rA^+[u_0^{\pm1},\ldots,u_n^{\pm1}]\left\{\frac{u_0-1}{[p]_{q_0}},\ldots,\frac{u_n-1}{[p]_{q_0}}\right\}_{\delta}.\]
As in Lemma \ref{lem:qdr_actual_cech_envelope} we use $u_0=q_1/q_0$ and $u_i=\T_{i,1}/\T_{i,0}$, with the subscript denoting which of the maps $p_0,p_1: \rA^+ \to \rA^{+,(1)}$ we apply to $q$ or $\T_i$. The map $\theta_1$ can be made explicit integrally. For $g\in\Gamma$, using the coordinates $g=(x_0,\ldots, x_n)$ of Definition \ref{def:Cq}, the evaluation of $\theta_1(a) \in \Map(\Gamma,\rA^+)$\footnote{We haven't yet checked the output is a continuous function.} for $a\in \rA^{+,(1)}$ at $g$ is obtained by applying the $\rA^+$-algebra map determined by $u_i\mapsto(q^p)^{x_i(g)}$ for $0\le i\le n$.

Using \cite[Theorem 4.4]{GrosLeStumQuiros2023}, $\rA^+[u_0^{\pm1}]\left\{\frac{u_0-1}{[p]_{q_0}}\right\}^{\wedge}_{\delta}$ is, up to a change of variables, the $(p,q-1)$-completion of the free $\rA^+$-module on their $\Prism$-divided powers $\omega_{0,\Prism}^{\{n\}}$ (see \cite[Definition 1.7]{GrosLeStumQuiros2023}) as a $(p,q-1)$-adic ring. Indeed from the discussion preceding Theorem 4.4 in \cite{GrosLeStumQuiros2023}, they identify $\widehat{\bigoplus}_{n\ge 0} \rA^+ \omega_{0,\Prism}^{\{n\}}$ with the prismatic envelope of the $\delta$-pair $(\rA^+ \llbracket \xi \rrbracket,(\xi))$. Upon checking universal properties this is indeed our ring (one has $\xi \mapsto q_0(u_0-1)$).

Now we see where $\theta_1$ sends this basis. Using \cite[Remark 1.8(5)]{GrosLeStumQuiros2023} there is a ``blowing up'' map $\mathrm{Bl}$, yielding
\[\rA^+ \langle \xi \rangle_{q_0^p} \xlongrightarrow{\mathrm{Bl}} \widehat{\bigoplus}_{n\ge 0} \rA^+ \omega_{0,\Prism}^{\{n\}} \simeq \rA^+[u_0^{\pm1}]\left\{\frac{u_0-1}{[p]_{q_0}}\right\}^{\wedge}_{\delta} \xlongrightarrow{\theta_1} \Map(\Gamma,\rA^+)\]
such that the $q$-divided power $\xi^{[m]_{q_0^p}}$ is sent to $[p]_{q_0}^m \omega_{0,\Prism}^{\{m\}}$ under $\mathrm{Bl}$. This is helpful, because we know where $\xi^{[m]_{q_0^p}}$ is sent to under the composite map $\theta_1 \circ \mathrm{Bl}$: $\xi$ is sent by $\theta_1$ to
\[q((q^{p})^{x_0}-1) \in \Map(\Gamma,\rA^+).\]
Thus we may deduce
\[[p]_{q_0}^m \theta_1(\omega^{\{m\}}_{0,\Prism}) = q^m\gamma_{m,q^p}(x_0),\]
so in particular using the alternative formula in the footnote for Definition \ref{def:Cq} we see $\theta_1(\omega^{\{m\}}_{0,\Prism})=\eta_{0,m}$ using $[p]_{q_0}$-torsionfreeness. 

Next, we handle the remaining factors corresponding to $u_i$ for $i>0$. The argument is more or less the same. Using \cite[Definition 1.6 and Corollary 5.6]{GrosLeStumQuirosNegative}, the corresponding factor $\rA^+[u_i^{\pm1}]\left\{\frac{u_i-1}{[p]_{q_0}}\right\}^{\wedge}_{\delta}$ is the completion of the free $\rA^+$-module on divided powers $\omega^{\{m\}}_{i,\Prism}$, in their notation. The analogous map in \cite[Proposition 1.7]{GrosLeStumQuirosNegative} lets us deduce the formula
\[[p]_q^m\theta_1(\omega^{\{m\}}_{i,\Prism})=\gamma_{m,q^p}(x_i).\]
Similarly we find $\theta_1(\omega^{\{m\}}_{i,\Prism}) = \eta_{i,m}$. Thus $\theta_1$ induces a surjection onto $\left(\bigoplus_{\alpha\in\N^{n+1}}\rA^+\eta_\alpha\right)^{\wedge}_{(p,q-1)}$. One sees $\theta_1$ lands in $\rC(\Gamma,\rA^+)$ by the discussion after Definition \ref{def:Cq}. Injectivity into continuous functions $\rC(\Gamma,\rA^+)$ can be checked modulo $p$ by $p$-torsionfreeness (and $p$-adic separatedness of the source), where it follows from the change of basis in Lemma \ref{lem:Cq_criterion} using uniqueness of Mahler expansions. Hence the integral map is, as desired, an isomorphism onto $\left(\bigoplus_\alpha\rA^+\eta_\alpha\right)^{\wedge}_{(p,q-1)}$ regarded as a subring of $\rC(\Gamma,\rA^+)$; the isomorphism is of $(p,q-1)$-adic rings, so the desired claim about $\AnSpf$ follows too after Laurentizing.

For general $s$, $\rA^{+,(s)}$ is the $s$-fold iterated tensor product $(\rA^{+,(1)})^{\widehat{\otimes}^s_{\rA^+}}$. Using this one may check $\rA^{+,(s)}\simeq\left(\bigoplus_{\alpha\in\N^{s(n+1)}}\rA^+\eta_\alpha\right)^{\wedge}_{(p,q-1)}$, so by the same argument the general claim follows.
\end{proof}

Thus we can now produce a functor
\[\rho_{\rA}^{\Gamma,*}: \D_{\Box}((\Spf\R,\mathcal{M})_\Prism,\Laur) \longrightarrow \D_{\Box}(\rA)^{h\Gamma}.\]
We will be interested in proving this functor is fully faithful, so it is important to be able to understand the $\RHom$ in the target category as group cohomology. In our situation there is a straightforward analogue of the Koszul resolution for group cohomology (see \cite[Proposition B.3 and Lemma B.4]{BoscoDrinfeld}). If $\rG\simeq\Z_p^r$, let $\Z_p\llbracket\rG\rrbracket^{\Box}$ be the solid algebra associated to $\Z_p\llbracket\rG\rrbracket$. For a choice of topological generators $\gamma_1,\ldots,\gamma_r$, we use the resolution in solid $\Z_p\llbracket\rG\rrbracket^{\Box}$-modules
\begin{equation}
0 \longrightarrow (\Z_p\llbracket\rG\rrbracket^{\Box})^{\binom{r}{r}}\longrightarrow \cdots \longrightarrow (\Z_p\llbracket\rG\rrbracket^{\Box})^{\binom{r}{1}} \xlongrightarrow{(\gamma_1-1,\ldots,\gamma_r-1)}\Z_p\llbracket\rG\rrbracket^{\Box} \longrightarrow (\Z_p)_{\Box} \longrightarrow 0.
\label{eq:solid_group_resolution}
\end{equation}
In Lemma \ref{lem:nuclear_continuous_action} we will use this to give a Koszul complex description of the $\RHom$ in $\D_{\Box}(\rA)^{h\rG}$ under a hypothesis we call $p$-complete nuclearity. The term \emph{$p$-completely nuclear} just means $\M \dquot p \in \D_{\Box}(\Spa(\overline{\rA}^+[1/[p]_q], \overline{\rA}^+))$ is nuclear (where $\overline{\rA}=\rA/p$ and similarly for $\rA^+$). Recall that the nuclear objects $\M \in \D_{\Box}(\overline{\rA})^{\mathrm{nuc}}$ consist of the objects $\M$ for which, for every compact $\rP\in\D_{\Box}(\overline{\rA})^\omega$, the canonical map
\[\rP^\vee \otimes_{\overline{\rA}}^{\L,\Box}\M \longrightarrow \RcalHom_{\overline{\rA}}(\rP,\M)\]
is an equivalence. For a discrete analytic ring this condition is automatic for discrete modules; for example, $\D_{\Box}(\Z_{\mathrm{disc}})^{\mathrm{nuc}}=\D(\Z)$. We note also that the tensor products we used are implicitly $p$-completed, which is why we put $\widehat{(-)}$.

\begin{lemma} \label{lem:nuclear_continuous_action}
Let $\rG \simeq \Z_p^r$, and assume $\M$ is $p$-completely nuclear. Then
\[\RHom_{\D_{\Box}(\rA)^{h\rG}}(\rA,\M) \simeq \left[ 0 \longrightarrow (\M)^{\binom{r}{r}} \longrightarrow \cdots \longrightarrow (\M)^{\binom{r}{1}}  \xlongrightarrow{(\gamma_1-1,\ldots,\gamma_r-1)} \M \longrightarrow 0 \right](\ast).\]
\end{lemma}

\begin{proof}
On the left side, by definition of $\D_{\Box}(\rA)^{h\rG}$ we may compute the result as 
\[\RHom_{\D_{\Box}(\rA)^{h\rG}}(\rA,\M) \simeq \lim_{s\in \Delta} (\rC(\rG^s, \rA) \widehat{\otimes}_{\rA}^{\Box} \M)(\ast).\]
Note $\rC(\rG^s, \rA) = (\rA\llbracket\rG^s\rrbracket^{\Box})^{\vee}$ with $\rA\llbracket\rG^s\rrbracket^{\Box}$ compact modulo $p$. Thus by $p$-complete nuclearity of $\M$ this becomes $\left(\lim_{s \in \Delta} \RcalHom_{\rA}(\rA\llbracket\rG^s\rrbracket^{\Box},\M)\right)(\ast) = \RHom_{\rA\llbracket\rG\rrbracket^{\Box,t}}(\rA, \M)$. Here we use $\rA\llbracket\rG\rrbracket^{\Box,t}$ in the sense of \cite[Definition 1.46]{MikamiRelativelyDiscrete}. We may then use the resolution in \eqref{eq:solid_group_resolution} to finish.
\end{proof}

Iterating this lemma to apply it for $\Gamma$, we learn that if $\gamma_0\in 1+p\Z_p$ is the usual generator, then taking $\fib(\gamma_0-1)$ on the Koszul complex of Lemma \ref{lem:nuclear_continuous_action} for $\Z_p(1)^n$ computes $\RHom_{\D_{\Box}(\rA)^{h\Gamma}}(\rA,\M)$ when $\M$ is $p$-completely nuclear. We may therefore consider this as solid group cohomology, and write $\RGamma_{\Box}(\Gamma,\M) := \RHom_{\D_{\Box}(\rA)^{h\Gamma}}(\rA,\M)$ when $\M$ is $p$-completely nuclear.

To prove derived full faithfulness of $\rho_{\rA}^{\Gamma,*}$ we will reduce to computing $\RGamma_{\Box}(\Gamma,\rC_{q}(\Gamma,\rA))$. It is easy to see $\rC_{q}(\Gamma,\rA)$ is $p$-completely nuclear using the presentation in Definition \ref{def:Cq}.

One can carry out the needed cohomology calculation using the $\eta_{i,m}$ basis of Definition \ref{def:Cq}, but it is actually slightly easier to work modulo $p$ where a more natural characterization is available. Via the coordinates $x=(x_0,\ldots,x_n):\Gamma\to\Z_p^{n+1}$ of Definition \ref{def:Cq}\footnote{The map $x$ is only a homeomorphism, and does not respect the group structures.}, a continuous function $f:\Gamma \to \M$ has a unique \emph{Mahler expansion} $f(x)=\sum_\alpha a_\alpha\binom{x}{\alpha}$. Its Mahler coefficient of degree $\alpha$ is $a_\alpha=(\Delta_0^{\alpha_0}\cdots\Delta_n^{\alpha_n}f)(0)$, where $\Delta_i f(x)=f(x+e_i)-f(x)$ and $e_i$ is the $i$th standard basis vector of $\Z_p^{n+1}$. We will use this in the following to give an alternative characterization of $\rC_q$ using Mahler expansions, at least modulo $p$, that we use in Lemma \ref{lem:Cech_cohomology}.

\begin{lemma} \label{lem:Cq_criterion}
Let $\overline{\rA}=\rA/p$ (similarly for $\overline{\rA}^+$) and consider $f\in \rC(\Gamma,\overline{\rA})$. We take its Mahler expansion
\[f(x)=\sum_{\alpha\in \N^{n+1}} a_\alpha\binom{x}{\alpha},\]
where $\binom{x}{\alpha}=\prod_i\binom{x_i}{\alpha_i}$. Then $f(x)\in \rC_q(\Gamma,\overline{\rA})$ if and only if $(q-1)^{-|\alpha|} a_\alpha \longrightarrow 0$ as $|\alpha| \longrightarrow \infty$.
\end{lemma}

\begin{proof}
We will prove
\[\left(\bigoplus_{\alpha\in \N^{n+1}}\overline{\rA}^+\eta_\alpha\right)^\wedge_{(q-1)} = \left(\bigoplus_{\alpha\in \N^{n+1}}\overline{\rA}^+ (q-1)^{|\alpha|}\binom{x}{\alpha}\right)^\wedge_{(q-1)},\]
which clearly suffices.

We first compute the Mahler expansions of the $\eta_{\alpha}$ from Definition \ref{def:Cq}. As discussed above, to compute the coefficients if $\alpha=(\alpha_0,\ldots,\alpha_n)$ we must work out the coefficients
\[(\Delta_0^{\beta_0} \cdots \Delta_n^{\beta_n}\eta_{\alpha})(0)\]
of the expansion of $\eta_\alpha$ into binomial coefficients.

Recall that for $i>0$ we have $\eta_{i,m}(x_i)=(q-1)^m \cdot q^{p\binom{m}{2}}\binom{x_i}{m}_{q^p}$. Because we work in $\rA=\R_0^{\sq}\laur{q-1}_p^{\wedge}$ the leading factor can be ignored. The Mahler expansion of a $q$-binomial coefficient can be computed using the recurrence
\[\Delta_i \binom{x_i}{m}_{q^p} = \binom{x_i+1}{m}_{q^p} - \binom{x_i}{m}_{q^p} = (q^{pm}-1)\binom{x_i}{m}_{q^p} + \binom{x_i}{m-1}_{q^p}.\]
Bringing back the scaling factors $(q-1)^m \cdot q^{p\binom{m}{2}}$, this yields $\Delta_i \eta_{i,m} = (q^{pm}-1)\eta_{i,m} + (q-1)q^{p(m-1)}\eta_{i,m-1}$. For $i=0$ the result is similar, since these differ only by a factor of $q^m$ (which is irrelevant in the argument that follows).

This immediately tells us some of the Mahler coefficients. Namely, it tells us that the coefficient of $\binom{x_i}{j}$ for $\eta_{i,m}$ vanishes for $j<m$ and is $(q-1)^m \cdot \mathrm{unit}$ for $j=m$. The first claim is simply because $\eta_{i,m}(j)=0$ when $j<m$. For the second, note from the previous recurrence
\[\Delta^m_i \eta_{i,m}(0) = (q-1)q^{p(m-1)} \Delta^{m-1}_i \eta_{i,m-1}(0).\]
Thus, upon iterating this, we get $(q-1)^m \cdot \mathrm{unit}$ as the coefficient of $\binom{x_i}{m}$. For higher terms, all that matters is that a similar argument shows when $j>m$ the coefficient of $\binom{x_i}{j}$ in $\eta_{i,m}$ is in $(q-1)^{m+p(j-m)} \cdot (\R_0^{\sq}\llbracket q-1 \rrbracket/p)$. 

In summary, we have shown that in $\overline{\rA}^+ = (\R_0^{\sq}\llbracket q-1 \rrbracket/p)$ we have
\[\eta_{i,m} = u_{i,m}(q-1)^m \binom{x_i}{m} + \sum_{j>m}c_{i,m,j} (q-1)^j \binom{x_i}{j}.\]
Here $u_{i,m}\in(\overline{\rA}^+)^\times$, and $c_{i,m,j}\in(q-1)^{(p-1)(j-m)}\overline{\rA}^+$, relating the two proposed bases of $\rC_q(\Gamma,\overline{\rA})$. Since this is an upper triangular change of basis where the off diagonal part is a topologically nilpotent operator, the claim follows. More precisely this shows the $(q-1)^{|\alpha|}\binom{x}{\alpha}$ and the $\eta_\alpha$ bases are related by an operator $\D(1+\rN)$ where $\D$ is diagonal and invertible and
\[\rN^r \left(\bigoplus_{\alpha\in \N^{n+1}}\overline{\rA}^+ (q-1)^{|\alpha|}\binom{x}{\alpha}\right)^\wedge_{(q-1)} \subset (q-1)^{(p-1)r} \left(\bigoplus_{\alpha \in \N^{n+1}} \overline{\rA}^+ (q-1)^{|\alpha|}\binom{x}{\alpha}\right)^\wedge_{(q-1)},\]
hence $\rN$ is topologically nilpotent. Thus we deduce the claimed isomorphism.
\end{proof}

\begin{remark}
Interestingly, this shows that modulo $p$ we have $\rC_q(\Gamma,\overline{\rA}) = \rC_{\lambda-\mathrm{an}}(\Gamma,\overline{\rA})$ where $\lambda = - \log_p(p-1)$ as in \cite{PoratSolidLocallyAnalytic}. Strictly speaking Porat does not allow $\lambda < 0$, but as noted in \cite[Definition 4.1, footnote 1]{HauckLA} the definition still makes sense.
\label{rem:la_Cq}
\end{remark}

We now carry out the needed cohomology calculation, using Lemma \ref{lem:nuclear_continuous_action}. Right translation defines a continuous $\Gamma$-action on $\rC_q(\Gamma,\rA)$ by $(\gamma f)(g)=f(g\gamma)$, which we use below.

\begin{lemma} \label{lem:Cech_cohomology}
There is a natural equivalence
\begin{display}
\rA \ar{r}{\sim} & \RGamma_{\Box}\left(\Gamma,\rC_q(\Gamma,\rA)\right).
\end{display}
\end{lemma}

\begin{proof}
By derived Nakayama it is enough to work modulo $p$. The proof of Lemma \ref{lem:Cq_criterion} shows there is a lattice
\[\left(\bigoplus_{\alpha\in\N^{n+1}} \overline{\rA}^+(q-1)^{|\alpha|}\binom{x}{\alpha}\right)^{\wedge}_{q-1} \subset \rC_q(\Gamma,\overline{\rA})\]
Let $\Fil^\bullet$ be the $q-1$-adic filtration on this lattice; this filtration is complete and separated.

Let $\gamma_0,\ldots, \gamma_n$ be the usual generators of $\Gamma$, in the coordinates $x_i$ of $\Gamma$ chosen in Definition \ref{def:Cq}. Writing down the action $(\gamma f)(g) := f(g\gamma)$ on $\rC_q(\Gamma,\rA)$ explicitly, we see right multiplication by $\gamma_i$ acts on $g=(x_0,\ldots,x_n)$ by sending $x_i \mapsto x_i + 1 + px_0$.

Using Lemma \ref{lem:nuclear_continuous_action}, the desired claim will follow by showing that modulo $p$ the complex
\[\left[\rC_q(\Gamma,\overline{\rA}) \xlongrightarrow{\gamma_i-\id} \rC_q(\Gamma,\overline{\rA})\right]\]
is naturally quasi-isomorphic to the functions in $\rC_q(\Gamma,\overline{\rA})$ which are independent of $x_i$ (placed in degree zero). Thus if we iterate this, using \eqref{eq:solid_group_resolution} shows that we get $\overline{\rA}$ as $\RGamma_{\Box}(\Gamma,\rC_q(\Gamma,\overline{\rA}))$; one uses Hochschild-Serre to work with $\Gamma$ rather than $\Z_p^n$, by dealing with $\Z_p^n$ first and then $1+p\Z_p$.

This can be done explicitly. The idea is to argue that $\gamma_i-\id$ is a perturbation of the finite difference operator $\Delta_i$ which acts on $\rC_q(\Gamma,\rA)$ (for the coordinate $x_i$). Here we use that $\rC_q(\Gamma,\rA)$ is a ring and stable under the coordinate change $x_i \mapsto x_{i}+1$. For $\Delta_i$ the claim is clear, because in
\[\left[\rC_q(\Gamma,\overline{\rA})\xlongrightarrow{\Delta_i}\rC_q(\Gamma,\overline{\rA})\right]\]
the kernel consists of functions independent of $x_i$, and it is surjective since we can define
\[h_i \left((q-1)^{|\alpha|}\binom{x}{\alpha}\right) = (q-1)^{|\alpha|}\binom{x_i}{\alpha_i+1}\prod_{j\ne i}\binom{x_j}{\alpha_j}\]
giving a continuous endomorphism of $\rC_q(\Gamma,\overline{\rA})$ such that $\Delta_i h_i=\id$. 

Now we observe that 
\[(\gamma_i-\id) \cdot f(x_0,\ldots,x_n) = f(x_0,\ldots,x_i + 1 + px_0, \ldots, x_n) - f(x_0,\ldots,x_n).\]
In particular, $\gamma_i-\id$ is nearly the finite difference operator $\Delta_i$ except for the $px_0$ term. We have a commutative diagram
\begin{display}
\rC_q(\Gamma,\overline{\rA}) \ar{r}{\gamma_i-\id} & \rC_q(\Gamma,\overline{\rA}) \\
\rC_q(\Gamma,\overline{\rA}) \ar{u}{\id} \ar{r}{\Delta_i} & \rC_q(\Gamma,\overline{\rA}) \ar{u}{(\gamma_i-\id)h_i}.
\end{display}
We claim $(\gamma_i-\id)h_i$ differs from $\id$ by an operator that carries $\Fil^r$ into $\Fil^{r+p-2}$. Indeed note that 
\[(\gamma_i-\id)h_i - \id = ((\gamma_i-\id)-\Delta_i)h_i.\]
Following the effect of $(\gamma_i-\id)-\Delta_i$ on $\binom{x_i}{\alpha}$, we note the identity 
\[((\gamma_i-\id)-\Delta_i) \cdot \binom{x_i}{\alpha_i} = \binom{(x_i+1)+px_0}{\alpha_i}- \binom{x_i+1}{\alpha_i} = \sum_{\ell \ge 1} \binom{x_i+1}{\alpha_i-p\ell} \binom{x_0}{\ell}\]
modulo $p$. In particular, using this for $(q-1)^{|\alpha|}\binom{x}{\alpha}$ shows $(\gamma_i-\id)-\Delta_i$ sends $\Fil^r$ into $\Fil^{r+p-1}$ (for the $q-1$-adic filtration). Inspecting the formula for $h_i$, we then deduce $((\gamma_i-\id)-\Delta_i)h_i$ sends $\Fil^r$ into $\Fil^{r+p-2}$. As $p>2$, this is filtration raising and $(\gamma_i-\id)h_i$ has an inverse.

We conclude in the previous diagram that there is a quasi-isomorphism between the rows, so the desired claim follows.
\end{proof}

\begin{remark}
The previous argument is very short and direct, but there is likely a more indirect (but more insightful) argument that draws upon an analogy to Shapiro's lemma. Using Remark \ref{rem:la_Cq}, the claim in question is equivalently that $\RGamma_{\Box}\left(\Gamma,\rC^{\lambda-\mathrm{an}}(\Gamma,\overline{\rA})\right) = \overline{\rA}$.
\end{remark}

\begin{proposition} \label{prop:Gamma_mod_ff}
Evaluation on the log $q$-de Rham prism attached to $\rA^+=\R_0^{\sq}\llbracket q-1\rrbracket$ induces a fully faithful functor
\[\rho_{\rA}^{\Gamma,*}: \D_{\Box}((\Spf\R,\mathcal{M})_\Prism,\Laur) \longrightarrow \D_{\Box}(\rA)^{h\Gamma}.\]
\end{proposition}

\begin{proof}
Using Proposition \ref{prop:stack_site_compare}, we may freely move between working with $\D_{\Box}(-)$ on the Laurent stacks $(\X,\mathcal{M})^{\mathcal{L},\Box}$ from \S 6.1, or the site-theoretic version $\D_{\Box}((\Spf\R,\mathcal{M})_\Prism,\Laur)$. Working with stacks, the functor $\rho_{\rA}^{\Gamma,*}$ is induced on $\D_{\Box}(-)$ by pullback along
\[\rho^{\Gamma}_{\rA}: [\AnSpf \rA/\Gamma] \longrightarrow (\Spf \R,\mathcal{M})^{\mathcal{L},\Box}.\]
See Definition \ref{def:anspf} for the meaning of $\AnSpf$. Taking the presentation $\colim_{[n]\in \Delta^\op} \AnSpf \rA^{(n)}$ of the target given by Proposition \ref{prop:log_BK_analytic_presentation} after base change to $(\X,\mathcal{M})^{\mathcal{L},\Box}$ and the natural presentation $\colim_{[n]\in \Delta^\op} \AnSpf \rC(\Gamma^n,\rA)$ of the source one obtains the description of the functor appearing before Lemma \ref{lem:Cech_function_calc}. Here, both presentations are obtained by taking the \v{C}ech nerve of the cover by $\AnSpf \rA$. We note that one should not construct this map by claiming only that the map $\AnSpf \rA \to (\Spf \R,\mathcal{M})^{\mathcal{L},\Box}$ obtained from the fpqc epimorphism of formal stacks obtained in Lemma \ref{lem:log_prism_stack_map} is $\Gamma$-equivariant, because a priori the action might not be continuous. But once one uses presentations of both stacks and Lemma \ref{lem:Cech_function_calc} to check the action is continuous, as in the site-theoretic construction it becomes well-defined.

The assertion reduces to cohomological properness of $\rho^{\Gamma}_{\rA}$ in the sense of \cite[Definition 2.1.3.1]{Kipp}, which we prove in Lemma \ref{lem:rho_cohom_prop}. We use the same conventions as \cite{Kipp}, and work with the six functor formalism on analytic stacks. In our case the analytic stacks all lie over $\AnSpec \Z_{p,\Box}$ with its solid analytic ring structure.

The claim that $\rho_{\rA}^{\Gamma,\ast}$ is fully faithful amounts to the claim that the map
\[\RHom(\mathcal{E}',\mathcal{E}) \longrightarrow \RHom(\rho_{\rA}^{\Gamma,\ast} \mathcal{E}', \rho_{\rA}^{\Gamma,\ast} \mathcal{E}) \simeq \RHom(\mathcal{E}', \R \rho^{\Gamma}_{\rA,*} \rho_{\rA}^{\Gamma,\ast} \mathcal{E})\]
is an isomorphism. It therefore is equivalent to check $\R \rho^{\Gamma}_{\rA,*} \rho_{\rA}^{\Gamma,\ast} \mathcal{E} \simeq \mathcal{E}$. As $\rho^{\Gamma}_{\rA}$ is cohomologically proper, we have access to the projection formula and can compute
\[\R \rho^{\Gamma}_{\rA,*} \rho_{\rA}^{\Gamma,\ast} \mathcal{E} \simeq \R \rho^{\Gamma}_{\rA,!}(\O \otimes^{\L} \rho_{\rA}^{\Gamma,\ast} \mathcal{E}) \simeq \R \rho^{\Gamma}_{\rA,\ast} \O \otimes^{\L} \mathcal{E}.\]
In particular we must show the canonical map $\O_{(\Spf \R,\mathcal{M})^{\mathcal{L},\Box}} \to \R \rho^{\Gamma}_{\rA,\ast} \O_{[\AnSpf \rA/\Gamma]}$ is an isomorphism.

This task can be simplified further. Again as $\rho^{\Gamma}_{\rA}$ is cohomologically proper we also have proper base change. We will use the square
\begin{display}
{[\AnSpf \rA^{(1)}/\Gamma]} \ar{d}{f} \ar{r} & {[\AnSpf \rA/\Gamma]} \ar{d}{\rho^{\Gamma}_{\rA}} \\
\AnSpf \rA \ar{r}{\rho_{\rA}} & (\Spf \R,\mathcal{M})^{\mathcal{L},\Box}.
\end{display}
This square is Cartesian, using that $\AnSpf \rA \times_{(\Spf \R,\mathcal{M})^{\mathcal{L},\Box}} \AnSpf \rA = \AnSpf \rA^{(1)}$ from the \v{C}ech nerve presentation of $(\Spf \R,\mathcal{M})^{\mathcal{L},\Box}$ in Proposition \ref{prop:log_BK_analytic_presentation}. Note that the $\Gamma$ action is nontrivial in $[\AnSpf \rA^{(1)}/\Gamma]$ as it is acting by right translation. For the unit $\O \in \D_{\Box}([\AnSpf \rA/\Gamma])$ properness gives a base change isomorphism
\begin{display}
\rho_{\rA}^* \R\rho^{\Gamma}_{\rA,*} \O \ar{r}{\sim} & \R f_* \O.
\end{display}
Formally, showing $\O_{(\Spf \R,\mathcal{M})^{\mathcal{L},\Box}} \to \R \rho^{\Gamma}_{\rA,\ast} \O_{[\AnSpf \rA/\Gamma]}$ is an isomorphism then reduces to the cohomology computation that the canonical map $\rA \to \R\Gamma_{\Box}(\Gamma, \rC_q(\Gamma,\rA))$ is an equivalence. It is here that we use that the unit is $p$-completely nuclear, so $\R f_* \O$ can be computed as solid group cohomology. This calculation was done in Lemma \ref{lem:Cech_cohomology}.
\end{proof}

The following lemma proves the properness result we needed.

\begin{lemma} \label{lem:rho_cohom_prop}
The map 
\[\rho_{\rA}^{\Gamma}: [\AnSpf \rA/\Gamma] \longrightarrow (\Spf \R,\mathcal{M})^{\mathcal{L},\Box}\]
constructed in Proposition \ref{prop:Gamma_mod_ff} is cohomologically proper.
\end{lemma}

\begin{proof}
Proposition \ref{prop:log_BK_analytic_presentation} and base change to the Laurent stack $(\X,\mathcal{M})^{\mathcal{L},\Box}$ show that $\rho_{\rA}: \AnSpf \rA \to (\Spf \R,\mathcal{M})^{\mathcal{L},\Box}$ is a representable proper effective epimorphism. Thus we only need to deal with the quotient by $\Gamma$.

Since properness is local on the target for such covers by \cite[Remark 1.1.3.7]{Kipp}, it suffices to check the claim for
\[[\AnSpf \rA^{(1)}/\Gamma] \longrightarrow \AnSpf \rA\]
where we again use the square
\begin{display}
{[\AnSpf \rA^{(1)}/\Gamma]} \ar{d}{f} \ar{r} & {[\AnSpf \rA/\Gamma]} \ar{d}{\rho^{\Gamma}_{\rA}} \\
\AnSpf \rA \ar{r}{\rho_{\rA}} & (\Spf \R,\mathcal{M})^{\mathcal{L},\Box}.
\end{display}
This property is stable under base change, so this can be also checked using $\rA^+$ and $\rA^{+,(1)}$. The right translation action involved in $[\AnSpf \rA^{(1)}/\Gamma]$ fixes the image of $\rA^+$, giving a factorization of the map as
\[[\AnSpec \rA^{+,(1)}/\Gamma]\longrightarrow \AnSpec \rA^+\times \B\Gamma\longrightarrow \AnSpec \rA^+.\]
Note we are using $(p,[p]_q)$-adic rings now, so we also use the notation $\Spa(\rA^{+,(1)},\rA^{+,(1)})$ elsewhere for $\AnSpec \rA^{+,(1)}$ as used here. After base change, pulling the first map back along $\AnSpec \rA^+ \to \AnSpec \rA^+\times \B\Gamma$ gives $\AnSpec \rA^{+,(1)}\to\AnSpec \rA^+$, which is proper by Lemma \ref{lem:cech_descendable}. It remains to show that
\[\B \Gamma \longrightarrow \AnSpec \Z_{p,\Box}\]
is proper for the six functor formalism on analytic stacks. 

We will use an argument similar to \cite[Proposition 5.3.2]{HeyerMannSixFunctors}. The map $g:\AnSpec \Z_{p,\Box}\to\B\Gamma$ is representable and proper, while $\B \Gamma \to \AnSpec \Z_{p,\Box}$ is $1$-truncated and its composite with $g$ is the identity. Then applying \cite[Corollary 4.7.5(iii)]{HeyerMannSixFunctors}, it remains only to verify that $g_*\O=\rC(\Gamma,\Z_p)$ is descendable. Using \eqref{eq:solid_group_resolution} applied to the normal subgroup $\Z_p(1)^n$ and quotient $1+p\Z_p\simeq\Z_p$, this shows that the trivial representation $\Z_p$ is perfect over $\Z_p\llbracket\Gamma\rrbracket$. Thus upon taking derived solid $\Z_p$-duals, we get
\[\Z_p\in\mathrm{Thick}(\rC(\Gamma,\Z_p))\subset\D_{\Box}(\B\Gamma)\]
as desired since $\Z_p\llbracket\Gamma\rrbracket^\vee\simeq\rC(\Gamma,\Z_p)$.
\end{proof}

\subsection{Roots in Laurent crystals}

We now turn our attention to proving Theorem \ref{thm:laur_root_crys}. As we now impose a perfectness or coherence condition, condensed mathematics is not particularly relevant in this subsection.

In Appendix \ref{app:solid_laur} we show in Corollary \ref{cor:laur_t_struct} that perfect Laurent crystals have a canonical $t$-structure in the normal crossings setting, allowing us to define $\Coh((\Spf \R,\mathcal{M})_\Prism,\Laur)$.

For a framed log affine, put $\rA = \R_0^{\sq} \laur{q-1}_p^{\wedge}$ as before. We use $\rA$ to mean this throughout this subsection, and use $\overline{\rA} := \rA/p$. One can understand the category of perfect $\rA$-complexes with continuous $\Gamma$ action or $\Perf(\rA)^{h\Gamma}$ as the totalization
\begin{equation}
\lim \left(
\begin{tikzcd}[shift left/.default=+0.8ex]
\Perf(\rA)
  \ar[shift left]{r}{p_0}
  \ar[shift right,swap]{r}{p_1}
&
\Perf(\rC(\Gamma,\rA))
  \ar{l}
  \ar[shift right, shift right]{r}
  \ar{r}
  \ar[shift left, shift left]{r}
&
\Perf(\rC(\Gamma^2,\rA))
  \ar[shift right]{l}
  \ar[shift left]{l}
  \ar[shift left, shift left, shift left]{r}
  \ar[shift left]{r}
  \ar[shift right]{r}
  \ar[shift right, shift right, shift right]{r}
&
\cdots
  \ar[shift left, shift left]{l}
  \ar{l}
  \ar[shift right, shift right]{l}
\end{tikzcd}
\right)
\label{eq:cobar}
\end{equation}
where $\rA$ has the topology specified after Lemma \ref{lem:qdr_actual_cech_envelope}. This is compatible with the previous very general description $\lim_{[s]\in\Delta}\D_{\Box}(\rC(\Gamma^s,\rA))$ by taking dualizable objects, which will only depend on the underlying ring (i.e. evaluation on $*$ is an equivalence with ordinary perfect complexes).

In this subsection we will be interested in coherent objects killed by $p$. In this case, to write down the data we may truncate and only keep the terms up to $\rC(\Gamma^2,\overline{\rA})$ in \eqref{eq:cobar}. The descent data can naturally be encoded as a comodule over the Hopf algebroid $(\overline{\rA},\rC(\Gamma,\overline{\rA}))$. One can also equivalently encode this data in terms of a stratification and the cocycle condition; for example the coaction is recovered as $\M \to \rC(\Gamma,\rA) \otimes_{p_1,\rA} \M \simeq \rC(\Gamma,\rA) \otimes_{p_0,\rA} \M$ where the second map uses the stratification isomorphism. The identity and cocycle conditions translate to the counit and coassociativity axioms on the comodule side.

Upon restricting to $\Coh$ the description is also compatible with classical theory describing $\Gamma$-modules, where one equips finitely presented modules with their weak topology as in \cite[Lemmas 8.2 and 8.22]{Schneider} (which discusses the case of a point). The author thanks Daishi Kiyohara for bringing this point to their attention.

In particular, one may check that there is an equivalence between $\Coh(\rA)^{h\Gamma}$ (the heart of the cobar description described above) and the category of finitely generated $\rA$-modules with a continuous $\Gamma$-action, with the $\rA$-module given its weak topology. Over a topological ring $\rA$ the weak topology on a finitely presented module $\M$ is the quotient topology induced by
\[\rA^n \twoheadrightarrow \M\]
which is independent of the choice of generators. The following lemma clarifies the situation -- we have tried to work this out in maximum generality so that the general principle can be applied to other types of base rings for $(\varphi,\Gamma)$-modules not used in the article.

\begin{lemma} \label{lem:weak_top}
Let $\rA$ be a topological ring, and let $\Gamma$ be a profinite group acting continuously on $\rA$. We make two topological assumptions on $\rA$:
\begin{itemize}
  \item The topology on $\rA$ is complete Hausdorff metrizable.
  \item For a finitely generated module $\M \subseteq \rA^n$ the weak and subspace topologies agree.
\end{itemize} 
Let $\Mod^{\mathrm{fp}}(\rA)$ be the category of finitely presented $\rA$-modules. Then there is an equivalence
\[\mathcal{F}: \Mod^{\mathrm{fp}}(\rA)^{h\Gamma}\simeq \left\{\begin{gathered}\M\in \Mod^{\mathrm{fp}}_{\rA} \text{ with a continuous semilinear } \\ \Gamma\text{-action for the weak topology on } \M\end{gathered} \right\}.\]
\end{lemma}

\begin{proof}
The main content of the argument is actually just checking $\mathcal{F}$ is well-defined: through the cobar description we immediately get some $\Gamma$-module $\M$, but there is no immediate reason for the action to be continuous for the weak topology.

What we must check to establish the lemma is that the coaction map
\[\M \longrightarrow \rC(\Gamma,\rA) \otimes_{p_0, \rA} \M \simeq \rC(\Gamma,\M)\]
is continuous for the weak topology on $\M$. Then the inverse to $\mathscr{F}$ is just forgetting the topological data. Since the first map is evidently continuous, the claim amounts to verifying that second isomorphism is an isomorphism of topological modules where $\rC(\Gamma,\rA)$ has the compact-open topology and $\M$ has the weak topology.

We will use the following general property of $\rC(\Gamma,-)$ with respect to exact sequences.

\begin{lemma} \label{lem:strict_exact}
Let $\rA$ be complete Hausdorff metrizable, and let
\[0 \longrightarrow \rK \longrightarrow \rA^n \longrightarrow \M \longrightarrow 0\]
be an exact sequence of left modules where $\rA^n$ has its natural product topology, $\rK$ has the subspace topology and is closed, and $\M$ the quotient topology. Then
\[0 \longrightarrow \rC(\Gamma,\rK) \longrightarrow \rC(\Gamma,\rA^n) \longrightarrow \rC(\Gamma,\M) \longrightarrow 0\]
is a strictly exact sequence of topological modules over $\rC(\Gamma,\rA)$ where each is equipped with the compact-open topology.
\end{lemma}

Here the term ``strict exact'' means that the leftmost term has the subspace topology and the rightmost term has the quotient topology.

\begin{proof}
The strict exactness claims are more or less clear using the subbasis of opens $\V_{\Gamma',\U} = \{f \in \rC(\Gamma,\rK): f(\Gamma')\subset \U\}$ where $\Gamma'\subseteq \Gamma$ is compact and $\U$ open for the compact-open topology, except for right strict exactness $\rC(\Gamma,\rA^n) \longrightarrow \rC(\Gamma,\M) \longrightarrow 0$. We must show both surjectivity, and that the compact-open topology on $\rC(\Gamma,\M)$ becomes the quotient topology induced from $\rC(\Gamma,\rA^n)$.

Using \cite[Theorem M, p. 86]{Wigner1973}, every continuous $f:\Gamma\to\M$ lifts to $\tilde{f}:\Gamma\to\rA^n$ viewing both as additive topological groups. This proves surjectivity. We must also show that $\rC(\Gamma,\M)$ acquires the quotient topology, which is what we need for strictness. The map is continuous, so we must check $\rC(\Gamma,\rA^n) \longrightarrow \rC(\Gamma,\M)$ is an open map.

Letting $\pi: \rA^n \to \M$ denote the quotient map, by abuse of notation we also use $\pi$ for the map $\rC(\Gamma,\rA^n) \longrightarrow \rC(\Gamma,\M)$. By virtue of being a homomorphism we must check $\pi_*(\V_{\Gamma,\U})=\V_{\Gamma,\pi(\U)}$ (again using the usual subbasis of the compact-open topology; this is a basis around zero). More precisely we must check $f: \Gamma \to \pi(\U) \subset \M$ lifts to $\tilde{f}: \Gamma \to \U \subset \rA^n$. Viewing $f$ as a map to $\M$ the previous result gives a lift, which need not land in $\U$. Around $\gamma \in \Gamma$ we may always correct this lift in an open neighborhood of $\gamma$ by taking $\tilde{f}-k_\gamma$, where $k_\gamma\in \rK$. Since $\Gamma$ is profinite we can take a finite clopen partition into $\Gamma_i$ (not as subgroups!), allowing us to correct $\tilde{f}$ piecewise by $k_{\gamma_i}$ to land in $\U$. This shows the claim.
\end{proof}

Granting the lemma the proof is straightforward. When $\M=\rA^n$ the claim that $\mathcal{F}$ is well-defined is clear. For general $\M$, consider the presentation $0 \to \rK \to \rA^n \twoheadrightarrow \M$. We can see $\rK \to \rA^n$ is closed, using the weak topology to identify $\rK \simeq \rA^m/\ker f$ for the map $\rA^m \xrightarrow{f} \rA^n$ giving the presentation of $\M$. By virtue of $\rA$ being Hausdorff, $\ker f$ is closed. Thus all Cauchy sequences in $\rK$ converge, and knowing the weak and subspace topologies agree we see $\rK$ is a complete additive subgroup of $\rA^n$, therefore closed.

By the lemma we therefore get an exact sequence
\[0 \longrightarrow \rC(\Gamma,\rK) \longrightarrow \rC(\Gamma,\rA^n) \longrightarrow \rC(\Gamma,\M) \longrightarrow 0\]
and applying it again with $\rA^m \twoheadrightarrow \rK$ we get a right exact sequence
\[\rC(\Gamma,\rA^m) \longrightarrow \rC(\Gamma,\rA^n) \longrightarrow \rC(\Gamma,\M) \longrightarrow 0.\]
For the first two we have $\rC(\Gamma,\rA) \otimes_{p_0, \rA} \rA^n \simeq \rC(\Gamma,\rA^n)$ and similarly for $\rA^m$, so keeping track of topologies we get an isomorphism of topological modules $\rC(\Gamma,\rA) \otimes_{p_0, \rA} \M \simeq \rC(\Gamma,\M)$ where $\M$ has its weak topology.
\end{proof}

Both assumptions are easy to verify for $\R_0^{\sq}\laur{q-1}_p^{\wedge}$.

Lemmas \ref{lem:qdr_laur_cech} and \ref{lem:perf_comparison} give the same comodule description for Laurent crystals that we just discussed for semilinear $\Gamma$-modules. Put $\overline{\rA}^{(m)}=\rA^{(m)}/p$. In the coaction formulation, coherent Laurent crystals are described by replacing $\rC(\Gamma,\overline{\rA})$ and $\rC(\Gamma^2,\overline{\rA})$ with $\overline{\rA}^{(1)}$ and $\overline{\rA}^{(2)}$. For the following lemma, we note there is a natural identification $\rC(\Gamma^s,\rN) \simeq \rC(\Gamma^s,\overline{\rA}) \otimes_{\overline{\rA},p_0}^{\Box} \rN$ for a finite $\overline{\rA}$-module $\rN$ equipped with its weak topology, as we showed in Lemma \ref{lem:weak_top}. We keep the solid decoration in what follows to remind the reader that we are working with modules with their weak topology or canonical solid structure as dualizable objects. In the coherent setting the distinction does not matter as they are identified after passing to the associated solid module.

As usual let $\overline{\rA}=\rA/p$, and let $p_0:\overline{\rA}\to\overline{\rA}^{(s)}$ be the map in the \v{C}ech nerve corresponding to the canonical inclusion $[0] \hookrightarrow [s]$. Recall the maps $\theta_s$ were introduced in Lemma \ref{lem:Cech_function_calc} to identify $\rA^{(s)}$ as subrings of $\rC(\Gamma^s,\rA)$.

\begin{lemma} \label{lem:Cech_finite_coefficients}
For every finitely generated $\overline{\rA}$-module $\M$, the map $\theta_s$ induces an injection
\[\overline{\rA}^{(s)}\otimes_{\overline{\rA},p_0}^{\Box}\M \longrightarrow \rC(\Gamma^s,\M).\]
If $\rN \subset \M$ is a submodule, then inside $\rC(\Gamma^s,\M)$ one further has
\[\overline{\rA}^{(s)}\otimes_{\overline{\rA},p_0}^{\Box}\rN = (\overline{\rA}^{(s)}\otimes_{\overline{\rA},p_0}^{\Box}\M) \cap\rC(\Gamma^s,\rN).\]
\end{lemma}

\begin{proof}
The first claim follows from a variant of Lemma \ref{lem:Cq_criterion} for coefficients, which we will briefly explain. For $f\in \rC(\Gamma^s,\M)$, take its Mahler expansion $f(x)=\sum_{\alpha\in \N^{s(n+1)}} m_\alpha\binom{x}{\alpha}$ where $\binom{x}{\alpha}=\prod_i\binom{x_i}{\alpha_i}$ (here allowing $x_i$ to be coordinates of $\Gamma^s$). Indeed the same argument as Lemma \ref{lem:Cq_criterion} yields $\rC_q(\Gamma^s,\overline{\rA}^+) \simeq \left(\bigoplus_{\alpha\in \N^{s(n+1)}}\overline{\rA}^+ (q-1)^{|\alpha|}\binom{x}{\alpha}\right)^\wedge_{(q-1)}$, so after tensoring with a lattice $\M^+\subset \M$ and Laurentizing we see that $\overline{\rA}^{(s)}\otimes_{\overline{\rA},p_0}^{\Box}\M$ is identified with functions having Mahler expansions of the form $m_{\alpha} = m'_{\alpha} (q-1)^{|\alpha|}$ for $m'_{\alpha} \to 0$ as $|\alpha| \to \infty$. In particular, it is injective.

We now attack the second claim, namely that if $\rN \subset \M$ is a submodule, then
\[\overline{\rA}^{(s)}\otimes_{\overline{\rA},p_0}^{\Box}\rN = (\overline{\rA}^{(s)}\otimes_{\overline{\rA},p_0}^{\Box}\M ) \cap\rC(\Gamma^s,\rN).\] 
It is useful to observe that the map $p_0:\rA\to\rA^{(s)}$ is $p$-completely flat: before applying $(-)[1/[p]_q]_p^{\wedge}$ we see $(p,[p]_q)$-complete flatness follows from flatness of prismatic envelopes (see \cite[Proposition 5.52]{MaoDerivedCrystalline}). Since $\rA^+$ is Noetherian, \cite[Lemma 5.15]{BhattCohenMacaulayAIC} makes $p_0:\rA\to\rA^{(s)}$ $p$-completely flat as it is the $p$-completion of a flat map. Its reduction modulo $p$, the map $p_0:\overline{\rA}\to\overline{\rA}^{(s)}$, is then flat.

Consider now the commutative diagram
\begin{display}
0 \ar{r} & \overline{\rA}^{(s)}\otimes_{\overline{\rA},p_0}^{\Box} \rN \ar{r} \ar[hook]{d} & \overline{\rA}^{(s)}\otimes_{\overline{\rA},p_0}^{\Box}\M \ar{r} \ar[hook]{d} & \overline{\rA}^{(s)}\otimes_{\overline{\rA},p_0}^{\Box}(\M/\rN) \ar[hook]{d} \ar{r} & 0 \\
0 \ar{r} & \rC(\Gamma^s,\rN) \ar{r} & \rC(\Gamma^s,\M) \ar{r} & \rC(\Gamma^s,\M/\rN) \ar{r} & 0.
\end{display}
The upper row is exact as we just explained above, and injectivity follows from the first claim in the lemma. The lower row is clearly exact. Applying the nine lemma (or snake lemma), there is then an exact sequence of cokernels of the vertical maps. Injectivity of the first map
\[\frac{\rC(\Gamma^s,\rN)}{\overline{\rA}^{(s)}\otimes_{\overline{\rA},p_0}^{\Box} \rN} \longrightarrow \frac{\rC(\Gamma^s,\M)}{\overline{\rA}^{(s)}\otimes_{\overline{\rA},p_0}^{\Box} \M}\]
of this exact sequence is equivalent to the desired claim.
\end{proof}

\begin{lemma} \label{lem:Gamma_Frob_improvement}
Let $\M$ be a finitely generated $\overline{\rA}$-module with continuous semilinear $\Gamma$-action. For $k\gg0$, the coaction of $\varphi^{*k}\M$ factors through
\[\rC_q(\Gamma,\overline{\rA}) \otimes_{\overline{\rA},p_0}^{\Box}\varphi^{*k}\M \subset \rC(\Gamma,\varphi^{*k} \M).\]
The same holds with $\Gamma^2$ as needed for the cocycle condition.
\end{lemma}

\begin{proof}
We explain the proof for $\Gamma$, as the argument for $\Gamma^2$ is similar. In Lemma \ref{lem:Cech_finite_coefficients} we extended Lemma \ref{lem:Cq_criterion} to coefficients to characterize $\rC_q(\Gamma,\overline{\rA}) \otimes_{\overline{\rA},p_0}^{\Box} \M \subset \rC(\Gamma,\M)$, which we will use again here.

Now we consider our $\M$ with a continuous semilinear $\Gamma$-action, and view it as a module $\M$ with a coaction 
\[\M \longrightarrow \rC(\Gamma,\rA) \otimes_{p_0,\rA} \M\]
satisfying counit and coassociativity axioms. Here we are using the perspective introduced at the start of this section that one may think of such a module as a comodule over the Hopf algebroid $(\overline{\rA},\rC(\Gamma,\overline{\rA}))$. In Lemma \ref{lem:weak_top} we explained that we can also think of this as a continuous map $\M \to \rC(\Gamma,\M)$ where we equip $\M$ with the weak topology.

What we wish to show is that upon taking a Frobenius pullback we have $\varphi^{*k} \M \to \rC(\Gamma,\varphi^{*k}\M)$ factoring through $\rC_q(\Gamma,\overline{\rA}) \otimes_{\overline{\rA},p_0}^{\Box} \varphi^{*k}\M$. Conceptually the remaining argument is very simple using the observation in Remark \ref{rem:la_Cq}: one observes that such $\Gamma$ modules always have a (lower) $\lambda$-analytic $\Gamma$-action for some $\lambda$, and taking Frobenius pullbacks improves the $\lambda$ until the factorization holds. 

Making this precise of course requires some more work. The notion of $\lambda$-analytic vectors is defined in \cite[Definition 2.10]{Porat_locally_analytic} but we will use the slightly different convention in \cite[Definition 4.1]{HauckLA}. In this convention we put 
\[\rC_{\lambda-\mathrm{an}}(\Gamma,\M) = \left\{\sum_\alpha m_\alpha \binom{x}{\alpha}: \nu_{q-1}(m_\alpha)-\frac{|\alpha|}{p^\lambda(p-1)} \to \infty\right\}.\]
Smaller $\lambda$ is a stronger condition. One also has upper $\lambda$-analytic vectors, replacing $\frac{|\alpha|}{p^\lambda(p-1)}$ with $\sum \nu_p\left(\lfloor \frac{\alpha_i}{p^\lambda} \rfloor!\right)$. One can check
\[\rC_{\lambda-\mathrm{an}}(\Gamma,\M)\subseteq \rC^{\lambda-\mathrm{an}}(\Gamma,\M) \subseteq \rC_{(\lambda+\varepsilon)-\mathrm{an}}(\Gamma,\M)\]
for every $\varepsilon > 0$ (see the similar argument in \cite[Lemma 4.4(2)]{PoratSolidLocallyAnalytic}).

In this specific case the results in \cite{PoratSolidLocallyAnalytic} do not literally apply, but \cite[Remark 5.10]{PoratSolidLocallyAnalytic} actually already contains everything needed for the argument. There Porat observes that for some small $\Gamma_0 \subset \Gamma$ we have $(g-1) (q-1)^i \overline{\rA}^+ \subset (q-1)^{i+1} \overline{\rA}^+$. Picking $\Gamma_0$ small enough to be uniform and picking a coordinate $x: \Z_p^{n+1} \to \Gamma_0$ (we remark this is not a group homomorphism in our case) specified by generators $\gamma_i$ Porat computes for $g = \prod_i \gamma_i^{x_i}$ and $x=(x_0,\ldots,x_n)$ that
\[g(a) = \sum_{\alpha \in \N^{n+1}} \binom{x}{\alpha} (\gamma_0-1)^{\alpha_0} \cdots (\gamma_n-1)^{\alpha_n}(a),\]
concluding that any element of $\overline{\rA}$ belongs to $\rC^{\lambda-\mathrm{an}}(\Gamma_0,\overline{\rA})^{\Gamma_0}$ for $\lambda \gg 0$, using the fact that $(g-1)$ is $q-1$-adic filtration raising. Of course, the same argument applies with finitely generated coefficients $\M$ equipped with their weak topology, with the elements being generators. One must be slightly careful here to run the first step of the argument as we must pick a finitely generated $\Gamma$-stable lattice $\M^+$ over $\overline{\rA}^+$ of $\M$, but this follows by continuity of the $\Gamma$ action for the weak topology on $\M$. Given a finite set of generators $\rS$ for any such lattice $\M^+$, by continuity and compactness of $\Gamma$ we obtain $\Gamma \rS \subset (q-1)^{-\rN} \M^+$. Now we simply replace $\M^+$ by $\overline{\rA}^+ \cdot \Gamma \rS \subset (q-1)^{-\rN} \M^+$. The $\overline{\rA}^+$-module $\overline{\rA}^+ \cdot \Gamma \rS$ is finitely generated as it's a submodule of $(q-1)^{-\rN} \M^+$, and by construction it is $\Gamma$-stable. The same argument gives $(g-1) (q-1)^i \M^+ \subset (q-1)^{i+1} \M^+$.

Thus we obtain that, before even taking a Frobenius pullback, the coaction factors as
\[\M \longrightarrow \rC^{\lambda-\mathrm{an}}(\Gamma,\M).\]
The previous argument gives this for $\rC^{\lambda-\mathrm{an}}(\Gamma_0,\M)$, but since $\Gamma_0\le \Gamma$ is a finite index subgroup we can obtain the desired claim by possibly increasing $\lambda$. The desired result then follows by observing that a Frobenius pullback improves $\lambda$ to $\lambda-k$ in the factorization of the coaction until Remark \ref{rem:la_Cq} applies. Note however that in Remark \ref{rem:la_Cq} we used lower locally analytic vectors, but up to further modifying $\lambda$ by taking Frobenius pullbacks we can obtain the same claim for $\rC_{\lambda-\mathrm{an}}(\Gamma,\M)$ by the earlier inclusions. 
\end{proof}

\begin{proposition} \label{prop:Gamma_to_crys}
Let $\M$ be a finitely generated $\overline{\rA}$-module with continuous semilinear $\Gamma$-action. For $k\gg0$, the module $\varphi^{*k}\M$ arises from a coherent Laurent crystal killed by $p$.
\end{proposition}

\begin{proof}
Choose $k$ large enough for Lemma \ref{lem:Gamma_Frob_improvement} to apply in both degrees one and two, and consider $\varphi^{*k}\M$. Applying Lemmas \ref{lem:Cech_function_calc} and \ref{lem:Cech_finite_coefficients}, we therefore get injective vertical maps in the diagram

\begin{display}[shift left/.default=+0.8ex]
\varphi^{*k}\M \ar[shift left]{r} \ar[shift right]{r} & \rC(\Gamma,\varphi^{*k}\M) \ar{l} \ar[shift right, shift right]{r} \ar{r} \ar[shift left, shift left]{r} & \rC(\Gamma^2,\varphi^{*k}\M) \ar[shift right]{l} \ar[shift left]{l} \\
\varphi^{*k}\M \ar[equal]{u} \ar[shift left]{r} \ar[shift right]{r} & \overline{\rA}^{(1)}\otimes_{\overline{\rA},p_0}^{\Box}\varphi^{*k}\M \ar[hook]{u} \ar{l} \ar[shift right, shift right]{r} \ar{r} \ar[shift left, shift left]{r} & \overline{\rA}^{(2)}\otimes_{\overline{\rA},p_0}^{\Box}\varphi^{*k}\M \ar[hook]{u} \ar[shift right]{l} \ar[shift left]{l}
\end{display}

We are given a full descent datum for $\varphi^{*k}\M$ on the upper row, and what we want to show is that it actually comes from the lower row. The coaction $\varphi^{*k}\M \to \rC(\Gamma,\varphi^{*k}\M)$ factors uniquely through the lower degree one term $\overline{\rA}^{(1)}\otimes_{\overline{\rA},p_0}^{\Box}\varphi^{*k}\M$ by Lemma \ref{lem:Gamma_Frob_improvement}. Note also that the diagram commutes using the description of the vertical maps from Lemma \ref{lem:Cech_function_calc}.

Finally, the counit and coassociativity identities hold in the upper row, hence in the lower row since the vertical maps are injective. Thus the lower row defines a descent datum on $\varphi^{*k}\M$ as a coherent Laurent crystal and the claim follows.
\end{proof}

\begin{lemma} \label{lem:Cech_submodule}
Let $\M$ be a finitely generated $\overline{\rA}$-module with a continuous semilinear $\Gamma$-action. If $\M$ is the evaluation of a coherent Laurent crystal killed by $p$ and $\rN \subset \M$ is $\Gamma$-stable, then $\rN$ also arises from a coherent Laurent crystal.
\end{lemma}

\begin{proof}
Denote by $\mathcal{E}$ the coherent Laurent crystal from the lemma statement. The coaction of $\mathcal{E}$, using the log $q$-de Rham prism and the comodule description from the start of this subsection, is given by
\[\M\longrightarrow \overline{\rA}^{(1)}\otimes_{\overline{\rA},p_0}^{\Box}\M \hookrightarrow \rC(\Gamma,\M).\]
We are given $\rN \subset \M$ equipped with a coaction $\rN \to \rC(\Gamma,\rN)$ as a $\Gamma$-module (using Lemma \ref{lem:weak_top}). The fact that $\rN$ is $\Gamma$-stable amounts to the claim that the diagram
\begin{display}
\M \ar{r} & \rC(\Gamma,\M) \\
\rN \ar[hook]{u} \ar{r} & \rC(\Gamma,\rN) \ar[hook]{u}
\end{display}
commutes. By the second claim of Lemma \ref{lem:Cech_finite_coefficients} we see the coaction restricts to $\rN$, giving
\[\rN \longrightarrow \overline{\rA}^{(1)}\otimes_{\overline{\rA},p_0}^{\Box}\rN.\]
The cocycle and counit conditions follow by a similar argument. Thus $\rN$ has the structure of a coherent Laurent crystal.
\end{proof}

\begin{proof}[Proof of Theorem \ref{thm:laur_root_crys}]
We let $\beta: \M \hookrightarrow \varphi^* \M$ be the map of $\Gamma$-modules we are given. By repeatedly applying the root map, we obtain $\Gamma$-equivariant injections
\[\M \hookrightarrow \varphi^{*k}\M\]
for all $k\ge 1$. By Proposition \ref{prop:Gamma_to_crys} eventually $\varphi^{*k}\M$ comes from a coherent Laurent crystal, and then applying Lemma \ref{lem:Cech_submodule} we see $\M=\rho_{\rA}^{\Gamma,\ast}(\mathcal{E})$ comes from a coherent Laurent crystal. The map $\beta$ also upgrades to a map of coherent Laurent crystals as desired, by applying Proposition \ref{prop:Gamma_mod_ff} (note that the functor there is also compatible with $\varphi^*$). Note that the result applies to coherent crystals by passing to dualizables in Proposition \ref{prop:Gamma_mod_ff} and then using Corollary \ref{cor:laur_t_struct} and the fact that $\rho_{\rA}^{\Gamma,\ast}$ is $t$-exact.
\end{proof}

As promised, we now investigate other bases.

\begin{corollary} \label{cor:coh_laur_other_bases}
Both Proposition \ref{prop:Gamma_mod_ff} and Theorem \ref{thm:laur_root_crys} hold over $\W(k)[\zeta_{p^\ell}]$ for $\ell > 0$, by changing to $\Gamma = (1+p^\ell \Z_p) \ltimes \Z_p(1)^n$ and keeping the same base $\rA=\R_0^{\sq}\laur{q-1}^{\wedge}_p$.

For $\ell=0$, again keeping the same base $\rA=\R_0^{\sq}\laur{q-1}^{\wedge}_p$ but using $\Gamma = \Z_p^\times \ltimes \Z_p(1)^n$, Proposition \ref{prop:Gamma_mod_ff} holds and Theorem \ref{thm:laur_root_crys} also holds.
\end{corollary}

\begin{proof}
For $\ell > 0$ the modification is straightforward, and can be done by simply setting the prismatic ideal in the $q$-de Rham prism to be $\varphi^{\ell-1}([p]_q)$. The toric factor acts by $\T_i\mapsto q^{p^\ell b_i}\T_i$, so this action is trivial modulo $\varphi^{\ell-1}([p]_q)$.

For $\ell=0$ and $p$ odd, one can evade redoing the calculations in Proposition \ref{prop:Gamma_mod_ff} in favor of a simpler calculation. We will first show that the pullback along
\[f: [\W(k)[\zeta_p]^{\mathcal{L},\Box}/\mu_{p-1}] \xlongrightarrow{} \W(k)^{\mathcal{L},\Box}\]
is fully faithful on $\D_{\Box}(-)$. Here the $\mu_{p-1}$ action comes from functoriality and the Galois action on $\W(k)[\zeta_p]/\W(k)$. One can use the same argument as in Proposition \ref{prop:Gamma_mod_ff} to reduce to checking if the pushforward of the structure sheaf on the source is $\O_{\W(k)^{\mathcal{L},\Box}}=\Laur$ via the projection formula after we verify $f$ is prim (for the six functor formalism on analytic stacks). If we know this, then it is formal to deduce the claim in the general case for $[(\R_0[\zeta_p],\mathcal{M})^{\mathcal{L},\Box}/\mu_{p-1}] \to (\R_0,\mathcal{M})^{\mathcal{L},\Box}$: one may use base change for prim maps to deduce the required isomorphism. Verifying $f$ is prim turns out to be formal from our results so far. One can use the diagram
\begin{display}
{[\AnSpf \B/\Z_p^\times]} \ar{dd} \ar{r} & {[\AnSpf \W(k)\laur{q-1}^{\wedge}_p/\Z_p^\times]} \ar{d}{\rho^{\Gamma}_{\W(k)\laur{q-1}^{\wedge}_p}} \\
 & {[\W(k)[\zeta_p]^{\mathcal{L},\Box}/\mu_{p-1}]} \ar{d}{f} \\
{\AnSpf \W(k)\laur{\tilde{p}}^{\wedge}_p} \ar{r} & {\W(k)^{\mathcal{L},\Box}}
\end{display}
where the outer rectangle is Cartesian (with $\tilde{p}$ as in \cite[Notation 3.8.9]{APC}). We may compute $\B$ as a prismatic envelope using Lemma \ref{lem:mixed_BK_fiber_product}; the analytic realization $(-)^{\Box}$ of Definition \ref{def:box_analytic_realization} and passage to Laurent stacks preserve finite limits.

The upper right arrow $\rho^{\Gamma}_{\W(k)\laur{q-1}^{\wedge}_p}$ is proper, as in the argument of Proposition \ref{prop:Gamma_mod_ff}. Moreover the argument also shows $\rho^{\Gamma}_{\W(k)\laur{q-1}^{\wedge}_p,*}\O \simeq \O$. The left vertical map $[\AnSpf \B/\Z_p^\times] \longrightarrow \AnSpf \W(k)\laur{\tilde{p}}^{\wedge}_p$ is proper by the same argument as in Proposition \ref{prop:Gamma_mod_ff}, using Lemma \ref{lem:cech_descendable} and properness of $\B \Z_p^\times$. Hence the composite of the right vertical maps is proper, since the bottom map is an effective epimorphism of analytic stacks by Proposition \ref{prop:ordinary_BK_analytic_presentation}. Now \cite[Corollary 4.7.5(ii)]{HeyerMannSixFunctors} shows that $f$ is prim, which by \cite[Corollary 4.5.11(ii) and Lemma 4.5.13(ii)]{HeyerMannSixFunctors} gives the required projection formula and base change properties for the argument to work.

Thus we check $\O_{\W(k)^{\mathcal{L},\Box}} \simeq \R f_* \O_{[\W(k)[\zeta_p]^{\mathcal{L},\Box}/\mu_{p-1}]}$. As previously noted, we will use the same argument as in Proposition \ref{prop:Gamma_mod_ff} to reduce to checking if the pushforward of the structure sheaf on the source is $\Laur$. Using Proposition \ref{prop:stack_site_compare}, we may work site-theoretically, so one may reduce to checking the stronger claim
\[\W(k)\llbracket \tilde{p} \rrbracket \xlongrightarrow{\sim} \Prism_{\W(k)[\zeta_p]/\W(k)\llbracket \tilde{p} \rrbracket}^{h \mu_{p-1}}.\]
Here we use conservativity of evaluation on $\W(k)\llbracket \tilde{p} \rrbracket$ to make this reduction. This can be checked directly, using the Hodge-Tate comparison. Indeed, $\L_{\W(k)[\zeta_p]/\W(k)} \simeq \widehat{\Omega}^1_{\W(k)[\zeta_p]/\W(k)}$ compatibly with the Galois action on $\W(k)[\zeta_p]/\W(k)$, but the target has no ordinary $\mu_{p-1}$ invariants; higher exterior powers also follow, as in this case they are just shifts of $\widehat{\Omega}^1_{\W(k)[\zeta_p]/\W(k)}$ with Breuil-Kisin twists. The $\mu_{p-1}$ invariants are exact so the desired result follows from
\[\W(k)\llbracket \tilde{p} \rrbracket = \W(k)\llbracket q-1 \rrbracket^{\mu_{p-1}}.\]

We now explain the second assertion, the proof of Theorem \ref{thm:laur_root_crys}. The \v{C}ech calculation carries over for the log prism $(\R_0^{\sq}\llbracket \tilde{p}\rrbracket, (\tilde{p}),\mathcal{M})$ -- this prism is constructed analogously to the log $q$-de Rham prism via a framing, but we note that one must be cautious as it does not possess a full $\Gamma$ action. Its first \v{C}ech term is
\[(\R_0^{\sq}\llbracket \tilde{p}\rrbracket)^{(1)} \simeq \R_0^{\sq}\llbracket \tilde{p}_0,\tilde{p}_1\rrbracket[u_1^{\pm1},\ldots,u_n^{\pm1}] \left\{ \frac{\tilde{p}_1}{\tilde{p}_0}, \frac{u_1-1}{\tilde{p}_0}, \ldots, \frac{u_n-1}{\tilde{p}_0} \right\}_{\delta,(p,\tilde{p}_0)}^{\wedge},\]
where $u_i=\T_{i,1}/\T_{i,0}$ and $\delta(u_i)=0$.
As we work modulo $p$ always, what is really needed is an analogue of Lemma \ref{lem:Cq_criterion}. This doesn't seem to quite work, but we can use the following trick. By a similar argument we may deduce
\[\left((\R_0^{\sq}\llbracket \tilde{p}\rrbracket)^{(1)} \otimes_{\R_0^{\sq}\llbracket \tilde{p}\rrbracket} \R_0^{\sq}\llbracket q-1\rrbracket\right)/p \simeq \left(\bigoplus_{\alpha \in\N^{n+1}} (\R_0^{\sq}/p)\llbracket q-1 \rrbracket \tilde{p}^{\alpha_0}(q-1)^{|(\alpha_1,\ldots,\alpha_n)|}\binom{x}{\alpha}\right)^{\wedge}_{(q-1)}.\]
Indeed for the toric directions we can repeat verbatim the same argument as already appearing in Lemma \ref{lem:Cech_function_calc} for $\eta_{i,m}$ when $i>0$, which then works identically through the method of Lemma \ref{lem:Cq_criterion}. For the remaining part one similarly observes that
\[\W(k)\llbracket \tilde{p}_0, \tilde{p}_1 \rrbracket\{\tilde{p}_1/\tilde{p}_0\}^{\wedge}_{\delta,(p,\tilde{p}_0)}/p \simeq \left(\bigoplus_{\alpha_0\in \N} k\llbracket \tilde{p}_0 \rrbracket \tilde{p}^{\alpha_0} \binom{x_0}{\alpha_0}\right)^{\wedge}_{\tilde{p}_0} \subset \rC(1+p\Z_p, k\llbracket \tilde{p}_0\rrbracket).\]
Combining these gives the claimed identification, noting that we need base extension primarily to deal with the toric part (otherwise $q-1$ is not an element of the ring, so the identification makes no sense).

This identification shows that a single Frobenius pullback makes 
\[\mathcal{E}\in \Coh((\R_0 \otimes_{\W(k)} \W(k)[\zeta_p],\mathcal{M})_\Prism,\Laur)_{p=0}^{h\mu_{p-1}}\]  
descend to $\Coh((\R_0,\mathcal{M})_\Prism,\Laur)$. Indeed comparing the previous identification to $(\R_0^{\sq}\llbracket q-1\rrbracket)^{(1)}/p$ shows a single Frobenius pullback lands the coaction of a Laurent crystal in the Laurentization of $\left((\R_0^{\sq}\llbracket \tilde{p}\rrbracket)^{(1)} \otimes_{\R_0^{\sq}\llbracket \tilde{p}\rrbracket} \R_0^{\sq}\llbracket q-1\rrbracket\right)/p$. But then taking $\mu_{p-1}$ invariants further lands the coaction in $(\R_0^{\sq}\llbracket \tilde{p}\rrbracket)^{(1)}/p [1/\tilde{p}]$ as desired, giving the crystal.

Thus it follows that given such a $\Z_p^\times \ltimes \Z_p(1)^n$-module over $\rA$ some high Frobenius pullback comes from a crystal. To finish the argument, we note that the description of 
\[\left((\R_0^{\sq}\llbracket \tilde{p}\rrbracket)^{(1)} \otimes_{\R_0^{\sq}\llbracket \tilde{p}\rrbracket} \R_0^{\sq}\llbracket q-1\rrbracket\right)/p\]
yields the corresponding claim in Lemma \ref{lem:Cech_finite_coefficients} as the required points (e.g. flatness of the maps in the \v{C}ech nerve) again hold. We also get the description for the second term of the \v{C}ech nerve, since this is $(\R_0^{\sq}\llbracket \tilde{p}\rrbracket)^{(1)} \widehat{\otimes}_{\R_0^{\sq}\llbracket \tilde{p}\rrbracket} (\R_0^{\sq}\llbracket \tilde{p}\rrbracket)^{(1)}$. Thus we obtain the analogue of Lemma \ref{lem:Cech_submodule}, and can conclude Theorem \ref{thm:laur_root_crys} in the same way.
\end{proof}

\section{Finitely generated unit Laurent crystals}

\subsection{Basic constructions}

Following \cite{EmertonKisin2} and the more general version in \cite{BhattLurieModpRH} (which explains how to remove finiteness hypotheses), we first generalize the notion of finitely generated unit modules.

\begin{definition} \label{def:module_fgu}
Let $(\rA^+,\I)$ be an orientable prism with $\I=(d)$. Equip $\rA=\rA^+[1/d]_p^{\wedge}$ with $\D_{\Box}(\rA)$ as defined in \ref{def:laur_category}. Note $\rA$ has a natural Frobenius $\varphi$ induced by the Frobenius on $\rA^+$ (since $\varphi(d)\equiv d^p \pmod{p}$).

Write $\rA\langle \rF\rangle$ for
\[\widehat{\bigoplus}_{i\ge 0} \rA \rF^i,\]
where $\rF x = \varphi(x) \rF$. We put $\D_{\Box}(\rA\langle \rF\rangle)$ for $p$-complete solid modules over this noncommutative analytic ring, with the analytic ring structure induced from $\rA$.

For $\M\in\D_{\Box}(\rA\langle \rF\rangle)$, we say $\M$ is \emph{unit} if
\[\varphi_{\rA}^*\M:=\rA\otimes^{\L,\Box}_{\rA,\varphi_{\rA}}\M\longrightarrow\M\]
given by $a\otimes m\mapsto a\rF_\M(m)$ is an equivalence. A unit object is \emph{finitely generated unit} if it lies in the thick subcategory generated by $\rA\langle \rF\rangle$.
\end{definition}

We remark that we can naturally view these as certain unit Frobenius modules, so $\D_{\fgu}(\rA)\subset \D_{\Box}(\rA)^{\varphi=1}$ as a full subcategory. In fact one can show the category on a fixed prism is insensitive to the solid structure, but we will still remember it so that our crystals live in an ambient category which is a stack.

\begin{lemma} \label{lem:solid_Frob_perf_comparison}
Let $(\rA^+,(d))$ be an orientable prism, and let $\rA=\rA^+[1/d]_p^{\wedge}$ with $\D_{\Box}(\rA)$ as in Definition \ref{def:laur_category}. Then 
\[\mathrm{Thick}_{\rA\langle \rF \rangle} \D_{\Box}(\rA\langle \rF \rangle) \simeq \mathrm{Thick}_{\rA\langle \rF \rangle} \D(\rA\langle \rF \rangle)\]
via $\M \mapsto \M(*)$.

The solid modules in the source thick subcategory are all $p$-completely nuclear.
\end{lemma}

\begin{proof}
The functor in question is exact and sends the solid variant of $\rA\langle \rF \rangle$ to $\rA\langle \rF \rangle$. It therefore suffices to compare derived endomorphisms of $\rA\langle \rF \rangle$ and its solid counterpart.

We can compute
\[\RHom_{\D_{\Box}(\rA\langle \rF\rangle)}(\rA\langle \rF\rangle,\rA\langle \rF \rangle)\simeq \RHom_{\D_{\Box}(\rA)}(\rA,\rA \langle \rF \rangle) \simeq \rA \langle \rF \rangle.\]
To justify this, by derived $p$-completeness it suffices to check equivalences modulo $p^n$, which makes the argument slightly easier as then $\rA\langle \rF \rangle \dquot p^n = (\rA \dquot p^n)[\rF]$. The first equivalence uses the free/forgetful adjunction, and the second equivalence is then clear. This then matches the corresponding calculation in ordinary modules, so it follows $\mathrm{Thick}_{\rA\langle \rF \rangle} \D(\rA\langle \rF \rangle) \simeq \mathrm{Thick}_{\rA\langle \rF \rangle} \D_{\Box}(\rA\langle \rF \rangle)$.

Now we address the nuclearity claim. The unit $\rA^+$ of $\D_{\Box}(\Spa(\rA^+,\rA^+))$ is nuclear, and since $\rA$ is the filtered colimit of multiplication by $d$ in derived $p$-complete solid modules it is also $p$-completely nuclear. The $p$-completely nuclear objects are closed under colimits, so it follows that $\rA\langle \rF\rangle=\widehat{\bigoplus}_{i\ge0}\rA \rF^i$ is $p$-completely nuclear. Nuclearity is preserved under the operations cutting out a thick subcategory, for example retracts or finite limits, so the claim follows.
\end{proof}

\begin{definition} \label{def:fgu}
Let $(\X,\mathcal{M})$ be a bounded fs log $p$-adic formal scheme. We say that a crystal
\[\mathcal{E}\in \D_{\Box}((\X,\mathcal{M})_\Prism,\Laur)^{\varphi=1}\]
is \emph{finitely generated unit} if for every object $(\rA,\I,\M_{\Spf \rA})^a$ of $(\X,\mathcal{M})_\Prism$ we have $\mathcal{E}((\rA,\I,\M_{\Spf \rA})^a)\in \D_{\fgu}(\rA[1/\I]^{\wedge}_p)$ in the sense of Definition \ref{def:module_fgu}. We write $\D_{\fgu}((\X,\mathcal{M})_\Prism,\Laur)$ for the full subcategory of finitely generated unit crystals, and $\Mod_{\fgu}((\X,\mathcal{M})_\Prism,\Laur)$ for those concentrated in degree zero.
\end{definition}

\begin{remark} \label{rem:fgu_mod_p}
Let $\B$ be a $p$-complete ring, not necessarily commutative, where $p$ is central. The typical example we will want to apply this to in this article is $\rA\langle \rF \rangle$. If $\M \in \widehat{\D}(\B)$ is derived $p$-complete, then we can show that
\[\M \in \Perf(\B) \Longleftrightarrow (\B/p)\otimes^{\L}_{\B}\M\in\Perf(\B/p).\]
This follows from a similar argument to \cite[Tag \href{https://stacks.math.columbia.edu/tag/09AW}{09AW}]{StacksProject} (which cannot be directly applied, as it works with commutative rings). We can even use this in our present setting, since passing through Lemma \ref{lem:solid_Frob_perf_comparison} we can run the same argument in the solid setting for $\rA\langle \rF \rangle$, thus detecting the finitely generated unit condition modulo $p$.
\end{remark}

When $(\X,\mathcal{M})$ satisfies Assumption \ref{ass:horizontal_nc} over $\W(k)[\zeta_p]$, Remark \ref{rem:qdR_test_sufficient} shows that it suffices to test Definition \ref{def:fgu} on the log $q$-de Rham prisms associated to framed log affines.

Our first main goal will be to give an explicit local description of this category, as well as deduce strict \'{e}tale descent. We note that while $\D_{\Box}((\X,\mathcal{M})_\Prism,\Laur)^{\varphi=1}$ is a strict \'{e}tale stack (by Corollary \ref{cor:laur_strict_etale}), it is not immediate that the finitely generated unit subcategory is also a strict \'{e}tale stack. While ordinary dualizable objects can be detected termwise in a limit, because we work over $\rA\langle \rF \rangle$ we use the thick subcategory generated by the unit (which satisfies no such property). Instead, descent of the fgu condition will use its characterization by coherent roots, which we now study.

\begin{definition} \label{def:roots}
Let $\rA$ be as in \ref{def:module_fgu}, and assume that $\rA$ is Noetherian regular. If $\M$ is a unit Frobenius module over $\rA$, a \emph{root} of $\M$ is a coherent $\rA$-submodule $\M_0\subset \M$ together with an injective map $\beta:\M_0 \hookrightarrow \varphi_{\rA}^*\M_0$ along with an identification
\[\M\simeq \colim_n \varphi_{\rA}^{*n}\M_0.\]
It is \emph{minimal} if it has no proper subroot, and it is distinguished if it is the unique minimal root. We call $\M$ the \emph{unitalization} of $\M_0$.

One can make the analogous definitions in $\QCoh_{\Box}((\X,\mathcal{M})_\Prism,\Laur)^{\varphi=1}$. When referring to a root of a module we mean the former, and when referring to a root of a crystal we mean the latter.
\end{definition}

The essential idea of the proof of strict \'{e}tale descent is to characterize the finitely generated unit condition in terms of distinguished minimal roots \emph{as crystals} in Proposition \ref{prop:root_crys_minimal}, and use that coherent crystals satisfy strict \'{e}tale descent to glue these. Once we have this, showing the relevant pullbacks of minimal roots are minimal forces them to agree on overlaps; this allows us to glue local distinguished minimal roots into a global one.

We now specialize to the situation where $(\Spf \R,\mathcal{M})$ is a framed log affine in the sense of Definition \ref{def:fr_log_aff} and put $\rA=\R_0^{\sq}\laur{q-1}^{\wedge}_p$ as before. In this situation we can extract a more practical criterion for verifying the finitely generated unit condition.

\begin{lemma} \label{lem:fgu_actually_fgu}
For $\rA = \R_0^{\sq}\laur{q-1}^{\wedge}_p$, a unit $\rA\langle \rF \rangle$-complex $\M$ is perfect over $\rA\langle \rF \rangle$ if and only if $\M$ is $t$-bounded and each $\H^i(\M)$ is finitely generated as a left $\rA\langle \rF \rangle$-module.
\end{lemma}

\begin{proof}
Supposing $\M$ is perfect over $\rA\langle \rF \rangle$, then so is $\M \dquot p$. Then \cite[Proposition 11.3.9]{BhattLurieModpRH} shows all $\H^i(\M \dquot p)$ are finitely generated over $\overline{\rA}[\rF]$. Now
\[\H^i(\M)/p \hookrightarrow \H^i(\M \dquot p)\]
is a unit submodule, hence finitely generated. Then lifting generators to get a map $\rA\langle \rF \rangle^{\oplus n} \twoheadrightarrow \H^i(\M)$ and applying derived Nakayama to the derived $p$-complete cokernel for surjectivity shows $\H^i(\M)$ is finitely generated.

Conversely, suppose all $\H^i(\M)$ are finitely generated over $\rA\langle \rF \rangle$. By Remark \ref{rem:fgu_mod_p} we need the analogous claim for $\M \dquot p$. Using the exact sequences
\begin{display}
0 \ar{r} & \H^i(\M)/p \ar{r} & \H^i(\M \dquot p) \ar{r} & \H^{i+1}(\M)[p] \ar{r} & 0.
\end{display}
we reduce to proving the assertion that if $\M$ is a finitely generated unit left module over $\rA\langle \rF \rangle$ then $\M[p]$ is too. This requires some work, as we work in a noncommutative setting (and the claim is false without the unit condition).

Modulo $p^\ell$, this is already known. By \cite[Lemma 1.2]{EmertonCoherentRings} we see $(\M/p^\ell)$ is finitely presented over $\rA/p^\ell[\rF]$, and moreover $\rA/p^\ell[\rF]$ is left coherent by \cite[Proposition 1.3]{EmertonCoherentRings}. It follows that $(\M/p^\ell)[p] \in \Mod^{\mathrm{fg}}(\rA/p^\ell[\rF])$ as it is the kernel of multiplication by $p$ on a finitely presented module, and it is clearly unit. For the actual claim we want, it will then suffice to show that $\M[p] \hookrightarrow (\M/p^\ell)[p]$ for some $\ell$: the cokernel of this map is finitely generated unit, thus showing $\M[p]$ arises as a kernel of a map of finitely generated unit modules over $\rA/p^\ell[\rF]$. Then left coherence shows $\M[p]$ is finitely generated over $\rA/p^\ell[\rF]$, which suffices.

To obtain $\M[p] \hookrightarrow (\M/p^\ell)[p]$, we show $\M[p^\infty]=\M[p^\ell]$ for some $\ell$. As $\M$ is derived $p$-complete it contains no $p$-divisible submodules, so showing $\M[p^\ell] \hookrightarrow \M[p^\infty]$ is surjective modulo $p$ suffices (as the cokernel would consist of $p$-divisible elements of the derived $p$-complete module $\M/\M[p^\ell]$). Considering the modules
\[\mathrm{im}(\M[p^k] \longrightarrow \M/p\M) \subset (\M/p\M)\]
we can see that these are all finitely generated unit, since they are unit submodules of the finitely generated unit module $\M/p\M$. But then taking distinguished minimal roots, since $\overline{\rA}$ is Noetherian this chain necessarily stabilizes at some $\ell$, upon which $\M[p^\ell] \hookrightarrow \M[p^\infty]$ is surjective modulo $p$. The claim of the lemma then follows.
\end{proof}

In Theorem \ref{thm:qdR_crys_Gamma_mod}, we will identify finitely generated unit crystals with certain $(\varphi,\Gamma)$-modules.

A useful consequence of this is that we may completely dispense with any condensed mathematics, at least locally. We will find that instead we can work with $\D_{\fgu}(\rA)^{h\Gamma}$, by which we will mean the category of unit objects in
\[\lim_{[s]\in \Delta} \mathrm{Thick}_{\rC(\Gamma^s,\rA\langle \rF \rangle)} \D(\rC(\Gamma^s,\rA\langle \rF \rangle))\]
with no solid structures involved.

With this, we can prove the following explicit description of $\D_{\fgu}((\Spf \R,\mathcal{M})_\Prism, \Laur)$.

\begin{theorem} \label{thm:qdR_crys_Gamma_mod}
Let $(\Spf \R,\mathcal{M})$ be a framed log affine in the sense of Definition \ref{def:fr_log_aff} and put $\rA=\R_0^{\sq}\laur{q-1}^{\wedge}_p$. Evaluation on the log $q$-de Rham prism $(\R_0^{\sq}\llbracket q-1\rrbracket,[p]_q,\mathcal{M})^a$ induces an equivalence
\[\rho_{\rA}^{\Gamma,*}:\D_{\fgu}((\Spf \R,\mathcal{M})_\Prism, \Laur) \xlongrightarrow{\sim} \D_{\fgu}(\rA)^{h\Gamma}.\]
\end{theorem}

\begin{proof}
First, we note that the functor is fully faithful. This follows from Proposition \ref{prop:Gamma_mod_ff}, but for the solid version of the target given by 
\[\D_{\Box}(\rA)^{h\Gamma, \varphi=1} = \lim_{[s]\in \Delta} \D_{\Box}(\rC(\Gamma^s,\rA))^{\varphi=1}.\]
The essential image, by definition of $\D_{\fgu}((\Spf \R,\mathcal{M})_\Prism, \Laur)$, satisfies the property that upon evaluation on the log prism $\rC(\Gamma^s,\rA^+)$ we get a unit Frobenius module in the category $\mathrm{Thick}_{\rC(\Gamma^s,\rA\langle \rF \rangle)} \D_{\Box}(\rC(\Gamma^s,\rA\langle \rF \rangle))$. By Lemma \ref{lem:solid_Frob_perf_comparison}, it follows that equivalently the functor lands in $\D_{\fgu}(\rA)^{h\Gamma}$, by which we mean unit objects in
\[\lim_{[s]\in \Delta} \mathrm{Thick}_{\rC(\Gamma^s,\rA\langle \rF \rangle)} \D(\rC(\Gamma^s,\rA\langle \rF \rangle)).\]

Thus we obtain the desired fully faithful functor
\[\rho_{\rA}^{\Gamma,*}:\D_{\fgu}((\Spf \R,\mathcal{M})_\Prism, \Laur) \longrightarrow \D_{\fgu}(\rA)^{h\Gamma}.\]
Now we turn to essential surjectivity. We first prove essential surjectivity onto the heart of the target modulo $p$. Put $\overline{\rA}=\rA/p$, and let $\M$ be a finitely generated unit $\overline{\rA}$-module with continuous semilinear $\Gamma$-action.

Since $\overline{\rA}$ is regular and $\rF$-finite, \cite[Corollary 2.26]{BlickleGamma} gives a distinguished minimal root
\[\beta:\M_{\min} \hookrightarrow \varphi^*\M_{\min}.\]
Every $\Gamma$-translate is again a root. Since it is a distinguished minimal root, it becomes $\Gamma$-equivariant and Theorem \ref{thm:laur_root_crys} lifts it to a coherent Laurent crystal $\mathcal{E}_0$. Its unitalization $\mathcal{E}:= \colim_k \varphi^{*k} \mathcal{E}_0$ is necessarily finitely generated unit\footnote{\label{fn:root_unitalization_fgu}If $\B=\B^+[1/\I]_p^{\wedge}$ for a log prism $(\B^+,\I, \M_{\Spf \B^+})^a$, let $\M_0,\M$ be the evaluations of $\mathcal{E}_0,\mathcal{E}$. We have
\[\M\simeq\bigl[\B \langle \rF \rangle \otimes_{\B}^{\L,\Box}\M_0\xlongrightarrow{1-\widetilde{\beta}}\B \langle \rF \rangle \otimes_{\B}^{\L,\Box}\M_0\bigr]\]
where $\widetilde{\beta}$ is induced by $\beta:\M_0\to\varphi^*\M_0$, showing $\mathcal{E}$ is finitely generated unit.}, and satisfies $\rho_{\rA}^*\mathcal{E}\simeq\M$, giving the desired preimage for objects in the heart of the target $\D_{\fgu}(\rA)^{h\Gamma}$ modulo $p$. 

This suffices to prove essential surjectivity. Since evaluation on the log prism $\rA^+$ is exact, in particular compatible with shifts and fibers, a Postnikov d\'{e}vissage (see \cite[Lemma 5.4.3]{HauckStacky}) on the target upgrades essential surjectivity on the heart modulo $p$ to essential surjectivity on the derived categories. Thus $\rho_{\rA}^{\Gamma,*}$ is an equivalence modulo $p$. Lemma \ref{lem:categorical_nakayama} now gives the claim.
\end{proof}

\begin{remark} \label{rem:weak_top_fgu}
We claim that objects in $\Mod_{\fgu}(\overline{\rA})^{h\Gamma}$ can be understood as $\overline{\rA}[\rF]$ modules with their weak topology and a continuous $\Gamma$ action. Moreover, modulo $p$ for a minimal root the $\Gamma$ action is continuous for the weak topology.

One would like to use a variant of Lemma \ref{lem:weak_top} for noncommutative rings. Despite the hypotheses failing, the argument goes through if we establish Lemma \ref{lem:strict_exact} for finitely generated unit modules, which we briefly sketch here. Modulo $p$ we have access to minimal roots. Given
\[0 \longrightarrow \overline{\rK} \longrightarrow \overline{\rA}[\rF]^\ell \longrightarrow \overline{\M} \longrightarrow 0\] 
we may observe that $\overline{\rA}[\rF]^\ell = \colim_n \bigoplus_{i\le n} \overline{\rA}^\ell \rF^i$ and $\overline{\M} = \colim_n \varphi^{\ast n} \overline{\M}_0$ for a minimal root $\overline{\M}_0$. We use the presentation induced by the root $\overline{\M}_0$ to obtain the rightmost surjection. Taking the kernel for each $n$ gives a colimit presentation of $\overline{\rK}$ as well. Using the root presentations, one can verify the weak and subspace topologies agree for $\overline{\rK}$ by using the claim for roots. For strict exactness after applying $\rC(\Gamma,-)$, in the nontrivial lifting step of Lemma \ref{lem:strict_exact} a continuous function $\Gamma \to \overline{\M}$ factors through some term $\varphi^{\ast n} \overline{\M}_0$ where the argument of Lemma \ref{lem:strict_exact} still applies; the same holds for the other colimit presentations. With these two claims established the remainder of Lemma \ref{lem:weak_top} goes through.

Finally, we note that by combining Remark \ref{rem:fgu_mod_p} and Lemma \ref{lem:fgu_actually_fgu} we get $\Mod_{\fgu}(\rA)^{h\Gamma} \simeq \lim_n \Mod_{\fgu}(\rA/p^n)^{h\Gamma}$, and the claim modulo $p^n$ can be deduced similarly as roots also exist in this setting. We omit the full details here, as we will only ever need these claims modulo $p$.
\end{remark}

Showing $t$-exactness is more subtle, because when working in the solid setting it is unclear if $\rho_{\rA}^{\Gamma,*}$ is $t$-exact when we use the larger source category $\D_{\Box}((\Spf \R,\mathcal{M})_\Prism, \Laur)$. However the finitely generated unit condition is sufficient to deal with this issue.

\begin{lemma} \label{lem:fgu_local_tstruct}
In Theorem \ref{thm:qdR_crys_Gamma_mod}, the $t$-structure on solid Laurent crystals restricts to the source $\D_{\fgu}((\Spf \R,\mathcal{M})_\Prism, \Laur)$ and $\rho_{\rA}^{\Gamma,*}$ is $t$-exact.
\end{lemma}

\begin{proof}
We note that in the proof of essential surjectivity in Theorem \ref{thm:qdR_crys_Gamma_mod}, implicitly the preimages constructed are all in the heart of $\D_{\Box}((\Spf \R,\mathcal{M})_\Prism, \Laur)$. Indeed mod $p$ the preimage is a unitalization of a coherent root, and coherent crystals are in the heart of the solid $t$-structure by Corollary \ref{cor:laur_t_struct} (note also that $\varphi^*\mathcal{E}_0$ remains coherent even for the derived Frobenius pullback, since we can test coherence after evaluation on $\rA$ where $\varphi$ is flat). Working generally, we have $\Mod_{\fgu}(\rA)^{h\Gamma} = \lim_\ell \Mod_{\fgu}(\rA/p^\ell)^{h\Gamma}$, so it suffices to construct compatible preimages modulo $p^\ell$. Leveraging derived full faithfulness, we may do so by taking iterated extensions and using the mod $p$ preimages.

The desired claim is then completely formal using the inverse equivalence in Theorem \ref{thm:qdR_crys_Gamma_mod}. Indeed, we have seen the inverse gives a $t$-exact functor
\[\D_{\fgu}(\rA)^{h\Gamma} \longrightarrow \D_{\Box}((\Spf \R,\mathcal{M})_\Prism, \Laur)^{\varphi=1}\]
since the preimages constructed lie in the heart. The essential image consists precisely of the finitely generated unit Laurent crystals, and on this subcategory the inverse $\rho_{\rA}^{\Gamma,*}$ must then be $t$-exact.
\end{proof}

By applying Theorem \ref{thm:laur_root_crys} we can prove a root criterion for the finitely generated unit condition.

\begin{proposition} \label{prop:root_crys_minimal}
Let $(\Spf \R,\mathcal{M})$ be a framed log affine. A crystal 
\[\mathcal{E}\in \D_{\Box}((\Spf \R,\mathcal{M})_\Prism,\Laur)^{\varphi=1}\] 
is finitely generated unit if and only if it is $t$-bounded and each crystal $\H^i(\mathcal{E}\dquot p)$ admits a distinguished minimal root in $\Coh((\Spf \R,\mathcal{M})_\Prism,\Laur)_{p=0}$.
\end{proposition}

\begin{proof}
Combining Theorem \ref{thm:qdR_crys_Gamma_mod}, Lemma \ref{lem:fgu_actually_fgu} and Remark \ref{rem:fgu_mod_p} to work modulo $p$, we see that each $\H^i(\mathcal{E}\dquot p)$ admitting a distinguished minimal root in $\Coh((\Spf \R,\mathcal{M})_\Prism,\Laur)_{p=0}$ is sufficient. Indeed we see the unitalization of a coherent crystal root is a finitely generated unit Laurent crystal, by the same argument as in the footnote \ref{fn:root_unitalization_fgu}.

Conversely, take a finitely generated unit Laurent crystal. Reducing modulo $p$ using Theorem \ref{thm:qdR_crys_Gamma_mod} and Lemma \ref{lem:fgu_actually_fgu} we reduce to the case where $\mathcal{E}$ corresponds to $\M \in \Mod_{\fgu}(\rA)_{p=0}^{h\Gamma}$. Here we also implicitly use Lemma \ref{lem:fgu_local_tstruct}. By the previous argument in Theorem \ref{thm:qdR_crys_Gamma_mod}, by \cite[Corollary 2.26]{BlickleGamma} there is a distinguished $\Gamma$-equivariant minimal root $\M_0$ of $\M$. Using Theorem \ref{thm:laur_root_crys} this root arises from a crystal $\mathcal{E}_0$. Together with Proposition \ref{prop:Gamma_mod_ff}, this theorem identifies roots as Laurent crystals with $\Gamma$-equivariant roots of $\M$. Since $\M_0$ has the stronger property of being a distinguished minimal root, the claim follows.
\end{proof}

\subsection{Strict \'{e}tale descent}

In Proposition \ref{prop:root_crys_minimal}, we have already produced the crucial observation we will need for descent: the strategy will be to show that for a strict \'{e}tale cover, if we know the Laurent crystals glue then so do distinguished minimal roots, as uniqueness forces them to agree on overlaps. This glosses over an important point though, which is that we need pullback stability for minimal roots (by virtue of everything being regular $\rF$-finite, \cite{BlickleGamma} will again make these minimal roots automatically distinguished minimal).

\begin{lemma} \label{lem:completion_minimal}
Let $\rA$ be a regular local $\rF$-finite $\F_p$-algebra with completion $\widehat{\rA}$, and let $\rN_0\subset \rN$ be a root of a finitely generated unit $\rA[\rF]$-module. Then $\rN_0$ is minimal if and only if $\widehat{\rA}\otimes_\rA \rN_0 \subset \widehat{\rA}\otimes_\rA \rN$ is minimal.
\end{lemma}

\begin{proof}
The forward implication is \cite[Corollary 2.11]{BlickleIH}. Conversely, assume the completed root is minimal and let $\rN_1\subset \rN_0$ be a subroot. Its completion is a subroot, hence $\widehat{\rA}\otimes_\rA \rN_1=\widehat{\rA}\otimes_\rA \rN_0$.

Since $\rA\to\widehat{\rA}$ is faithfully flat, it follows that $\rN_1=\rN_0$. Thus $\rN_0$ has no proper subroot, hence is minimal.
\end{proof}

\begin{lemma} \label{lem:minimal_root_etale_pullback}
Let $\Spa \B \to \Spa \rA$ be an \'{e}tale morphism of regular affinoids over $k\laur{q-1}$. Then pullback along this map sends minimal roots of finitely generated unit modules to minimal roots.
\end{lemma}

\begin{proof}
We first note that $\rA$ and $\B$ are both regular $\rF$-finite, so all of the hypotheses of \cite{BlickleGamma} are satisfied. For testing minimality, the results in \cite{BlickleGamma} give us several useful tools. Let us suppose that $\M$ is a finitely generated unit module over say $\rA$, with root $\beta: \M_0 \to \varphi^* \M_0$.

Following \cite[Lemma 2.16]{BlickleGamma}, Blickle constructs $\sigma^{-1}_{\beta}(\M_0)\subset\M_0$ which can be characterized as the smallest submodule $\M_0'$ such that $\beta(\M_0)\subset \varphi_{\rA}^* \M_0'$. By \cite[Lemma 2.17]{BlickleGamma} this construction commutes with localization, and by \cite[Proposition 2.18]{BlickleGamma} a root is minimal in our sense precisely when $\sigma^{-1}_{\beta}(\M_0)=\M_0$.

This construction is given by an explicit formula, which we now explain. We let $\D_{\rA}^{(1)} = \End_{\rA}(\rA^{(1)})$ be the algebra of right $\rA$-linear endomorphisms of the Frobenius twist $\rA^{(1)}$. Blickle gives a Frobenius descent functor
\[\sigma^{-1}: \Mod_{\D_{\rA}^{(1)}} \longrightarrow \Mod_{\rA}\]
as the inverse to $\rA^{(1)} \otimes_{\rA} -$ (which is an equivalence since $\rA^{(1)}$ is finite locally free over $\rA$). We may compute $\sigma^{-1}_\beta(\M_0)$ as $\sigma^{-1}(\D_{\rA}^{(1)}\cdot \beta(\M_0)) \subseteq \sigma^{-1} \sigma^\ast \M_0 = \M_0$.

Now we prove the claim. By localization compatibility, a root $\M_0$ over $\rA$ is minimal if and only if its localization at every maximal ideal is again minimal. We use here that $\sigma_{\beta}^{-1}(\M_0)=\M_0$ is an open condition, see \cite[Proposition 2.18]{BlickleGamma}. The same then also holds for $\B$. Combining this with Lemma \ref{lem:completion_minimal}, we see minimality may be checked on the \emph{completed} local rings at the classical points of the rigid spaces.

It is straightforward to reduce to the finite \'{e}tale case using this observation. The morphism $\Spa \B \to \Spa \rA$ factors near a classical point $x$ into an affinoid open immersion and a finite \'{e}tale morphism. However since affinoid open immersions induce isomorphisms on completed local rings by \cite[7.3.2/3]{BoschGuntzerRemmert}, after localizing at $x$ we have now reduced to showing the claim for the map
\[\widehat{\rA}_{\mathfrak{n}} \longrightarrow \widehat{\B}_{\m}\]
which is now finite \'{e}tale with source and target regular $\rF$-finite. Here the maximal ideal $\m$ corresponds to the rigid point $x$, and $\mathfrak{n} = \m \cap \rA$.

Thus, let us assume the ring map $\rA \to \B$ is finite \'{e}tale. We would be done if we showed the identity
\[\sigma^{-1}_{\beta_{\B}} (\B \otimes_{\rA} \M_0) = \B \otimes_{\rA} \sigma^{-1}_{\beta}(\M_0),\]
as using minimality the right side simplifies to $\B \otimes_{\rA} \M_0$, meaning $\B \otimes_{\rA} \M_0$ is minimal. Using the explicit formula, 
\[\sigma^{-1}_{\beta_{\B}} (\B \otimes_{\rA} \M_0) = \sigma^{-1}(\D_{\B}^{(1)} \beta_{\B} (\B \otimes_{\rA} \M_0)).\]
Now substitute $\B\otimes_{\rA}\D_{\rA}^{(1)}\simeq\D_{\B}^{(1)}$. This follows here since the relative Frobenius is an isomorphism, so $\B \otimes_{\rA} \rA^{(1)} \simeq \B^{(1)}$, and taking $\End$ gives the claim.

Applying this substitution yields 
\[\sigma^{-1}_{\beta_{\B}} (\B \otimes_{\rA} \M_0) = \sigma^{-1}(\B \otimes_{\rA} \D_{\rA}^{(1)} \cdot \beta_{\B} (\B \otimes_{\rA} \M_0)) = \sigma^{-1}(\D_{\rA}^{(1)} \cdot \beta(\M_0)) \otimes_{\rA} \B.\]  
In the last step we have used that $\sigma^{-1}$ commutes with \'{e}tale base change, which again can be seen from $\B\otimes_{\rA}\D_{\rA}^{(1)}\simeq\D_{\B}^{(1)}$. Recalling the definition of $\sigma^{-1}_{\beta}(\M_0)$ we have proved $\sigma^{-1}_{\beta_{\B}} (\B \otimes_{\rA} \M_0) = \B \otimes_{\rA} \sigma^{-1}_{\beta}(\M_0)$ as desired.
\end{proof}

\begin{corollary} \label{cor:minimal_root_framed_pullback}
Let $f:(\Spf \rS,\mathcal{M}_{\Spf \rS})\to(\Spf \R,\mathcal{M}_{\Spf \R})$ be a strict \'{e}tale map where the target is a framed log affine. If $\mathcal{E}\in \D_{\fgu}((\Spf \R,\mathcal{M}_{\Spf \R})_\Prism,\Laur)_{p=0}$ with distinguished minimal root $\mathcal{E}_{0}$, then $f^*\mathcal{E}$ is finitely generated unit and $f^*\mathcal{E}_{0}$ is its distinguished minimal root.
\end{corollary}

\begin{proof}
If $(\Spf \R,\mathcal{M}_{\Spf \R})$ has framing $\sq$, give $(\Spf \rS,\mathcal{M}_{\Spf \rS})$ the framing $\sq \circ f$. The claim follows by combining Lemma \ref{lem:minimal_root_etale_pullback} and Proposition \ref{prop:root_crys_minimal}, using the description in Theorem \ref{thm:qdR_crys_Gamma_mod} to apply the first lemma. 
\end{proof}

Now the proof of strict \'{e}tale descent for $\D_{\fgu}((\X,\mathcal{M})_\Prism,\Laur)$ follows.

\begin{theorem} \label{thm:fgu_etale_desc}
The assignment
\[(\X,\mathcal{M}) \mapsto \D_{\fgu}((\X,\mathcal{M})_\Prism, \Laur)\]
is a stack for the strict \'{e}tale topology on $(\X,\mathcal{M})$ satisfying Assumption \ref{ass:horizontal_nc} over $\W(k)[\zeta_p]$.
\end{theorem}

\begin{proof}
By Corollary \ref{cor:laur_strict_etale}, $\D_{\Box}((\X,\mathcal{M})_\Prism,\Laur)^{\varphi=1}$ satisfies strict \'{e}tale descent. Let $\{\rho_i:\U_i\to(\X,\mathcal{M})\}_{i\in\I}$ be a strict \'{e}tale cover such that every $\rho_i^*\mathcal{E}$ is finitely generated unit.

Using Remark \ref{rem:fgu_mod_p} and Proposition \ref{prop:root_crys_minimal} we may reduce the claim to showing
\[\left\{\rho_i^*\mathcal{E} \in \Mod_{\fgu}((\U_i)_\Prism, \Laur)_{p=0}\right\}_{i\in \I} \implies \left\{ \mathcal{E}\in \Mod_{\fgu}((\X,\mathcal{M})_\Prism, \Laur)_{p=0}\right\}\]
for $\mathcal{E} \in \QCoh_{\Box}((\X,\mathcal{M})_\Prism,\Laur)^{\varphi=1}$ killed by $p$, or that the condition can be detected locally on a strict \'{e}tale cover. Note that the converse is obvious, because any test log prism for the $\U_i$ is also one for $(\X,\mathcal{M})$, so the finitely generated unit condition on the right side already includes the conditions imposed on the left hand side.

By Lemma \ref{lem:strict_etale_local_framing}, after a strict \'{e}tale refinement we may assume without loss of generality that every $\U_i = \Spf\R_i$ is a framed log affine. Put $\U_{ij} = \U_i\times_{(\X,\mathcal{M})}\U_j$ with projections $p_i,p_j$; each affine open of $\U_{ij}$ acquires a framing, and hence an associated logarithmic $q$-de Rham prism, through each of $p_i$ and $p_j$. Proposition \ref{prop:root_crys_minimal} gives distinguished minimal roots $\mathcal{E}_{0,i} \subset \rho_i^*\mathcal{E}$. Corollary \ref{cor:minimal_root_framed_pullback}, applied with each framing, makes $p_i^*\mathcal{E}_{0,i}$ and $p_j^*\mathcal{E}_{0,j}$ distinguished minimal roots on any affine open of $\U_{ij}$. We may further cover $\U_{ij}$ by framed log affines.

Their uniqueness in Proposition \ref{prop:root_crys_minimal}, or what we called `distinguished', gives an identification $p_i^*\mathcal{E}_{0,i} = p_j^*\mathcal{E}_{0,j}$ as subcrystals of $\mathcal{E}|_{\U_{ij}}$. Indeed by the remark after Corollary \ref{cor:laur_t_struct} one has descent for coherent Laurent crystals $\Coh^{\varphi}((\X,\mathcal{M})_\Prism,\Laur)$ with a lax Frobenius. Thus we can check agreement locally on the chosen open cover of $\U_{ij}$ by framed log affines, which also shows agreement of the root maps. These equalities also satisfy the cocycle condition on triple overlaps, which is all we need to check as we work with coherent Laurent crystals.

Strict \'{e}tale descent for $\Coh^{\varphi}((\X,\mathcal{M})_\Prism,\Laur)_{p=0}$ and $\D^{\varphi}_{\Box}((\X,\mathcal{M})_\Prism,\Laur)^\heart_{p=0}$ produces a global coherent root $\mathcal{E}_0$. We can check this global root is in fact a distinguished minimal root by restricting any subroot or any other root to the cover. It is also clear the unitalization remains $\mathcal{E}$ (as the claim is local and pullbacks commute with colimits). By the argument in Footnote \ref{fn:root_unitalization_fgu} we then deduce $\mathcal{E}$ is finitely generated unit from the existence of a root and the claim follows.
\end{proof}

\begin{remark} \label{rem:qdR_test_sufficient}
Combined with Theorem \ref{thm:qdR_crys_Gamma_mod}, this shows that the finitely generated unit condition can equivalently be tested upon evaluation by a covering family of log $q$-de Rham prisms as claimed in the introduction.
\end{remark}

Let us record an important fact we showed in the proof of the theorem.

\begin{corollary} \label{cor:global_root_fgu_condition}
Let $(\X,\mathcal{M})$ satisfy Assumption \ref{ass:horizontal_nc} over $\W(k)[\zeta_p]$. A crystal 
\[\mathcal{E}\in \D_{\Box}((\X,\mathcal{M})_\Prism,\Laur)^{\varphi=1}\]
is finitely generated unit if and only if it is $t$-bounded and, for every $\ell$, the unit Laurent crystal $\H^\ell(\mathcal{E}\dquot p)$ admits a unique global minimal coherent root.
\end{corollary}

\begin{proof}
The forward implication is the construction in the proof of Theorem \ref{thm:fgu_etale_desc}, applied to the cohomology sheaves of $\mathcal{E}\dquot p$ using Proposition \ref{prop:root_crys_minimal} locally. The converse follows from the first paragraph of the proof of Proposition \ref{prop:root_crys_minimal}, which explains why it suffices, on a framed log affine, to know each $\H^\ell(\mathcal{E}\dquot p)$ admits a root. Then as we know the property of being finitely generated unit is strict \'{e}tale local by Theorem \ref{thm:fgu_etale_desc} the general claim follows.
\end{proof}

\begin{remark}
This is an analogue of the root criterion for finitely generated unit modules on a smooth quasicompact $k$-scheme $\X$. By \cite[Theorem 6.1.3]{EmertonKisin2}, a finitely generated unit $\O_{\rF,\X}$-module $\M$ admits a global coherent root $\beta: \M_0 \to \varphi^* \M_0$. In this setting it is well-known that the unit Frobenius on $\M$ gives it a canonical connection with nilpotent $p$-curvature, upgrading it to a crystal on $\X_{\cris}$.

We let $\varphi$ denote the absolute Frobenius. Now $\M_0$ need not have the structure of a crystal, but $\varphi^* \M_0$ does have the structure of a crystal by Cartier descent. Recall that Cartier descent (see \cite[Theorem 5.1]{KatzNilpMon}) tells us that pullback along the relative Frobenius $\mathrm{Fr}: \X \to \X^{(1)}$ gives a fully faithful functor
\[\mathrm{Fr}^*: \QCoh(\X^{(1)}) \longrightarrow \QCoh((\X/k)_\cris).\]
In particular we send $\mathscr{F}$ to $(\mathrm{Fr}^*\mathscr{F},\nabla_{\mathrm{can}})$. The observation is then that we have a factorization
\begin{display}
\varphi: \X \ar{r}{\mathrm{Fr}} & \X^{(1)} \ar{r}{\pi} & \X
\end{display}
of the absolute Frobenius. In particular, when we compute $\varphi^* \M_0$ the underlying crystal has a canonical connection $(\mathrm{Fr}^*(\pi^* \M_0),\nabla_{\mathrm{can}})$. Moreover since we applied a Frobenius pullback, and $\X$ is smooth (making $\varphi$ flat) we see $\varphi^* \beta: \varphi^* \M_0 \to \varphi^{*2} \M_0$ remains a root of $\M$. By construction, this root map is also a morphism of crystals where we use the canonical Cartier connection.
\end{remark}

\begin{corollary} \label{cor:pullback}
Let $f:(\Y,\mathcal{N})\to(\X,\mathcal{M})$ be a morphism between pairs satisfying Assumption \ref{ass:horizontal_nc} over $\W(k)[\zeta_p]$. Then pullback preserves finitely generated unit crystals.
\end{corollary}

\begin{proof}
It suffices to check the claim modulo $p$ for a finitely generated unit Laurent crystal $\mathcal{E}$ in the heart, by applying derived Nakayama and a Postnikov d\'{e}vissage. Upon writing $\mathcal{E}\simeq \colim_n \varphi^{\ast n} \mathcal{E}_0$ using the root from Corollary \ref{cor:global_root_fgu_condition} and noting the pullback preserves colimits, the pullback is finitely generated unit.
\end{proof}

\subsection{The Kummer-\'{e}tale realization and constructibility}

In this section we construct a functor
\[\T_\et: \D_{\fgu}((\X,\mathcal{M})_\Prism,\Laur)^\op \longrightarrow \D((\X,\mathcal{M})_{\eta,\proket}, \Z_p),\]
and show that its essential image is constrained to certain $p$-Kummer constructible sheaves that we will define shortly.

To make the definition of $\T_\et$, it is convenient to work on a basis of log perfectoids for the pro-Kummer-\'{e}tale site.

\begin{lemma}
\label{lem:perfectoid_basis}
Let $(\X,\mathcal{M})$ be a locally Noetherian fs log $p$-adic formal scheme. Then the log affinoid perfectoid objects form a basis $\B$ of $(\X,\mathcal{M})_{\eta,\proket}$ stable under fiber products.
\end{lemma}

\begin{proof}
This is \cite[Propositions 5.3.11 and 5.3.12]{DiaoLanLiuZhuLogAdic}.
\end{proof}

For an object $\U=\lim_i\U_i$ in the basis from the previous lemma, we write $\widehat{\U}=\Spa(\rS,\rS^+)$ for its associated perfectoid space, equipped with the natural saturated log structure as in \cite[Remark 5.3.5]{DiaoLanLiuZhuLogAdic}.

\begin{definition} \label{def:Tet}
Let $(\X,\mathcal{M})$ be a locally Noetherian fs log $p$-adic formal scheme, and let $\mathcal{E}$ be a finitely generated unit Laurent crystal. Let $\U\in \B$ be as in Lemma \ref{lem:perfectoid_basis}, with associated log perfectoid space
\[\widehat{\U} = (\Spa(\rS,\rS^+),\mathcal{M}_{\rS}).\]
We let $\Prism_{\rS^+}$ be the log prism with underlying prism $(\W(\rS^{+,\flat}),\ker \theta)$ and log structure induced by the tilted log structure on $\rS^\flat$. We define
\[\T_\et(\mathcal{E})(\U) := \R\Hom_{\D^{\varphi}_{\Box}(\W(\rS^\flat))}(\mathcal{E}(\Prism_{\rS^+}), \W(\rS^\flat)) \in \widehat{\D}(\Z_p).\]
These assignments define a presheaf on $\B$, which by ``unfolding'' (or right Kan extension) yields a presheaf on $(\X,\mathcal{M})_{\eta,\proket}$ we denote by $\T_\et(\mathcal{E})$.
\end{definition}

\begin{remark}
Our convention so far has been to use the strict site, and unfortunately these perfect prisms do not live in the strict site. This is not an issue because Remark \ref{rem:crystal_indep} allows us to use the saturated log prismatic site and obtain the same category of crystals.
\end{remark}

Here 
\[\D_{\Box}(\Ainf(\rS)) = \D_{\Box}(\W(\rS^\flat)) := \lim_n \D_{\Box}(\Spa(\Ainf(\rS)/p^n,\Ainf(\rS^+)/p^n)),\] 
with $\Spa(\Ainf(\rS)/p^n,\Ainf(\rS^+)/p^n)$ given its natural analytic ring structure from the Huber pair as in Definition \ref{def:laur_category}. One can see $\Ainf(\rS) = \Ainf(\rS^+)[1/\xi]_p^{\wedge}$ where $\ker \theta = (\xi)$ as both are strict $p$-rings.

Before proceeding further, we show that the output is actually a sheaf. 

\begin{lemma} \label{lem:T_et_proket_sheaf}
Let $(\X,\mathcal{M})$ be a locally Noetherian fs log $p$-adic formal scheme. Definition \ref{def:Tet} yields a functor
\[\T_\et: \D_{\fgu}((\X,\mathcal{M})_\Prism,\Laur)^\op \longrightarrow \D((\X,\mathcal{M})_{\eta,\proket},\Z_p).\]
\end{lemma}

\begin{proof}
We will actually check hyperdescent, which is made easier as hypersheaves on the basis $\B$ from Lemma \ref{lem:perfectoid_basis} are the same as hypersheaves on $(\X,\mathcal{M})_{\eta,\proket}$.

For $\U \in \B$ let $\U_\bullet\to\U$ be a hypercover by basis objects in $\B$, with associated perfectoid spaces $\widehat{\U}=\Spa(\rS,\rS^+)$ and $\widehat{\U}_m=\Spa(\rS_m,\rS_m^+)$. On the associated perfectoid spaces one gets a $v$-hypercover
\[\widehat{\U}_\bullet \longrightarrow \widehat{\U}.\]
Here we use that the assignment $\U \mapsto \widehat{\U}$ commutes with fiber products, and the easier fact that pro-Kummer-\'{e}tale covers induce $v$-covers.

Put $\M=\mathcal{E}(\Prism_{\rS^+})$ and $\M_m=\mathcal{E}(\Prism_{\rS_m^+})$. The desired descent claim is that 
\[\lim_{[m]\in\Delta} \RHom_{\D_{\Box}^\varphi(\Ainf(\rS_m))}(\M_m,\Ainf(\rS_m)) \simeq \RHom_{\D_{\Box}^\varphi(\Ainf(\rS))}(\M,\Ainf(\rS)).\] 
We can use the extension and restriction adjunction to make all $\RHom$s in $\D_{\Box}^\varphi(\Ainf(\rS))$, simplifying the left side to 
\[\lim_{[m]\in\Delta} \RHom_{\D_{\Box}^\varphi(\Ainf(\rS))}(\M,\Ainf(\rS_m)).\]
Thus the needed claim is actually that $\Ainf(\rS) \simeq\R\lim_{[m]\in\Delta}\Ainf(\rS_m)$. The claim can be checked modulo $p$, where it reduces to the $v$-acyclicity of the structure sheaves on $\Spa(\rS_m^\flat, \rS_m^{+,\flat})$ shown in \cite[Theorem 1.2]{etCohDiamonds}. Here we must also keep track of Frobenius and solid structures, but this much is clear.
\end{proof}

It is easy to see that for every morphism $f:(\Y,\mathcal{N})\to(\X,\mathcal{M})$ between pairs satisfying Assumption \ref{ass:horizontal_nc}, the construction is compatible with pullback:
\[\T_{\et,\Y}(\L f^*\mathcal{E})\simeq\L f_\eta^*\T_{\et,\X}(\mathcal{E}).\]

As previously claimed, the functor $\T_\et$ actually lands in the full subcategory of $p$-Kummer sheaves which we now define. Recall that for a log geometric point $\overline{x} \in (\X,\mathcal{M})_\eta$ we have an exact sequence
\begin{display}
0 \ar{r} & \I_{\overline{x}}^{\log} \ar{r} & \pi_1^{\proket}(\overline{x}_{\log}) \ar{r} & \pi_1^{\proet}(\overline{x}) \ar{r} & 0
\end{display}
and one can compute $\I_{\overline{x}}^{\log} = \Hom(\overline{\mathcal{M}}^{\mathrm{gp}}_{\overline{x}}, \widehat{\Z}(1))$. The prime-to-$p$ part $\I_{\overline{x}}^{\log, (p)'}$ is given by replacing $\widehat{\Z}(1)$ with its prime-to-$p$ part. If we set 
\[\rho(\overline{x})=\mathrm{rk}_{\Z} \overline{\mathcal{M}}^{\mathrm{gp}}_{\overline{x}}\] 
for geometric points $\overline{x}$ of the adic generic fiber $\X_\eta$, the logarithmic strata in $(\X, \mathcal{M})_\eta$ are the (connected components of) loci where $\rho=a$. We will write $\D^{(b)}_{\cons,\D}((\X,\mathcal{M}_\D)_{\eta},\Z_p)$ for Kummer-\'{e}tale sheaves which are derived lisse on these strata. The following definition extracts the full subcategory of sheaves which are $p$-Kummer.

\begin{definition} \label{def:p-Kum}
Let $(\X,\mathcal{M})$ be as in Assumption \ref{ass:horizontal_nc}. We call $\mathscr{F}\in\D^{(b)}_{\cons,\D}((\X,\mathcal{M})_\eta,\Z_p)$ \emph{$p$-Kummer} if prime-to-$p$ logarithmic inertia $\I_{\overline{x}}^{\log, (p)'}$ acts trivially on every cohomology stalk of $\mathscr{F}\dquot p^n$ at every log geometric point, for every $n\ge1$. 
\end{definition}

For a derived lisse sheaf, this agrees with \cite[Definition 2, and Proposition 5]{KoshikawaYaoCorrigendum} (which is where the name ``$p$-Kummer'' comes from). The remainder of this section will prove that $\T_\et(\mathcal{E})$ is constructible for the stratification induced by $\D$, and as a byproduct that the result lands in the subcategory of $p$-Kummer sheaves.

We first need the following descent result. When dealing with these results on $p$-Kummer sheaves we allow a bit more generality, namely we can also allow $K/\W(k)[1/p]$ ramified.

\begin{lemma} \label{lem:constructible_strict_etale_descent}
On pairs satisfying Assumption \ref{ass:horizontal_nc} over $\O_K$, the full subcategories of $p$-Kummer sheaves in
\[\D^{(b)}_{\cons,\D}((\X,\mathcal{M})_\eta,\Z_p)\]
form a stack for the strict \'{e}tale topology.
\end{lemma}

We will denote these by $\D^{(b),\mathrm{pKum}}_{\cons,\D}((\X,\mathcal{M})_\eta,\Z_p)$ later.

\begin{proof}
We first check the claim without the $p$-Kummer condition imposed. This can be checked fairly directly. Let $f:(\U,\mathcal{N})\to(\X,\mathcal{M})$ be a strict \'{e}tale cover on the log formal schemes, and let $\U^{(\bullet)}$ be its \v{C}ech nerve. 

Passing to adic generic fibers, this is compatible with fiber products so we can simply take log adic generic fibers of $\U^{(\bullet)}$ to obtain the corresponding \v{C}ech nerve on the log adic generic fiber. On the larger category of pro-Kummer-\'{e}tale sheaves, trivially we have
\[\D^{(b)}((\X,\mathcal{M})_{\eta,\proket},\Z_p) \simeq \lim_{[n]\in\Delta} \D^{(b)}((\U^{(n)},\mathcal{N}^{(n)})_{\eta,\proket},\Z_p).\]
We also have such descent for $\D^{(b)}_{\mathrm{lisse}}$. By strictness, if $\X_\eta^{(a)}$ and $\U_\eta^{(a)}$ denote the loci where $\mathrm{rk}_{\Z} \overline{\mathcal{M}}_{\bar{x}}^{\mathrm{gp}} = a$ (respectively for $\mathcal{N}$), then
\[\U_\eta^{(a)}=\U_\eta\times_{\X_\eta}\X_\eta^{(a)}.\]
The logarithmic strata in the \v{C}ech nerve are therefore compatible, and we also have descent for lisse sheaves on each stratum (the strata are obtained by base change, so the corresponding \v{C}ech nerve again consists of strict \'{e}tale maps). As we are interested in constructible sheaves for this stratification the claim follows.

For the categories $\D^{(b),\mathrm{pKum}}_{\cons,\D}((\X,\mathcal{M})_\eta,\Z_p)$ with the $p$-Kummer condition imposed, note that strictness identifies the logarithmic inertia groups so the claim follows.
\end{proof}

This has a useful implication: by Lemmas \ref{lem:strict_etale_local_framing} and \ref{lem:constructible_strict_etale_descent}, constructibility for the logarithmic stratification induced by $\D$ may be checked on a framed log affine (as well as the $p$-Kummer condition). Thus for the desired result that $\T_\et$ lands in the full subcategory $\D^{(b),\mathrm{pKum}}_{\cons,\D}((\X,\mathcal{M})_\eta,\Z_p)$ we can now immediately assume that $(\X,\mathcal{M})$, rather than simply satisfying Assumption \ref{ass:horizontal_nc}, is a framed log affine.

In this case there is a natural perfectoid cover we will use. Let $(\X,\mathcal{M})$ be a framed log affine as in Definition \ref{def:fr_log_aff}. Then the log formal scheme $\mathbb{T}_r^n$ induces a natural log perfectoid cover
\[(\Spf \W(k)\langle \zeta_{p^\infty}\rangle\langle \T_1^{1/p^\infty}, \ldots, \T_r^{1/p^\infty}, \T_{r+1}^{\pm 1/p^\infty}, \ldots,  \T_{n}^{\pm 1/p^\infty} \rangle, \N[1/p]^r)_\eta \longrightarrow (\mathbb{T}_r^n)_\eta.\]
Pulling this back along $\sq_\eta$, we get a log perfectoid cover
\[(\X,\mathcal{M})_{\infty,\eta} \longrightarrow (\X,\mathcal{M})_\eta\]
which is a quasi-pro-Kummer-\'{e}tale torsor for $\Gamma$. When the framed log affine is over $\W(k)[\zeta_p]$ as we usually assume, $\Gamma = (1+p\Z_p) \ltimes \Z_p(1)^n$ as usual. We change the cyclotomic factor $1+p\Z_p$ depending on the base, for example it is $\Z_p^\times$ when working over $\W(k)$.

Using this perfectoid cover, we can give another useful criterion to check the $p$-Kummer condition.

\begin{lemma} \label{lem:pKum_framed_criterion}
Let $(\X,\mathcal{M})$ be a framed log affine with $\X=\Spf \R$ and let $\mathscr{F}$ be derived $p$-complete in $\D^{(b)}((\Spf \R,\mathcal{M})_{\eta,\proket},\Z_p)$. For $\Lambda\subseteq\{1,\ldots,r\}$, put
\[(\X,\mathcal{M})_{\infty,\eta,\Lambda}=(\X,\mathcal{M})_{\infty,\eta}\times_{(\X,\mathcal{M})_\eta}(\X,\mathcal{M})_{\eta,\Lambda},\]
and let $\epsilon_\infty:((\X,\mathcal{M})_{\infty,\eta},\mathcal{M}_\infty)_{\proket}\to((\X,\mathcal{M})_{\infty,\eta})_{\proet}$ be the natural morphism of sites. The following are equivalent:
\begin{enumerate}
\item $\mathscr{F}$ is constructible for the stratification induced by $\D$ and is $p$-Kummer;
\item Its pullback to $(\X,\mathcal{M})_{\infty,\eta}$ is $\epsilon_\infty^*\mathcal{L}_\infty$ for some $\mathcal{L}_\infty\in\D^{(b)}(((\X,\mathcal{M})_{\infty,\eta})_{\proet},\Z_p)$ which is constructible for the finite stratification $(\X,\mathcal{M})_{\infty,\eta}=\bigsqcup_\Lambda (\X,\mathcal{M})_{\infty,\eta,\Lambda}$.
\end{enumerate}
\end{lemma}

\begin{proof}
Constructibility is clear, because the pro-Kummer-\'{e}tale Galois cover $(\X,\mathcal{M})_{\infty,\eta}\to (\X,\mathcal{M})_\eta$ has Galois group $\Gamma$, whose action respects the strata $(\X,\mathcal{M})_{\infty,\eta}=\bigsqcup_\Lambda(\X,\mathcal{M})_{\infty,\eta,\Lambda}$. In particular we may test lissity on the strata induced by $\D$ after pullback to $(\X,\mathcal{M})_{\infty,\eta,\Lambda}$, and then constructibility along $\D$ after pullback to $(\X,\mathcal{M})_{\infty,\eta}=\bigsqcup_\Lambda(\X,\mathcal{M})_{\infty,\eta,\Lambda}$ by testing lissity on each stratum.

For the $p$-Kummer condition, we check stalkwise. The point is that the logarithmic inertia on the perfectoid cover $(\X,\mathcal{M})_{\infty,\eta}$ extracts precisely the prime-to-$p$ part on the base $(\X,\mathcal{M})_{\eta}$: for $\overline{y}\in (\X,\mathcal{M})_{\infty,\eta}$ a log geometric point above $\overline{x}\in (\X,\mathcal{M})_{\eta}$,
\[\I_{\infty,\overline{y}}^{\log} \simeq \Hom(\Z[1/p]^{|\Lambda|}, \widehat{\Z}(1)) = \I_{\overline{x}}^{\log, (p)'}\]
since $\I_{\overline{x}}^{\log, (p)'}$ is the prime-to-$p$ part of $\Hom(\overline{\mathcal{M}}_{\overline{x}}^{\mathrm{gp}},\widehat{\Z}(1))$. The condition we ask, that the pullback is of the form $\epsilon_\infty^*\mathcal{L}_\infty$, is equivalent to saying this acts trivially.
\end{proof}

We let $\T_\et^{\perf}$ denote the \emph{schematic} solution functor introduced in Theorem \ref{thm:EK_perfect}. We let $(\R_0^{\sq})_{\infty,\eta}$ denote the associated perfectoid ring for $(\X,\mathcal{M})_{\infty,\eta}$ over $\Q_p$, and in what follows we will apply $\T_\et^{\mathrm{perf}}$ to $\Spec (\R_0^{\sq})^\flat_{\infty,\eta}$. For a finitely generated unit module $\M$ over $\R_0^{\sq}\laur{q-1}^{\wedge}_p$, we use the convention
\[\T_{\et}^{\perf}(\M) := \T_{\et}^{\perf}\left(\M\otimes_{\R_0^{\sq}\laur{q-1}^{\wedge}_p}^{\L}\W((\R_0^{\sq})_{\infty,\eta}^{\flat})\right),\]
as it does not literally make sense as an input before base changing.

The basic structure of the argument for constructibility will be to produce a factorization of $\T_\et$ in Lemma \ref{lem:factorization} on a framed log affine involving $\T_\et^{\perf}$. With the input of $\T_\et$ regarded as an object in $\D_{\fgu}(\R_0^{\sq}\laur{q-1}_p^{\wedge})^{h\Gamma}$ by Theorem \ref{thm:qdR_crys_Gamma_mod}, we can ultimately reduce to objects $\M$ in the heart killed by $p$ and argue that they have $\Gamma$-equivariant minimal roots $\M_0$ as in the proof of Theorem \ref{thm:qdR_crys_Gamma_mod}. Then looking at the root map
\[\beta: \M_0 \hookrightarrow \varphi^* \M_0\]
one sees $\M_0, \varphi^*\M_0$, and the cokernel of $\beta$ all are finitely generated $\Gamma$-equivariant modules. We constrain these by studying their Fitting ideals (which are then $\Gamma$-stable ideals), and classifying the possibilities enough to force the root map to be an equivalence on each stratum. We do this by forcing $\M_0, \varphi^* \M_0$ and $\coker \beta$ to be vector bundles through the usual Fitting criterion in \cite[Tag \href{https://stacks.math.columbia.edu/tag/07ZD}{07ZD}]{StacksProject}. To conclude, one then deduces constructibility using the factorization in Lemma \ref{lem:factorization}, using Katz's theorem to show the result on each stratum is lisse at the perfectoid level (whereupon the criterion in Lemma \ref{lem:pKum_framed_criterion} is useful).

We begin with a classification of stable ideals.

\begin{lemma} \label{lem:toric_stable_ideal}
Let $(\Spf \R,\mathcal{M})$ be a framed log affine, and put $\overline{\R}=\R_0/p$. If $\mathfrak{a}\subset \overline{\R}^{\sq}\laur{q-1}$ is an ideal stable under $\Z_p(1)^n\subset\Gamma$, then, on each connected component of
\[\Spa \overline{\R}^{\sq}\laur{q-1}\setminus \V(\T_1\cdots \T_r)\]
the image of $\mathfrak{a}[1/\T_1,\ldots,1/\T_r]$ is either zero or the unit ideal.

Moreover, the same conclusion holds for the scheme $\Spec \overline{\R}^{\sq}\laur{q-1}[1/\T_1,\ldots,1/\T_r]$.
\end{lemma}

\begin{proof}
Put $\X=\Spa \overline{\R}^{\sq}\laur{q-1}$ and $\X^\circ=\X\setminus \V(\T_1\cdots \T_r)$. Due to the framing $\sq$, as rigid analytic spaces over $k\laur{q-1}$ we have an \'{e}tale map 
\[\X\longrightarrow \Spa k\laur{q-1} \langle \T_1, \ldots, \T_r, \T_{r+1}^{\pm 1}, \ldots, \T_n^{\pm 1} \rangle.\]
The target is reduced, hence so is $\X$. Through this \'{e}tale map, we also get an action of $\Z_p(1)^n$ on $\X$ which preserves $\X^{\circ}$. On the base, elements $x\in \Z_p(1)_i$ of the $i$th copy of $\Z_p(1)$ act by $\T_i\mapsto q^{px}\T_i$ and fix the other coordinates. We first show the claim for the adic space using a similar idea to \cite[Lemma 3.9]{syncrisloc}.

Put $\mathcal{I}=\mathfrak{a}\O_{\X^\circ}$, and let $x$ be a classical point for which $\mathcal{I}_x$ is a proper ideal. Since $\O_{\X^\circ,x}$ is a local ring, this means $\mathcal{I}_x$ is contained in $\m_x \O_{\X^\circ,x}$ and so all $f\in \mathfrak{a}$ vanish at $x$. Up to replacing $k\laur{q-1}$ by the residue field $\kappa(x)/k\laur{q-1}$, the rigid analytic inverse function theorem will identify some neighborhood of $x$ with a polydisk. Whether or not $\mathcal{I}_x$ vanishes can be checked after this extension, so without loss of generality we can assume $x$ is a rational point with a polydisk neighborhood.

Now we use $\Z_p(1)^n$ stability of $\mathfrak{a}$. Considering the elements $\gamma_\alpha = (p^{\alpha_1}, \ldots, p^{\alpha_n})\in \Z_p(1)^n$, stability means for $f\in \mathfrak{a}$ that $\gamma_\alpha \cdot f \in \mathfrak{a}$. They all vanish at $x$, so we will also get vanishing at $q^{p^{\alpha_i+1}}\T_i(x)$. For $\min_i \alpha_i \gg 0$ these live in the polydisk around $x$, by continuity of the action of $\Z_p(1)^n$. Since we work in positive characteristic we therefore get
\[q^{p^{\alpha+1}}-1=(q-1)^{p^{\alpha+1}} \longrightarrow 0,\]
so we have produced an accumulating grid of zeroes in a polydisk around $x$. Applying the single variable rigid analytic identity theorem repeatedly shows that a finite set of generators of $\mathfrak{a}$ vanishes on a neighborhood of $x$. Thus $\mathcal{I}_x=0$ if $\mathcal{I}_x$ is a proper ideal. In particular either $\mathcal{I}_x=\O_{\X^\circ,x}$ or it is zero, and we have seen the locus where it is zero is open (i.e. it propagates to a polydisk around $x$). Thus on connected components it is either zero or the unit ideal.

This proves the first assertion. To deduce the claim for the scheme, the same claims hold for $\mathfrak{a}_\m$ for maximal ideals $\m$ of $\Spec \overline{\R}^{\sq}\laur{q-1}[1/\T_1,\ldots,1/\T_r]$. Indeed the claim that $\mathcal{I}_x$ is zero or the unit ideal is simply about the base change to the completed local ring, which is faithfully flat over $\overline{\R}^{\sq}\laur{q-1}[1/\T_1,\ldots,1/\T_r]_\m$ in our situation. Thus by faithfully flat descent, we deduce $\mathfrak{a}_\m$ is zero or the unit ideal for maximal ideals in $\Spec \overline{\R}^{\sq}\laur{q-1}[1/\T_1,\ldots,1/\T_r]$. 

This is enough to deduce the desired claim that $\mathfrak{a}$ is $0$ or $1$ on connected components for the scheme. To see this, note that connected components are integral, and so reducing to the affine case we can use the fact that for a domain $\D$ we have
\[\D = \bigcap_{\m \in \operatorname{mSpec} \D} \D_\m.\]
The ideal $\mathfrak{a}$ is zero or one in each $\D_{\m}$, hence on $\D$. The claim then follows.
\end{proof}

\begin{lemma} \label{lem:derived_stratum_fgu}
Let $\overline{\rA}=\overline{\R}^{\sq}\laur{q-1}$ for a framed log affine, and let $\rJ_{\Lambda}=(\T_i)_{i\in \Lambda}$ for some $\Lambda\subset\{1,\ldots,r\}$. If $\M \in \D_{\fgu}(\overline{\rA})^{h\Gamma}$, then
\[(\overline{\rA}/\rJ_{\Lambda})\otimes^\L_{\overline{\rA}}\M\]
is in $\D_{\fgu}(\overline{\rA}/\rJ_{\Lambda})^{h\Gamma}$.
\end{lemma}

\begin{proof}
There is an induced $\Gamma$ action because the ideal $\rJ_{\Lambda}$ is stable under $\Gamma$. The finitely generated unit claim is clear because
\[(\overline{\rA}/\rJ_{\Lambda})[\rF]\otimes^{\L}_{\overline{\rA}[\rF]}\M\simeq(\overline{\rA}/\rJ_{\Lambda})\otimes^{\L}_{\overline{\rA}} \M,\]
showing the result is in $\mathrm{Thick}((\overline{\rA}/\rJ_{\Lambda})[\rF])$. It is unit because $\rJ_{\Lambda}$ is Frobenius stable (as is any ideal in positive characteristic).
\end{proof}

In the following proposition recall that $(\R_0^{\sq})_{\infty,\eta}$ is the ring of functions of the affinoid perfectoid space $(\Spf \R,\mathcal{M})_{\infty,\eta}$.

\begin{proposition} \label{prop:perfectoid_constructibility}
Let $(\Spf \R,\mathcal{M}_{\Spf \R})$ be a framed log affine with framing $\sq$, and let $\mathcal{E}$ be an object of $\D_{\fgu}((\Spf \R,\mathcal{M}_{\Spf \R})_\Prism,\Laur)$. Then $\T_{\et}^{\perf}(\mathcal{E}(\R_0^{\sq}\laur{q-1}^{\wedge}_p))$ is constructible for the stratification of $\Spec (\R_0^{\sq})_{\infty,\eta}^{\flat}$ cut out by the divisor $\D_\infty=\bigcup_i \V(\T_i^\flat)$.
\end{proposition}

\begin{proof}
Because constructibility can be detected after reduction modulo $p$, from the start we can assume $\mathcal{E}\in \D_{\fgu}((\Spf \R,\mathcal{M}_{\Spf \R})_\Prism,\Laur)_{p=0}$. As usual we put $\overline{\rA} = \overline{\R}_0^{\sq}\laur{q-1}$ for the mod $p$ reduction of $\R_0^{\sq}\laur{q-1}^{\wedge}_p$.

Let $\M=\mathcal{E}(\R_0^{\sq}\laur{q-1}^{\wedge}_p)$, which we may regard as an object $\D_{\fgu}(\overline{\rA})^{h\Gamma}$ using Theorem \ref{thm:qdR_crys_Gamma_mod}. We defined $\T_{\et}^{\perf}(\M)$ to be 
\[\T_{\et}^{\perf}\left(\M\otimes_{\R_0^{\sq}\laur{q-1}^{\wedge}_p}^{\L} \W((\R_0^{\sq})_{\infty,\eta}^{\flat})\right).\]
The input is in $\D_{\fgu}(\W((\R_0^{\sq})_{\infty,\eta}^{\flat}))$. For testing the claim, we may reduce to the case where $\M$ is a module (note the base change is flat, by flatness of Frobenius). To simply show that the output of $\T_\et^{\perf}$ in this case is lisse on the open stratum, or complement of $\D_\infty$, it would suffice to show that $\M[1/\T_1,\ldots,1/\T_r]$ is a vector bundle over $\overline{\rA}[1/\T_1,\ldots,1/\T_r]$ -- then the contravariant version of Katz's theorem from Remark \ref{rem:katz_contra} shows the output is lisse.

Let us dispense with this case first. We regard everything in what follows as a module over $\overline{\rA}[1/\T_1,\ldots,1/\T_r]$ (i.e. we replace $\M$ with its base change to this ring). As we work modulo $p$, we can pick a $\Gamma$-equivariant minimal root $\beta: \M_0 \to \varphi^* \M_0$ for $\M$ (just as in the proof of Theorem \ref{thm:qdR_crys_Gamma_mod} using \cite{BlickleGamma}). Now Lemma \ref{lem:toric_stable_ideal} provides some serious constraints: each $\Fitt_i(\M_0)$ is $\Z_p(1)^n$-stable as are Fitting ideals of $\varphi^* \M_0$ and $\coker \beta$. Pick the minimal $i$ so that $\Fitt_i(\M_0)\neq 0$. One may test the desired constructibility claim locally, so we may reduce to the connected case and apply Lemma \ref{lem:toric_stable_ideal} to deduce this is the unit ideal. Finally applying \cite[Tag \href{https://stacks.math.columbia.edu/tag/07ZD}{07ZD}]{StacksProject} makes $\M_0$ locally free of rank $i$. We note $\varphi^*\M_0$ has the same rank and is necessarily a vector bundle.

By the same argument, we also find $\coker \beta$ is a vector bundle of some possibly different rank. But then we have the situation where we have an exact sequence
\begin{display}
0 \ar{r} & \M_0 \ar[hook]{r}{\beta} & \varphi^* \M_0 \ar{r} & \coker \beta \ar{r} & 0
\end{display}
consisting entirely of vector bundles, but the first two have the same rank. Thus $\coker \beta$ has rank zero, and $\beta$ is an isomorphism. We conclude $\M=\M_0$, and the contravariant version of Katz's theorem tells us that $\T_{\et}^{\perf}(\M)$ is lisse away from $\D_\infty$.

One can run exactly the same argument for any other stratum in the stratification induced by $\D_\infty$ to see the result is lisse on this stratum. What we must do is consider $\overline{\R}_0^{\sq}\laur{q-1}/\rJ_{\Lambda}$ as in Lemma \ref{lem:derived_stratum_fgu} and use this lemma to deduce the derived base change of $\M$ is finitely generated unit, and then we can invert any surviving $\T_i$ if desired. The result follows by making the same reduction to the heart and repeating the argument.
\end{proof}

\begin{remark}
For the trivial log structure there is only the open stratum, so the same argument shows we obtain a derived lisse complex. In particular without log structures, the theory degenerates to the ordinary theory of $(\varphi,\Gamma)$-modules.
\end{remark}

Let $\nu:(\X_{\infty,\eta}^\flat)_\et\to\Spec((\R_0^{\sq})_{\infty,\eta}^{\flat})_\et$ be the support map, as used in \cite[Notation 6.2.1]{BhattPadicHodgeNotes}, and let $\lambda:((\X,\mathcal{M})_{\infty,\eta}^\flat)_\proet\to((\X,\mathcal{M})_{\infty,\eta}^\flat)_\et$ be the natural morphism of sites. For a constructible $\Z_p$-complex $\mathscr{F}$, set
\[\mathrm{an}^*\mathscr{F}:=\R\lim_m\lambda^*\nu^*(\mathscr{F}\dquot p^m)\in\D((\X_{\infty,\eta}^\flat)_\proet,\Z_p).\]
For a $\Gamma$-equivariant constructible complex $\mathscr{F}$, set
\[\mathrm{an}_\Gamma^*\mathscr{F}:=\left(\epsilon_\infty^*((\mathrm{an}^*\mathscr{F})^{\#})\right)^{h\Gamma}\in\D((\Spf \R,\mathcal{M})_{\eta,\proket},\Z_p),\]
where $(-)^{\#}$ denotes transport across the tilting equivalence (to characteristic zero) and $h\Gamma$ denotes descent along $(\X,\mathcal{M})_{\infty,\eta}\to(\X,\mathcal{M})_\eta$. Here everything is also $p$-completed.

\begin{lemma} \label{lem:factorization}
Let $(\Spf \R,\mathcal{M})$ be a framed log affine with framing $\sq:(\Spf \R,\mathcal{M})\to\mathbb{T}_r^n$. Evaluation on $\Prism_{(\R_0^{\sq})_\infty}$ gives a commutative diagram
\begin{display}[column sep=small,row sep=large]
\D_{\fgu}((\Spf \R,\mathcal{M})_\Prism,\Laur)^\op \ar{rrr}{\T_\et} \ar[swap]{d}{\mathcal{E}\mapsto\mathcal{E}(\Prism_{(\R_0^{\sq})_\infty})} & & & \D((\Spf \R,\mathcal{M})_{\eta,\proket},\Z_p) \\
\left(\D_{\fgu}(\W((\R_0^{\sq})_{\infty,\eta}^{\flat}))^{h\Gamma}\right)^\op \ar{rrr}{\T_\et^{\perf}} & & & \D^{(b)}_{\cons}(\Spec (\R_0^{\sq})_{\infty,\eta}^{\flat},\Z_p)^{h\Gamma} \ar[swap]{u}{\mathrm{an}_\Gamma^*}
\end{display}
\end{lemma}

\begin{proof}
We showed already that the formula for $\T_\et$ outputs a sheaf in $\D((\Spf \R,\mathcal{M})_{\eta,\proket},\Z_p)$ in Lemma \ref{lem:T_et_proket_sheaf}, so we may compute its output using the \v{C}ech nerve of $(\X,\mathcal{M})_{\infty,\eta}$ which is a pro-Kummer-\'{e}tale $\Gamma$-torsor over $(\Spf \R,\mathcal{M})_\eta$. In particular, once we use this specific cover to compute $\T_\et(\mathcal{E})$ we get a factorization
\[\T_\et: \D_{\fgu}((\Spf \R,\mathcal{M})_\Prism,\Laur)^\op \longrightarrow \D((\X,\mathcal{M})_{\infty,\eta,\proket},\Z_p)^{h\Gamma}.\]
Setting $\M_{\infty}=\mathcal{E}(\Prism_{(\R_0^{\sq})_\infty})$, the goal is therefore a $\Gamma$-equivariant identification
\[\T_\et(\mathcal{E})(\U) \xlongrightarrow{\sim} \epsilon_\infty^*((\mathrm{an}^*\T_\et^{\perf}(\M_\infty))^{\#})(\U)\]
for log perfectoids $\U=(\Spa(\rS,\rS^+),\mathcal{M}_{\rS}) \to (\X,\mathcal{M})_{\infty,\eta}$. The equivariance is automatic once we do this functorially in $\U$, as the \v{C}ech nerve of $(\X,\mathcal{M})_{\infty,\eta}$ over $(\X,\mathcal{M})_\eta$ is a bar construction in $\Gamma$.

To construct the map is straightforward. Let $\mathscr{F}=\T_\et^{\perf}(\M_{\infty})$, and write $\mathscr{F}_{\rS}$ for its pullback to $\Spec \rS^\flat$. The formula for $\T_\et(\mathcal{E})(\U)$ identifies it directly with $\R\Gamma(\Spec \rS^{\flat}_{\proet}, \mathscr{F}_{\rS})$. The right hand side $\epsilon_\infty^*((\mathrm{an}^*\T_\et^{\perf}(\M_\infty))^{\#})(\U)$ identifies with $\R\Gamma(\Spa(\rS^{\flat},\rS^{+,\flat})_{\proet}, \mathrm{an}^* \mathscr{F}_{\rS})$. Thus analytification, by functoriality, induces a comparison morphism.

By derived Nakayama, derived $p$-completeness reduces us to working with mod $p$ coefficients. The comparison map between $\R\Gamma(\Spa(\rS^{\flat},\rS^{+,\flat})_{\proet}, \mathrm{an}^* \mathscr{F}_{\rS})$ and $\R\Gamma(\Spec \rS^{\flat}_{\proet}, \mathscr{F}_{\rS})$ is then an equivalence by \cite[Proposition 5.1]{AnschutzExtendingTorsors}\footnote{We note this comparison can be reproved using $\mathrm{arc}_{\varpi}$ descent results for torsion \'{e}tale cohomology from \cite[Theorem 5.4 and Corollary 6.17]{BhattMathewArc}, where $\varpi$ is the tilted pseudouniformizer.}, applied to $(\rS^\flat,\rS^{+,\flat})$, together with the \'{e}tale versus pro\'{e}tale comparisons in \cite[Corollary 5.1.6]{proet} and \cite[Proposition 8.5(i)]{etCohDiamonds}.
\end{proof}

\begin{remark}
There is an instructive warning explained in \cite[Warning 6.2.8]{BhattPadicHodgeNotes} related to the previous proof, which is that $\nu^*$ is \emph{not} fully faithful on constructible sheaves. The issue with full faithfulness comes from the fact that distinct constructible sheaves can have the same analytification, and in particular there is a nonzero constructible $\mathscr{F}$ with $\nu^* \mathscr{F}=0$ (as can be constructed from the warning). However this particular example also has no cohomology, so it is perfectly compatible with the previous argument.
\end{remark}

\begin{remark}
\label{rem:Tet_sym_mon}
The factorization in Lemma \ref{lem:factorization} shows that $\T_\et$ from Definition \ref{def:Tet} is symmetric monoidal: for $\T_{\et}^{\perf}$, we showed it is symmetric monoidal in Theorem \ref{thm:EK_perfect}, and the remaining functors are symmetric monoidal. Note that one must also check that finitely generated unit modules are closed under $\otimes$, which can be done modulo $p$ by Remark \ref{rem:fgu_mod_p} and then using the roots provided by Corollary \ref{cor:global_root_fgu_condition} and the fact that $\otimes$ commutes with colimits.
\end{remark}

The same factorization proved in Lemma \ref{lem:factorization} controls the essential image of $\T_\et$.

\begin{proposition} \label{prop:T_et_cons}
Every object in the essential image of $\T_\et$ on $\D_{\fgu}((\X,\mathcal{M})_\Prism,\Laur)^\op$ is a $p$-Kummer sheaf in $\D^{(b)}_{\cons,\D}((\X,\mathcal{M})_\eta,\Z_p)$.
\end{proposition}

\begin{proof}
The result is automatically $p$-Kummer by Lemma \ref{lem:factorization} and Lemma \ref{lem:pKum_framed_criterion}. 

For constructibility, we again use Lemma \ref{lem:factorization}. Inspecting the diagram, it suffices to know that the result in $\D^{(b)}_{\cons}(\Spec (\R_0^{\sq})_{\infty,\eta}^{\flat},\Z_p)^{h\Gamma}$ is constructible for the stratification induced by $\V(\T_i^\flat)$. Indeed, analytification and untilting respect these strata, as does the $\Gamma$-action as the strata on $(\X,\mathcal{M})_{\infty,\eta}$ are stable under the $\Gamma$ action, and descend to the strata $(\X,\mathcal{M})_{\Lambda,\eta}$ induced by $\D_\eta$ on $(\X,\mathcal{M})_\eta$. 

Proposition \ref{prop:perfectoid_constructibility} shows exactly this result, so the desired claim follows.
\end{proof}

\section{The main equivalence}

In this section we will prove the main theorem, which is that
\[\T_\et: \D_{\fgu}((\X,\mathcal{M})_\Prism,\Laur)^\op \longrightarrow \D^{(b),\mathrm{pKum}}_{\cons,\D}((\X,\mathcal{M})_\eta,\Z_p)\]
is an equivalence. 

In the first subsection, we focus on the open stratum $\X_\eta \setminus \D_\eta$. We will always denote the canonical inclusion map by $j: \X_\eta \setminus \D_\eta \to \X_\eta$.

We first construct a log prismatic crystal we call $\Laur(\ast \D)$ which corresponds to $j_! \underline{\Z_p}$, and argue that
\[\Perf((\X,\mathcal{M})_\Prism,\Laur(\ast \D))^{\varphi=1,\op}\]
is sent by $\T_\et$ to derived lisse sheaves on $\X_\eta \setminus \D_\eta$. More precisely, the essential image consists of these sheaves but embedded in $\D^{(b),\mathrm{pKum}}_{\cons,\D}((\X,\mathcal{M})_\eta,\Z_p)$ via the fully faithful functor $j_!$. This appears as Theorem \ref{thm:open_stratum}.

While superficially looking similar to \cite{KoshikawaYaoCorrigendum}, this result on the open stratum is actually the most nontrivial part of the main theorem because the main technique of passing to completed perfections (as used in \cite{BhattScholzePrismaticFCrystals} and \cite{KoshikawaYaoCorrigendum}) fails.\footnote{See Remark \ref{rem:perf_invar_counter}.} Briefly, the approach is to first argue strict \'{e}tale descent to work on a framed log affine with framing $\sq$ and associated log $q$-de Rham prism $(\R_0^{\sq}\llbracket q-1 \rrbracket, [p]_q, \mathcal{M})^a$. Once we are in this situation, we can use Lemma \ref{lem:factorization} and an analogue of Theorem \ref{thm:qdR_crys_Gamma_mod} to reduce full faithfulness to Proposition \ref{prop:open_decompletion_hom}, equivalently full faithfulness of the base change functor from $\Perf(\R_0^{\sq}\laur{q-1}^{\wedge}_p[1/\T_1,\ldots,1/\T_r]_p^{\wedge})^{\varphi=1,h\Gamma}$ to 
\[\Perf(\Ainf((\R_0^{\sq})_\infty)[1/\xi,1/[\T_1^\flat],\ldots,1/[\T_r^\flat]]^{\wedge}_p)^{\varphi=1,h\Gamma}.\]
We reduce this to a calculation in $\Gamma$-cohomology which is then carried out using decompletion techniques from \cite{KedlayaLiu2016}. We must also use a result on full faithfulness of analytification from \cite{BhattPadicHodgeNotes}.

For essential surjectivity, we make use of \cite{KoshikawaYaoCorrigendum}. Dualizing their correspondence gives an equivalence 
\[\Perf((\X,\mathcal{M})_\Prism,\Laur)^{\varphi=1,\op} \xlongrightarrow{\sim} \D^{(b),\mathrm{pKum}}_{\mathrm{lisse}}((\X,\mathcal{M})_\eta,\Z_p).\]
The idea is then to use the main result of \cite{lutkebohermert}, rephrased in \cite{DiaoLanLiuZhuLogAdic} to assert that a local system on the open stratum admits an extension to a Kummer-\'{e}tale one. After a finite prime-to-$p$ Kummer cover to make this a $p$-Kummer local system we may then produce a crystal $\mathcal{E}$ through \cite{KoshikawaYaoCorrigendum}, at which point we can tensor with $\Laur(\ast \D)$ and descend to a preimage in $\Perf((\X,\mathcal{M})_\Prism,\Laur(\ast \D))^{\varphi=1,\op}$.

The remainder of the argument to deduce the main theorem uses a strategy to reduce to thick generators in the form of Lemma \ref{lem:thick_equiv_criterion}. We prove a variant of Kashiwara's lemma to handle closed strata, and deduce a hollow variant of the theorem on the open stratum (for technical reasons). We can then prove the theorem by matching thick generators. More precisely, we know by this point that we have a functor
\[\T_\et: \D_{\fgu}((\X,\mathcal{M})_\Prism,\Laur)^\op \longrightarrow \D^{(b),\mathrm{pKum}}_{\cons,\D}((\X,\mathcal{M})_\eta,\Z_p)\]
and it is easy to see that it is exact. The target category is thickly generated by extensions by zero of local systems on locally closed log strata, so for essential surjectivity one must provide corresponding crystals on the source. We show these corresponding objects are also thick generators of the source, and moreover for maps between generators $\T_\et$ is fully faithful. The claim then follows.

\subsection{The open stratum}

We begin by constructing $\Laur(\ast \D)$, as explained in the overview. For an $\O_\X^\times$-torsor, let $\mathcal{L}_{\O_\X^\times}(-)$ denote passage to the associated line bundle over $\O_\X$. Recall also that part of the definition of the log structure $\mathcal{M}$ is a map
\[\alpha: \mathcal{M} \longrightarrow (\O_\X, \cdot),\]
which we use in the following construction. We remark that this is similar to \cite{TsujiLogCrystallinePD}, where such a crystal is constructed for the log crystalline site.

\begin{definition}
Let $\mathcal{M}(-\D)$ be the $\O_\X^\times$-torsor of local sections $m\in\mathcal{M}$ such that $\alpha_\X(m)\O_\X=\O_\X(-\D)$. Let $(\rA,\I,\M_{\Spf \rA})^a$ be a bounded log prism in $(\X,\mathcal{M})_\Prism^\op$, with structure map
\[f:(\Spf \rA/\I,\M_{\Spf \rA/\I})\to (\X,\mathcal{M}).\]
Let $i:\Spf \rA/\I\to\Spf \rA$ be the closed immersion. On $(\Spf \rA/\I)_\et$, set $f^*\mathcal{M}(-\D)$ to be obtained from $f^{-1}\mathcal{M}(-\D)$ as the $\O^\times_{\Spf \rA/\I}$-torsor given by change of structure group along $f^{-1} \O_\X^\times \to \O^\times_{\Spf \rA/\I}$. The map $f^\flat: f^{-1}\mathcal{M} \to \M_{\Spf \rA/\I}$ obtained from the structure map $f$ induces a morphism of \'{e}tale sheaves $f^*\mathcal{M}(-\D)\to\M_{\Spf \rA/\I}$. On $(\Spf \rA)_\et$, set $\widetilde{\mathcal{M}}_{\rA}(-\D):=\M_{\Spf \rA}\times_{i_*\M_{\Spf \rA/\I}}i_*f^*\mathcal{M}(-\D)$. We then finally take
\[\O_{\Prism}(-\D)(\rA,\I,\M_{\Spf \rA}):=\Gamma\left(\Spf \rA,\mathcal{L}_{\O_{\Spf \rA}^\times}(\widetilde{\mathcal{M}}_{\rA}(-\D))\right),\]
which one can think of as taking an $\O_{\Spf \rA}$-span of the generators in $\widetilde{\mathcal{M}}_{\rA}(-\D)$ we have extracted.
\end{definition}

The previous formulation was useful for working globally in an easy manner, but we point out that when working locally one can pick $m_\D\in \mathcal{M}$ so $\alpha(m_\D)$ is a local defining equation for $\D$. What we called $\mathcal{M}(-\D)$ is a globalization of this generated by such sections. Then we have
\[\widetilde{\mathcal{M}}_{\rA}(-\D) = \left\{m\in \M_{\rA}: \mathrm{red}(m) = u f^{\flat}(m_\D) \text{ locally for } u\in \O_{\rA/\I}^\times \right\},\]
where $f^\flat: f^{-1}\mathcal{M} \to \M_{\Spf \rA/\I}$ is obtained from the structure map $f$ and $\mathrm{red}: \M_{\rA} \to \M_{\rA/\I}$ is the quotient by $1+\I\O_{\rA}$.

We can now define $\Laur(\ast \D)$.

\begin{definition} \label{def:open_crystal}
For $n\ge 0$, set
\[\O_\Prism(n\D):=\underline{\Hom}_{\O_\Prism}\left(\O_\Prism(-\D)^{\otimes^{\Box} n},\O_\Prism\right).\]
The map $\O_\Prism(-\D)\to \O_\Prism$ gives transition maps $\O_\Prism(n\D)\to \O_\Prism((n+1)\D)$. We define
\[\Laur(\ast\D):=\left(\left(\colim_{n\ge0}\O_\Prism(n\D)\right)^{\wedge}_p[1/\I_\Prism]\right)^{\wedge}_p.\]
\end{definition}

We use the solid tensor product here, as we want to emphasize that the colimit happens in the ambient category $\D_{\Box}((\X,\mathcal{M})_{\Prism},\Laur)$. From this description we also see Frobenius on $\O_\Prism$ extends to $\Laur(\ast\D)$. Strict \'{e}tale locally on a framed log affine, if $\D$ is cut out by parameters $\T_1,\ldots,\T_r$, then $\Laur(\ast\D)$ evaluates on the corresponding log $q$-de Rham prism to $\R_0^{\sq}\laur{q-1}^{\wedge}_p[1/\T_1,\ldots,1/\T_r]^{\wedge}_p$. Lemma \ref{lem:factorization} then allows us to compute $\T_\et(\Laur(\ast\D))\simeq j_!\underline{\Z_p}$.

As this evaluates to a ring on every log prism, we may consider $\Perf((\X,\mathcal{M})_\Prism,\Laur(\ast \D))^{\varphi=1}$. This may be described explicitly by a \v{C}ech nerve for a framed log affine $(\Spf \R,\mathcal{M})$ with framing $\sq$. Let $\rA^{(m)}$ be the Laurentized \v{C}ech nerve of the log $q$-de Rham prism $\rA^+ = \R_0^{\sq}\llbracket q-1 \rrbracket$. Lemma \ref{lem:qdr_laur_cech} together with compatibility of taking module objects with limits yields
\[\D_{\Box}((\X,\mathcal{M})_\Prism,\Laur(\ast\D)) \simeq \lim_{[m]\in\Delta}\D_{\Box}(\rA^{\circ,(m)}),\]
where $\rA^{\circ,(m)}$ are defined as the evaluations of $\Laur(\ast\D)$ on the \v{C}ech nerve of the log $q$-de Rham prism. Thus taking dualizable objects and applying an equalizer $(-)^{\varphi=1}$ we obtain 
\[\Perf((\X,\mathcal{M})_\Prism,\Laur(\ast\D))^{\varphi=1} \simeq \lim_{[m]\in\Delta}\Perf(\rA^{\circ,(m)})^{\varphi=1}.\]
Here, as in Lemma \ref{lem:perf_comparison} in the Appendix, we appeal to the fact that (in our situation) one can ignore solid structures when restricting to dualizables.

\begin{lemma} \label{lem:open_localized_fgu}
We have $\Laur(\ast \D)\in \D_{\fgu}((\X,\mathcal{M})_\Prism,\Laur)$.

Moreover, any crystal in $\Perf((\X,\mathcal{M})_\Prism,\Laur(\ast \D))^{\varphi=1}$ is also in $\D_{\fgu}((\X,\mathcal{M})_\Prism,\Laur)$.
\end{lemma}

\begin{proof}
For the first claim that $\Laur(\ast \D)\in \D_{\fgu}((\X,\mathcal{M})_\Prism,\Laur)$, as a first step note that the finitely generated unit condition can be checked modulo $p$ by Remark \ref{rem:fgu_mod_p}. Upon doing so, the colimit presentation in Definition \ref{def:open_crystal} can be rewritten as $\left(\colim_{n\ge0}\varphi^{*n}\O_\Prism(\D)/p \right)[1/\I_\Prism]$. This colimit presentation using the root $\O_\Prism(\D)/p$ makes it clear that after evaluation on any log prism $(\rA^+,\I,\mathcal{M}_{\rA^+})^a$ the result is in $\mathrm{Thick}(\overline{\rA}[1/\I][\rF])$ and unit: unit is clear, and the fact that it lies in the desired thick subcategory can be seen by the argument in the footnote in Theorem \ref{thm:qdR_crys_Gamma_mod}.

For the second claim, we use Theorem \ref{thm:qdR_crys_Gamma_mod} as a shortcut. By Remark \ref{rem:qdR_test_sufficient} we can test membership in $\D_{\fgu}((\X,\mathcal{M})_\Prism,\Laur)$ using log $q$-de Rham prisms, so we pick a log $q$-de Rham prism $(\R_0^{\sq}\llbracket q-1 \rrbracket, [p]_q,\M)^a$ and wish to show for $\rA=\R_0^{\sq}\laur{q-1}^{\wedge}_p$ that
\[\Perf(\rA[1/\T_1,\ldots,1/\T_r]_p^{\wedge})^{\varphi=1} \hookrightarrow \D_{\fgu}(\rA).\]
Again using Remark \ref{rem:fgu_mod_p}, this can be tested modulo $p$. On the target, in our present situation it amounts to finite generation of the cohomology modules over $\overline{\rA}[\rF]$, so we can check on the heart. But there it is clear, as we can choose a coherent $\overline{\rA}$-submodule $\M_0\subset \M$ with $\M_0[1/\T]=\M$, and let $\beta:\M \to\varphi^*\M$ be inverse to the unit Frobenius on $\M$. But then for some $\rN \ge 0$ we have
\[\beta(\M_0)\subseteq\T^{-\rN}\varphi^*\M_0.\]
It follows that $\T^{-k} \M_0$ is the desired root once $k$ is large; the claim follows.
\end{proof}

Thus, we can now make sense of the main theorem's statement for this subsection:

\begin{theorem} \label{thm:open_stratum}
Let $(\X,\mathcal{M})$ be as in Assumption \ref{ass:horizontal_nc} over $\W(k)[\zeta_p]$. Then $\T_\et$ fits in a commutative diagram
\begin{display}
\Perf((\X,\mathcal{M})_\Prism,\Laur(\ast \D))^{\varphi=1,\op} \ar[hook]{d} \ar{r}{\sim} & \D^{(b)}_{\mathrm{lisse}}((\X_\eta \setminus \D_\eta)_{\proet}, \Z_p) \ar[hook]{d}{j_!} \\
\D_{\fgu}((\X,\mathcal{M})_\Prism,\Laur)^\op \ar{r}{\T_\et} & \D^{(b)}_{\cons,\D}((\X,\mathcal{M})_\eta,\Z_p).
\end{display}
The top horizontal arrow is an equivalence.
\end{theorem}

\begin{remark}
\label{rem:strata}
Although \cite[Corollary 3.7]{BhattScholzePrismaticFCrystals} describes the target using ordinary Laurent $\rF$-crystals, it is not directly suited to the functor $j_!$. Among several issues is that the complement $\X_\eta\setminus\D_\eta$ need not be quasicompact; for example, consider $\P^1_{\Q_p}\setminus\{\infty\}=\A^{1,\an}_{\Q_p}$. Moreover, a formal model of this complement need not map to the fixed smooth model $\X$ with generic fiber $\X_\eta$, and even worse Raynaud models require noncanonical admissible blowups to acquire such a map and need not remain smooth. We fix these problems by using log formal schemes to encode the stratification.
\end{remark}

A second remark worth making concerns a comparison to variants of \cite{yamashita} and \cite{Faltings1989}. As we can access more general constructible sheaves it seems one can likely pursue this further to work out comparison theorems with constructible coefficients, rather than just $j_! \mathbb{L}$ as in the remark. 

\begin{remark} \label{rem:comparison}
Let $\pi=\zeta_p-1$ be a uniformizer of $\W(k)[\zeta_p]$. If one does not invert $\I_\Prism$ and instead takes the unitalization of $\O_\Prism(\D) \to \varphi^* \O_\Prism(\D)$, one can calculate that the crystalline realization (i.e. the unitalization $\mathcal{E}$ of $\O_\cris(\D) \to \varphi^* \O_\cris(\D)$) satisfies
\[\R\Hom_{\D^{\varphi}(((\X,\mathcal{M})_{\pi=0}/\W(k))_\cris,\O_\cris)} (\mathcal{E},\O_\cris) \simeq \left(\R\Gamma_{\cris,c}((\X,\mathcal{M})_{\pi=0}/\W(k))\right)^{\varphi=1}.\]
Here compactly supported log crystalline cohomology is as defined in \cite{TsujiLogCrystallinePD}. This shows the crystalline realization of the natural non-Laurentized variant of $\Laur(\ast \D)$ is compatible with the expectation one would have from the open $\mathrm{C}_{\mathrm{st}}$ conjecture.
\end{remark}

Before we begin the proof, we must establish some additional basic facts about the category $\Perf((\X,\mathcal{M})_\Prism,\Laur(\ast \D))^{\varphi=1}$.

\begin{lemma} \label{lem:open_stratum_desc}
Let $(\X,\mathcal{M})$ be as in Assumption \ref{ass:horizontal_nc} over $\W(k)[\zeta_p]$. The category 
\[\Perf((\X,\mathcal{M})_\Prism,\Laur(\ast \D))^{\varphi=1,\op}\]
is a stack for the strict \'{e}tale topology on such $(\X,\mathcal{M})$.
\end{lemma}

\begin{proof}
This follows from Theorem \ref{thm:fgu_etale_desc}, because we can extract the desired category 
\[\Perf((\X,\mathcal{M})_\Prism,\Laur(\ast \D))^{\varphi=1}\]
as dualizable modules over the idempotent algebra $\Laur(\ast \D)$ in $\D_{\fgu}((\X,\mathcal{M})_\Prism,\Laur)$. Formally, both of these operations preserve the limit diagram we use to test descent along a strict \'{e}tale cover $(\U,\mathcal{M}_\U) \to (\X,\mathcal{M})$, namely that we have a limit diagram
\[\D_{\fgu}((\X,\mathcal{M})_\Prism,\Laur) \simeq \lim_{[n]\in \Delta} \D_{\fgu}((\U^{(n)},\mathcal{M}_{\U^{(n)}})_\Prism,\Laur)\]
for the \v{C}ech nerve. Because $\Laur(\ast \D)$ is an idempotent algebra in $\D_{\fgu}((\X,\mathcal{M})_\Prism,\Laur)$, taking modules over it gives a limit diagram, and then the same holds for taking dualizable objects. In fact this last step is actually not necessary, since the Fitting ideal argument in Proposition \ref{prop:perfectoid_constructibility} shows it is redundant.
\end{proof}

\begin{lemma} \label{lem:open_crys_gamma_mod}
Let $(\Spf \R,\mathcal{M})$ be a framed log affine. Then
\[\Perf((\Spf \R,\mathcal{M})_\Prism, \Laur(\ast \D))^{\varphi=1} \simeq \Perf(\R_0^{\sq}\laur{q-1}^{\wedge}_p[1/\T_1, \ldots, 1/\T_r]^{\wedge}_p)^{\varphi=1, h\Gamma}.\]
\end{lemma}

\begin{proof}
Put $\rA=\R_0^{\sq}\laur{q-1}^{\wedge}_p$ and $\rA^{\circ}=\rA[1/\T_1,\ldots,1/\T_r]^{\wedge}_p$, so that evaluation on the logarithmic $q$-de Rham prism carries $\Laur(\ast\D)$ to $\rA^\circ$ with its natural $\Gamma$-action. Using Lemma \ref{lem:open_localized_fgu} we see $\Perf((\Spf \R,\mathcal{M})_\Prism, \Laur(\ast \D))^{\varphi=1}$ is a full subcategory of finitely generated unit Laurent crystals, and due to Theorem \ref{thm:qdR_crys_Gamma_mod} we see $\Perf((\Spf \R,\mathcal{M})_\Prism, \Laur(\ast \D))^{\varphi=1} \xlongrightarrow{} \Perf(\rA^{\circ})^{\varphi=1, h\Gamma}$ is fully faithful.

Conversely, let $\M\in \Perf(\rA^{\circ})^{\varphi=1,h\Gamma}$. This category clearly sits inside $\D_{\fgu}(\rA)^{h\Gamma}$, thus Theorem \ref{thm:qdR_crys_Gamma_mod} gives an fgu crystal $\mathcal{E}$ with evaluation $\M$. The localization map
\[\mathcal{E}\longrightarrow\mathcal{E}\otimes\Laur(\ast\D)\]
is an equivalence after evaluation on the log $q$-de Rham prism, hence is an equivalence by conservativity. Thus $\mathcal{E}$ is an $\Laur(\ast\D)$-module. Using the description 
\[\Perf((\Spf \R,\mathcal{M})_\Prism,\Laur(\ast\D))\simeq \lim_{[m]\in\Delta}\Perf(\rA^{\circ,(m)}),\]
since dualizability is detected termwise in a limit of symmetric monoidal categories $\mathcal{E}$ is perfect over $\Laur(\ast\D)$. The claim follows. Here we have also used the crystal property, which shows the evaluations on $\rA^{\circ,(m)}$ are perfect since they are base changes of the perfect evaluation on $\rA^{\circ}$.
\end{proof}

We now tackle full faithfulness. Our strategy will be to show this in the perfectoid setting (using the formulation of $\T_\et$ in Lemma \ref{lem:factorization}). However unlike Laurent $\rF$-crystals in perfect complexes, it is not possible to use only the unit Frobenius to reduce to the perfectoid setting. We explain this further in Remark \ref{rem:perf_invar_counter}.

\begin{proposition} \label{prop:open_decompletion_hom}
Let $(\Spf \R,\mathcal{M})$ be a framed log affine. Then evaluating on the perfect log prism $\Ainf((\R_0^{\sq})_\infty)$ gives a fully faithful functor from $\Perf((\Spf \R,\mathcal{M})_\Prism,\Laur(\ast \D))^{\varphi=1}$ to 
\[\Perf(\Ainf((\R_0^{\sq})_\infty)[1/\xi,1/[\T_1^\flat],\ldots,1/[\T_r^\flat]]^{\wedge}_p)^{\varphi=1,h\Gamma}.\]
\end{proposition}

\begin{proof}
As usual, associated to the framed log affine we have a logarithmic $q$-de Rham prism, so we let $\rA=\R_0^{\sq} \laur{q-1}_p^{\wedge}$. For this particular proof we only care about a localization, so put $\rA^{\circ} := \rA[1/\T_1,\ldots,1/\T_r]_p^{\wedge}$ and $\rA_{\perf}^{\circ} = \Ainf((\R_0^{\sq})_\infty)[1/\xi,1/[\T_1^\flat],\ldots,1/[\T_r^\flat]]^{\wedge}_p$. We can express the source as $\Perf(\rA^{\circ})^{\varphi=1,h\Gamma}$ using Lemma \ref{lem:open_crys_gamma_mod}; the target is already $\Perf(\rA_\perf^{\circ})^{\varphi=1,h\Gamma}$.

By the tensor/Hom adjunction, it suffices to show that
\[\RGamma_{\mathrm{cont}}(\Gamma,\M)^{\varphi=1}\longrightarrow\RGamma_{\mathrm{cont}}(\Gamma,\M\otimes_{\rA^{\circ}} \rA^{\circ}_{\perf})^{\varphi=1}\]
is an isomorphism for $\M \in \Perf(\rA^{\circ})^{\varphi=1,h\Gamma}$. We can also use the ordinary tensor product $\M\otimes_{\rA^{\circ}} \rA^{\circ}_{\perf}$ in the target because the relevant map $\rA^{\circ} \to \rA^{\circ}_{\perf}$ is faithfully flat. Indeed, Kunz's theorem shows $\rA^+\to\rA^+_{\perf}$ is $p$-completely faithfully flat, and the map $\rA^{\circ}\to\rA^{\circ}_{\perf}$ is its $p$-completed base change after inverting $[p]_q,\T_1,\ldots,\T_r$, so it is $p$-completely faithfully flat; since $\rA^{\circ}$ is Noetherian, it is faithfully flat.

By derived Nakayama we can work with mod $p$ coefficients, which we now do. The source and target both have standard $t$-structures modulo $p$ obtained by restricting the $t$-structures on $\D(\rA^{\circ}/p)$ and $\D(\rA^{\circ}_{\perf}/p)$. To check these are well-defined one must check truncations remain perfect -- this is clear for the source because $\rA^{\circ}/p$ is Noetherian regular. For the target the $t$-structure exists before applying $(-)^{h\Gamma}$, which is implicit in the proof of \cite[Proposition 3.4]{BhattScholzePrismaticFCrystals}: in particular they show that the top cohomology of a perfect unit Frobenius module over an $\F_p$-algebra is a vector bundle. Using a Postnikov d\'{e}vissage (\cite[Lemma 5.4.3]{HauckStacky}) we reduce to testing the claim for $\M$ in the heart killed by $p$ (which must lie in $\Vect(\rA^{\circ}/p)^{\varphi=1,h\Gamma}$).

We now prove $\RGamma_{\mathrm{cont}}(\Gamma,\M)^{\varphi=1}\simeq \RGamma_{\mathrm{cont}}(\Gamma,\M\otimes_{\rA^{\circ}} \rA^{\circ}_{\perf})^{\varphi=1}$ for $\M$ in $\Vect(\rA^{\circ}/p)^{\varphi=1,h\Gamma}$. The faithfully flat map $\rA^\circ/p \longrightarrow \rA^\circ_{\perf}/p$ is injective, with cokernel $(\rA^\circ_{\perf}/\rA^\circ)/p$. Since $\M\in\Vect(\rA^\circ/p)$, tensoring this sequence with $\M$ is exact and we may reduce to showing
\[\RGamma_{\mathrm{cont}}(\Gamma,\M\otimes_{\rA^{\circ}}(\rA^\circ_{\perf}/\rA^{\circ}))^{\varphi=1}=0.\]
Now choose $m\gg0$ and let $(\R_0^{\sq})_m$ denote the ring obtained by $p$-completely adjoining $p^m$th roots of the $\T_i$. Define $\varphi^{-m}\rA^{\circ}=(\R_0^{\sq})_m\laur{q^{1/p^m}-1}[1/\T_1,\ldots,1/\T_r]_p^{\wedge}$ and $\varphi^{-m}\rA^{\circ,\cyc}=(\R_0^{\sq})_m\laur{q^{1/p^\infty}-1}^{\wedge}[1/\T_1,\ldots,1/\T_r]^{\wedge}_p$, viewed as intermediate subrings
\[\rA^{\circ} \subset \varphi^{-m}\rA^{\circ}\subset \varphi^{-m}\rA^{\circ,\cyc}\subset \rA^{\circ}_{\perf}.\]
Since $\varphi$ and $\Gamma$ commute, the unit Frobenius of $\M$ gives a $\Gamma$-equivariant equivalence
\[\RGamma_{\mathrm{cont}}(\Gamma,\M)^{\varphi=1}\simeq\RGamma_{\mathrm{cont}}(\Gamma,\M\otimes_{\rA^\circ}\varphi^{-m}\rA^\circ)^{\varphi=1}.\]

For $m \gg 0$, the exact triangles associated to 
\[\varphi^{-m}\rA^{\circ}/\rA^{\circ} \longrightarrow \rA^{\circ}_{\perf}/\rA^{\circ} \longrightarrow \rA^{\circ}_{\perf}/\varphi^{-m}\rA^{\circ}\]  
and to $\varphi^{-m}\rA^{\circ,\cyc}/\varphi^{-m}\rA^{\circ} \longrightarrow \rA^{\circ}_{\perf}/\varphi^{-m}\rA^{\circ} \longrightarrow \rA^{\circ}_{\perf}/\varphi^{-m}\rA^{\circ,\cyc}$ show
\[\RGamma_{\mathrm{cont}}(\Gamma,\M\otimes_{\rA^{\circ}}(\rA^{\circ}_{\perf}/\rA^{\circ}))^{\varphi=1}=0\]
by applying Lemmas \ref{lem:open_decompletion_toric_acyclic} and \ref{lem:open_decompletion_q_acyclic}. The claim then follows.
\end{proof}

We will next provide the needed technical lemmas on continuous $\Gamma$ cohomology. We retain the notation of Proposition \ref{prop:open_decompletion_hom}.

To make sense of the cohomology calculations, we must work out the relevant topologies on $\rA^{\circ}/p$ and related rings to compute continuous group cohomology with. Using Remark \ref{rem:weak_top_fgu} we see the ring
\[\overline{\rA}^{\circ}:=\rA^\circ/p=\bigcup_{i\ge0}\T^{-i}\overline{\R}^{\sq}\laur{q-1}\]
is equipped with the inductive limit topology and each $\T^{-i}\overline{\R}^{\sq}\laur{q-1}$ has the topology induced from $\overline{\R}^{\sq}\laur{q-1}$. We must give vector bundles $\M$ over this ring their weak topology, which is all we need to deal with. Note that the statement of Lemma \ref{lem:weak_top} also applies over $\overline{\rA}^{\circ}$ since we use vector bundles, where the proof is immediate. We use similar constructions for $\rA^{\circ,\cyc}/p$ and $\rA^{\circ}_{\perf}/p$.

Let $\overline{\rA}^+ = (\R_0^{\sq}/p)\llbracket q-1 \rrbracket$. By picking a $\Gamma$-stable lattice $\M^+\in \Mod_{\overline{\rA}^+}$ so $\M^+[1/(q-1), 1/\T_1, \ldots, 1/\T_r]=\M$, we see just as in the proof of Lemma \ref{lem:Gamma_Frob_improvement} that after shrinking $\Gamma$ to $\Gamma_0 \le \Gamma$, eventually $\gamma_i-1$ raises the $q-1$-adic filtration $\Fil^{\bullet} \T^{-\alpha} \M^+[1/(q-1)]$ for each multi-index $\alpha \in \N^r$. We specifically can choose $\Gamma_0=(1+p^m\Z_p)\ltimes(p^{m-1}\Z_p(1))^n \subset \Gamma$ for $m\gg 0$, which will ultimately be the choice of $m$ we need for Proposition \ref{prop:open_decompletion_hom}.

To prove the vanishing of $\R\Gamma(\Gamma,-)$ in Lemmas \ref{lem:open_decompletion_q_acyclic} and \ref{lem:open_decompletion_toric_acyclic}, we adapt the method of \cite[Lemmas 5.6.3 and 5.6.4]{KedlayaLiu2016}.\footnote{It also appears there are some gaps in \S 5, but the overall method is fine.} In the above, what we have been implicitly working with is an Ind-Banach $\F_p$-module. That is, $\M = \colim_\rN \M^{(\rN)}$, which can be taken in the previous discussion to be $\T^{-\rN} \M^+[1/(q-1)]$. Here $\T^{-\rN} := \T_1^{-\rN} \cdots \T_r^{-\rN}$, as this is a cofinal subsystem. 

We work with a more restrictive variant of Ind-Banach $\F_p$-modules. We will abbreviate this just to ``LB'' modules and rings, but our convention is that an LB module means a presentation $\M=\colim_{\rN\in\N}\M^{(\rN)}$ by Banach modules with closed contractive embeddings. We then say LB rings are algebra objects with LB underlying modules. We also write
\[\widehat{\bigoplus}^{\mathrm{LB}}_j\M_j:=\colim_\rN\widehat{\bigoplus}_j\M_j^{(\rN)},\]
where the completed sum is taken for the supremum norm; in applications, all norms are normalized by the same value of $|q-1|$.

For a closed subgroup $\Gamma_0\le\Gamma$, choose topological generators $g_1,\ldots,g_k$ of $\Gamma_0 \cap\Z_p(1)^n\simeq\Z_p^k$. If the image of $\Gamma_0$ in $1+p\Z_p$ is nontrivial, choose $g_0\in\Gamma_0$ lifting a topological generator of this image. Using the usual resolution for group cohomology of $\Z_p^k$ and Hochschild-Serre we may compute $\RGamma_{\mathrm{cont}}(\Gamma_0,\M)$ as the fiber of $g_0-1$ on
\[\mathrm{Kos}^\bullet(\M;g_1-1,\ldots,g_k-1),\]
where $g_0$ acts on the coefficients and by conjugation on the Koszul generators. If the image in $1+p\Z_p$ is trivial, the Koszul complex itself computes cohomology.

This cochain model of cohomology allows us to make the following definitions, subject to the choices of generators above: a complex is \emph{strictly contractible} if it admits a continuous contracting homotopy. For a family of complexes $\rC_j^\bullet$, such homotopies $h_j^i$ are \emph{uniformly bounded} if they preserve the $\rC_j^{i,(\rN)}$ and
\[\sup_{j,i} \bigl\| h_j^i:\rC_j^{i,(\rN)} \longrightarrow \rC_j^{i-1,(\rN)} \bigr\| < \infty\]
for every $\rN$.

The following simple observation will be used throughout.

\begin{lemma} \label{lem:kedlaya-liu}
Let $\Gamma_0 \le \Gamma$ be a closed subgroup and let $\B$ be an LB ring with a continuous action of $\Gamma_0$, and let $\V$ be an LB $\B$-module with a continuous $\Gamma_0$-action.

Suppose that we have a $\Gamma_0$-equivariant decomposition of LB modules
\[\widehat{\bigoplus}^{\mathrm{LB}}_{j\ge0} \V_j \xlongrightarrow{\sim} \V,\]
such that the cochain complexes $\R\Gamma(\Gamma_0,\V_j)$ admit uniformly bounded contracting homotopies. Then $\R\Gamma(\Gamma_0,\V)=0$.
\end{lemma}

\begin{proof}
This is immediate from the definitions: uniform boundedness extends the contracting homotopy to $\R\Gamma(\Gamma_0,\V)$.
\end{proof}

Note here we implicitly have made a choice of generators to produce the cochain models for cohomology. 

\begin{lemma} \label{lem:open_decompletion_q_acyclic}
Let $\M\in \Vect(\rA^{\circ}/p)^{\varphi=1,h\Gamma}$. For $m\gg0$,
\[\RGamma_{\mathrm{cont}}(\Gamma,\M \widehat{\otimes}_{\rA^{\circ}}(\varphi^{-m}\rA^{\circ,\cyc}/\varphi^{-m}\rA^{\circ}))=0.\]
\end{lemma}

\begin{proof}
It suffices to show the cohomology for $\Gamma_0=(1+p^m\Z_p)\ltimes (p^{m-1}\Z_p(1))^n$ vanishes, where $m$ is chosen as before so that $\gamma_i-1$ raises the $(q-1)$-adic filtration. We call the normal subgroup $\H=(p^{m-1}\Z_p(1))^n$.

To apply Lemma \ref{lem:kedlaya-liu}, we will decompose $\M\otimes_{\overline{\rA}^{\circ}}(\varphi^{-m}\overline{\rA}^{\circ,\cyc}/\varphi^{-m}\overline{\rA}^{\circ})$. To do this, we set $\B_0 = \varphi^{-m} \overline{\rA}^{\circ} = (\overline{\R}_0^{\sq})_m\laur{q^{1/p^{m}}-1}[1/\T_1,\ldots,1/\T_r]$. Recall we used the notation $(\overline{\R}_0^{\sq})_m$ to denote that we have adjoined $p^m$th roots of the $\T_i$. We then set $\B_j = (\overline{\R}_0^{\sq})_m \laur{q^{1/p^{m+j}}-1}[1/\T_1,\ldots,1/\T_r]$, so $\B_\infty$ is $\varphi^{-m}\overline{\rA}^{\circ,\cyc}$. We write $\M_{\B_j}$ for the base change of $\M$ to one of these LB rings. 

The decomposition $\B_{j+1}=\bigoplus_{a=0}^{p-1}\B_jq^{a/p^{m+j+1}}$ allows us to define $\Gamma$-equivariant $\B_j$-linear retractions $\psi_j$ by projection onto $a=0$. With the prescribed LB structures, one obtains
\[\widehat{\bigoplus}^{\mathrm{LB}}_j\ker(\psi_j) \simeq \B_\infty/\B_0.\]
Tensoring up with the vector bundle $\M$, we get a decomposition by $\M_{\B_j} \widehat{\otimes} \ker \psi_j$ of 
\[\M\widehat{\otimes}_{\overline{\rA}^{\circ}}(\varphi^{-m}\overline{\rA}^{\circ,\cyc}/\varphi^{-m}\overline{\rA}^{\circ}) = \M_{\B_{\infty}}/\M_{\B_0}.\]

Thus it only remains to produce the contracting homotopies. Lemma \ref{lem:kedlaya-liu} still applies here, as we still have produced an explicit cochain model for $\Gamma_0$-cohomology in the discussion preceding Lemma \ref{lem:kedlaya-liu} (despite $\Gamma_0$ not being isomorphic to $\Z_p(1)^n$). For this model, we recall we had to fix generators. We fix $\gamma_0 = 1 + p^{m}$ generating $1+p^m \Z_p$, and $\gamma_i$ generators of $\H=(p^{m-1}\Z_p(1))^n$.

It suffices actually to produce a contracting homotopy for the fiber of $\gamma_0-1$, but for the entire Koszul complex $\mathrm{Kos}_{\H}^\bullet(\M_{\B_j} \otimes_{\B_j} \ker \psi_j)$. Using the aforementioned decompositions $\B_{j+1}=\bigoplus_{a=0}^{p-1}\B_jq^{a/p^{m+j+1}}$, since the decomposition is $\H$-equivariant we can decompose further into
\[\mathrm{Kos}_{\H}^\bullet(\M_{\B_j} \otimes_{\B_j} \ker \psi_j) \simeq \bigoplus_{1\le a \le p-1} \mathrm{Kos}_{\H}^\bullet(\M_{\B_j} q^{a/p^{m+j+1}}).\]
To apply Lemma \ref{lem:kedlaya-liu} it will suffice to produce a contracting homotopy for the action of $\gamma_0-1$ on these summands. 

For each $a$, the summands $\mathrm{Kos}_{\H}^\bullet(\M_{\B_j} q^{a/p^{m+j+1}})$ in question can be written as
\[\M_{\B_j} q^{a/p^{m+j+1}} \xlongrightarrow{} \M_{\B_j} q^{a/p^{m+j+1}} \otimes_{\Z_p} \Lambda^1(\H) \xlongrightarrow{} \cdots \xlongrightarrow{} \M_{\B_j} q^{a/p^{m+j+1}} \otimes_{\Z_p} \Lambda^n(\H).\]
For $\I\subseteq\{1,\ldots,n\}$, write $e_\I=\bigwedge_{i\in \I}e_i$ for the corresponding basis element of the term for $\Lambda^{|\I|}(\H)$. Using $\gamma_0\gamma_i=\gamma_i^{1+p^m} \gamma_0$ we see
\[\gamma_0(- \cdot e_\I)=\left(\prod_{i\in \I}\frac{\gamma_i^{1+p^m}-1}{\gamma_i-1}\right)^{-1}\gamma_0(-)e_\I.\]
On the summand $\M_{\B_j} q^{a/p^{m+j+1}}$ we have $\gamma_0 - 1 = (q^{a/p^{j+1}}-1) \left[1+ (q^{a/p^{j+1}}-1)^{-1}(\gamma_0-q^{a/p^{j+1}})\right]$. Thus to give a contracting homotopy on the mapping fiber, it suffices to show the bound $\|(q^{a/p^{j+1}}-1)^{-1}(\gamma_0-q^{a/p^{j+1}})\|_{\rN} < 1$. We actually do not care if the bound depends on $j$, since in computing $\|(\gamma_0 - 1)^{-1}\|_{\rN}$ the inverse of $\left[1+ (q^{a/p^{j+1}}-1)^{-1}(\gamma_0-q^{a/p^{j+1}})\right]$ has norm $\le 1$ -- thus the bound is entirely dependent on $|(q^{a/p^{j+1}}-1)|^{-1} < |q-1|^{-1}$, making it uniform.

By using the action on $\B_j$ and using that $\gamma_0-1$ raises the $q-1$-adic filtration for coefficients this action eventually satisfies for $m \gg 0$
\[\left\|(q^{a/p^{j+1}}-1)^{-1}(\gamma_0-q^{a/p^{j+1}})\right\|_\rN \le |q-1|^{p^{-j}(1-p^{-m}-p^{-1})}<1.\]
To deduce this bound, first note that $\left\|(q^{a/p^{j+1}}-1)^{-1}(\gamma_0-q^{a/p^{j+1}})\right\|_\rN$ can be bounded by $\left\|(q^{a/p^{j+1}}-1)\right\|_{\rN}^{-1} \left\|(\gamma_0-q^{a/p^{j+1}})\right\|_\rN$, and right away we can bound the first factor by $|q-1|^{-p^{-j-1}}$. The second factor, since $\gamma_0-1$ raises the $q-1$-adic filtration on $\M_{\B_j}$ and we calculate the norm on $\M_{\B_j} q^{a/p^{m+j+1}}$, has norm bounded by $|q-1|^{p^{-j}-p^{-j-m}}$. 

This handles the fiber of $\gamma_0-1$ on $\M_{\B_j} q^{a/p^{m+j+1}}$. For the full Koszul complex, note that there are finitely many $\I\subseteq\{1,\ldots,n\}$ to deal with and $\left\|\left(\prod_{i\in \I}\frac{\gamma_i^{1+p^m}-1}{\gamma_i-1}\right)^{-1} - 1\right\|_\rN \le |q-1|^{p^m-1}$ (again crucially using that $\gamma_i-1$ raises the $q-1$-adic filtration). Thus we have produced the required homotopies and Lemma \ref{lem:kedlaya-liu} finishes the claim. Note that we only showed vanishing of $\Gamma_0$-cohomology, but by Hochschild-Serre this proves there is no $\Gamma$-cohomology as well.
\end{proof}

Fortunately, the remaining lemma is easier because we can deal with an abelian group.

\begin{lemma} \label{lem:open_decompletion_toric_acyclic}
Let $\M\in \Vect(\rA^{\circ}/p)^{\varphi=1,h\Gamma}$. For $m\gg0$,
\[\RGamma_{\mathrm{cont}}(\Gamma,\M \widehat{\otimes}_{\rA^{\circ}}(\rA^{\circ}_{\perf}/\varphi^{-m}\rA^{\circ,\cyc}))=0.\]
\end{lemma}

\begin{proof}
The argument is very similar, so we simply indicate the decomposition we use and contracting homotopies to apply the same strategy as the previous lemma.

In this case we put $\B_j = \varphi^{-(m+j)}\rA^{\circ,\cyc}/p$ and write $\M_{\B_j}$ for the base changes to these rings, and $\Gamma_0=(p^{m-1}\Z_p(1))^n$. We must produce a $\Gamma_0$-equivariant decomposition of
\[\M\otimes_{\rA^{\circ}}(\rA^{\circ}_{\perf}/\varphi^{-m}\rA^{\circ,\cyc}) = \M_{\B_{\infty}}/\M_{\B_0}.\]
Now for $\alpha=(a_1,\ldots,a_n)$, there is a decomposition $\B_{j+1}=\bigoplus_{\alpha\in\{0,\ldots,p-1\}^n}\B_j\prod_{\nu=1}^n\T_\nu^{a_\nu/p^{m+j+1}}$. Through this, we may define $\psi_j$ as the projection onto the $\alpha = \vec{0}$ component. With the obvious LB structures, these give
\[\widehat{\bigoplus}^{\mathrm{LB}}_{j\ge0}\ker(\psi_j)\simeq\B_\infty/\B_0\]
as LB modules, which we can similarly tensor with $\M$ to get a decomposition of $\M_{\B_{\infty}}/\M_{\B_0}$. We must then produce contracting homotopies for the individual summands $\M_{\B_j} \widehat{\otimes} \ker \psi_j$.

Using the generators $\gamma_1,\ldots,\gamma_n$ of $\Gamma_0$, the preceding $\Gamma_0$-equivariant decomposition reduces us to contracting the complexes
\[\mathrm{Kos}_{\Gamma_0}^\bullet\left(\M_{\B_j}\prod_\nu\T_\nu^{a_\nu/p^{m+j+1}}\right)\]
for $\alpha\neq0$. Choose the least $i$ with $a_i\neq0$. We have 
\[\gamma_i-1 = (q^{a_i/p^{j+1}}-1)\left[1+(q^{a_i/p^{j+1}}-1)^{-1}(\gamma_i-q^{a_i/p^{j+1}})\right].\] 
To contract the fiber of $\gamma_i-1$ uniformly, as before it suffices to show 
\[\left\| (q^{a_i/p^{j+1}}-1)^{-1}(\gamma_i-q^{a_i/p^{j+1}}) \right\|_{\rN}<1.\]
By using the action on $\B_j$ and that $\gamma_i-1$ is $q-1$-adic filtration raising for coefficients (by our choice of $m\gg 0$) we similarly obtain
\[\left\| (q^{a_i/p^{j+1}}-1)^{-1}(\gamma_i-q^{a_i/p^{j+1}}) \right\|_{\rN} \le |q-1|^{(p-1)p^{-j-1}}<1.\]
This gives $(\gamma_i-1)^{-1}$ of norm at most $|q-1|^{-1}$ on each of the summands for $\M_{\B_j}$ where $\alpha \neq 0$. For the full Koszul complex, we can use a similar argument to the previous lemma for $\gamma_i-1$. These homotopies are compatible and have the same bound, so Lemma \ref{lem:kedlaya-liu} then gives the vanishing for continuous $\Gamma_0$-cohomology, and Hochschild-Serre gives it for $\Gamma$.
\end{proof}

\begin{remark} \label{rem:perf_invar_counter}
It is tempting to try to prove Lemma \ref{lem:open_decompletion_q_acyclic} only using $\varphi$, as this can be done for ordinary $(\varphi,\Gamma)$-modules using invariance of \'{e}tale $\varphi$-modules under completed perfections. However in our situation this does not suffice. Consider the map
\[\rA=k\laur{q-1}\langle \T\rangle[1/\T]\to \B=k \laur{q^{1/p^\infty}-1}^{\wedge}\langle (\T^\flat)^{1/p^\infty}\rangle[1/\T^\flat],\]
with $\T\mapsto \T^\flat$. If full faithfulness followed just from invariance under this completed perfection, then $\Perf(\rA)^{\varphi=1}\to \Perf(\B)^{\varphi=1}$ would be fully faithful. By Katz's theorem this would force $\H^1_\et(\rA,\underline{\F_p})\to \H^1_\et(\B,\underline{\F_p})$ to be an isomorphism. We use the Artin-Schreier class on $\B$ defined by
\[f=\sum_{n\ge 1}(q-1)^n(\T^\flat)^{-1+1/p^n}.\]
Two Artin-Schreier classes $f$ and $g$ are equivalent if $f-g=h^p-h$ for some $h$. We get a counterexample if it is impossible to have such a relation where $g\in \rA$. Here, the example works because the exponents $-1+1/p^n$ lie in distinct orbits in $\Z[1/p]/p^{\Z}$ and the first integral representatives $1-p^n$ in the orbit grow arbitrarily negative. Since $g\in \rA$, in contrast there is a lowest exponent $\T^{-\ell}$, and we choose $n$ so $1-p^n$ is smaller. Then when solving $f-g=h^p-h$, past this point one may work out $h$ exponent by exponent. We find the coefficient of $(\T^\flat)^{-1+1/p^n}$ is necessarily zero (otherwise $h$ is forced to have unbounded poles), but for $(\T^\flat)^{-1/p^r+1/p^{n+r}}$ it is necessarily $(q-1)^{n/p^r}$. This shows there is no such convergent power series $h$.

This implies the \'{e}tale cohomology groups are not isomorphic via the natural base change map. The issue that occurred here is roughly that passage to $\B$ mixes completed and decompleted perfections: if we used the usual completed perfection, there would be no issue.
\end{remark}

\begin{proposition} \label{prop:open_fullfaithful}
The top horizontal functor in Theorem \ref{thm:open_stratum} is fully faithful.
\end{proposition}

\begin{proof}
We have strict \'{e}tale descent of the source by Lemma \ref{lem:open_stratum_desc} and for the target $\D^{(b)}_{\mathrm{lisse}}((\X_\eta\setminus\D_\eta)_{\proet},\Z_p)$ it is clear. Thus using Lemma \ref{lem:strict_etale_local_framing}, we may reduce the claim to a framed log affine. We will use Lemma \ref{lem:factorization} to analyze full faithfulness.

To begin, the source category can be identified with $\Perf(\R_0^{\sq}\laur{q-1}^{\wedge}_p[1/\T_1,\ldots,1/\T_r]^{\wedge}_p)^{\varphi=1,h\Gamma}$ by Lemma \ref{lem:open_crys_gamma_mod}, and Proposition \ref{prop:open_decompletion_hom} makes the base change functor to 
\[\Ainf((\R_0^{\sq})_\infty)[1/\xi,1/[\T_1^\flat],\ldots,1/[\T_r^\flat]]^{\wedge}_p\]
fully faithful. Lemma \ref{lem:factorization} supplies a formula for $\T_\et$ where this is the first step, and the remaining steps are obviously fully faithful (for example tilting or $\Gamma$-descent) with the exception of analytification.

Fortunately this claim about analytifications follows from results in \cite{BhattPadicHodgeNotes}. Put
\[\U = \Spec (\R_0^{\sq})_{\infty,\eta}^{\flat}[1/\T_1^\flat,\ldots,1/\T_r^\flat].\]
For $\mathbb{L},\mathbb{L}'\in \D^{(b)}_{\mathrm{lisse}}(\U_{\et},\F_p)$, depending on the rank $m_i$ of $\H^i\RcalHom(\mathbb{L},\mathbb{L}')$ the image of the associated representation of $\pi_1^{\et}(\U)$ is contained in a finite group $\GL_{m_i}(\F_p)$. Because of this, we therefore have a single trivializing \'{e}tale cover $\U' \to \U$ on which $\H^i\RcalHom(\mathbb{L},\mathbb{L}') \simeq \underline{\F_p}^{m_i}$ for every $i$. On $\U'$, \cite[Remark 6.2.7]{BhattPadicHodgeNotes} further gives
\[\R\Gamma(\U',\underline{\F_p}) \xlongrightarrow{\sim} \R\Gamma(\U'^{\an},\underline{\F_p}).\]
Finally, a Postnikov d\'{e}vissage combined with finite \'{e}tale descent to trivialize coefficients to $\underline{\F_p}^{\oplus \ell}$ gives full faithfulness modulo $p$; derived Nakayama then gives full faithfulness of analytification.
\end{proof}

We next prove essential surjectivity after passing, on a framed log affine, to a finite prime-to-$p$ Kummer cover. The following proposition allows us to descend back down from this cover.

\begin{proposition} \label{prop:open_strat_Kum_desc}
Let $(\Spf \R,\mathcal{M})$ be a framed log affine with coordinates $\T_1,\ldots,\T_r$, and let $\rN$ be prime to $p$. Put
\[\R_\rN=\R\langle \rS_1,\ldots,\rS_r\rangle/(\rS_1^\rN-\T_1,\ldots,\rS_r^\rN-\T_r),\]
and equip $\Spf \R_\rN$ with the log structure $\mathcal{M}_{\rN}$ induced by $\V(\rS_1\cdots \rS_{r})$. Then
\[\left(\Perf((\Spf \R_{\rN},\mathcal{M}_{\rN})_\Prism,\Laur(\ast \D_{\rN}))^{\varphi=1}\right)^{h\mu_{\rN}^r} \simeq \Perf((\Spf \R,\mathcal{M})_\Prism,\Laur(\ast \D))^{\varphi=1}.\]
\end{proposition}

In applications, we simply pass to opposite categories on both sides.

\begin{proof}
Let $f_{\rN}:(\Spf \R_{\rN},\mathcal{M}_{\rN})\to(\Spf\R,\mathcal{M})$ be the natural map of log formal schemes, where as above $\R_\rN=\R\langle \rS_1,\ldots,\rS_r\rangle/(\rS_1^\rN-\T_1,\ldots,\rS_r^\rN-\T_r)$. Since $f_{\rN}^*\D=\rN\D_{\rN}$, there is a canonical equivalence
\[f_{\rN}^*\Laur(\ast\D) \simeq \Laur(\ast(\rN\D_{\rN})) \simeq \Laur(\ast\D_{\rN}).\]
Thus pullback induces a functor on the corresponding module categories
\[\Perf((\Spf \R,\mathcal{M})_\Prism,\Laur(\ast \D))^{\varphi=1} \longrightarrow \Perf((\Spf \R_{\rN},\mathcal{M}_{\rN})_\Prism,\Laur(\ast \D_{\rN}))^{\varphi=1},\]
which we will use to induce the desired equivalence.

The claim becomes more transparent when we write out the \v{C}ech nerves for the source and target. Base changing the framing for $(\Spf \R,\mathcal{M})$ along $\T_i\mapsto\rS_i^{\rN}$ gives $(\Spf \R_{\rN},\mathcal{M}_{\rN})$ a framing $\sq_{\rN}$, so it also has an associated $q$-de Rham prism. Let $\rA = \R_0^{\sq}\laur{q-1}_p^{\wedge}$, and $\rA_{\rN} = \R_{0,\rN}^{\sq_{\rN}}\laur{q-1}_p^{\wedge}$ be their corresponding (Laurentized) log $q$-de Rham prisms. As usual we write $\rA^{(s)}$ and $\rA_{\rN}^{(s)}$ for the Laurentized \v{C}ech nerves, and we write $\rA^{(s)}(\ast \D)$ and $\rA_{\rN}^{(s)}(\ast \D)$ for the evaluations of $\Laur(\ast \D)$ on these prisms. 

We then obtain the following diagram:
\begin{display}[shift left/.default=+0.8ex]
\Perf(\rA_{\rN}(\ast \D_{\rN}))^{\varphi=1} \ar[shift left]{r} \ar[shift right]{r} & \Perf(\rA_{\rN}^{(1)}(\ast \D_{\rN}))^{\varphi=1} \ar{l} \ar[shift right, shift right]{r} \ar{r} \ar[shift left, shift left]{r} & \Perf(\rA_{\rN}^{(2)}(\ast \D_{\rN}))^{\varphi=1} \ar[shift right]{l} \ar[shift left]{l} \ar[shift left, shift left, shift left]{r} \ar[shift left]{r} \ar[shift right]{r} \ar[shift right, shift right, shift right]{r} & \cdots \ar[shift left, shift left]{l} \ar{l} \ar[shift right, shift right]{l} \\

\Perf(\rA(\ast \D))^{\varphi=1} \ar{u} \ar[shift left]{r} \ar[shift right]{r} & \Perf(\rA^{(1)}(\ast \D))^{\varphi=1} \ar{u} \ar{l} \ar[shift right, shift right]{r} \ar{r} \ar[shift left, shift left]{r} & \Perf(\rA^{(2)}(\ast \D))^{\varphi=1} \ar{u} \ar[shift right]{l} \ar[shift left]{l} \ar[shift left, shift left, shift left]{r} \ar[shift left]{r} \ar[shift right]{r} \ar[shift right, shift right, shift right]{r} & \cdots \ar[shift left, shift left]{l} \ar{l} \ar[shift right, shift right]{l} 
\end{display}
The rings on top all have a grading by $\alpha \in (\Z/\rN \Z)^r$, making $\Perf$ acquire an action of $\mu_{\rN}^r$\footnote{We also need the condition $\R_{\alpha} \cdot \R_\beta \subseteq \R_{\alpha + \beta}$ for $\alpha,\beta \in (\Z/\rN\Z)^r$ for a general ring $\R$, which is satisfied here.}. For example, $\R_{0,\rN}^{\sq_{\rN}}=\R_0^{\sq}\langle \rS_1,\ldots,\rS_r\rangle/(\rS_1^\rN-\T_1,\ldots,\rS_r^\rN-\T_r)$, so 
\[\rA_{\rN}(\ast \D_{\rN}) = \bigoplus_{\alpha \in (\Z/\rN \Z)^r} \rA(\ast \D) \cdot \rS^{\alpha},\]
with the $\alpha=0$ component being the natural inclusion of $\rA(\ast \D)$. For the first \v{C}ech term, the analogous grading induced from the one on $\R_{0,\rN}^{\sq_{\rN}}$ again holds, even at the level of $\rA_{\rN}^{+,(1)}$ and $\rA^{+,(1)}$. We see this by noting that the elements $v_i = \rS_{i,1}/\rS_{i,0}$ we adjoin (in the notation of Lemma \ref{lem:qdr_actual_cech_envelope}; for $\rA_{\rN}^{+,(1)}$ we use $v_i$ rather than $u_i$ to distinguish them) satisfy $v_i=(\rS_{i,1}/\rS_{i,0})=(\T_{i,1}/\T_{i,0})^{1/\rN} = u_i^{1/\rN} \in \rA^{+,(1)}$ since $\rN$ is prime to $p$. Indeed
\[u_i^{1/\rN} = \sum_{k\ge 0} \binom{1/\rN}{k} (u_i-1)^k,\]
which converges $(p,[p]_q)$-adically since we force $u_i-1\in [p]_q \rA^{+,(1)}$. For higher \v{C}ech terms, the claim is clear because they are obtained from iterated tensor products of $\rA^{+,(1)}$ over $\rA^+$ (and similarly for $\rA_{\rN}$).

Even better, here
\[\rA^{(s)}(\ast \D) \longrightarrow \rA_{\rN}^{(s)}(\ast \D_{\rN})\]
is a finite \'{e}tale $\mu_{\rN}^r$ torsor. Indeed, since we have $p$-completely inverted $\T_1,\ldots, \T_r$ and $\rS_1,\ldots,\rS_r$ the map becomes finite \'{e}tale.

It therefore follows that $\Perf(\rA^{(s)}(\ast \D))^{\varphi=1} \simeq \Perf(\rA_{\rN}^{(s)}(\ast \D_{\rN}))^{\varphi=1,h\mu_{\rN}^r}$ by finite \'{e}tale descent, noting that the $\mu_{\rN}^r$ action commutes with $\varphi$. The claim then follows as these equivalences are compatible with the maps in the \v{C}ech nerve.
\end{proof}

We note that concretely one may descend a crystal down by taking the projection into $\alpha=0$ components, which is indeed $\Gamma$-equivariant. We will use this shortly.

\begin{corollary}
The top horizontal functor in Theorem \ref{thm:open_stratum} is an equivalence.
\end{corollary}

\begin{proof}
By Proposition \ref{prop:open_fullfaithful}, it only remains to prove essential surjectivity. Since $\X_\eta$ is quasicompact and smooth, any $\mathbb{L}\in \D^{(b)}_{\mathrm{lisse}}((\X_\eta\setminus\D_\eta)_{\proet},\Z_p)$ is actually globally bounded rather than simply locally bounded. Thus for the purpose of showing essential surjectivity, since $\mathbb{L}\in\D^{(b)}_{\mathrm{lisse}}((\X_\eta\setminus\D_\eta)_{\proet},\Z_p)^\heart$ generate bounded objects under shifts and fibers and the functor is exact we only need essential surjectivity on the heart. Moreover, by strict \'{e}tale descent we may also reduce to testing the equivalence when $(\X,\mathcal{M})$ is a framed log affine.

By the previous discussion, it will then suffice to directly construct a preimage 
\[\mathcal{E} \in \Perf((\X,\mathcal{M})_\Prism,\Laur(\ast \D))^{\varphi=1,\op}\] 
so for such an $\mathbb{L}$ we get $\T_\et(\mathcal{E})=j_! \mathbb{L}$. By \cite[Corollary 6.3.4]{DiaoLanLiuZhuLogAdic}, $\mathbb{L}^{\mathrm{ext}}:=j_{\mathrm{k\acute{e}t},*}\mathbb{L}$ is a Kummer-\'{e}tale $\Z_p$-local system on the log adic generic fiber. We want to argue there is a single prime-to-$p$ Kummer cover $f_{\rN}:(\X_{\rN},\mathcal{M}_{\rN})\to(\X,\mathcal{M})$ as in Proposition \ref{prop:open_strat_Kum_desc} such that $f_{\rN}^*\mathbb{L}^{\mathrm{ext}}$ is $p$-Kummer. Without loss of generality we assume $(\Spf \R,\mathcal{M})_{\eta}$ is connected, as we can take the least common multiple of such $\rN$ on each component. By \cite[Definition 4]{KoshikawaYaoCorrigendum} we want to show that if $\overline{x}$ is a log geometric point of $(\Spf \R,\mathcal{M})_\eta$ then
\[\pi_1^{\mathrm{k}\et}((\Spf \R,\mathcal{M})_{\eta,\overline{x}}, \overline{x}) \xlongrightarrow{\rho} \mathrm{Aut}_{\Z_p}(\mathbb{L}^{\mathrm{ext}}_{\overline{x}})\]
has pro-$p$ image after pullback along $f_{\rN}$ for some $\rN$. On the left hand side, we use the pullback log structure and $\rho$ denotes the associated representation to $\mathbb{L}^{\mathrm{ext}}$. The key observation is that the kernel of $\mathrm{Aut}_{\Z_p}(\mathbb{L}^{\mathrm{ext}}_{\overline{x}}) \to \GL_n(\F_p)$ is pro-$p$, where $n$ is the rank of the mod $p$ reduction of $\mathbb{L}$ (which makes sense by connectedness). Thus to obtain a pro-$p$ image we can take $\rN$ to be the prime-to-$p$ part of $|\GL_n(\F_p)|$, as the effect of pullback along $f_{\rN}$ takes $\rN$th powers in the image (see the proof of Proposition 5 in \cite{KoshikawaYaoCorrigendum}).

Thus we may assume $f_{\rN}^*\mathbb{L}^{\mathrm{ext}}$ is $p$-Kummer. By the contravariant version of \cite[Theorem 7]{KoshikawaYaoCorrigendum} (see Remark \ref{rem:katz_contra}) we have an equivalence
\[\Perf((\X_{\rN},\mathcal{M}_{\rN})_\Prism,\Laur)^{\varphi=1,\op} \simeq \D^{(b),\mathrm{pKum}}_{\mathrm{lisse}}((\X_{\rN},\mathcal{M}_{\rN})_\eta,\Z_p),\] 
therefore $f_{\rN}^*\mathbb{L}^{\mathrm{ext}}$ has a corresponding perfect Laurent $\rF$-crystal $\mathcal{E}$ on $(\X_{\rN},\mathcal{M}_{\rN})_\Prism$, satisfying $\T_\et(\mathcal{E})\simeq f_{\rN}^*\mathbb{L}^{\mathrm{ext}}$.

Let $j_{\rN}:\X_{\rN,\eta}\setminus\D_{\rN,\eta}\hookrightarrow(\X_{\rN},\mathcal{M}_{\rN})_\eta$ denote the inclusion. Since $\T_\et(\Laur(\ast \D_{\rN}))\simeq j_{\rN,!}\underline{\Z_p}$ using Lemma \ref{lem:factorization} and $\T_\et$ is symmetric monoidal (from Remark \ref{rem:Tet_sym_mon}) we deduce
\[\T_\et(\mathcal{E}\otimes\Laur(\ast \D_{\rN}))\simeq j_{\rN,!}f_{\rN}^*\mathbb{L}.\]
Now since the prime-to-$p$ Kummer cover we use is a $\mu_{\rN}^r$-torsor, it follows that $f_{\rN}^*\mathbb{L}^{\mathrm{ext}}$ is equipped with a $\mu_{\rN}^r$ action. Translating this across the (contravariant version of) the Koshikawa-Yao correspondence from \cite{KoshikawaYaoCorrigendum}, we get a $\mu_{\rN}^r$ action on the corresponding Laurent $\rF$-crystal $\mathcal{E}$, and therefore such an action on $\mathcal{E}\otimes\Laur(\ast \D_{\rN})$ (noting $\Laur(\ast \D_{\rN})$ has its own $\mu_{\rN}^r$-action). Finally, Proposition \ref{prop:open_strat_Kum_desc} gives us a crystal downstairs on $(\Spf \R,\mathcal{M})_{\Prism}$ which realizes $j_!\mathbb{L}$. We may check this by verifying that before descending the $\mu_{\rN}^r$ actions correspond on each side of the correspondence.
\end{proof}

\subsection{Proof of the main theorem}

By this point, we have seen that both the source and target of 
\[\T_\et: \D_{\fgu}((\X,\mathcal{M})_\Prism,\Laur)^\op \longrightarrow \D^{(b),\mathrm{pKum}}_{\cons,\D}((\X,\mathcal{M})_\eta,\Z_p)\]
satisfy descent for strict \'{e}tale maps. Therefore we will localize to the situation of a framed log affine for the majority of this section, using Lemma \ref{lem:strict_etale_local_framing}.

Additionally we now know that there is an equivalence on the open stratum, and that the functor lands in the correct target by Proposition \ref{prop:T_et_cons}. The strategy we will carry out here is to match thick generators of both sides, using the following standard lemma.

\begin{lemma} \label{lem:thick_equiv_criterion}
Let $\rF: \mathcal{C} \longrightarrow \mathcal{D}$ be an exact functor between idempotent complete stable $\infty$-categories which are thickly generated by objects $\{\rC_{i}\}_{i\in \I}$ and $\{\D_i\}_{i\in \I}$, such that $\rF(\rC_i)=\D_i$ and the morphism 
\[\RHom(\rC_i,\rC_j) \longrightarrow \RHom(\D_i,\D_j)\]  
induced by $\rF$ is an equivalence for any $i,j\in \I$. Then $\rF$ is an equivalence.
\end{lemma}

\begin{proof}
The objects $\rC\in \mathcal{C}$ such that $\RHom(\rC_i,\rC) \longrightarrow \RHom(\D_i,\rF(\rC))$ is an equivalence form a thick subcategory containing the $\rC_i$, hence all of $\mathcal{C}$. Similarly we can run this argument in the other variable, deducing $\rF$ is fully faithful. Showing essential surjectivity is then straightforward, as the essential image is a thick subcategory containing the $\D_i$ hence all of $\mathcal{D}$. Here we have used idempotent completeness of $\mathcal{C}$ to ensure the essential image is closed under retracts.
\end{proof}

In our situation with $\T_\et$ on a framed log affine, it is easy to see all the hypotheses on the source and target categories are satisfied. Idempotent completeness is clear on the source (which is actually all we used in the proof!) by the description in Theorem \ref{thm:qdR_crys_Gamma_mod}. We can also see it directly, because $\D_{\Box}((\Spf \R,\mathcal{M})_\Prism,\Laur)^{\varphi=1}$ is idempotent complete and the finitely generated unit definition is, by construction, closed under retracts as we ask for membership in $\mathrm{Thick}(\rA[1/\I]_p^{\wedge}\langle \rF \rangle)$ for every log prism.

The following toy model explains a little more about the architecture of the proof, using Lemma \ref{lem:thick_equiv_criterion}.

\textbf{A toy model.} Consider the disk $\mathbb{D}_{\Q_p(\zeta_p)}$ punctured at the origin; we use the coordinate $\T$ so $\mathbb{D}=\Spa \Q_p(\zeta_p)\langle \T \rangle$, equipped with the log structure $\mathcal{M}_{\V(\T)}$. We use the usual formal model $(\Spf \Z_p[\zeta_p]\langle \T \rangle, \mathcal{M}_{\V(\T)})$.

Let $j:\mathbb{D}^\times\hookrightarrow\mathbb{D}$ and $i:\V(\T)\hookrightarrow\mathbb{D}$ denote the two strata. If we set $\rA=\Z_p\langle \T \rangle \laur{q-1}_p^{\wedge}$ and $\Gamma=(1+p\Z_p)\ltimes \Z_p(1)$, we may summarize the setup with the following diagram:
\begin{display}[column sep=large,row sep=large]
{\left\{\substack{\text{open stratum unit}\\\text{Laurent crystals}}\right\}^\op} \ar[hook]{r} \ar{d}{\T_\et} & {\left\{\substack{\text{fgu Laurent crystals}\\\text{ on }  (\Spf \Z_p[\zeta_p]\langle \T \rangle, \mathcal{M}_{\V(\T)})_\Prism}\right\}^\op} \ar{d}{\T_\et} & {\left(\D_{\fgu}(\rA/(\T))^{h\Gamma}\right)^\op} \ar[dashed,hook',swap]{l}{(i_+(-))[-1]} \ar{d}{\T_\et} \\
{\left\{\substack{p\text{-Kummer local systems}\\\text{on }\mathbb{D}^\times}\right\}} \ar[hook]{r}{j_!} & {\left\{\substack{p\text{-Kummer constructible sheaves}\\\text{on }(\mathbb{D},\mathcal{M}_{\V(\T)})}\right\}} & {\left\{\substack{p\text{-Kummer local systems}\\\text{ on }(\V(\T),i^*\mathcal{M}_{\V(\T)})}\right\}} \ar[hook',swap]{l}{i_*}.
\end{display}
The left vertical arrow is an equivalence by Theorem \ref{thm:open_stratum}, and tensoring with $\Laur(\ast\V(\T))$ identifies with $j_!$ under $\T_\et$. The dashed arrow in the diagram is supplied by Lemma \ref{lem:kashiwara} (Kashiwara's lemma).

We will ultimately also end up checking that the right vertical arrow is also an equivalence. In this particular case, it is substantially easier than the open stratum result because it is a hollow log point: one can run the usual proof of the $(\varphi,\Gamma)$-module equivalence for a point, except with the perfectoid cover $(\Spf \Z_p\langle \zeta_{p^\infty}\rangle)_\eta$ we must equip it with the hollow log structure $\N[1/p] \to 0$ to be a pro-Kummer-\'{e}tale $\Gamma$ torsor over a hollow log point. Running a tilting and deperfection argument (noting the analogue of the $p$-Kummer criterion using Lemma \ref{lem:pKum_framed_criterion} holds here) one obtains the equivalence.

One must then show that on either the top or bottom row, objects coming from the hollow point or open stratum thickly generate. For the bottom row the strategy is to use localization triangles
\[j_! j^* \mathscr{F} \longrightarrow \mathscr{F} \longrightarrow i_* i^* \mathscr{F}\]
and on the top row we will eventually do the same. Thus to use Lemma \ref{lem:thick_equiv_criterion}, we only need the $\RHom$ equivalences. There are no maps from sheaves supported on the open stratum to the closed stratum, and analogously on the finitely generated unit module side, so apart from the equivalences we already know the only point to check is whether or not we have full faithfulness for maps
\[\mathcal{E}_{\mathbb{D}^\times} \longrightarrow \mathcal{E}'_{\V(\T)},\]
with the arrow reversed after applying $\T_\et$. Here the source is an open stratum Laurent crystal and the target lies in $\D_{\fgu}(\rA/(\T))^{h\Gamma}$.

Put $j_! \mathbb{L}=\T_\et(\mathcal{E}_{\mathbb{D}^\times})$ and $i_*\mathbb{L}'=\T_\et(\mathcal{E}'_{\V(\T)})$.
We will show that we may assume without loss of generality $\mathbb{L}$ extends across the origin to the $p$-Kummer derived local system $\R j_*\mathbb{L}$: after replacing the disk by a prime-to-$p$ Kummer cover $\mathbb{D}_{\rN}$ this will be true, and ultimately we can descend the full faithfulness conclusion from this case.

Thus we assume there is 
\[\mathcal{E}_{\mathbb{D}}\in \Perf((\Spf \Z_p[\zeta_p]\langle \T \rangle,\mathcal{M})_\Prism,\Laur)^{\varphi=1}\]
such that $\T_\et(\mathcal{E}_{\mathbb{D}}) \simeq \R j_*\mathbb{L}$ and $\mathcal{E}_{\mathbb{D}^{\times}} \simeq \mathcal{E}_{\mathbb{D}}(*\V(\T))$. On the Kummer-\'{e}tale side, one has
\[\RHom(i_* \mathbb{L}', \R j_* \mathbb{L}) = 0.\]
Indeed, the adjunction between $j^*$ and $\R j_*$ identifies this with $\RHom(j^* i_* \mathbb{L}', \mathbb{L}) = 0$. Applying $\RHom(i_*\mathbb{L}',-)$ to the triangle
\[j_! \mathbb{L} \longrightarrow \R j_* \mathbb{L} \longrightarrow i_* i^* \R j_* \mathbb{L}\]
and using the adjunction between $i^*$ and $i_*$ gives $\RHom(i_* \mathbb{L}', j_! \mathbb{L}) \simeq \RHom(\mathbb{L}', i^* \R j_* \mathbb{L})[-1]$. This now consists of sheaves supported at the origin $\V(\T)\subset \mathbb{D}$, so producing an analogous isomorphism on the crystal side would reduce the needed result to a hollow point where the claim is known. The general proof of the theorem will proceed along similar lines to this toy model, but when there are multiple irreducible components the combinatorics is more complicated.

\textbf{The general proof}. We now observe there is a natural class of thick generators for the target category that we already know. In the situation of a framed log affine $(\Spf \R,\mathcal{M}) \xlongrightarrow{\sq} \mathbb{T}_r^n$ we have log strata for $\Lambda \subseteq \{1,\ldots,r\}$
\[(\Spf \R,\mathcal{M})_{\Lambda,\eta} = \bigcap_{\lambda \in \Lambda} \D_{\lambda,\eta} \setminus \bigcup_{\lambda \notin \Lambda} \D_{\lambda,\eta}\]
where the open stratum corresponds to the empty set. Here $\D_{\lambda}$ is the zero locus of the \'{e}tale lift of $\T_\lambda$, or one of the irreducible components of the divisor.

\begin{lemma} \label{lem:ket_thick_gens}
Let $(\Spf \R,\mathcal{M})$ be a framed log affine. The category $\D^{(b),\mathrm{pKum}}_{\cons,\D}((\Spf \R,\mathcal{M})_\eta,\Z_p)$ is thickly generated by $j_{\Lambda,!}\mathbb{L}$ for $\mathbb{L}$ a $p$-Kummer sheaf in $\D^{(b)}_{\mathrm{lisse}}((\Spf \R,\mathcal{M})_{\Lambda,\eta},\Z_p)$, where
\[j_{\Lambda}: (\Spf \R,\mathcal{M})_{\Lambda,\eta} \hookrightarrow (\X,\mathcal{M})_\eta\]
is the inclusion of the stratum for $\Lambda \subseteq \{1,\ldots,r\}$.
\end{lemma}

\begin{proof}
Inductively apply the localization triangle
\[j_! j^* \mathscr{F} \longrightarrow \mathscr{F} \longrightarrow i_* i^* \mathscr{F}\]
for an open/closed stratification ($j$ being the open inclusion, $i$ the closed inclusion). Repeating this produces $\mathscr{F}$ in terms of the thick generators. It is worth pointing out that in this procedure $p$-Kummer sheaves are closed under extensions, since invariants for $\I_{\overline{x}}^{\log,(p)'}$ are exact (due to this being the prime-to-$p$ part).
\end{proof}

Thus far, in $\D_{\fgu}((\Spf \R,\mathcal{M})_\Prism,\Laur)$ we have obtained only one class of these generators for the open stratum. In the situation of a framed log affine, we will now introduce the categories needed to handle other strata.

\phantomsection
\label{notation:main_thm}
\textbf{Notation for this section}. Consider a framed log affine $(\Spf \R,\mathcal{M})$ with framing $\sq$ as in Definition \ref{def:fr_log_aff}, and put $\rA=\R_0^{\sq}\laur{q-1}^{\wedge}_p$. For $\Lambda\subseteq\{1,\ldots,r\}$, put
\[\rA_\Lambda=\rA/(\T_i)_{i\in \Lambda}[1/\T_j]^{\wedge}_{p,j \notin \Lambda}.\]
We introduce the category
\[\D_{\fgu}(\rA_{\Lambda})^{h\Gamma},\]
which actually by the Fitting ideal argument in Proposition \ref{prop:perfectoid_constructibility} is simply $\Perf(\rA_{\Lambda})^{\varphi=1,h\Gamma}$.

Ultimately we will argue these categories correspond to $j_{\Lambda,!}\D^{(b),\mathrm{pKum}}_{\mathrm{lisse}}((\Spf \R,\mathcal{M})_{\Lambda,\eta},\Z_p)$, up to a shift which only matters for perversity. Note that unlike the open stratum, we do actually need to pay attention to the $p$-Kummer condition on these more general strata. They are equipped with their pullback log structures, so there will be a hollow component to the log structure on which the $p$-Kummer condition is not vacuous. To properly embed these categories as full subcategories into $\D_{\fgu}(\rA)^{h\Gamma}$ we need a variant of Kashiwara's lemma.

We write $\D_{\fgu}((\X,\mathcal{M})_\Prism,\Laur)_{\rZ}$ for the full subcategory of finitely generated unit Laurent crystals supported on $\rZ$. Let us explain the general globalized notion of this when $\D=\bigcup_i \D_i$ is a simple normal crossings divisor with irreducible components $\D_i$. For any nonempty intersection $\rZ=\D_1\cap\cdots\cap\D_c$ of distinct irreducible components, define
\[\D_{\fgu}((\X,\mathcal{M})_\Prism,\Laur)_{\rZ}:= \bigcap_{a=1}^c \ker\left(-\otimes^{\L,\Box}\Laur(\ast\D_a)\right).\]
Similarly one may make a definition for connected components. This covers all $\rZ$ which are smooth closed log strata of the stratification induced by $\D$.

Under Theorem \ref{thm:qdR_crys_Gamma_mod}, we may translate the previous definition on a framed log affine to obtain the analogous category $\D_{\fgu}(\rA)^{h\Gamma}_{\rJ_{\Lambda}}$ of finitely generated unit modules supported on $\T_i$ for $i\in \Lambda\subseteq \{1,\ldots,r\}$, where $\rJ_{\Lambda}=(\T_i)_{i\in \Lambda}$. We then can prove the following:

\begin{lemma}[Kashiwara's lemma] \label{lem:kashiwara}
Let $(\Spf \R,\mathcal{M})$ be a framed log affine, and as usual put $\rA=\R_0^{\sq}\laur{q-1}^{\wedge}_p$ and fix $\Lambda \subset \{1,\ldots,r\}$. We put $\rJ_{\Lambda} = (\T_i)_{i\in \Lambda}$.

Then there is an equivalence
\[i^!:=\left((\rA/\rJ_{\Lambda})\otimes^{\L}_{\rA}-\right)[-|\Lambda|]: \D_{\fgu}(\rA)^{h\Gamma}_{\rJ_{\Lambda}}\xlongrightarrow{\sim} \D_{\fgu}(\rA/\rJ_{\Lambda})^{h\Gamma}.\]
\end{lemma}

\begin{proof}
Put $\overline{\rA}=\rA/p$ and $\overline{\rJ}_{\Lambda}=\rJ_{\Lambda}\overline{\rA}$. The categories and functor above satisfy the hypotheses of Lemma \ref{lem:categorical_nakayama} and we may check the support modulo $p$. It is therefore enough to prove that
\[i^!:=\left((\overline{\rA}/\overline{\rJ}_{\Lambda})\otimes^{\L,\Box}_{\overline{\rA}}-\right)[-|\Lambda|]:\D_{\fgu}(\overline{\rA})_{\overline{\rJ}_{\Lambda}}^{h\Gamma}\longrightarrow\D_{\fgu}(\overline{\rA}/\overline{\rJ}_{\Lambda})^{h\Gamma}\]
is an equivalence. The functor is well-defined: one can use Lemma \ref{lem:derived_stratum_fgu}, and in particular $i^!$ is $\Gamma$-equivariant.

In this case, the functor $i^!$ is an equivalence \emph{before} we apply $(-)^{h\Gamma}$. One may use the proof of \cite[Corollary 5.11.3]{EmertonKisin2} in the present regular $\rF$-finite situation. Thus one obtains the adjoint equivalences
\[i_+:\D_{\fgu}(\overline{\rA}/\overline{\rJ}_{\Lambda}) \rightleftarrows \D_{\fgu}(\overline{\rA})_{\overline{\rJ}_{\Lambda}}: i^!.\]
This survives after applying a cobar construction $(-)^{h\Gamma}$ to both sides, because evaluation on the degree zero term is conservative. Thus one may test $i^! i_+ \simeq \id$ and $i_+ i^! \simeq \id$ after evaluation on the degree zero term, at which point it reduces to the ordinary Kashiwara lemma for finitely generated unit modules.
\end{proof}

We use $\rJ_{\Lambda}$ as in Lemma \ref{lem:kashiwara}. Applying the \emph{inverse} $i_+$ in Kashiwara's lemma yields a more general fully faithful functor
\[\D_{\fgu}(\rA/\rJ_{\Lambda})^{h\Gamma,\op} \xlongrightarrow{i_+} \D_{\fgu}(\rA)^{h\Gamma,\op},\]
and there is an obvious inclusion, corresponding to $!$-extension, 
\[\Perf(\rA_{\Lambda})^{h\Gamma,\varphi=1,\op} \hookrightarrow \D_{\fgu}(\rA/\rJ_{\Lambda})^{h\Gamma,\op}.\]
We denote the composite functor, again fully faithful, by $j_{\Lambda,!}$.

The category $\Perf(\rA_{\Lambda})^{h\Gamma,\varphi=1,\op}$ is also equipped with a functor
\[\T_\et: \Perf(\rA_{\Lambda})^{h\Gamma,\varphi=1,\op} \longrightarrow \D^{(b),\mathrm{pKum}}_{\mathrm{lisse}}((\Spf \R,\mathcal{M})_{\Lambda,\eta},\Z_p)\]
arising in the same way as Lemma \ref{lem:factorization}.

\begin{corollary} \label{cor:thick_gen_embedding}
The following diagram commutes:
\begin{display}
\Perf(\rA_{\Lambda})^{h\Gamma,\varphi=1,\op} \ar{d}{\T_\et} \ar{r}{j_{\Lambda,!}} & \D_{\fgu}(\rA)^{h\Gamma,\op} \ar{d}{\T_\et} \\
\D^{(b),\mathrm{pKum}}_{\mathrm{lisse}}((\Spf \R,\mathcal{M})_{\Lambda,\eta},\Z_p) \ar{r}{j_{\Lambda,!}[-|\Lambda|]} & \D^{(b),\mathrm{pKum}}_{\cons,\D}((\Spf \R,\mathcal{M})_\eta,\Z_p)
\end{display}
On the right $\T_\et$ is passed through the equivalence of Theorem \ref{thm:qdR_crys_Gamma_mod}, which may then be explicitly computed via Lemma \ref{lem:factorization}.
\end{corollary}

\begin{proof}
This is more or less clear by construction from Lemma \ref{lem:factorization}, by expanding what the square means:
\begin{display}
\Perf(\rA_{\Lambda})^{h\Gamma,\varphi=1,\op} \ar[out=180,in=180,looseness=2,start anchor=west,end anchor=west,swap]{dd}{\T_\et} \ar{d}{\T_\et^{\perf}} \ar{r}{j_{\Lambda,!}} & \D_{\fgu}(\rA)^{h\Gamma,\op} \ar{d}{\T_\et^{\perf}} \ar[out=0,in=0,looseness=2,start anchor=east,end anchor=east]{dd}{\T_\et} \\ 
\D^{(b)}_{\mathrm{lisse}}(\Spec (\R_0^{\sq})_{\infty,\Lambda,\eta}^{\flat},\Z_p)^{h\Gamma} \ar{d}{\mathrm{an}_\Gamma^*} \ar{r}{j_{\Lambda,!}[-|\Lambda|]} & \D^{(b)}_{\cons, \D_\infty}(\Spec (\R_0^{\sq})_{\infty,\eta}^{\flat},\Z_p)^{h\Gamma} \ar{d}{\mathrm{an}_\Gamma^*} \\
\D^{(b),\mathrm{pKum}}_{\mathrm{lisse}}((\Spf \R,\mathcal{M})_{\Lambda,\eta},\Z_p) \ar{r}{j_{\Lambda,!}[-|\Lambda|]} & \D^{(b),\mathrm{pKum}}_{\cons,\D}((\Spf \R,\mathcal{M})_\eta,\Z_p).
\end{display}
Here $\Spec (\R_0^{\sq})_{\infty,\Lambda,\eta}^{\flat}$ similarly means that we take the reduced ring where $\T_i^{\flat} = 0$ for $i\in \Lambda$, and invert them for $i \notin \Lambda$ (note this is not a perfectoid ring). By abuse of notation we base change to the perfection first and then apply $\T_\et^\perf$ -- more precisely, in Lemma \ref{lem:factorization} we have denoted the composite of the first two arrows by simply $\T_\et^{\perf}$.

The bottom square is commutative, so we reduce to checking commutativity of the top square. This is clear for the first component of
\[j_{\Lambda,!}: \Perf(\rA_{\Lambda})^{h\Gamma,\varphi=1,\op} \hookrightarrow \D_{\fgu}(\rA/\rJ_{\Lambda})^{h\Gamma,\op} \xlongrightarrow{i_+} \D_{\fgu}(\rA)^{h\Gamma,\op}\]
so we reduce to studying compatibility for $i_+$.

One can further reduce to checking this modulo $p$ by derived Nakayama, as there is already a morphism we can test is an equivalence to witness commutativity. Indeed for the closed part $i_+$ of our construction of $j_{\Lambda,!}$ there is a map
\[\T^{\perf}_\et(i_+ \M) \longrightarrow i_* i^* \T^{\perf}_\et(i_+ \M) \longrightarrow i_* \T^{\perf}_\et(\rA/\rJ_{\Lambda} \otimes^{\L}_{\rA} i_+ \M) = i_* \T^{\perf}_\et(\M)[-|\Lambda|].\]
This is then not too difficult, because one can see that for the solution functor of \cite{EmertonKisin2} one has for $\M \in \Perf(\overline{\rA}_{\Lambda})^{h\Gamma,\varphi=1}$ that
\[\mathrm{RSol}(i_+ \M) = i_* \mathrm{RSol}(\M)[-|\Lambda|],\]
compatibly with the $\Gamma$-action. On the left we use the solution functor for $\D_{\fgu}(\overline{\rA})$ and on the right for $\D_{\fgu}(\overline{\rA}/\rJ_{\Lambda})^{h\Gamma}$. By adjunction one can reduce the check to one for $i^!$ (revealing the source of the shift).

For the comparison we actually needed to work with $\T_\et^{\perf}$, but base changing to $\Spec (\R_0^{\sq})_{\infty,\eta}^{\flat}$ preserves this identity.
\end{proof}

In fact, the left vertical arrow in the diagram of Corollary \ref{cor:thick_gen_embedding} is an equivalence. This can be proved directly.

\begin{proposition} \label{prop:thick_gen_equiv}
The functor
\[\T_\et: \Perf(\rA_{\Lambda})^{h\Gamma,\varphi=1,\op} \longrightarrow \D^{(b),\mathrm{pKum}}_{\mathrm{lisse}}((\Spf \R,\mathcal{M})_{\Lambda,\eta},\Z_p)\]
is an equivalence.
\end{proposition}

\begin{proof}
We first check the case $|\Lambda|=r$, which we will explain in full detail. This case is a bit easier because the corresponding stratum on the log adic generic fiber is the adic generic fiber of an ordinary framed affine $\rZ=\bigcap_{i=1}^r \D_i$ over a rank $n-r$ torus with a rank $r$ hollow log structure. What will matter in the argument is that we have not inverted any of the coordinates $\T_i$.

Equipping the target torus with the hollow chart $\N^r\to0$, we get an induced strict \'{e}tale framing
\[\sq:(\rZ,\mathcal{M}|_{\rZ})\longrightarrow\left(\Spf \W(k)[\zeta_p]\langle \T_{r+1}^{\pm1},\ldots,\T_n^{\pm1}\rangle,(\N^r\to0)^a\right).\]
The torus has a natural perfectoid cover $\Spf \W(k)^{\cyc} \langle \T_{r+1}^{\pm 1/p^\infty}, \ldots, \T_n^{\pm 1/p^\infty} \rangle$ with hollow log structure $\N[1/p]^r\to 0$, so using $\sq$ we obtain an associated perfectoid cover $(\R_0^{\sq})_{\infty,\rZ}[1/p]$ again with the chart $\N[1/p]^r\to 0$; we let $(\rA/\rJ_{\Lambda})_\perf$ denote $\Ainf((\R_0^{\sq})_{\infty,\rZ}[1/p])$. Katz's theorem, tilting, and pro-Kummer-\'{e}tale descent for the same full group $\Gamma$ give
\[\left(\Perf((\rA/\rJ_{\Lambda})_\perf)^{\varphi=1,h\Gamma}\right)^\op \xlongrightarrow{\sim} \D^{(b),\mathrm{pKum}}_{\mathrm{lisse}}\left(\left(\rZ,\N^r\right)_\eta,\Z_p\right) \xhookrightarrow{i_*} \D^{(b)}_{\cons,\D}((\Spf \R,\mathcal{M})_\eta, \Z_p).\]
In this case the finitely generated unit condition is equivalent to being perfect with unit Frobenius by the proof of Proposition \ref{prop:perfectoid_constructibility}. Moreover,
\[\D_{\fgu}(\rA/\rJ_{\Lambda})^{h\Gamma}\simeq\Perf((\rA/\rJ_{\Lambda})_\perf)^{\varphi=1,h\Gamma}\]
by \cite[Proposition 3.6]{BhattScholzePrismaticFCrystals} and Lemma \ref{lem:categorical_nakayama}; thus the desired claim follows. We note Lemma \ref{lem:factorization} identifies this construction with $\T_\et$.

In the remaining cases where $|\Lambda|$ is smaller we must involve part of the open stratum theorem, but fortunately the hard work has already been done. The result we need is that
\[\Perf(\rA_{\Lambda})^{h\Gamma,\varphi=1,\op} \xlongrightarrow{\sim} \D^{(b),\mathrm{pKum}}_{\mathrm{lisse}}((\Spf \R,\mathcal{M})_{\Lambda,\eta},\Z_p),\]
where the target has its natural hollow pullback log structure. The equivalence on the open stratum tells us that there is an equivalence
\[\Perf(\rA_{\Lambda})^{h\Gamma',\varphi=1,\op} \xlongrightarrow{\sim} \D_{\mathrm{lisse}}^{(b)}((\Spf \R)_{\Lambda,\eta},\Z_p)\]
where on the right we use the trivial log structure and $\Gamma'$ is the quotient $(1+p\Z_p) \ltimes \Z_p(1)^{n-|\Lambda|}$. We can check the equivalence modulo $p$ by Lemma \ref{lem:categorical_nakayama}. 

Full faithfulness follows using the same argument as Theorem \ref{thm:open_stratum}. The only nontrivial point is to note that adding hollow directions to Proposition \ref{prop:open_decompletion_hom} does not affect the reduction argument. Moreover in the helper lemmas, since we work modulo $p$ by increasing the parameter $m$ we can make the hollow part act trivially after shrinking it. This follows because the action of $\Z_p(1)^{|\Lambda|}$ factors through a finite $p$-group modulo $p$. This follows on $\Vect(\overline{\rA}_{\Lambda})^{h\Gamma,\varphi=1,\op}$ by using Katz's correspondence from \cite[Proposition 3.4]{BhattScholzePrismaticFCrystals}, so it amounts to the result that an $\F_p$-local system on a regular affine $\F_p$-scheme has a finite dimensional vector space of endomorphisms (and $\Z_p(1)^{|\Lambda|}$ acts by such endomorphisms). The hollow directions will then not affect the remainder of the argument, which can be carried out in the same way.

With full faithfulness in hand, we will next show essential surjectivity. We first observe that there is a commutative diagram
\begin{display}
\Perf((\rA/(\T_i)_{i\in \Lambda}))^{h\Gamma,\varphi=1,\op} \ar{d} \ar{r}{\T_\et} & \D^{(b),\mathrm{pKum}}_{\mathrm{lisse}}(\rZ_{\Lambda,\eta},\Z_p) \ar{d}{j_{\Lambda}^\ast} \\
\Perf(\rA_{\Lambda})^{h\Gamma,\varphi=1,\op} \ar[hook]{r}{\T_\et} & \D^{(b),\mathrm{pKum}}_{\mathrm{lisse}}((\Spf \R,\mathcal{M})_{\Lambda,\eta},\Z_p)
\end{display}
To see the top arrow is an equivalence, because we have not inverted any $\T_i$ the proof is the same argument as the $|\Lambda|=r$ case of Proposition \ref{prop:thick_gen_equiv} but with the log perfectoid cover of $\rZ_{\Lambda,\eta}$ being the zero locus of $\T_i^{1/p^\infty}$ for $i\in \Lambda$ in $(\X,\mathcal{M})_{\infty,\eta}$ with its pullback log structure.

To show essential surjectivity of the bottom arrow, as in the proof of Theorem \ref{thm:open_stratum} the first step is to check that there is some prime-to-$p$ Kummer cover $f: (\Spf \R_{\rN},\mathcal{M}_{\rN}) \to (\Spf \R,\mathcal{M})$ so that after pullback along this cover some $\mathbb{L} \in \D^{(b),\mathrm{pKum}}_{\mathrm{lisse}}((\Spf \R,\mathcal{M})_{\Lambda,\eta},\Z_p)$ has $f^*\mathbb{L}$ admit an extension to $\rZ_{\Lambda,\eta} \times_{(\Spf \R,\mathcal{M})_\eta} (\Spf \R_{\rN},\mathcal{M}_{\rN})_{\eta}$. Here we have abused notation slightly: we use $f^*$ just to mean the induced pullback functor on the log stratum for $\Lambda$.

The general situation we must now consider is a framed log affine $(\Spf \R,\mathcal{M})$ where we allow a hollow component $\N^{|\Lambda|}$ in the log structure $\mathcal{M}$ and denote by $(\Spf \R,\mathcal{M})_\eta^\circ$ the open stratum. It suffices to know the general claim that if $\mathbb{L}\in \D^{(b)}_{\mathrm{lisse}}((\Spf \R,\mathcal{M})_\eta^\circ, \Z_p)$ there is a prime-to-$p$ Kummer cover $f: (\Spf \R_{\rN},\mathcal{M}_{\rN}) \to (\Spf \R,\mathcal{M})$ so $f_\eta^* \mathbb{L}$ admits a $p$-Kummer extension to $(\Spf \R_{\rN},\mathcal{M}_{\rN})_\eta$. Upon keeping track of hollow log inertia, we see \cite[Corollary 6.3.4]{DiaoLanLiuZhuLogAdic} still holds in our setting. The derived version of this log purity result, namely the equivalence
\[\R j_*: \D^{(b)}_{\mathrm{lisse}}((\Spf \R,\mathcal{M})_\eta^\circ, \Z_p) \xlongrightarrow{\sim} \D^{(b)}_{\mathrm{lisse}}((\Spf \R,\mathcal{M})_\eta, \Z_p),\]
follows by applying Lemma \ref{lem:categorical_nakayama} to first reduce to $\F_p$ coefficients. For full faithfulness, as we did in Proposition \ref{prop:open_decompletion_hom} we can use the tensor/Hom adjunction and a Postnikov d\'{e}vissage to reduce to derived full faithfulness in the heart, which then follows from $\mathbb{L} \simeq \R j_* j^* \mathbb{L}$ for ordinary local systems. Essential surjectivity, again by a Postnikov d\'{e}vissage, also then follows from essential surjectivity for ordinary local systems. Given this equivalence, the claim is then clear from boundedness that a single prime-to-$p$ Kummer cover makes a derived lisse sheaf $p$-Kummer (using the same argument about the finite image of prime-to-$p$ log inertia on each cohomology sheaf as in the proof of Theorem \ref{thm:open_stratum}).

Just as in the proof of Theorem \ref{thm:open_stratum}, for essential surjectivity the final step is to descend. The analogue of Proposition \ref{prop:open_strat_Kum_desc} is actually not too difficult in this setting, because instead of the \v{C}ech nerve one can also run the same argument with the cobar complex $\rC(\Gamma^n, \rA_{\Lambda})$. The key point is again that for the prime-to-$p$ Kummer cover, the corresponding $\rA_{\Lambda,\rN}$ has a grading by $(\Z/\rN \Z)^{r-|\Lambda|}$ and this extends to $\rC(\Gamma^n, \rA_{\Lambda,\rN})$ compatibly with the maps in the cobar complex. Thus the claim follows.
\end{proof}

\begin{remark}
We can also argue essential surjectivity through a thick generation argument modulo $p$: it turns out the derived lisse sheaves with trivial hollow log inertia thickly generate modulo $p$, and through the existing equivalence in Theorem \ref{thm:open_stratum} we can produce these.
\end{remark}

In Lemma \ref{lem:thick_equiv_criterion}, we must now consider the different types of generators we need to compare. We also have not yet shown that $\Perf(\rA_{\Lambda})^{h\Gamma,\varphi=1}$, properly embedded via $j_{\Lambda,!}$, thickly generate $\D_{\fgu}(\rA)^{h\Gamma}$. We establish this first; the proof method is actually nearly identical to Lemma \ref{lem:ket_thick_gens}.

For a closed subset $i:\rZ_\Lambda\hookrightarrow\Spf\R$ cut out by $\rJ_{\Lambda}=(\T_i)_{i\in \Lambda}$ and $\M\in\D_{\fgu}(\rA/\rJ_{\Lambda})^{h\Gamma}$, the equivalence above in Lemma \ref{lem:kashiwara} gives
\[\T_\et(i_+\M)\simeq i_*\T_\et(\M)[-|\Lambda|].\]
It is also easy to give a formula for $i_+ i^! (-)$ even when the input is not supported on $\rZ$. In this situation local cohomology along $\rZ$ is computed by
\[\RGamma_{\rZ}(\M)\simeq \M\otimes^{\L}_{\rA}\bigotimes_{i\in\Lambda}^{\L} \left[\rA\longrightarrow(\rA[1/\T_i])^{\wedge}_p\right].\]
The complex $\RGamma_{\rZ}(\M)$ is equipped with an action by $\Gamma$, since each $\gamma\in\Gamma$ carries every $\T_i$ to a unit multiple of itself. Checking $\RGamma_{\rZ}(\M)$ is finitely generated unit is straightforward, and using Kashiwara's lemma we can check
\[\RGamma_{\rZ}(\M)\simeq i_+i^!\M.\]

\begin{lemma} \label{lem:crys_thick_gens}
Let $(\Spf \R,\mathcal{M})$ be a framed log affine. The category $\D_{\fgu}(\rA)^{h\Gamma}$ is thickly generated by $j_{\Lambda,!}\M_\Lambda$ for $\M_\Lambda\in \Perf(\rA_{\Lambda})^{h\Gamma,\varphi=1}$.
\end{lemma}

\begin{proof}
Similar to Lemma \ref{lem:ket_thick_gens}, we need to produce analogous triangles extracting open and closed parts. Let $\M \in \D_{\fgu}(\rA)^{h\Gamma}$.

For an irreducible component $\D_i$ in $\D$, there is an exact triangle
\[\RGamma_{\D_i}(\M) \longrightarrow \M \longrightarrow \M(\ast \D_i),\]
extracting $\M(\ast \D_i)$ supported away from $\D_i$ and $\RGamma_{\D_i}(\M)$ supported on $\D_i$ (we implicitly use Lemma \ref{lem:kashiwara} here to see that this part comes from $\D_{\fgu}(\rA/\T_i)^{h\Gamma}$). Repeating this for collections $\{\D_i\}_{i\in \Lambda}$ we can eventually show that $\M$ is in the thick subcategory generated by $j_{\Lambda,!}\M_\Lambda$ for $\M_\Lambda\in \Perf(\rA_{\Lambda})^{h\Gamma,\varphi=1}$.

Once we reach an individual stratum, a priori we only have $\M_\Lambda\in \D_{\fgu}(\rA_{\Lambda})^{h\Gamma}$. Fortunately the Fitting ideal argument in Proposition \ref{prop:perfectoid_constructibility} shows that these are necessarily perfect. Then we run the same argument: the $\Gamma$-invariant ideals are either zero or the unit ideal on connected components of the scheme. We know an object $\M \in \Mod_{\fgu}(\overline{\rA}_{\Lambda})^{h\Gamma}$ admits a $\Gamma$-equivariant minimal root $\beta: \M_0 \to \varphi^* \M_0$. Taking Fitting ideals, we see the source and target are vector bundles of the same rank and $\coker \beta$ is a vector bundle. Using ranks the root map is an equivalence, so all of the objects in the heart are vector bundles modulo $p$. By a Postnikov d\'{e}vissage, this proves perfectness modulo $p$. We can then test perfectness modulo $p$ since $\rA_{\Lambda}$ is still derived $p$-complete (as well as $\M_{\Lambda}$), so the desired result follows.
\end{proof}

We have now established everything in Lemma \ref{lem:thick_equiv_criterion} that we need to establish the main theorem by comparing thick generators, except the comparisons of $\RHom$s on generators. For these, the general strategy will be a reduction to one of the known equivalences in Proposition \ref{prop:thick_gen_equiv}. For this reduction, we will need $\RHom$ comparisons on both sides of the correspondence. We work out the Kummer-\'{e}tale side first, using the new notation
\[\rZ_{\Lambda,\eta} = \bigcap_{\lambda \in \Lambda} \D_{\lambda,\eta} \subset (\Spf \R,\mathcal{M})_{\eta}.\]

\begin{lemma} \label{lem:orthog_red_ket}
Consider a framed log affine $(\Spf \R,\mathcal{M})$, and choose (possibly empty) subsets $\Lambda \subseteq \Lambda'\subseteq \{1,\ldots,r\}$. Letting $\mathbb{L}'\in \D^{(b)}_{\mathrm{lisse}}((\Spf \R,\mathcal{M})_{\Lambda',\eta},\Z_p)$ and $\mathbb{L}$ be derived lisse on $\rZ_{\Lambda,\eta} \subset (\Spf \R,\mathcal{M})_\eta$, we have
\[\RHom_{\rZ_{\Lambda,\eta}}(j_{\Lambda',!} \mathbb{L}', j_{\Lambda,!} j_{\Lambda}^* \mathbb{L}) \simeq \RHom_{\rZ_{\Lambda',\eta}}(j_{\Lambda',!}\mathbb{L}', i_{\rZ_{\Lambda',\eta}}^* \mathbb{L})[-|\Lambda'\setminus \Lambda|].\]
Here $j_{\Lambda},j_{\Lambda'}$ are the inclusions of the strata $(\Spf \R,\mathcal{M})_{\Lambda,\eta}$ and $(\Spf \R,\mathcal{M})_{\Lambda',\eta}$ with their pullback log structures. On the right side we regard $j_{\Lambda'}$ as the inclusion of the open stratum in $\rZ_{\Lambda',\eta}$.
\end{lemma}

\begin{proof}
We may assume $\Lambda \subset \Lambda'$ is a strict inclusion, as the claim is trivial otherwise. We can reframe this by working in $\rZ_{\Lambda,\eta} \setminus \bigcup_{j \notin \Lambda'} \D_{j,\eta}$. The point of this restriction is that the inclusion
\[(\Spf \R,\mathcal{M})_{\Lambda',\eta} \xlongrightarrow{i_{\Lambda'}} \rZ_{\Lambda,\eta} \setminus \bigcup_{j \notin \Lambda'} \D_{j,\eta}\]
is a strict closed immersion, while the inclusion $j_\Lambda$ of $(\Spf \R,\mathcal{M})_{\Lambda,\eta}$ is a strict open immersion.

We make use of a resolution similar to the one appearing in \cite[\S 2]{yamashita}, namely that $j_{\Lambda,!} j_\Lambda^* \mathbb{L}$ can be resolved by
\[\mathbb{L} \longrightarrow \bigoplus_{i} \mathbb{L}|_{\D_{i,\eta}} \longrightarrow \bigoplus_{i < j} \mathbb{L}|_{\D_{i,\eta} \cap \D_{j,\eta}} \longrightarrow \cdots \longrightarrow i_{\Lambda',*} i_{\Lambda'}^* \mathbb{L},\]
starting in cohomological degree zero. The direct sums are indexed over members of $\Lambda' \setminus \Lambda$, and we also abuse notation slightly: when we write restriction to some intersection of $\D_{i,\eta}$, we mean pullback to the intersection followed by extension by zero to regard it as an object in the appropriate category. This resolution can be seen, for example, by checking stalks on both sides. 

We wish to compute $\RHom_{\rZ_{\Lambda,\eta}}(j_{\Lambda',!} \mathbb{L}', j_{\Lambda,!} j_{\Lambda}^* \mathbb{L})$, and upon using this resolution for the right term $j_{\Lambda,!} j_{\Lambda}^* \mathbb{L}$ we claim all terms except $i_{\Lambda',*} i_{\Lambda'}^* \mathbb{L}$ have no contribution. This would then give exactly the desired result, with the observation that $i_{\Lambda',*} i_{\Lambda'}^* \mathbb{L}$ appears last in the resolution of length $|\Lambda' \setminus \Lambda|$ accounting for the shift.

Vanishing of $\RHom(i_{\Lambda',*} \mathbb{L}', -)$ for other terms of the resolution can be checked modulo $p$ and after a Postnikov d\'{e}vissage to the heart, where we can use logarithmic purity as in \cite[Corollary 4.6.7]{DiaoLanLiuZhuLogAdic} to deduce it. The present hollow case can again be deduced from this by keeping track of hollow log inertia. Once logarithmic purity is established, the vanishing claim for all but the last term in the resolution of $j_{\Lambda,!} j_\Lambda^* \mathbb{L}$ follows immediately by adjunction. For example, with the first term in the resolution we compute by directly applying logarithmic purity that $\RHom(i_{\Lambda',*} \mathbb{L}', \mathbb{L}) = \RHom(i_{\Lambda',*} \mathbb{L}', \R j_{\Lambda,*} j_{\Lambda}^* \mathbb{L}) = 0$ by adjunction, since $j_{\Lambda}^* i_{\Lambda',*}=0$. The same follows for all terms in the resolution except the final one $i_{\Lambda',*} i_{\Lambda'}^* \mathbb{L}$.
\end{proof}

This reduces to an $\RHom$ between local systems on the stratum corresponding to $\Lambda'$, because now the inclusion of $(\X,\mathcal{M})_{\Lambda',\eta} \hookrightarrow \rZ_{\Lambda',\eta}$ is a strict open immersion. It is worth pointing out that the case handled above where $\Lambda \subseteq \Lambda'$ is the only interesting case, as in other cases the $\RHom$ is zero.

Recall we have already fully faithfully embedded $\Perf(\rA_{\Lambda})^{h\Gamma,\varphi=1}$ for all $\Lambda$ into $\D_{\fgu}(\rA)^{h\Gamma}$ using the construction preceding Corollary \ref{cor:thick_gen_embedding}. The next lemma repeats Lemma \ref{lem:orthog_red_ket} but on the finitely generated unit module side. We write $\D_{\Lambda}$ for $\bigcup_{i \notin \Lambda} \D_i$.

\begin{lemma} \label{lem:orthog_red}
Consider a framed log affine $(\Spf \R,\mathcal{M})$, and choose (possibly empty) subsets $\Lambda \subseteq \Lambda'\subseteq \{1,\ldots,r\}$. Let $\M_{\Lambda} \in \Perf(\rA/(\T_i)_{i\in \Lambda})^{h\Gamma,\varphi=1}$, and let $\M_{\Lambda'} \in \Perf(\rA_{\Lambda'})^{h\Gamma,\varphi=1}$, and regard both as objects in $\D_{\fgu}(\rA/(\T_i)_{i\in \Lambda})^{h\Gamma}$ using Kashiwara's lemma (Lemma \ref{lem:kashiwara}). Then 
\[\RHom(\M_{\Lambda}(\ast \D_{\Lambda}), \M_{\Lambda'}) \simeq \RHom(\RGamma_{\bigcap_{j \in \Lambda'\setminus \Lambda} \D_j} \M_{\Lambda}, \M_{\Lambda'})[-|\Lambda'\setminus \Lambda|].\]
\end{lemma}

\begin{proof}
The argument is more or less identical, but translated to the finitely generated unit side. We no longer have access to logarithmic purity so we will need to prove the required vanishing in a different way. 

For convenience of notation put $|\Lambda'\setminus \Lambda|=c$. Without affecting any of the results, we may instead work over $\rA_{\Lambda,\Lambda'} := \rA/(\T_i)_{i\in \Lambda}[1/\T_j]_{p,j \notin \Lambda'}^{\wedge}$, since $\M_{\Lambda'}$ is already localized at these variables. This is the analogue in the previous proof of passing to $\rZ_{\Lambda,\eta} \setminus \bigcup_{j \notin \Lambda'} \D_{j,\eta}$.

By tensoring with the ($\Gamma$-equivariant) localization triangles
\[\RGamma_{\D_i}(\rA_{\Lambda,\Lambda'}) \longrightarrow \rA_{\Lambda,\Lambda'} \longrightarrow \rA_{\Lambda,\Lambda'}(\ast \D_i)\]
for $i\in \Lambda'\setminus \Lambda$ we see that $\M_{\Lambda}(\ast \bigcup_{i\in \Lambda'\setminus \Lambda} \D_i)$ has a resolution
\[\RGamma_{\bigcap_{i\in \Lambda'\setminus \Lambda} \D_i} \M_{\Lambda} \longrightarrow \bigoplus_{|\I| = c-1} \RGamma_{\bigcap_{i\in \I} \D_i} \M_{\Lambda} \longrightarrow \cdots \longrightarrow \bigoplus_i \RGamma_{\D_i} \M_{\Lambda} \longrightarrow \M_{\Lambda}\]
where we always index over $\I \subseteq \Lambda' \setminus \Lambda$, and the resolution is in degrees $-c$ to zero. For the purposes of mapping out to $\M_{\Lambda'}$, this can be used in place of $\M_{\Lambda}(\ast \D_{\Lambda})$.

What remains is an orthogonality result: the desired result follows if we show only the first term has a nonzero $\RHom(-,\M_{\Lambda'})$.

Abstractly, for $\RGamma_{\bigcap_{i\in \I} \D_i} \M_{\Lambda}$ one sees by Kashiwara that
\[\RHom(\RGamma_{\bigcap_{i\in \I} \D_i} \M_{\Lambda},\M_{\Lambda'})\]
can be computed as an $\RHom$ in $\rA_{\Lambda,\Lambda'}/(\T_i)_{i\in \I}$ where the first argument is perfect and the second argument is supported on $\D_{j}$ for $j  \notin \I$ (here we assume $\I$ is proper). This follows from Lemma \ref{lem:orthogonality}, after reducing to working modulo $p$.
\end{proof}

We will use the following lemma only on $\overline{\rA}_{\Lambda}$, but it is convenient to handle the argument in full generality.

\begin{lemma} \label{lem:orthogonality}
Let $\overline{\rA}$ be a regular $\rF$-finite $k\laur{q-1}$ algebra, with a continuous action of $\Gamma = (1+p\Z_p) \ltimes \Z_p(1)^n$. Let $\T\in \overline{\rA}$ be a nonzerodivisor such that all generators of $\Gamma$ act trivially on $\T$, except for $\gamma_i$ generating the $i$th copy of $\Z_p(1)_i$ which acts by
\[\gamma_i \cdot \T = q^{p} \T\]
and fixes $\overline{\rA}/\T$. Let $\M\in\Perf(\overline{\rA})^{\varphi=1,h\Gamma}$ and $\M_i\in\D_{\fgu}(\overline{\rA}/\T)^{h\Gamma}$. Then
\[\RHom(\M,i_+\M_i)=0.\]
\end{lemma}

We note the proof of Kashiwara's lemma still provides $i_+$, as we used nothing more than these hypotheses.

\begin{proof}
Let $\gamma_i$ be a generator of $\Z_p(1)_i$. Because $\M$ is dualizable, using the tensor/Hom adjunction we may assume $\M=\overline{\rA}$ by absorbing it into $i_+ \M_i$. Passing through the Emerton-Kisin correspondence (established in the present regular setting in \cite{MathewModpRH}), upon setting $\X=\Spec \overline{\rA}$ and $\rZ=\Spec \overline{\rA}/\T$ it will suffice to show
\[\gamma_i-1\in \mathrm{Aut}(\RHom(i_* \mathscr{F}, \F_p)),\]
where $i_* \mathscr{F}$ corresponds to $i_+ \M_i$ under $\mathrm{RSol}$ (using Kashiwara's lemma). Indeed once we show this claim, transporting it to the corresponding $\RHom$ in $\D_{\fgu}(\overline{\rA})^{h\Gamma}$ shows the vanishing for derived $\Gamma$-invariants by normality of $\Z_p(1)_i \le \Gamma$ using Hochschild-Serre. We can check the claim that $\gamma_i-1$ acts by an automorphism after passing through the Emerton-Kisin correspondence, which on the module side we know to be sufficient for the desired vanishing.

We will prove this claim for any constructible $\F_p$-sheaf $\mathscr{F}$ on $\rZ$; we can assume it is an ordinary constructible sheaf by first applying a Postnikov d\'{e}vissage to the heart. Since $\Z_p(1)_i$ is pro-$p$ and $i_* \mathscr{F}$ is constructible and supported on the locus $\rZ$ which is fixed by $\Z_p(1)_i$, it will act on $\mathscr{F}$ through a finite $p$-group quotient.\footnote{We have used here that the vector space of endomorphisms is finite dimensional over $\F_p$. This is not obvious on the finitely generated unit module side, which is why we pass through \cite{MathewModpRH} to run the argument.} Thus upon replacing $\Z_p(1)_i$ with $p^\ell \Z_p(1)_i$ (and $\gamma_i$ with $p^\ell \gamma_i$), we can assume without loss of generality that the action of $\gamma_i$ is trivial on $\mathscr{F}$.

Switching to $\rZ$ via the $\Z_p(1)_i$-equivariant adjunction $i_* \leftrightarrows i^!$ to exploit triviality of the action on $\mathscr{F}$,
\[\RHom_{\X}(i_* \mathscr{F}, \underline{\F}_p)^{h(p^\ell \Z_p(1)_i)} = \RHom_{\rZ}(\mathscr{F}, (i^! \underline{\F}_p)^{h(p^{\ell} \Z_p(1)_i)}).\]
We then aim to argue $(i^! \underline{\F}_p)^{h p^{\ell} \Z_p(1)_i}=0$. We may check the claim on stalks, since we are testing if a finite limit ($\fib(\gamma_i-1)$) vanishes. For a geometric point $\overline{z}$ of $\rZ$, using Artin-Schreier we have
\[(i^!\underline{\F_p})_{\overline{z}} \simeq \fib((\O^{\mathrm{sh}}_{\X,\overline{z}}[1/\T]/\O^{\mathrm{sh}}_{\X,\overline{z}})[-1] \xlongrightarrow{\varphi-1} (\O^{\mathrm{sh}}_{\X,\overline{z}}[1/\T]/\O^{\mathrm{sh}}_{\X,\overline{z}})[-1]).\]
Now we note $p^{\ell}\Z_p(1)_i$ fixes $\O^{\mathrm{sh}}_{\X,\overline{z}}/(\T)$ and the generator $p^\ell \gamma_i$ sends $\T$ to $q^{p^{\ell+1}}\T$. Using this hypothesis one can see the map $\gamma_i^{p^\ell}-1: \O^{\mathrm{sh}}_{\X,\overline{z}}[1/\T]/\O^{\mathrm{sh}}_{\X,\overline{z}} \to \O^{\mathrm{sh}}_{\X,\overline{z}}[1/\T]/\O^{\mathrm{sh}}_{\X,\overline{z}}$ is invertible by checking the action on the associated graded of the filtration $\T^{-m}\O^{\mathrm{sh}}_{\X,\overline{z}}/\O^{\mathrm{sh}}_{\X,\overline{z}}$ (indeed it acts as multiplication by $q^{-p^{\ell+1}m}-1$ on the $m$th graded piece), and then the claim follows.
\end{proof}

With Lemma \ref{lem:orthog_red_ket} and Lemma \ref{lem:orthog_red} in hand, we can now prove the required $\RHom$ comparisons by making a reduction to the known comparison on $\Perf(\rA_{\Lambda'})^{\varphi=1,h\Gamma}$ (with $\Lambda'$ as in the notation of the theorem below).

\begin{theorem} \label{thm:RHom_comparisons}
Choose (possibly empty) subsets $\Lambda, \Lambda'\subseteq \{1,\ldots,r\}$. Let $\M_{\Lambda} \in \Perf(\rA_{\Lambda})^{h\Gamma,\varphi=1}$, and let $\M_{\Lambda'} \in \Perf(\rA_{\Lambda'})^{h\Gamma,\varphi=1}$. Then $\T_\et$ induces an equivalence
\[\RHom_{\D_{\fgu}(\rA)^{h\Gamma}}(j_{\Lambda,!}\M_{\Lambda},j_{\Lambda',!}\M_{\Lambda'}) \simeq \RHom(\T_\et(j_{\Lambda',!}\M_{\Lambda'}),\T_\et(j_{\Lambda,!}\M_{\Lambda})).\]
\end{theorem}

\begin{proof}
The comparisons we need fit into three cases: where $\Lambda=\Lambda'$, where $\Lambda \subset \Lambda'$, and where $\Lambda \not \subseteq \Lambda'$. In the third case, on both sides it is easy to see that the $\RHom$ is zero. In the first case we have already proven this by combining Corollary \ref{cor:thick_gen_embedding} and Proposition \ref{prop:thick_gen_equiv}.

Thus, all of the remaining work is the case $\Lambda \subseteq \Lambda'$. Geometrically, on the Kummer-\'{e}tale side both $\T_\et(j_{\Lambda',!}\M_{\Lambda'})$ and $\T_\et(j_{\Lambda,!}\M_{\Lambda})$ are supported on the closed locus
\[\rZ_{\Lambda,\eta} = \bigcap_{\lambda \in \Lambda} \D_{\lambda,\eta} \subset (\Spf \R,\mathcal{M})_{\eta}\]
equipped with the pullback log structure. Inside of this closed locus on the Kummer-\'{e}tale side of the correspondence, the $\RHom$ we are computing involves maps from $j_{\Lambda',!} \mathbb{L}' \to j_{\Lambda,!} \mathbb{L}$ (up to some shifts). In other words, maps from an arbitrary deeper stratum to $\mathbb{L}$ supported on the open stratum $(\X,\mathcal{M})_{\Lambda,\eta} = \rZ_{\Lambda,\eta}^{\circ} \subset \rZ_{\Lambda,\eta}$. Via Kashiwara's lemma, we may therefore drop $j_{\Lambda',!}$ and $j_{\Lambda,!}$ and simply consider objects $\M_{\Lambda'}, \M_{\Lambda}$ in $\D_{\fgu}(\rA/(\T_i)_{i\in \Lambda})^{h\Gamma}$ supported on the deeper stratum $(\X,\mathcal{M})_{\Lambda',\eta} \subset \rZ_{\Lambda,\eta}$ and the open stratum respectively.

We first argue in a special case where $\M_{\Lambda}$ admits an extension to 
\[\M_{\Lambda}^{\mathrm{ext}} \in \Perf(\rA/(\T_i)_{i\in \Lambda})^{\varphi=1,h\Gamma}.\]
This means $\M_{\Lambda}$ is recovered as $\M_{\Lambda}^{\mathrm{ext}} \otimes \rA_{\Lambda} = \M_{\Lambda}^{\mathrm{ext}}(\ast \D_{\Lambda})$, i.e. we $p$-completely invert the remaining coordinates. Geometrically, this asks on the Kummer-\'{e}tale side that it extends to a derived lisse sheaf on $\rZ_{\Lambda,\eta}$ with its pullback log structure. 

In this situation, in $\D_{\fgu}(\rA/(\T_i)_{i\in \Lambda})^{h\Gamma}$ we have the equivalence
\[\RHom(\M_{\Lambda}^{\mathrm{ext}}(\ast \D_{\Lambda}), \M_{\Lambda'}) \simeq \RHom(\RGamma_{\bigcap_{j \in \Lambda'\setminus \Lambda} \D_j} \M_{\Lambda}^{\mathrm{ext}}, \M_{\Lambda'})[-|\Lambda'\setminus \Lambda|]\]
provided by Lemma \ref{lem:orthog_red}. On the right, using Kashiwara again gives an $\RHom$ in $\D_{\fgu}(\rA/(\T_i)_{i\in \Lambda'})^{h\Gamma}$, corresponding to a smaller closed locus. But then the second argument of the $\RHom$ lies in 
\[\Perf(\rA_{\Lambda'})^{\varphi=1,h\Gamma} \subset \D_{\fgu}(\rA/(\T_i)_{i\in \Lambda'})^{h\Gamma},\] 
and in particular using the localization adjunction we obtain the $\RHom$
\[\RHom(\RGamma_{\bigcap_{j \in \Lambda'\setminus \Lambda} \D_j} \M_{\Lambda}^{\mathrm{ext}}(\ast \D_{\Lambda'}), \M_{\Lambda'})[-|\Lambda'\setminus \Lambda|],\]
now taken in $\Perf(\rA_{\Lambda'})^{\varphi=1,h\Gamma}$ (so both arguments are being implicitly passed through Lemma \ref{lem:kashiwara}).

On the Kummer-\'{e}tale side we can run the same reduction. By Lemma \ref{lem:orthog_red_ket} we obtain the analogous statement
\[\RHom_{\rZ_{\Lambda,\eta}}(j_{\Lambda',!} \mathbb{L}', j_{\Lambda,!} j_{\Lambda}^* \mathbb{L}^{\mathrm{ext}}) \simeq \RHom_{\rZ_{\Lambda',\eta}}(j_{\Lambda',!}\mathbb{L}', i_{\rZ_{\Lambda',\eta}}^* \mathbb{L}^{\mathrm{ext}})[-|\Lambda'\setminus \Lambda|],\]
where subscripts denote the log adic space we use when taking $\RHom$ of pro-Kummer-\'{e}tale sheaves. Further restricting to the open locus $(\X,\mathcal{M})_{\Lambda',\eta} \subset \rZ_{\Lambda',\eta}$, since $j_{\Lambda',!}$ is left adjoint to $j_{\Lambda'}^*$ regarded as the inclusion of the open stratum $(\X,\mathcal{M})_{\Lambda',\eta} \subset \rZ_{\Lambda',\eta}$ this simplifies to 
\[\RHom_{(\X,\mathcal{M})_{\Lambda',\eta}}(\mathbb{L}', j_{\Lambda'}^* i_{\rZ_{\Lambda',\eta}}^* \mathbb{L}^{\mathrm{ext}})[-|\Lambda'\setminus \Lambda|]\]
where the source and target are derived lisse. These two final $\RHom$s we obtained, on either side, have their reduction steps exchanged by $\T_\et$. Once the reduction is made they now involve derived lisse sheaves or objects in $\Perf(\rA_{\Lambda'})^{\varphi=1,h\Gamma}$, and using Proposition \ref{prop:thick_gen_equiv} they can be seen to be identified by $\T_\et$.

The general case then follows from Lemma \ref{lem:primetop_Kum_red}, which reduces the claim to the situation where $\M_{\Lambda}^{\mathrm{ext}}$ exists. Indeed on the Kummer-\'{e}tale side the same formula for $\mu^{r}_{\rN}$ invariants holds, precisely because
\[(\Spf \R_{\rN},\mathcal{M}_{\rN})_\eta \longrightarrow (\Spf \R,\mathcal{M})_\eta\]
is a torsor in the pro-Kummer-\'{e}tale topology for this group.
\end{proof}

We owe an additional lemma, which lets us assume that $\M_{\Lambda}$ admits an extension. We let $\rA_{\rN}$ be the ring of the (Laurentized) log $q$-de Rham prism of the Kummer cover as in Proposition \ref{prop:open_strat_Kum_desc}. We also use the same notation for the Kummer cover.

\begin{lemma} \label{lem:primetop_Kum_red}
Choose subsets $\Lambda \subseteq \Lambda'\subseteq \{1,\ldots,r\}$, and suppose we are given objects $\M_{\Lambda'}$ and $\M_{\Lambda}$ in $\D_{\fgu}(\rA/(\T_i)_{i\in \Lambda})^{h\Gamma}$ supported on the strata for $\Lambda'$ and $\Lambda$ respectively. There is a prime-to-$p$ Kummer cover $f: (\Spf \R_{\rN},\mathcal{M}_{\rN}) \to (\Spf \R,\mathcal{M})$ such that 
\[f^* \M_{\Lambda} = \M_{\Lambda,\rN}^{\mathrm{ext}}(\ast \D_{\Lambda,\rN})\]
where $\M_{\Lambda,\rN}^{\mathrm{ext}}\in \Perf((\rA_{\rN}/(\T^{1/\rN}_i)_{i\in \Lambda}))^{h\Gamma,\varphi=1}$. Moreover,
\[\RHom(f^* \M_{\Lambda}, f^* \M_{\Lambda'})^{h\mu^r_{\rN}} = \RHom(\M_{\Lambda}, \M_{\Lambda'}).\]
\end{lemma}

\begin{proof}
The first assertion was already explained in the proof of Proposition \ref{prop:thick_gen_equiv} in the essential surjectivity step.

For the second assertion, as in Proposition \ref{prop:open_strat_Kum_desc}, we have 
\[\rA_{\rN} = \bigoplus_{\alpha \in (\Z/\rN\Z)^r} \rA \cdot \rS^{\alpha},\] 
with a natural $\mu_{\rN}^r$ action associated to the grading, which gives $\RHom(f^* \M_{\Lambda}, f^* \M_{\Lambda'})$ a $\mu_{\rN}^r$-action. In our specific case, the components for $i\in \Lambda$ act trivially -- in what follows we will simply prove a more general assertion, regarding $\M_{\Lambda}$ and $\M_{\Lambda'}$ as arbitrary finitely generated unit modules in $\D_{\fgu}(\rA)^{h\Gamma}$ (using Kashiwara's lemma).

Since $\rN$ is invertible in $\Z_p$, for any $\M \in \D_{\fgu}(\rA)^{h\Gamma}$ the homotopy fixed points of its pullback are computed by taking the degree zero summand in the decomposition $\rA_{\rN} = \bigoplus_{\alpha \in (\Z/\rN\Z)^r} \rA \cdot \rS^{\alpha}$. We also note that the projection onto this summand is both $\Gamma$ and Frobenius equivariant. We therefore obtain
\begin{align*}
\RHom(f^*\M_{\Lambda},f^*\M_{\Lambda'})^{h\mu_{\rN}^r}
&\simeq \RHom\left(\M_{\Lambda},(\rA_{\rN}\widehat{\otimes}_{\rA}^{\L}\M_{\Lambda'})^{h\mu_{\rN}^r}\right) \\
&\simeq \RHom(\M_{\Lambda},\M_{\Lambda'}).
\end{align*}
In the first line we use adjunction, and that $\RHom(\M_{\Lambda},-)$ preserves limits to move $h\mu_{\rN}^r$ inside. This proves the second claim. Here all $\RHom$s are calculated in $\D_{\fgu}(\rA)^{h\Gamma}$.
\end{proof}

Combining our previous results, we may deduce the main theorem using Lemma \ref{lem:thick_equiv_criterion}.

\begin{theorem} \label{thm:RH_main}
Let $(\X,\mathcal{M})$ satisfy Assumption \ref{ass:horizontal_nc} over $\W(k)[\zeta_p]$. Let $\D^{(b),\mathrm{pKum}}_{\cons,\D}((\X,\mathcal{M})_\eta, \Z_p)$ denote the full subcategory of $p$-Kummer constructible sheaves in $\D^{(b)}_{\cons,\D}((\X,\mathcal{M})_\eta, \Z_p)$. Then
\begin{display}
\T_\et: \D_{\fgu}((\X,\mathcal{M})_\Prism,\Laur)^\op \ar{r}{\sim} & \D^{(b),\mathrm{pKum}}_{\cons,\D}((\X,\mathcal{M})_\eta, \Z_p)
\end{display}
is a symmetric monoidal equivalence.
\end{theorem}

\begin{proof}
By Lemma \ref{lem:strict_etale_local_framing} and strict \'{e}tale descent (for the source Theorem \ref{thm:fgu_etale_desc} and for the target Lemma \ref{lem:constructible_strict_etale_descent}), we may first work on a framed log affine $(\Spf\R,\mathcal{M})$.

In this case, we apply Lemma \ref{lem:thick_equiv_criterion}. We use the generators from Lemma \ref{lem:ket_thick_gens} and Lemma \ref{lem:crys_thick_gens}, which correspond (up to dimension shifts, which do not matter for the argument) by combining Corollary \ref{cor:thick_gen_embedding} and Proposition \ref{prop:thick_gen_equiv}. 

We have now checked all the hypotheses:
\begin{itemize}
  \item Both source and target are idempotent complete, and $\T_\et$ is exact. This is explained already immediately after Lemma \ref{lem:thick_equiv_criterion}.
  \item The generators of the source and target correspond under $\T_\et$, as explained above.
  \item The comparisons on $\RHom$s are done by Theorem \ref{thm:RHom_comparisons}.
\end{itemize}
Thus $\T_\et$ is an equivalence.

We further claimed that it is a symmetric monoidal equivalence. The fact that the equivalence is symmetric monoidal can also be tested locally, where via Lemma \ref{lem:factorization} it reduces to symmetric monoidality in Theorem \ref{thm:EK_perfect}.
\end{proof}

We must also understand perverse $t$-exactness. We use a definition of the $p^+$-perverse $t$-structure on the Kummer-\'{e}tale side of the correspondence which is similar to \cite{BhattHansenCons}.

\begin{definition} \label{def:perverse_pKum}
For each logarithmic stratum $j_\Lambda:(\X,\mathcal{M})_{\Lambda,\eta} \hookrightarrow (\X,\mathcal{M})_\eta$, put $d_\Lambda=\dim((\X,\mathcal{M})_{\Lambda,\eta})$. We define the $p$-perverse $t$-structure on $\D^{(b)}_{\cons,\D}((\X,\mathcal{M})_\eta,\Z_p)$ by
\[{}^p\D^{\le0} = \{\mathscr{F} \mid j_\Lambda^*\mathscr{F} \in \D^{\le-d_\Lambda}((\X,\mathcal{M})_{\Lambda,\eta},\Z_p)\text{ for every }\Lambda\},\]
\[{}^p\D^{\ge0} = \{\mathscr{F} \mid j_\Lambda^!\mathscr{F} \in \D^{\ge-d_\Lambda}((\X,\mathcal{M})_{\Lambda,\eta},\Z_p)\text{ for every }\Lambda\}.\]
We write ${}^p\H^i$ for its cohomology functors, and note that a similar argument to \cite[Proposition 4.6(1)]{BhattHansenCons} shows it is well-defined (here in their argument pick $\U$ to be the open stratum, and use recollement to work down the stratification).

The $p^+$-perverse $t$-structure is then defined by
\[{}^{p^+}\D^{\le0}=\{\mathscr{F}\in{}^p\D^{\le1} \mid {}^p\H^1(\mathscr{F})\text{ is locally bounded }p\text{-power torsion}\},\]
with the $\ge 0$ part simply given by
\[{}^{p^+}\D^{\ge0}=\{\mathscr{F} \mid \mathscr{F}\dquot p\in{}^p\D^{\ge0}((\X,\mathcal{M})_\eta,\F_p)\}.\]
One may check these $t$-structures restrict to the $p$-Kummer subcategory (i.e. that perverse truncations of $p$-Kummer sheaves are $p$-Kummer).
\end{definition}

\begin{proposition} \label{prop:Tet_perverse_exact}
Assume that $(\X,\mathcal{M})_\eta$ is pure of dimension $d$. The equivalence of Theorem \ref{thm:RH_main}, shifted by $[d]$, is $t$-exact for the standard $t$-structure on the source and the $p^+$-perverse $t$-structure of Definition \ref{def:perverse_pKum} on the target.
\end{proposition}

\begin{proof}
The source $t$-structure is strict \'{e}tale local as we have already shown in Lemma \ref{lem:fgu_local_tstruct} that the $t$-structure on finitely generated unit Laurent crystals is restricted from the $t$-structure on $\D_{\Box}((\X,\mathcal{M})_\Prism,\Laur)$, and the same holds for the $t$-structure on the target in Definition \ref{def:perverse_pKum}. We may therefore work on a framed log affine.

We first prove the claim modulo $p$ and therefore use the base $\overline{\rA}=\overline{\R}^{\sq}\laur{q-1}$. We may consider $\M\in\D_{\fgu}(\overline{\rA})^{h\Gamma}$, using Theorem \ref{thm:qdR_crys_Gamma_mod}. For $\Lambda\subseteq\{1,\ldots,r\}$, as before put
\[\overline{\rA}_\Lambda = \left(\overline{\rA}/(\T_i)_{i\in \Lambda}\right)[1/\T_j]_{j \notin \Lambda, j\le r}\]
and $\M_{\Lambda} := \M \otimes_{\overline{\rA}}^{\L} \overline{\rA}_\Lambda$. This is a perfect complex by the argument from Proposition \ref{prop:perfectoid_constructibility}.

Now let $j_\Lambda:\X_{\eta,\Lambda} \hookrightarrow (\X,\mathcal{M})_\eta$ denote inclusion of the corresponding stratum (which is of dimension $d-|\Lambda|$). Using Lemma \ref{lem:factorization} to compute $\T_\et$ (and its derived reduction $\T_{\et,\F_p}:=\T_\et \dquot p$), we deduce from the contravariant version of Katz's correspondence that
\[j_\Lambda^*\T_{\et,\F_p}(\M)[d]\in\D^{\le-(d-|\Lambda|)} \Longleftrightarrow \M_\Lambda\in\D^{\ge-|\Lambda|}.\]
We now claim that $\M\in\D^{\ge0}$ if and only if $\M_\Lambda\in\D^{\ge-|\Lambda|}$ for every $\Lambda$, where we use the ordinary $t$-structure on $\D_{\fgu}(\overline{\rA})^{h\Gamma}$. Applying \cite[Lemma 3.1]{AlonsoJeremiasSaorinTStructures} along the natural stratification induced by $\V(\T_i)$ (using $i_\Lambda^!\M \simeq \M_\Lambda[-|\Lambda|]$) gives the claim.

Definition \ref{def:perverse_pKum} now gives
\[\T_{\et,\F_p}(\M)[d] \in {}^p\D^{\le0} \Longleftrightarrow \M \in \D^{\ge0}.\]
Since $\T_{\et,\F_p}$ is contravariant, it will exchange coconnectivity on the source with perverse connectivity on the target. The previous claim also gave the converse, that perverse connectivity on the target implies coconnectivity on the source. Since $\T_{\et,\F_p}$ is also an antiequivalence we then get perverse $t$-exactness modulo $p$.

To deduce the general claim, we see $\M \in \D^{\le 0}$ if and only if $\M \dquot p \in \D^{\le 0}$ by derived Nakayama. On the Kummer-\'{e}tale side we have analogous reduction criteria from Definition \ref{def:perverse_pKum}, so for the $\le 0$ part on the source we may reduce to the mod $p$ assertion. Orthogonality then proves the general claim.
\end{proof}

It is also possible to extract ordinary constructible sheaves from the correspondence. We warn that extracting ordinary constructible sheaves in the derived sense cannot work without altering the crystal side of the correspondence, as these are not a full subcategory of the derived category of Kummer-\'{e}tale sheaves (simply because Kummer-\'{e}tale cohomology and \'{e}tale cohomology differ).

\begin{remark} \label{rem:eff}
In \cite[Theorem 7.5.6]{KedlayaLiu2016}, there is a notion of an ``effective'' $(\varphi,\Gamma)$-module which corresponds to an ordinary local system. We may isolate a similar subcategory in our setting, although it is slightly more complex because our correspondence is perverse $t$-exact. In particular, this means that the finitely generated unit $(\varphi,\Gamma)$ modules corresponding to ordinary constructible sheaves are typically not in the heart of the ordinary $t$-structure (otherwise they would correspond to perverse sheaves). Nevertheless, we may isolate a similar subcategory: for a framed log affine and $\rA=\R_0^{\sq}\laur{q-1}^{\wedge}_p$, we say $\M \in \D_{\fgu}(\rA)^{h\Gamma}$ is effective if for every $\T_i$ the $i$th $\Z_p(1)$-factor in $\Gamma$ acts trivially on
\[(\rA/(\T_i))\otimes^{\L}_{\rA}\M\]
and $\T_\et(\M)$ lies in the heart of the ordinary $t$-structure.

In the derived setting, it seems plausible that one could take the ``effective part'' of an $\R\Hom$ to obtain the derived category of ordinary (i.e. without any log structures) constructible sheaves along $\D$.
\end{remark}

We remark that trivial modifications of the previous arguments allow one to deduce the theorem over a slightly more general base. In particular our restriction to $\W(k)[\zeta_p]$ up until this point is mainly for expository reasons to avoid tracking additional variables and cases throughout.

\begin{corollary} \label{cor:other_bases}
Theorem \ref{thm:RH_main} and all preceding results stated over $\W(k)[\zeta_p]$ admit corresponding variants over $\W(k)[\zeta_{p^\ell}]$ for every $\ell\ge0$.
\end{corollary}

\begin{proof}
Corollary \ref{cor:coh_laur_other_bases} is the only nontrivial change needed in the argument, which covers these cases. With this, one can deduce Theorem \ref{thm:fgu_etale_desc} verbatim (as it uses nothing about $\Gamma$). The following arguments are insensitive to the change in prismatic ideal or $\Gamma$, except for Proposition \ref{prop:open_decompletion_hom}; however there the reduction is to vanishing of a cohomology group, so changing $\Gamma$ from $(1+p\Z_p) \ltimes \Z_p(1)^n$ to $(\mu_{p-1} \times (1+p\Z_p)) \ltimes \Z_p(1)^n$ does not change the result by using Hochschild-Serre.
\end{proof}

\begin{remark}
The entire local theory using $(\varphi,\Gamma)$-modules seems to have no obstruction to working over $\A_K$, but unfortunately the globalization step does not appear to work. One would want to have a prism $(\A_K^+, \varphi^n([p]_q))$ which satisfies 
\[\A_K^+[1/\varphi^{n}([p]_q)]_p^{\wedge} \simeq \A_K.\]
However it seems such a prism cannot exist in general without dropping some desirable properties. See \cite[Appendix 2, Proposition 23]{Watanabe} for more details. The argument there showing nonexistence uses that the mod $p$ reduction should be $\E_K^+$, the ring of integers of the field of norms. It is still possible that $\A_K^+$ exists by dropping this condition, but this seems unnatural to do.
\end{remark}

The following explicit $(\varphi,\Gamma)$-module gives a constructible sheaf with nontrivial gluing between the two strata.

\begin{example} \label{ex:nontrivial_cons}
Consider the formal disk $(\X,\mathcal{M})=(\Spf \Z_p\langle \T \rangle,\N)$ with its log structure induced by $\D=\{\T=0\}$. We will denote by $j:\U_\eta = \X_\eta\setminus\D_\eta \hookrightarrow (\X_\eta,\N)$ the inclusion of the punctured disk, with its trivial log structure. For the $(\varphi,\Gamma)$-modules in the following example, we will set $\overline{\R} = \F_p\langle \T \rangle$ and $\rA = \overline{\R}^{\sq}\laur{q-1}$.

Following \cite{Brinon}, let $f:\mathcal{E} \to \mathbb{D} = \X_\eta$ be the projectivization of the Tate family with parameter $\tilde{q}=p^p\T$, given by
\[v^2+uv=u^3-b_2(\tilde{q})u-b_3(\tilde{q}),\]
where $b_2$ and $b_3$ are the standard coefficients for the Tate curve in $\tilde{q}\Z_p\llbracket \tilde{q} \rrbracket$.

It is smooth away from $\T=0$, where it becomes a nodal cubic. We consider the ordinary Zariski constructible sheaf
\[\mathscr{F} = \R^1 f_* \underline{\F_p}.\]
On the open $\U_\eta$, we have $\mathscr{F}|_{\U_\eta}=\underline{\Hom}_{\F_p}(\mathcal{E}[p],\underline{\F_p})$. As explained in \cite[\S 3.1]{Brinon} one may describe $\mathcal{E}[p]$ on $\U_\eta$ as an extension of $\underline{\F_p}$ by $\mu_p$ via the Kummer class $[\tilde{q}]$, and since $\tilde{q}=p^p\T$ this is also the nonzero Kummer class of $\T$. That is, $[\tilde{q}]=[\T]$ in $\H^1_\et(\U_\eta,\mu_p)$ as classes obtained via the Kummer map $\O(\U_\eta)^\times \to \H^1_\et(\U_\eta,\mu_p)$ since the kernel consists of $p$th powers. Brinon also explains the explicit formula for the Galois action, using this Kummer class $[\tilde{q}]$ for the Tate parameter.

Globally, the sheaf $\mathscr{F}$ is a nonsplit extension
\[0 \longrightarrow \underline{\F_p} \longrightarrow \mathscr{F} \longrightarrow j_!\underline{\F_p}(-1) \longrightarrow 0.\]
Here we deduce it is nonsplit by dualizing the previous description of $\mathcal{E}[p]$ as a nonsplit extension corresponding to $[\tilde{q}] \neq 0$ on the open stratum.

Under $\T_{\et,\F_p}$, the corresponding $(\varphi,\Gamma)$-module has underlying module
\[\M=\rA[1/\T] e_1 \oplus \rA e_2.\]
The submodule $\rA[1/\T]e_1$ and quotient $\rA e_2$ correspond respectively to $j_!\underline{\F_p}(-1)$ and $\underline{\F_p}$, in accordance with the contravariance of $\T_\et$. Frobenius fixes both basis vectors; this can be seen since the Kummer class $[\tilde{q}]$ trivializes on the perfectoid cover $(\X,\mathcal{M})_{\infty,\eta}$, because $\T$ admits a compatible system of $p$-power roots. For $\gamma=(a_\gamma,b_\gamma) \in \Gamma=\Z_p^\times \ltimes \Z_p(1)$, we use the matrix
\[\R_{\gamma} = \begin{pmatrix} \overline{a}_\gamma & \overline{b}_\gamma \\ 0 & 1 \end{pmatrix}\]
where the bars denote reduction modulo $p$ to give the action of $\Gamma$. This means for $a_\gamma$, in the decomposition $\Z_p^\times = \F_p^\times \times (1+p\Z_p)$ we extract only the $\F_p^\times$ component as $\overline{a}_\gamma$. Note $\R_{\gamma} \cdot \R_{\delta} = \R_{\gamma \delta}$ as needed. One can compute this action of $\Gamma$ on trivialization of the Kummer class governing the extension.
\end{example}

\section{Pushforward compatibility}

\subsection{Six functors and Laurent stacks}

The category $\D_{\Box}((\X,\mathcal{M})_\Prism,\O_\Prism)$ is large enough that obtaining a six functor formalism in this category formally implies the Laurent analogue. For ordinary formal schemes (i.e. non-logarithmic), such a six functor formalism is developed in \cite{Kipp,Kipp2} by constructing an analytic stack $\X^{\Prism,\Box}$ using the formalism of analytic stacks developed by Clausen and Scholze in \cite{clausenscholze}. We will mainly use \cite{Kipp2} adapted to a logarithmic setting, as the conventions there will exactly recover the category of solid Laurent crystals that we have thus far been treating site-theoretically.

By and large, adapting these results to the logarithmic setting is straightforward as many of the important results in \cite{Kipp2} are stable under arbitrary base change. Let $(\X,\mathcal{M})$ satisfy Assumption \ref{ass:horizontal_nc} over $\O_K$, where $K/\W(k)[1/p]$ is a finite totally ramified extension. As explained in more detail in Definition \ref{def:Olsson_stack}, the log prismatic stack is defined via the Cartesian square
\begin{display}
(\X,\mathcal{M})^{\Prism} \ar{r} \ar{d} & \X^{\Prism} \ar{d} \\
\Z_p^{\Prism} \times \Log \ar{r} & \Log^{\Prism}
\end{display}
where $\Log$ is the $p$-completion of Olsson's log stack from \cite{OlssonLogStacks}. This is a relative prismatic stack over $\Log$, viewed as a ``$\delta$-stack''. We use this as a starting point to transfer Kipp's construction of $\X^{\Prism,\Box}$ to one for $(\X,\mathcal{M})^{\Prism,\Box}$, but we will need to upgrade the remaining terms (besides $\X^{\Prism}$, which is already handled by Kipp) in the Cartesian square to analytic stacks.

In \cite[Construction 2.2.7]{Kipp2}, the \'{e}id topology is defined. Roughly, one takes the topology generated by \'{e}tale covers and integral descendable covers; as explained in \cite[2.2.5]{Kipp2}, integrality for formal schemes is tested on ordinary schemes, so for formal affines one checks it on rings modulo an ideal of definition. We will always use integrality in this sense.

\begin{definition} \label{def:box_analytic_realization}
For a formal prestack $\mathcal{Y}$ over $\Z_p^\Prism$, let $\mathrm{sh}_{\mathrm{\acute{e}id}}\mathcal{Y}$ denote the \'{e}id sheafification of the associated presheaf on $\mathrm{fSch}^{\mathrm{aff}}_{/\Z_p^\Prism}$ from \cite[Construction 2.2.7]{Kipp2}. We define
\[\mathcal{Y}^{\Box}:=(-)^{\Box}_{\mathrm{Kipp}}(\mathrm{sh}_{\mathrm{\acute{e}id}}\mathcal{Y}),\]
which we call informally the \emph{analytic realization} where
\[(-)^{\Box}_{\mathrm{Kipp}}:(\Z_p^\Prism)_{\mathrm{\acute{e}id}}\longrightarrow\mathrm{AnStack}_{/\Z_p^{\Prism,\Box}}\]
is the left exact functor of \cite[Theorem 2.3.3(A) and \S3.6.1]{Kipp2}. We note $(-)^{\Box}$ preserves finite limits and colimits, and $(\Spf\B)^{\Box}=\Spa(\B,\B)$ for a prism $\B$.
\end{definition}

Strictly speaking, this construction is not obviously the same as in \cite{Kipp2}; however this definition will be more convenient for us, as the log prismatic stack is a type of relative stack so we need a version of analytification which also works there. It turns out, however, that this is the same as Kipp's construction. In \cite[3.1.9]{Kipp2} the topos $\mathsf{Sh}(\Spf(\Z_p)_{\sqrt{\infty},\et})$ is introduced, or sheaves for the topology generated by infinite root covers and \'{e}tale covers. The functor $\Spf \R \mapsto \Spa(\Prism_\R,\Prism_\R)$ sends covers in this topology to effective epimorphisms of analytic stacks. It is shown that semiperfectoids form a basis of the topology $\Spf(\Z_p)_{\sqrt{\infty},\et}$, and so one obtains a colimit preserving functor
\[\mathsf{Sh}(\Spf(\Z_p)_{\sqrt{\infty},\et}) \longrightarrow \mathrm{AnStack}_{/\Z_p^{\Prism,\Box}}\]
induced by $\Spf \R \mapsto \Spa(\Prism_\R,\Prism_\R)$ on the semiperfectoid basis. This is fine in the non-logarithmic setting, but unfortunately in the logarithmic setting semiperfectoids are problematic (they may not even have maps to the log prismatic stack, as they often lack fine log structures).

To show the two constructions are the same and yield $\X^{\Prism,\Box}$, it will suffice to show that for an infinite root/\'{e}tale cover $\Spf \R \to \X$
\[\Spf \Prism_\R \longrightarrow \X^{\Prism}\]
is an effective epimorphism in the \'{e}id topos. Then we may use the corresponding colimit presentation coming from the \v{C}ech nerve to present $\X^{\Prism}$ in the \'{e}id topos, and then since both functors are colimit preserving we get the same colimit presentation in analytic stacks.

In fact we claim more generally that for an infinite root/\'{e}tale cover $\Spf \rA \to \Spf \B$ after prismatization and \'{e}id sheafification we get an effective epimorphism in the \'{e}id topos. It suffices to check this on generating covers; for \'{e}tale covers it's clear, so one must show 
\[\mathrm{sh}_{\mathrm{\acute{e}id}} \V_\infty^\Prism \longrightarrow \mathrm{sh}_{\mathrm{\acute{e}id}} \V^\Prism\]
is an effective epimorphism where $\V=\Spf \Z_p\langle t_i\rangle_{i\in \I}$ and $\V_\infty = \Spf \Z_p\langle t_i^{1/p^\infty}\rangle_{i\in \I}$. In the affine case it is easy to see the \'{e}id sheafification is in fact redundant, because $\W(-)$ is an \'{e}id sheaf in spectra by \cite[3.2.1]{Kipp2} (in particular one may apply the argument in the easy affine case of \cite[Lemma 7.3]{APC2}). The argument of \cite[Lemma 6.3]{APC2} applies to see we get an effective epimorphism. Consider an $\R$-point of $\V^{\Prism}$. The explicit cover $\R \to \rS$ constructed there (giving an $\rS$-point lift to $\V_\infty^\Prism$) is in fact an \'{e}id cover. Upon unwinding their proof one first constructs, after orienting the Cartier-Witt divisor as $(d)\to \W(\R)$ by localizing,
\[\B=\W(\R)\langle x_i^{1/p^\infty} \rangle_{i\in \I} \left\{\frac{x_i-\tilde{t}_i}{d}\right\}^{\wedge}_{(p,d)}\]
where $\tilde{t}_i$ are lifts of $t_i$ to $\W(\R)$ and $\delta(x_i^{1/p^n})=0$. Then $\rS=\B \otimes_{\W(\R)}\R$ (underived by their flatness result). Descendability follows from the proof of \cite[Proposition 3.2.1.8]{Kipp}. For integrality, Lemma \ref{lem:prismatic_envelope_integral} shows that $\pi_0(\B\dquot(p,d))$ is integral over the algebra obtained from $\pi_0(\W(\R)\dquot(p,d))$ by adjoining compatible $p$-power roots of the images of the $\tilde{t}_i$. This algebra is itself integral over $\pi_0(\W(\R)\dquot(p,d))$. Since the image of $(p,d)$ in $\R$ is nilpotent, we deduce $\R\to\rS$ is integral, hence an \'{e}id cover.

Thus, we use Definition \ref{def:box_analytic_realization} throughout. In particular, we may then safely define $(\X,\mathcal{M})^{\Prism,\Box}$ as the fiber product
\begin{display}
(\X,\mathcal{M})^{\Prism,\Box} \ar{r} \ar{d} & \X^{\Prism,\Box} \ar{d} \\
(\Z_p^{\Prism} \times \Log)^{\Box} \ar{r} & \Log^{\Prism,\Box}
\end{display}
where the notation on the bottom is interpreted according to Definition \ref{def:box_analytic_realization}, and thus includes \'{e}id sheafification. Finite (derived) limit preservation of $(-)^{\Box}$ shows one may also construct $(\X,\mathcal{M})^{\Prism,\Box}$ by applying analytic realization to the already existing formal stack over $\Z_p^{\Prism}$.

The logarithmic Hodge-Tate locus is the fiber product $(\X,\mathcal{M})^{\Prism,\Box} \times_{\Z_p^{\Prism,\Box}} \Z_p^{\mathrm{HT},\Box}$, where $\Z_p^{\mathrm{HT},\Box}$ is as in \cite{Kipp2}. Following \cite[Construction 4.3.7 and \S4.3.9]{Kipp2}, we write $(\X,\mathcal{M})^{\mathcal{L},\Box}$ for the closure of the open complement of the Hodge-Tate locus, further base changed along $\colim_n \Spa(\Z/p^n) \to \Spa \Z_p$. Thus, after pullback along an orientable prism $(\B,(d))$, it is given by
\[\colim_{n\ge 1}\Spa\left((\B\dquot p^n)[1/d],\B\dquot p^n\right).\]
We use this construction often, so we have given it the notation $\AnSpf$ in Definition \ref{def:anspf}.

Since this stack is obtained from $\X^{\Prism,\Box}$ by base change, the relevant results of \cite{Kipp2} apply to it. See \cite[\S1.1]{Kipp2} for background on cohomological smoothness and weak cohomological properness in six functor formalisms.

\begin{lemma} \label{lem:formal_f!}
Let $f: \X \to \Y$ be a morphism of analytic stacks, and suppose it is cohomologically smooth and weakly cohomologically proper. Then given a square
\begin{display}
\X \times_{\Y} \rZ \ar{r}{\iota_{\X}} \ar{d}{f_{\rZ}} & \X \ar{d}{f} \\
\rZ \ar{r}{\iota_\Y} & \Y
\end{display}
the map $f_{\rZ}: \X \times_\Y \rZ \to \rZ$ is cohomologically smooth and weakly cohomologically proper. Moreover, $f_{\rZ}^! \O_{\rZ}\simeq \iota_\X^* f^! \O_{\Y}$.
\end{lemma}

\begin{proof}
Weak cohomological properness and cohomological smoothness are stable under arbitrary base change by \cite[Remark 1.1.8]{Kipp2}, so $f_{\rZ}$ is again cohomologically smooth and weakly cohomologically proper. 

For the second claim about $f_{\rZ}^! \O_{\rZ}$, this follows from cohomological smoothness (see \cite[Definition 1.2.1.1(b)]{Kipp}).
\end{proof}

\begin{lemma} \label{lem:formal_pushforward_duality}
Let $f:\X\to\Y$ be a morphism of analytic stacks, and  suppose that it is cohomologically smooth and weakly cohomologically proper. If
\[\R f_! \simeq \R f_*\]
via the canonical map, $\R f_*$ preserves dualizable objects. Moreover, for every dualizable object $\mathcal{E} \in \D_{\Box}(\X)$ there is a natural equivalence
\[(\R f_*\mathcal{E})^\vee \simeq \R f_*(\mathcal{E}^\vee\otimes f^! \O_Y).\]
\end{lemma}

\begin{proof}
For every $\mathcal{E}'\in\D_{\Box}(\Y)$, one has
\[\RcalHom_\Y(\R f_*\mathcal{E},\mathcal{E}') \simeq \R f_*\RcalHom_\X(\mathcal{E},f^!\mathcal{E}') \simeq \R f_*(\mathcal{E}^\vee\otimes f^! \O_Y) \otimes \mathcal{E}'.\]
The first isomorphism follows from $\RcalHom_\Y(\R f_!\mathcal{E},\mathcal{E}') \simeq \R f_*\RcalHom_\X(\mathcal{E},f^!\mathcal{E}')$ (which one can check by showing maps into each are the same, and using the tensor/Hom adjunction and the projection formula) and $\R f_! \simeq \R f_*$, while the second is just using the formula $f^!(-)\simeq f^! \O_Y \otimes \L f^*(-)$.

Thus $\R f_*\mathcal{E}$ is dualizable with dual $\R f_*(\mathcal{E}^\vee\otimes f^! \O_Y)$.
\end{proof}

\begin{lemma} \label{lem:prismatic_stack_cohom_smooth_proper}
Let $(\X,\mathcal{M}_\X) \to (\Y,\mathcal{M}_\Y)$ be a strict log smooth proper morphism of log $p$-adic formal schemes both satisfying Assumption \ref{ass:horizontal_nc} over $\O_K$ of pure relative dimension $d$. Then the induced maps $f$ on $(\X,\mathcal{M}_\X)^{\Prism,\Box}$ and $(\X,\mathcal{M}_\X)^{\mathcal{L},\Box}$ are cohomologically smooth and weakly cohomologically proper, and on both stacks one has
\[f^!\O\simeq\O\{d\}[2d].\]
Moreover $\R f_! \simeq \R f_*$ via the canonical map.
\end{lemma}

\begin{proof}
Since we assumed strictness, the underlying morphism of formal schemes $\X \to \Y$ is smooth proper. We may first reduce all claims to the analogous ones for $\X^{\Prism,\Box}$ and $\Y^{\Prism,\Box}$. Using strictness of $f$,
\[(\X,\mathcal{M}_\X)^{\Prism,\Box} \simeq \X^{\Prism,\Box} \times_{\Y^{\Prism,\Box}}(\Y,\mathcal{M}_\Y)^{\Prism,\Box}.\]
The claimed result for $(\X,\mathcal{M}_\X)^{\Prism,\Box}$ then follows from Lemma \ref{lem:formal_f!} and compatibility of $\R f_! \simeq \R f_*$ with arbitrary base change. For the Laurent stack $(\X,\mathcal{M}_\X)^{\mathcal{L},\Box}$, the claim again follows by base change.

Thus, we must handle the non-logarithmic case. The map $f: \X^{\Prism,\Box} \longrightarrow \Y^{\Prism,\Box}$ is cohomologically smooth and weakly cohomologically proper by Theorem 4.1.3(C) in \cite{Kipp2} (using Example 4.3.2 of \cite{Kipp2} to deduce analogous claims for the prismatization). Using Theorem 4.1.3(F) in \emph{loc. cit.} we also get $f^!\O \simeq \O\{d\}[2d]$. The only nontrivial claim is therefore that the canonical map $\R f_! \longrightarrow \R f_*$ is an isomorphism, as weak cohomological properness a priori only gives a \emph{noncanonical} isomorphism. For this we will use results from \cite{MikamiDuality}.

Using the diagonal, we can reduce to the case of regular closed immersions. We use $\D_{\Box}$ to denote the six functor formalism on analytic stacks, and let $\pi_2$ denote the second projection from $\X^{\Prism,\Box}\times_{\Y^{\Prism,\Box}}\X^{\Prism,\Box}$. By weak cohomological properness we see $f$ is $\D_{\Box}$-prim, so \cite[Corollary 1.28]{MikamiDuality} gives a codualizing object
\[\delta_f \simeq \R\pi_{2,*}\R\Delta_{f,!}\O\]
and a natural equivalence $\R f_!(\delta_f\otimes-) \xlongrightarrow{\sim} \R f_*$. It therefore will suffice to trivialize $\delta_f$ canonically. Since $\pi_2\circ\Delta_f=\id$, it is enough to construct a canonical equivalence 
\[\R\Delta_{f,!} \xlongrightarrow{\sim} \R\Delta_{f,*}.\]
Since $\X \to \Y$ was smooth proper, the diagonal $\X \to \X \times_\Y \X$ is now a regular closed immersion. It is worth pointing out that $\R\Delta_{f,!}$ does make sense in this situation, as the induced map $\Delta_f$ on $(-)^{\Prism,\Box}$ is $!$-able by \cite[Theorem 4.1.3(A)]{Kipp2}.

It will then suffice to abstractly consider the case where $i: \X \to \Y$ is some regular closed immersion, with $i_{\Prism}: \X^{\Prism,\Box} \to \Y^{\Prism,\Box}$ the induced map on analytic stacks, and show that $\R i_{\Prism,!} \xlongrightarrow{\sim} \R i_{\Prism,*}$ via the canonical map. Here we do not need to impose any hypotheses on $\X$ or $\Y$, and in fact they can be arbitrary derived $p$-adic formal schemes. We will apply \cite[Lemma 1.50]{MikamiDuality}\footnote{See \cite[Definition 1.10]{MikamiDuality} for why it applies to analytic stacks.}, upon which the desired equivalence holds provided that $i_{\Prism}$ is, universally $\D_{\Box}(-)^*$-locally on its target, represented by affine proper morphisms (compare \cite[\S2.1.10]{Kipp2} and the class in \cite[Definition 1.5(2)]{MikamiDuality}).

We now check that this is true. There is a semiperfectoid cover $\Spa(\Prism_{\R},\Prism_{\R})$ of $\Y^{\Prism,\Box}$ from \cite[Construction 3.7.1]{Kipp2}, which is an effective epimorphism of analytic stacks and therefore a universal $\D_{\Box}^{*}$-cover (we may need a family of such covers indexed over some set $\I$, but without loss of generality we can work with a single one for the claim). As $\R$ is semiperfectoid it is also a quotient of some perfectoid ring $\rS$; if one chooses local equations $\overline{f}_1,\ldots,\overline{f}_d$ for the regular closed immersion $i$ after base change along $\Spf \R$, we can lift these to $f_i\in \Prism_{\rS} \twoheadrightarrow \rS \twoheadrightarrow \R$, also giving $f_i \in \Prism_\R$ by applying the map $\Prism_{\rS} \to \Prism_\R$. Now upon base changing along this cover $\Spa(\Prism_{\R},\Prism_{\R}) \to \Y^{\Prism,\Box}$, the resulting situation is summarized by the Cartesian squares
\begin{display}
\Spa(\Prism_{\rS \dquot (f_1,\ldots,f_d)},\Prism_{\rS \dquot (f_1,\ldots,f_d)}) \ar{r} & \Spa(\Prism_{\rS},\Prism_{\rS}) \\
\Spa(\Prism_{\R \dquot (f_1,\ldots,f_d)},\Prism_{\R \dquot (f_1,\ldots,f_d)}) \ar{u} \ar{r} \ar[->>]{d} & \Spa(\Prism_{\R},\Prism_{\R}) \ar[->>]{d} \ar{u} \\
\X^{\Prism,\Box} \ar{r}{i_\Prism} & \Y^{\Prism,\Box}
\end{display}
On the bottom, we have used regularity of $i$ to obtain that the square is Cartesian (and compatibility of $(-)^{\Prism,\Box}$ with finite derived limits). Lemma \ref{lem:prismatic_envelope_integral} and \cite[Example 2.1.13]{Kipp2} then show that the top arrow is affine proper, and the middle arrow is a base change of the top arrow. Thus the claim follows by applying \cite[Lemma 1.50]{MikamiDuality}.
\end{proof}

\begin{remark} \label{rem:kipp_strengthening}
The previous argument can strengthen the results in \cite{Kipp2} to give a canonical identification $\R f_{\Prism,!} \simeq \R f_{\Prism,*}$ when $f$ is proper. We want to distinguish this from cohomological properness: the morphisms we consider need not be truncated, as is usually included in the definition of a cohomologically proper morphism. The previous argument establishes that for a regular closed immersion $i:\X \to \Y$ of derived formal schemes the induced map on prismatizations $i_{\Prism}: \X^{\Prism,\Box} \to \Y^{\Prism,\Box}$ satisfies this property, i.e. $\R i_{\Prism,!}\simeq \R i_{\Prism,\ast}$ canonically. As pointed out to the author by Niklas Kipp, using \cite[Theorem 1.1]{MaoDerivedCrystalline} one may show the claim holds for all closed immersions. Indeed by \'{e}tale descent we can reduce to the affine case. Then in the previous argument, we may first pass to a semiperfectoid cover $\Spf \R_0$ and then \cite[Theorem 1.1]{MaoDerivedCrystalline} yields a sifted colimit presentation in terms of regular closed immersions on derived $p$-complete animated rings over $\R_0$. Then using the fact that prismatization is compatible with derived limits and sifted colimits of integral maps are integral (using the integrality criterion for affine properness explained in \cite[Example 2.1.13]{Kipp2}), Lemma \ref{lem:prismatic_stack_cohom_smooth_proper} shows $\R i_{\Prism,!}\simeq \R i_{\Prism,\ast}$ canonically.

The claim for derived closed immersions suffices in general, as the strategy in \cite{Kipp2} is to show cosmoothness of $f_{\Prism}$ and weak cohomological properness of the diagonal, but we have now established in place of weak cohomological properness that $\R (\Delta_f)_{\Prism,!}\simeq \R (\Delta_f)_{\Prism,\ast}$. Thus the same argument establishes that if $f: \X \to \Y$ is proper then $f_{\Prism}: \X^{\Prism,\Box} \to \Y^{\Prism,\Box}$ satisfies $\R f_{\Prism,!}\simeq \R f_{\Prism,\ast}$.
\end{remark}

To actually make use of these stacks, we must show that they recover solid Laurent crystals defined using the site.

\begin{proposition} \label{prop:stack_site_compare}
For every log formal scheme $(\X,\mathcal{M})$ satisfying Assumption \ref{ass:horizontal_nc}, there is a natural equivalence
\[\mathrm{an}_{(\X,\mathcal{M})}: \D_{\Box}((\X,\mathcal{M})^{\mathcal{L},\Box}) \xlongrightarrow{\sim} \D_{\Box}((\X,\mathcal{M})_\Prism,\Laur).\]
\end{proposition}

\begin{proof}
First, the equivalence can be tested locally because both sides satisfy strict \'{e}tale descent. For the target this is shown in the appendix in Corollary \ref{cor:laur_strict_etale}. For the source, the stack itself sends strict \'{e}tale covers to effective epimorphisms of analytic stacks, using compatibility of $\X^{\Prism,\Box}$ with \'{e}tale covers from \cite[Theorem 4.1.3(D)]{Kipp2} (again using Example 4.3.2 of \cite{Kipp2} to deduce analogous claims for the prismatization) and base changing to obtain the claim for $(\X,\mathcal{M})^{\Prism,\Box}$.

Thus we may work on a framed log affine. In this case, there exists a choice of log Breuil-Kisin prism atlas $(\rA^+,\I,\mathcal{N})^a$ for this framed log affine such that we obtain the presentation
\[(\X,\mathcal{M})^{\Prism,\Box} = \colim_n \Spa(\rA^{+,(n)},\rA^{+,(n)})\]
proved in Proposition \ref{prop:log_BK_analytic_presentation}. Base changing to $(\X,\mathcal{M})^{\mathcal{L},\Box}$ we obtain an analogous presentation, but replace the terms $\Spa(\rA^{+,(n)},\rA^{+,(n)})$ with $\colim_m \Spa(\rA^{+,(n)}/p^m[1/\I],\rA^{+,(n)}/p^m)$.

In Definition \ref{def:laur_category}, the category $\D_{\Box}(\colim_m \Spa((\rA^{+,(n)}/p^m)[1/\I],\rA^{+,(n)}/p^m))$ is precisely what we assign to the log prisms $\rA^{+,(n)}$ appearing in the \v{C}ech nerve of $(\rA^+,\I,\mathcal{N})^a$. Combining with Theorem \ref{thm:laur_cech}, we obtain
\begin{align*}
\D_{\Box}((\X,\mathcal{M})_\Prism,\Laur) &\simeq \lim_n \D_{\Box}(\colim_m \Spa((\rA^{+,(n)}/p^m)[1/\I],\rA^{+,(n)}/p^m)) \\
&\simeq \D_{\Box}((\X,\mathcal{M})^{\mathcal{L},\Box}).
\end{align*}
Thus we get an identification for each log Breuil-Kisin prism. This identification is also independent of choices. One can see this by either assembling together the identifications for every log Breuil-Kisin prism into a canonical one, or noting that a canonical functor exists already. Indeed, for some $(\rA,(d),\mathcal{N})^a$ in the strict log prismatic site, applying the functor $(-)^{\Box}$ (Definition \ref{def:box_analytic_realization}) to the map of Lemma \ref{lem:log_prism_stack_map} gives
\[\Spa(\rA,\rA) \longrightarrow (\X,\mathcal{M})^{\Prism,\Box},\]
so passing to Laurent stacks yields a canonical functor
\[\D_{\Box}((\X,\mathcal{M})^{\mathcal{L},\Box}) \longrightarrow \D_{\Box}((\X,\mathcal{M})_\Prism,\Laur).\]
This we call $\mathrm{an}_{(\X,\mathcal{M})}$.
\end{proof}

Again under Assumption \ref{ass:horizontal_nc}, it is straightforward to obtain a pushforward compatibility result.

\begin{lemma} \label{lem:kipp_perf_compare}
The equivalence of Proposition \ref{prop:stack_site_compare} restricts to $\Perf$. 

Moreover, for strict log smooth proper maps $f$ the site-theoretic pushforward $\R f_{\Prism,*}$ preserves $\Perf$, and
\[\mathrm{an}_{(\Y,\mathcal{N})} \R f_{\mathcal{L},*} \simeq \R f_{\Prism,*}\mathrm{an}_{(\X,\mathcal{M})}\]
for any input in $\D_{\Box}((\X,\mathcal{M})^{\mathcal{L},\Box})$ where $f_{\mathcal{L}}$ denotes the induced map on Laurent stacks.
\end{lemma}

\begin{proof}
This is formal from the previous results. For the first point, we constructed the functor through pullbacks of stacks so it will be symmetric monoidal, and thus the equivalence restricts to dualizables. For the second point, this follows by taking right adjoints and observing $\mathrm{an}$ has pullback compatibility.

Once we know this, preservation of $\Perf$ is true on the stack side by applying Lemma \ref{lem:prismatic_stack_cohom_smooth_proper} and Lemma \ref{lem:formal_pushforward_duality}.
\end{proof}

We also obtain Poincar\'{e} duality, where we can obtain $f^!$ on the site side by transporting through Proposition \ref{prop:stack_site_compare}.

\begin{corollary}
Let $f:(\X,\mathcal{M})\to(\Y,\mathcal{N})$ be as in Lemma \ref{lem:prismatic_stack_cohom_smooth_proper}, and suppose both source and target satisfy Assumption \ref{ass:horizontal_nc}. Then $f^!\mathcal{F}\simeq f^*\mathcal{F}\{d\}[2d]$ for $\mathcal{F}\in\Perf((\Y,\mathcal{N})_\Prism,\Laur)$, and
\[\R f_{\Prism,*}\mathcal{E} \simeq \RcalHom_{\Laur}(\R f_{\Prism,*}(\mathcal{E}^{\vee}\{d\}),\Laur[-2d])\]
for every $\mathcal{E}\in\Perf((\X,\mathcal{M})_\Prism,\Laur)$.
\label{cor:laurent_prismatic_pd}
\end{corollary}

\begin{proof}
Transporting Lemma \ref{lem:prismatic_stack_cohom_smooth_proper} through Lemma \ref{lem:kipp_perf_compare} gives $f^!\mathcal{F}\simeq f^*\mathcal{F}\{d\}[2d]$. Lemma \ref{lem:formal_pushforward_duality}, together with the pushforward comparison in Lemma \ref{lem:kipp_perf_compare}, gives
\[(\R f_{\Prism,*}\mathcal{E})^\vee\simeq\R f_{\Prism,*}(\mathcal{E}^\vee\{d\})[2d].\]
Dualizing, this equivalence gives the claim.
\end{proof}

We also note that the Poincar\'{e} duality pairing is compatible with Frobenius: the identification $f^!\Laur\simeq\Laur\{d\}[2d]$ is Frobenius equivariant since it is upgraded to an identification of $\rF$-gauges in \cite[Theorem 4.1.3(F)]{Kipp2}.

\subsection{Pushforward compatibility}

Let $f:(\X,\mathcal{M}_\X) \to (\Y,\mathcal{M}_\Y)$ be a strict log smooth proper morphism of pure relative dimension $d$ between log $p$-adic formal schemes whose adic generic fibers have pure dimensions $d_\X$ and $d_\Y$. We assume $(\X,\mathcal{M}_\X)$ and $(\Y,\mathcal{M}_\Y)$ satisfy Assumption \ref{ass:horizontal_nc} over one of the bases $\Spf \W(k)[\zeta_{p^\ell}]$.

Let $d=d_\X-d_\Y$. In such a setup, define
\[\R f_+(-) := \R f_*(-)\{d\}[d].\]
We will show that $\R f_+$ induces a functor
\[\R f_+: \D_{\fgu}((\X,\mathcal{M}_\X)_\Prism,\Laur) \longrightarrow \D_{\fgu}((\Y,\mathcal{M}_\Y)_\Prism,\Laur).\]
After this, the pushforward compatibility result is formal.

We first need some preliminaries on relative stacks. Let $(\Y,\mathcal{M}_\Y)$ be a framed log affine with associated logarithmic $q$-de Rham prism $(\rA^+,\I,\mathcal{N})^a$, so $\Y=\Spf(\rA^+/\I)$ (with the induced log structure), and let $(\X,\mathcal{M}_\X) \to (\Y,\mathcal{M}_\Y)$ be strict smooth. By Proposition \ref{prop:log_BK_analytic_presentation}, the map $\rho_{\rA^+}:\Spa(\rA^+,\rA^+)\to(\Y,\mathcal{M}_\Y)^{\Prism,\Box}$ is a representable proper effective epimorphism. We will also consider the relative variant $((\X,\mathcal{M}_\X)/\rA^+)^{\Prism,\Box}$ of these log prismatic stacks, defined via the Cartesian square
\begin{display} 
((\X,\mathcal{M}_\X)/\rA^+)^{\Prism,\Box} \ar{r} \ar{d}{f_{\rA^+}} & (\X,\mathcal{M}_\X)^{\Prism,\Box} \ar{d} \\
\Spa(\rA^+,\rA^+) \ar{r}{\rho_{\rA^+}} & (\Y,\mathcal{M}_\Y)^{\Prism,\Box}.
\end{display}
We define $((\X,\mathcal{M}_\X)/\rA)^{\mathcal{L},\Box}$ by further base change over $\Spa(\rA^+,\rA^+)$ to 
\[(\Spf(\rA^+/\I)/\rA)^{\mathcal{L},\Box}=\AnSpf\rA=\colim_n\Spa((\rA^+/p^n)[1/\I],\rA^+/p^n).\]
Equivalently, we replace all terms in the Cartesian square with their Laurent versions. The method of Proposition \ref{prop:stack_site_compare} identifies $\Perf(((\X,\mathcal{M}_\X)/\rA)^{\mathcal{L},\Box})$ with $\Perf(((\X,\mathcal{M}_\X)/\rA^+)_\Prism,\Laur)$, by locally using the \v{C}ech presentation of the stack induced by a Breuil-Kisin prism for $\X$ and identifying it with the site-theoretic one.

There is a natural Frobenius $\varphi$ on $((\X,\mathcal{M}_\X)/\rA)^{\mathcal{L},\Box}$ over $\varphi_{\rA}$ on $(\Spf(\rA^+/\I)/\rA)^{\mathcal{L},\Box}$. In what follows the superscript $(\varphi)$ denotes base change along this $\varphi_{\rA}$. If $f_{\rA}$ denotes the structure map to $(\Spf(\rA^+/\I)/\rA)^{\mathcal{L},\Box}$, we use $\mathrm{Fr}$ or \emph{relative Frobenius} to mean the induced map 
\[\mathrm{Fr}=(\varphi,f_{\rA}): ((\X,\mathcal{M}_\X)/\rA)^{\mathcal{L},\Box} \to ((\X,\mathcal{M}_\X)/\rA)^{\mathcal{L},\Box,(\varphi)}.\] 

The following proposition is an analogue of Cartier descent. There are some variants already existing which do not quite suffice for our applications, such as in \cite{GLQcartier}. See especially \cite[Corollary 6.8]{GLQcartier}, which is similar to the general approach taken below. The proof is roughly a $q$-analogue of a typical proof of Cartier descent.

\begin{proposition} \label{prop:laurent_cartier_descent}
In the above setting, relative Frobenius induces a fully faithful functor
\[\mathrm{Fr}^*:\Perf(((\X,\mathcal{M}_\X)/\rA)^{\mathcal{L},\Box,(\varphi)}) \longrightarrow \Perf(((\X,\mathcal{M}_\X)/\rA)^{\mathcal{L},\Box}).\]
\end{proposition}

\begin{proof}
We may work strict \'{e}tale locally by descent of the source and target (this is again automatic from the claim for $(\X,\mathcal{M}_\X)^{\Prism,\Box}$ by base change): thus we assume that $\X$ is a framed log affine relative to $(\Y,\mathcal{M}_\Y)$ with relative coordinates $\T_1,\ldots,\T_d$. More precisely, we have 
\[(\X,\mathcal{M}_\X) \longrightarrow (\Y,\mathcal{M}_\Y) \times_{\Spf \W(k)[\zeta_{p^\ell}]} \widehat{\mathbf{G}}_m^d \longrightarrow \mathbb{T}^{n+d}_r\]
where $\widehat{\mathbf{G}}_m^d$ has a trivial log structure. Thus $(\X,\mathcal{M}_\X)$ inherits a framing, and has a corresponding log $q$-de Rham prism with ring $\B^+$ and $\B=\B^+[1/(q-1)]_p^{\wedge}$. A straightforward variant of Proposition \ref{prop:Gamma_mod_ff} gives a proper morphism $\rho:[\AnSpf\B/\Z_p^d]\to((\X,\mathcal{M}_\X)/\rA)^{\mathcal{L},\Box}$ with a fully faithful pullback where $\Z_p^d$ fixes $\rA$ and $\gamma_i \cdot \T_j =q^{p^\ell \delta_{ij}}\T_j$ for the generators $\gamma_i$ of $\Z_p^d$. For $\ell=0$, we use the same chart with $\I=([p]_q)$, which still gives a cover over $\W(k)$. In this case $\gamma_i\cdot \T_j =q^{p\delta_{ij}}\T_j$ and $\varphi(\I)=([p]_{q^p})$; the proof is similar for this case so we focus on $\ell \ge 1$.

Now we observe that Proposition \ref{prop:Gamma_mod_ff} generalizes directly to the relative setting with $\Z_p^d$ (indeed one just omits the harder ``cyclotomic'' component of the calculation), so in particular we obtain a fully faithful functor
\[\rho^*: \Perf(((\X,\mathcal{M}_\X)/\rA)^{\mathcal{L},\Box}) \longrightarrow \Perf([\AnSpf\B/\Z_p^d]) \simeq \Perf(\B)^{h\Z_p^d}.\]
Since we deal with perfect complexes, we may also dispense with any condensed mathematics here using the argument of Lemma \ref{lem:perf_comparison}. One similarly obtains a fully faithful evaluation functor $\rho^{*,(\varphi)}: \Perf(((\X,\mathcal{M}_\X)/\rA)^{\mathcal{L},\Box,(\varphi)}) \longrightarrow \Perf(\B^{(\varphi)})^{h\Z_p^d}$ for the Frobenius twist $\B^{(\varphi)}$. These functors will reduce the comparison below to a group cohomology calculation.

We may first reduce to checking full faithfulness after derived reduction modulo $p$ by derived Nakayama; we use $\overline{(-)}$ to indicate that we work modulo $p$ as usual. Since we work with perfect coefficients, as in Proposition \ref{prop:open_decompletion_hom} the tensor/Hom adjunction combined with a Postnikov d\'{e}vissage reduces us to verifying
\[\RHom_{\Perf(((\X,\mathcal{M}_\X)/\rA)^{\mathcal{L},\Box,(\varphi)})}(\O,\overline{\mathcal{E}}) \simeq \RHom_{\Perf(((\X,\mathcal{M}_\X)/\rA)^{\mathcal{L},\Box})}(\O,\mathrm{Fr}^* \overline{\mathcal{E}})\]
for $\overline{\mathcal{E}}\in \Coh(((\X,\mathcal{M}_\X)/\rA)^{\mathcal{L},\Box,(\varphi)})$ annihilated by $p$. The relevant $t$-structure is produced by Corollary \ref{cor:laur_t_struct}. Associated to $\overline{\mathcal{E}}$ there is a coherent $\overline{\B}^{(\varphi)}$-module with an action by $\Z_p^d$, given by evaluation on the log $q$-de Rham prism.

For any coherent $\Z_p^d$-module $\overline{\M}$, by compactness and finite generation of $\overline{\M}$ we can pick a $\Z_p^d$-stable lattice $\overline{\M}^+\subset \overline{\M}$ which is finitely generated over $\overline{\B}^{+,(\varphi)}$. On the heart, using the argument of Lemma \ref{lem:Cq_criterion} with coefficients (as extended in the proof of Lemma \ref{lem:Cech_finite_coefficients}) we see that evaluation of $\overline{\mathcal{E}} \in \Coh(((\X,\mathcal{M}_\X)/\rA)^{\mathcal{L},\Box, (\varphi)})$ on the log $q$-de Rham prism yields a coherent $\Z_p^d$-module $\overline{\M}$ where the coaction $\overline{\M} \to \rC(\Z_p^d,\overline{\M})$ lands in 
\[\left(\bigoplus_{\alpha \in \N^d} \overline{\M}^+(q^{p^\ell}-1)^{|\alpha|}\binom{x}{\alpha}\right)^{\wedge}_{(q-1)}[1/(q-1)] \subset \rC(\Z_p^d, \overline{\M}).\]
We have used here that $\rC(\Z_p^d, \overline{\M}^+)[1/(q-1)] = \rC(\Z_p^d, \overline{\M})$, again using compactness to ensure bounded denominators. We may reduce the Proposition to checking whether
\[\RGamma(\Z_p^d,\mathrm{Fr}^* \overline{\M}/\overline{\M}) \simeq \bigoplus_{0\neq \beta \in \{0,\ldots,p-1\}^d} \RGamma(\Z_p^d,\T^\beta \overline{\M}) \simeq 0\]
since each summand $\T^\beta \overline{\M}$ of $\mathrm{Fr}^* \overline{\M}/\overline{\M}$ is stable under the action of $\Z_p^d$. By the earlier result about the coaction for $\overline{\M}$, the operator $\partial_{i,\overline{\M}} = \frac{\gamma_i-1}{q^{p^\ell}-1}$ satisfies $\partial_i^k \longrightarrow 0$ on $\overline{\M}$ (and moreover preserves each summand). For the summand corresponding to $\beta$ there is some $i$ so $\beta_i\neq 0$, and it suffices to check the induced action of $\partial_{i,\overline{\M}}$ is invertible for some $\beta_i\neq 0$ to deduce the cohomology vanishes (as then $\fib(\gamma_i-1)=0$). The induced action of $\partial_{i,\overline{\M}}$ on the summand $\T^\beta \overline{\M}$ for $\beta \neq 0$ is by $[\beta_i]_{q^{p^\ell}} + q^{p^\ell \beta_i} \partial_{i,\overline{\M}}$. But since $\partial_{i,\overline{\M}}^k \longrightarrow 0$, this operator is invertible and the claim follows.
\end{proof}

\begin{remark}
A similar approach to Proposition \ref{prop:Gamma_mod_ff} gives a slightly longer argument that shows full faithfulness by checking $\R \mathrm{Fr}_* \O \simeq \O$. The ultimate cohomology calculation is essentially the same.
\end{remark}

\begin{proposition} \label{prop:frob_bc}
Let $f:(\X,\mathcal{M}_\X)\to(\Y,\mathcal{M}_\Y)$ be a strict log smooth proper morphism between log formal schemes satisfying Assumption \ref{ass:horizontal_nc} over $\W(k)[\zeta_{p^\ell}]$ with $\ell\ge0$. For every $\mathcal{E}\in\Perf((\X,\mathcal{M}_\X)_\Prism,\Laur)$, the canonical map
\[\varphi^*\R f_*\mathcal{E} \longrightarrow \R f_*\varphi^*\mathcal{E}\]
is an isomorphism.
\end{proposition}

\begin{proof}
Using proper base change and strict \'{e}tale descent, we may reduce to the case where $(\Y,\mathcal{M}_\Y)$ is a framed log affine with an associated log $q$-de Rham prism. In the $\ell=0$ case this log $q$-de Rham prism $(\rA^+,\I,\M_{\Spf \rA^+})$ is not a Breuil-Kisin atlas in the sense that $\Spf \rA^+/\I$ with its induced log structure is not $(\Y,\mathcal{M}_\Y)$. However, in this case it is still a cover so the argument is similar and we will therefore focus on the case $\ell > 0$.

There is a canonical comparison map $\varphi^*\R f_*\mathcal{E} \longrightarrow \R f_*\varphi^*\mathcal{E}$ where the source and target are perfect by Lemma \ref{lem:kipp_perf_compare}. We may test if this is an isomorphism after pullback along the proper effective epimorphism $\AnSpf \rA \to (\Y,\mathcal{M}_\Y)^{\mathcal{L},\Box}$. Note that this pullback commutes with pushforwards by properness. Writing $f_{\rA}$ for the structure map $(\X/\rA)^{\mathcal{L},\Box} \to \AnSpf \rA = (\Spf \rA^+/\I / \rA)^{\mathcal{L},\Box}$, we have therefore reduced to checking the canonical map
\[\varphi_{\rA}^*\R f_{\rA,*}\mathcal{E} \longrightarrow \R f_{\rA,*}\varphi^*\mathcal{E}\]
is an isomorphism, where by abuse of notation we retain $\mathcal{E}$ for the pullback to the relative stack $((\X,\mathcal{M}_\X)/\rA)^{\mathcal{L},\Box}$. We note here that Lemma \ref{lem:prismatic_stack_cohom_smooth_proper} shows $\R f_{\rA,*}=\R f_{\rA,!}$.

Consider now the following commutative diagram, which has a Cartesian right square:
\begin{display}
((\X,\mathcal{M}_\X)/\rA)^{\mathcal{L},\Box} \ar{r}{\mathrm{Fr}} \ar[swap]{dr}{f_{\rA}} & ((\X,\mathcal{M}_\X)/\rA)^{\mathcal{L},\Box,(\varphi)} \ar{r}{\widetilde{\varphi}_{\rA}} \ar{d}{f_{\rA}^{(\varphi)}} & ((\X,\mathcal{M}_\X)/\rA)^{\mathcal{L},\Box} \ar{d}{f_{\rA}} \\
& \AnSpf \rA \ar{r}{\varphi_{\rA}} & \AnSpf \rA
\end{display}
On the top row, we also have $\widetilde{\varphi}_{\rA} \circ \mathrm{Fr}=\varphi$. As everything is in $\Perf$, it suffices to know for perfect test objects that $\RHom(-,\varphi_{\rA}^*\R f_{\rA,*}\mathcal{E}) \simeq \RHom(-,\R f_{\rA,*}\varphi^*\mathcal{E})$. By adjunction, the target is also $\RHom(f_{\rA}^* -, \varphi^*\mathcal{E})$.

Applying proper base change for the right square followed by adjunction we deduce
\[\RHom(-,\varphi_{\rA}^*\R f_{\rA,*}\mathcal{E}) \simeq \RHom(-, \R f_{\rA,*}^{(\varphi)} \widetilde{\varphi}_{\rA}^* \mathcal{E}) \simeq \RHom(f_{\rA}^{(\varphi),*} -, \widetilde{\varphi}_{\rA}^* \mathcal{E}).\]
We wanted to show this was $\RHom(f_{\rA}^* -, \varphi^*\mathcal{E})$. Upon inspecting the previous commutative diagram, the claim then follows from full faithfulness of $\mathrm{Fr}^*$ in Proposition \ref{prop:laurent_cartier_descent}.
\end{proof}

\begin{theorem} \label{thm:fgu_pushforward}
Let $f: (\X,\mathcal{M}_\X) \to (\Y,\mathcal{M}_\Y)$ be strict log smooth proper, where the source and target satisfy Assumption \ref{ass:horizontal_nc} over $\Spf \W(k)[\zeta_{p^\ell}]$. Then $\R f_+$ induces a functor
\[\R f_+: \D_{\fgu}((\X,\mathcal{M}_\X)_\Prism,\Laur) \longrightarrow \D_{\fgu}((\Y,\mathcal{M}_\Y)_\Prism,\Laur).\]
\end{theorem}

\begin{proof}
Since the finitely generated unit category is stable under shifts and twists we may test the assertion for $\R f_*$. We may also test the assertion modulo $p$ by derived Nakayama, and further reduce to the case where the input $\mathcal{E}$ lies in the heart and is killed by $p$. Then by Corollary \ref{cor:global_root_fgu_condition}, there is a coherent root $\mathcal{E}_0$ of $\mathcal{E}$ so that $\mathcal{E}\simeq \colim_n \varphi^{*n} \mathcal{E}_0$.

Since $\R f_* \simeq \R f_!$ canonically, by applying Proposition \ref{prop:stack_site_compare} we see the pushforward commutes with colimits. In particular, using the previous root presentation of $\mathcal{E}$ we obtain
\[\R f_* (\colim_n \varphi^{*n} \mathcal{E}_0) \simeq \colim_n \R f_* \varphi^{*n} \mathcal{E}_0.\]
By Proposition \ref{prop:frob_bc} we see this is the same as $\colim_n \varphi^{*n} \R f_* \mathcal{E}_0$. Finally, using Lemma \ref{lem:kipp_perf_compare} we see $\R f_* \mathcal{E}_0$ is perfect, which proves the result is finitely generated unit using the same argument as Footnote \ref{fn:root_unitalization_fgu}.
\end{proof}

\begin{theorem} \label{thm:fgu_pushforward_compat}
Let $f:(\X,\mathcal{M}_\X)\to(\Y,\mathcal{M}_\Y)$ be a strict log smooth proper morphism of pure relative dimension $d$ where both satisfy Assumption \ref{ass:horizontal_nc} over $\Spf \W(k)[\zeta_{p^\ell}]$ whose log adic generic fibers are pure of dimensions $d_\X$ and $d_\Y$ respectively. Let $\D_\X$ and $\D_\Y$ be horizontal normal crossings divisors inducing $\mathcal{M}_\X$ and $\mathcal{M}_\Y$. There is a commutative diagram
\begin{display}
\D_{\fgu}((\X,\mathcal{M}_\X)_\Prism, \Laur)^\op \ar{r}{\T_{\et,\X}[d_\X]} \ar{d}{\R f_* \{d\} [d]} & \D^{(b)}_{\cons,\D_\X}((\X,\mathcal{M}_\X)_\eta,\Z_p) \ar{d}{\R f_{\eta,*}} \\
\D_{\fgu}((\Y,\mathcal{M}_\Y)_\Prism, \Laur)^\op \ar{r}{\T_{\et,\Y}[d_\Y]} & \D^{(b)}_{\cons,\D_\Y}((\Y,\mathcal{M}_\Y)_\eta,\Z_p).
\end{display}
\end{theorem}

\begin{proof}
Recall we set $\R f_+(-) := \R f_*(-)\{d\}[d]$. Our strategy will be to deduce pushforward compatibility from pullback compatibility (from Corollary \ref{cor:pullback}) by showing $\R f_+[d]$ is left adjoint to $\L f^*$. Indeed we may work on Laurent stacks through Proposition \ref{prop:stack_site_compare}, and there one may calculate using Lemma \ref{lem:prismatic_stack_cohom_smooth_proper}
\[
\begin{aligned}
\RHom(\R f_+\mathcal{E}[d],\mathcal{F})
&\simeq\RHom(\R f_*(\mathcal{E}\{d\}[2d]),\mathcal{F})\\
&\simeq\RHom(\mathcal{E}\{d\}[2d],f^!\mathcal{F})
\simeq\RHom(\mathcal{E},\L f^*\mathcal{F}).
\end{aligned}
\]
Everything in the previous calculation is Frobenius-equivariant (see the remark after Corollary \ref{cor:laurent_prismatic_pd}), so the adjunction persists after applying $(-)^{\varphi=1}$ to all categories in sight.

Thus, we see that $\R f_+[d]$ is left adjoint to $\L f^*$. By pullback compatibility from Corollary \ref{cor:pullback}, we get
\[\T_{\et,\X}^{-1} \circ \L f_\eta^* \rightleftarrows \T_{\et,\Y}\circ(\R f_+[d])\]
by transporting the adjunction across the antiequivalence in Theorem \ref{thm:RH_main}. We implicitly use here that $\R f_+[d]$ preserves the crystal side from Theorem \ref{thm:fgu_pushforward}, and that $\L f_\eta^*$ preserves $p$-Kummer sheaves. Uniqueness of right adjoints therefore gives
\[\T_{\et,\Y}(\R f_+\mathcal{E}[d])\simeq\R f_{\eta,*}\T_{\et,\X}(\mathcal{E})\]
once we know $\R f_{\eta,*}$ preserves $\D^{(b),\mathrm{pKum}}_{\mathrm{cons},\D}(-,\Z_p)$. Since strict smooth proper pushforwards preserve $\D^{(b)}_{\mathrm{lisse}}((-)_{\Lambda,\eta},\Z_p)$ (i.e. derived lisse sheaves on strata), we only need to check the $p$-Kummer condition is preserved: for example one may use a spectral sequence to understand the prime-to-$p$ log inertia on the pushforward at a geometric stalk.\footnote{If one uses a spectral sequence to argue this, it is important to note invariants for the prime-to-$p$ log inertia are exact. In particular it will again act trivially on extensions, which is why on the abutment of the spectral sequence it will still act trivially if it does on the $\E_2$ page.} The desired theorem then follows after noting $d_\X = d + d_\Y$ and keeping track of shifts.
\end{proof}

While we have primarily been interested in bases such as $\W(k)[1/p]$, we remark that it is possible to also work over a fairly general perfectoid base. In fact, over a perfectoid base the main results can be done more generally in the semistable case.

\begin{remark} \label{rem:perfectoid_base}
The generality we work in can be improved over $\O_{\C_p}$. For this, we introduce a new prism $\A_{\inf,\R}^{\sq}$, similar to the construction in \cite{MinWangRelativePhiGamma}. For 
\[\R=\O_{\C_p} \langle \T_1,\ldots,\T_r, \T_{r+1}^{\pm1},\ldots,\T_n^{\pm1}\rangle/(\T_1 \cdots \T_s -p)\] 
we use the log prism
\[(\Ainf(\O_{\C_p}) \langle \T_1,\ldots,\T_r, \T_{r+1}^{\pm1},\ldots,\T_n^{\pm1} \rangle_{(p,\I)}/(\T_1 \cdots \T_s - [p^\flat]), \I, (\O_{\C_p^\flat} \setminus 0) \oplus_{\N} \N^r )^a\]
where $\I = \ker(\theta: \Ainf(\O_{\C_p}) \to \O_{\C_p})$. Here $s \le r$, and we use the canonical lift $[\cdot]: \O_{\C_p}^\flat \to \Ainf(\O_{\C_p})$ to let the multiplicative monoid $(\O_{\C_p^\flat} \setminus 0)$ define the desired prelog structure. Replacing the $q$-de Rham prism by this prism nearly proves all results in the article verbatim for a quasicompact semistable $(\X,\mathcal{M})$ over $(\Spf \O_{\C_p}, \O_{\C_p} \setminus 0)$ by replacing $\Gamma$ with $\Z_p^{n-1}$. Indeed, while even non-Noetherian in general before Laurentizing, after Laurentizing the corresponding ring again has a regular $\rF$-finite reduction modulo $p$.

The only point in the argument where we encounter an issue is with Olsson's stack of log structures $\Log$. Strictly speaking, to consider the log prismatic stack with this definition one must insist that the log structure is fine. If we instead use the moduli stack of relatively fine log structures $\Log_{(\Spf \Ainf(\O_{\C_p}), \O_{\C_p}^\flat \setminus 0)}$ as the base, there is no issue and all arguments proceed the same. By this, we mean maps $\X \to \Log_{(\Spf \Ainf(\O_{\C_p}), \O_{\C_p}^\flat \setminus 0)}$ classify integral log structures $(\X,\mathcal{M}) \xlongrightarrow{f} (\Spf \Ainf(\O_{\C_p}), \O_{\C_p}^\flat \setminus 0)$ such that there is \'{e}tale locally a surjection $f^* (\O_{\C_p}^\flat \setminus 0) \oplus \N^r \to \mathcal{M}$.
\end{remark}

\appendix

\section{Log prismatic stacks}

\subsection{The log prismatic stack}
\label{app:log_prismatic_stacks}

In this appendix by formal stacks we mean accessible fpqc sheaves in groupoids (or more generally spaces) on animated $p$-nilpotent rings.

We denote by $\Log$ the $p$-completion of Olsson's stack of fine log structures over $\Z_p$ introduced in \cite{OlssonLogStacks}. For every fine monoid $\rP$, we have an associated toric stack
\[\mathcal{A}_{\rP}=[\Spf\Z_p\langle \rP\rangle/\Spf\Z_p\langle \rP^{\mathrm{gp}}\rangle].\]
These $\mathcal{A}_{\rP}$ cover $\Log$ by \cite[Corollary 5.25]{OlssonLogStacks}. One may check that these, and in fact $\Log$, are formal $\delta$-stacks in the sense of \cite[Construction 10.20]{AKNRelativeDelta}.

The general philosophy behind the log prismatic stack is easy to motivate, the principal idea being that the data of a log formal scheme $\X$ with a fine log structure amounts to a map $\X \to \Log$. Moreover, in this interpretation many logarithmic cohomology theories $\H_{\log}(-)$ can often be produced by taking a relative version $\H(\X/\Log)$ using the ordinary cohomology theory $\H(-)$. For example, logarithmic de Rham cohomology can be recovered this way. In \cite{AKNRelativeDelta} they define prismatic stacks relative to $\delta$-rings, but more generally, if $\X\to\mathcal{S}$ is a morphism to a formal $\delta$-stack, we may write
\[(\X/\mathcal{S})^\Prism:=\X^\Prism\times_{\mathcal{S}^\Prism}(\mathcal{S}\times\Z_p^\Prism)\]
for its \emph{relative prismatization}. The point we explain next is what the map $\Z_p^\Prism \times \mathcal{S} \to \mathcal{S}^\Prism$ is.

We will use this construction with $\mathcal{S}=\Log$, which we explain in more detail. We must specify a map $\Log\times\Z_p^\Prism\to\Log^\Prism$, as the other map is simply induced by functoriality of prismatization. Let $\R$ be a $p$-nilpotent ring and let $\mathcal{M}_\R\in\Log(\R)$. Regarding Witt vectors as the cofree $\delta$-ring and using the $\delta$-stack structure on $\Log$ gives
\[\Log(\R)\simeq\Map_{\delta}(\Spf\W(\R),\Log)\longrightarrow\Map(\Spf\W(\R),\Log).\]
Let $\mathcal{M}_{\W(\R)}$ be the resulting log structure. If $\rP\to\mathcal{M}_\R$ is a fine chart with structure map $\alpha:\rP\to\R$, then locally
\[\mathcal{M}_{\W(\R)}=\left(\rP \xlongrightarrow{m\mapsto[\alpha(m)]} \W(\R)\right)^a.\]
We then define the map as
\[\pi_{\Log}(\mathcal{M}_{\R},\I\to\W(\R))=\left(\I\to\W(\R),\mathcal{M}_{\W(\R)}|_{\Spec(\W(\R)\dquot \I)}\right).\]
Then fpqc sheafification yields $\pi_{\Log}:\Log\times\Z_p^\Prism\to\Log^\Prism$. We remark here that $\W(\R)\dquot \I$ has a slightly different meaning, explained in \cite[\S 2.1]{APC2}.

\begin{definition} \label{def:Olsson_stack}
Let $(\X,\mathcal{M})$ be a bounded fs log formal scheme. Then there is a map $\X \to \Log$ induced by the log structure. We define $(\X,\mathcal{M})^{\Prism}$ to be the fiber product
\begin{display}
(\X,\mathcal{M})^{\Prism} \ar{d} \ar{r} & \X^{\Prism} \ar{d} \\
\Log \times \Z_p^\Prism \ar{r}{\pi_{\Log}} & \Log^{\Prism}.
\end{display}
This is the log prismatic stack introduced by Olsson in
\cite{OlssonLogPrismaticStack}.\footnote{We note that the definition of this stack was also considered by Koshikawa.}
\end{definition}

Throughout this appendix, we work in the situation of Assumption \ref{ass:horizontal_nc} over $\O_K$, where $K/\W(k)[1/p]$ is a finite totally ramified extension.

Strict \'{e}tale locally, every such log formal scheme is a framed log affine, for which there is a logarithmic Breuil-Kisin atlas. We define this terminology as follows.

\begin{definition}
Let $(\X,\mathcal{M})$ satisfy Assumption \ref{ass:horizontal_nc} over $\O_K$. A logarithmic Breuil-Kisin prism for $(\X,\mathcal{M})$ is a log prism $(\rA,\I,\M_{\Spf \rA})^a$ such that
\[(\Spf \rA/\I, \M_{\Spf \rA/\I}) \longrightarrow (\X,\mathcal{M})\]
is a strict open immersion of log $p$-adic formal schemes.
\end{definition}

A logarithmic Breuil-Kisin atlas for $(\X,\mathcal{M})$ means this is an isomorphism.

\begin{lemma} \label{lem:standard_log_BK_prism}
Let $(\X,\mathcal{M})=(\Spf\R,\mathcal{M})$ be a framed log affine over $\O_K$. Then there is a logarithmic Breuil-Kisin prism $(\rA,\I,\M_{\Spf\rA})^a$ which is an atlas for $(\X,\mathcal{M})$.
\end{lemma}

\begin{proof}
This construction is also explained in \cite[Construction 5.10]{Inoue}, but we explain it again here for convenience of the reader. Write $\R\simeq\R_0\widehat{\otimes}_{\W(k)}\O_K$ by the remark after Definition \ref{def:fr_log_aff}. 

We first review what happens for a point. Choose a uniformizer $\pi$ of $\O_K$ with Eisenstein polynomial $\E(u)$, and give $\mathfrak{S}=\W(k)\llbracket u \rrbracket$ the $\delta$ structure with Frobenius $u\mapsto u^p$ (with the Witt vector Frobenius on $\W(k)$). This makes $(\mathfrak{S},(\E(u)))$ a Breuil-Kisin prism for $\O_K$. When $\O_K=\W(k)[\zeta_{p^\ell}]$ as in our applications, we may use the $q$-de Rham prism $(\W(k)\llbracket q-1 \rrbracket,(\varphi^{\ell-1}([p]_q)))$ with $q\mapsto q^p$ (for $\ell \ge 1$); this is also a Breuil-Kisin prism under our definition.

To construct a logarithmic Breuil-Kisin prism, take a framed log affine and give
\[\mathbb{T}_r^n = \Spf \W(k)\langle \T_1,\ldots,\T_r,\T_{r+1}^{\pm 1},\ldots,\T_n^{\pm 1}\rangle\]  
the Frobenius lift induced by $\T_i \mapsto \T_i^p$. Since $\Spf \R_0$ is $p$-completely \'{e}tale over $\mathbb{T}_r^n$, we may lift this Frobenius uniquely to $\R_0$. We then set
\[\rA=(\R_0\otimes_{\W(k)}\mathfrak{S})^{\wedge}_{(p,\E(u))}.\]
This has a log prism structure, with prismatic ideal $\E(u)$ and the prelog chart $\N^r\to\rA$ sending the $i$th generator to $\T_i$ and $\delta_{\log}(e_i)=0$. One may also directly check
\[(\Spf\rA/\I,\M_{\Spf\rA/\I}) \simeq (\X,\mathcal{M}),\]
so it is an atlas. The argument can be run for the logarithmic $q$-de Rham prism.
\end{proof}

\subsection{Quasicoherent sheaves on the log prismatic stack}

We show that quasicoherent sheaves on $(\X,\mathcal{M})^{\Prism}$ recover log prismatic crystals. This ends up being slightly more subtle than the argument for $\X^{\Prism}$, because we must generally avoid using $\Prism_{\R}$ for $\R$ quasiregular semiperfectoid. The basic issue is that in the log setting these will generally not come with fine log structures, so they will not admit maps to $\Log$ induced by their log structures.

We first show strict log prisms give points of Olsson's log prismatic stack.

\begin{lemma} \label{lem:log_prism_stack_map}
Let $(\X,\mathcal{M})$ satisfy Assumption \ref{ass:horizontal_nc} over $\O_K$, and let $(\rA,\I,\mathcal{N})^a$ be a bounded log prism with fine log structure and a strict map $(\Spf\rA/\I,\mathcal{N}|_{\Spf\rA/\I})\longrightarrow(\X,\mathcal{M})$. Then there is a natural map 
\[\rho_{\X,\rA}:\Spf\rA \longrightarrow (\X,\mathcal{M})^\Prism.\]

If $(\X,\mathcal{M})$ is a framed log affine and $(\rA,\I,\mathcal{N})^a$ is the logarithmic Breuil-Kisin prism associated to its framing, then $\rho_{\X,\rA}$ is an fpqc epimorphism.
\end{lemma}

\begin{proof}
Construction 8.14 of \cite{APC2} gives a map $\Spf\rA \to \X^\Prism$, using the prism structure. Using the definition of $(\X,\mathcal{M})^{\Prism}$ as a relative prismatic stack with respect to $\Log$, it remains to complete the diagram
\begin{display}
\Spf \rA \ar{r} \ar[dashed]{d} & \X^{\Prism} \ar{d} \\
\Z_p^\Prism \times \Log \ar{r}{\pi_{\Log}} & \Log^{\Prism}.
\end{display}
Currently we are missing the dashed arrow as we have not used the log structure, but this is supplied because $\mathcal{M}_{\Spf \rA}$ gives a map to $\Log$, with the other component being the classifying map to $\Z_p^{\Prism}$; unwinding definitions and using that we have a log prism one may check the composites to $\Log^\Prism$ agree, so we get the desired map. We can also produce this map using the induced map of relative stacks, as used in Lemma \ref{lem:mixed_BK_fiber_product}.

Now suppose $(\X,\mathcal{M})$ is a framed log affine. In this situation, Lemma \ref{lem:chart_log_prism} supplies an fpqc epimorphism 
\[\rho_{\X,\B_{\rP}}:\Spf \B_{\rP} \longrightarrow (\X,\mathcal{M})^{\Prism}.\]
We may check whether $\rho_{\X,\rA}: \Spf \rA \to (\X,\mathcal{M})^{\Prism}$ is an fpqc epimorphism after base change along one, and indeed we obtain
\begin{display}
\Spf (\rA \coprod_* \B_{\rP}) \ar{r} \ar{d} & \Spf \B_{\rP} \ar{d} \\
\Spf \rA \ar{r}{\rho_{\X,\rA}} & (\X,\mathcal{M})^{\Prism}
\end{display}
using Lemma \ref{lem:mixed_BK_fiber_product} (which only uses the first part of this lemma). Here $(\rA \coprod_* \B_{\rP})$ is the coproduct in the category of log prisms $(\X,\mathcal{M})_\Prism^{\mathrm{op}}$ which exists by Lemma \ref{lem:prism_coproduct}. But one can directly check that 
\[\Spf (\rA \coprod_* \B_{\rP}) \longrightarrow \Spf \B_{\rP}\] 
is fpqc as it is a prismatic envelope of a $(p,\I)$-completely regular sequence, using \cite[Proposition 3.13]{BhattScholzePrisms} (as constructed in Lemma \ref{lem:prism_coproduct}). This proves the claim.
\end{proof}

Using these fpqc epimorphisms, we next want to relate fiber products $\Spf\rA\times_{(\X,\mathcal{M})^\Prism}\Spf\B$ to coproducts in the category of log prisms over $(\X,\mathcal{M})$. However since \cite[Proposition 4.13]{KoshikawaLogPrismaticI} works with the relative rather than absolute site there is not a direct reference. We use a similar method to construct the envelope in our particular case.

\begin{lemma} \label{lem:prism_coproduct}
Let $(\X,\mathcal{M})$ be a framed log affine over $\O_K$ with $\X=\Spf\R$, let $(\rA,\I,\M_{\Spf\rA})^a$ be the logarithmic Breuil-Kisin prism associated with its framing, and let $(\B,(d),\mathcal{N})^a\in (\X,\mathcal{M})^\op_\Prism$ be a strict log prism. The coproduct $\rA\coprod_*\B$ in $(\X,\mathcal{M})_\Prism^\op$ exists.
\end{lemma}

\begin{proof}
We recall that the construction in Lemma \ref{lem:standard_log_BK_prism} has several variants: the first uses the usual Breuil-Kisin prism with prismatic ideal $\E(u)$ generated by an Eisenstein polynomial of a uniformizer $\pi$, or we may use a log $q$-de Rham prism. The argument in each case is very similar, so we will focus on the first case (the answer in the $q$-de Rham case, at least for a self-coproduct, is already recorded in Lemma \ref{lem:qdr_actual_cech_envelope}). The general argument is analogous to \cite[Proposition 4.13]{KoshikawaLogPrismaticI}.

Write $\R=\R_0\widehat{\otimes}_{\W(k)}\O_K$, with the strict \'{e}tale framing
\[\sq: (\Spf \R_0,\N^r) \longrightarrow (\Spf \W(k)\langle\T_1,\ldots,\T_r,\T_{r+1}^{\pm1},\ldots,\T_n^{\pm1}\rangle, \N^r).\]
We keep $\E(u)$ as a choice of Eisenstein polynomial for a choice of uniformizer $\pi$. From Lemma \ref{lem:standard_log_BK_prism}, we construct $\rA=\R_0 \llbracket u \rrbracket$ and $\I=(\E(u))$, where $u, \T_i$ are sent to $0$ by $\delta$; it is equipped with a log structure with prelog chart $\N^r$ lifted along the framing, say with basis $e_i$ inherited from $\mathcal{M}=\N^r$ on $\Spf \R_0$.

We wish to construct the coproduct $\rA\coprod_*\B$. For the other log prism $(\B,(d),\mathcal{N})^a$, put $\overline{\B}=\B/d$. By strictness and \cite[Lemma 3.8]{KoshikawaLogPrismaticI} we get a chart $\N^r \to \Gamma(\Spf\B,\mathcal{N})$ giving generators $m_i \in \Gamma(\Spf\B,\mathcal{N})$ (using the lemma to lift the sections on $\Gamma(\Spf \B/d, \mathcal{N}|_{\B/d})$ which is pulled back from $\mathcal{M}$ by strictness).

The coproduct $\rC=\rA\coprod_*\B$ is of course characterized by its universal property, and we will construct directly a log prism satisfying this universal property. We will need to find a log prism $\rC$ with maps from $\rA$ and $\B$ such that, for $\rS \in (\X,\mathcal{M})_{\Prism}^\op$,
\[\Hom_{(\X,\mathcal{M})_\Prism^{\mathrm{op}}}(\rC,\rS) \simeq \Hom_{(\X,\mathcal{M})_\Prism^{\mathrm{op}}}(\rA,\rS) \times \Hom_{(\X,\mathcal{M})_\Prism^{\mathrm{op}}}(\B,\rS).\]
To unwind exactly what this means, consider a pair of maps $f:\rA\to\rS$ and $g:\B\to\rS$ corresponding to an element of the right side. We will work here with the strict log prismatic site with $f$ and $g$ being morphisms of prisms over $(\X,\mathcal{M})$. Moreover, since we work with the strict log prismatic site $g(m_i)$ and $f(e_i)$ are forced to have the same reduction modulo $d$ (which we may use everywhere by rigidity of prisms). Then Lemma 3.8 of \cite{KoshikawaLogPrismaticI} and its proof yield a unique unit $\lambda_i\in 1+d\rS$ such that $g(m_i)=\lambda_i f(e_i)$.

Now let us construct $\rC$ explicitly. We first adjoin variables corresponding to these units as
\[\D^{\mathrm{ex}} = \left((\R_0\llbracket u \rrbracket \widehat{\otimes}_{\W(k)} \B)[v_1^{\pm1},\ldots,v_r^{\pm1}]/(\alpha(m_i)-v_i\T_i)_{i=1}^r\right)^{\wedge}_{(p,\E,d)},\]
where $\alpha:\Gamma(\Spf\B,\mathcal{N})\to \B$ is the log structure map. We will give this ring the log structure pulled back from $\B$ and set the $\delta_{\log}$-structure to be given by $\delta(v_i)=v_i^p\delta_{\log}(m_i)$.\footnote{This construction is essentially exactification of a $\delta_{\log}$-triple as in \cite[Construction 2.17]{KoshikawaLogPrismaticI}.} There is a natural map $q: \D^{\mathrm{ex}} \longrightarrow \overline{\B}$: we send all $v_i \mapsto 1$, and on the $\R_0\llbracket u \rrbracket$ factor use the map
\[\R_0\llbracket u \rrbracket \longrightarrow \R \longrightarrow \overline{\B}\]
and the obvious map $\B \to \overline{\B}$ on the other factor. We then set
\[\rC := \left(\D^{\mathrm{ex}}\left\{\frac{\ker q}{d}\right\}_\delta\right)^{\wedge}_{(p,d)},\]
noting this has a simple form as $\ker(q: \D^{\mathrm{ex}} \to \overline{\B})$ is $(d,x_0,\ldots,x_n)$ where the $x_i$ form a $(p,d)$-completely regular sequence. To be more precise, locally one can take the $x_i$ to be $v_1-1,\ldots,v_r-1$ followed by $u-\tilde{q}(u)$ and $\T_{r+1}-\tilde{q}(\T_{r+1}), \ldots, \T_n-\tilde{q}(\T_n)$ where $\tilde{q}(-)$ denotes a choice of lift of $q(-)\in \overline{\B}$ to $\B$ in $\D^{\mathrm{ex}}$. The map $\B\to\rC$ makes it so $\rC$ is a strict log prism in $(\X,\mathcal{M})_\Prism^\op$. Since $\E(u)\in\ker(q)$ maps into $d\rC$, the map $\rA\to\rC$ is a log prism map.

Now we need to match the universal properties. Recall that the universal property of this prismatic envelope is that for each prism $(\rS,\I_{\rS})$ and $\delta$-map
\[h: \D^{\mathrm{ex}} \longrightarrow \rS\]
carrying $\ker(q)$ into $\I_{\rS}$ there is a unique extension of $h$ to a $\delta$-map $\rC \to \rS$. This is automatically a map of prisms, since $d\in \ker(q)$.

Now suppose we have for some test log prism $\rS$, two log prism maps over $(\X,\mathcal{M})$: $\rA \to \rS$ and $\B \to \rS$. Then the construction of $\D^{\mathrm{ex}}$ makes it so that there is a unique $\delta_{\log}$-map
\[\D^{\mathrm{ex}} \longrightarrow \rS\]
that carries $v_i$ to $\lambda_i$. The maps agree over $(\X,\mathcal{M})$, forcing this map to also carry $\ker(q)$ into $\I_{\rS}$. Now we use the universal property of the prismatic envelope: there is a unique compatible map $\rC \to \rS$, which is forced to be a map of log prisms. This gives an injection 
\[\Hom_{(\X,\mathcal{M})_\Prism^{\mathrm{op}}}(\rA,\rS) \times \Hom_{(\X,\mathcal{M})_\Prism^{\mathrm{op}}}(\B,\rS) \hookrightarrow \Hom_{(\X,\mathcal{M})_\Prism^{\mathrm{op}}}(\rC,\rS).\]
It is surjective, since we can restrict a log prism map $h: \rC \to \rS$ along the log prism maps $\rA \to \rC$ and $\B \to \rC$ to obtain $f$ and $g$. These must satisfy $g(m_i) = h(v_i) f(e_i)$ by construction, where $h(v_i)\in 1 + \I_{\rS}$. But then by uniqueness we get $h(v_i)=\lambda_i \in 1 + \I_{\rS}$, and it follows that the map 
\[\D^{\mathrm{ex}} \longrightarrow \rS\]
we then obtain is $h|_{\D^{\mathrm{ex}}}$, so by the previous procedure we recover $h$ from $f$ and $g$.
\end{proof}

We now show how to recover these coproducts via the log prismatic stack.

\begin{lemma} \label{lem:mixed_BK_fiber_product}
Let $(\X,\mathcal{M})$ be a framed log affine over $\O_K$ and let $(\rA,\I,\M_{\Spf\rA})^a$ be the logarithmic Breuil-Kisin prism associated with its framing. Let $(\B,(d),\mathcal{N})^a$ be a bounded log prism in $(\X,\mathcal{M})_\Prism^\op$. Then
\[\Spf(\rA\coprod_*\B)\simeq\Spf\rA\times_{(\X,\mathcal{M})^\Prism}\Spf\B.\]
The analogous assertion also holds with trivial log structures.
\end{lemma}

\begin{proof}
By Lemma \ref{lem:prism_coproduct}, the coproduct exists and therefore the statement makes sense. We will use the notation of Lemma \ref{lem:prism_coproduct}, with $\overline{\B}=\B/d$. 

In the situation of a framed log affine, the morphism $\X \to \Log$ induced by the log structure factors as $\X \to \mathcal{A}_{\N^r} \to \Log$. Here $\mathcal{A}_{\N^r}\to\Log$ is an \'{e}tale chart for the log stack, given by $\mathcal{A}_{\rP}:=[\Spf \Z_p\langle \rP\rangle/\Spf \Z_p\langle \rP^{\mathrm{gp}}\rangle]$. This fact is useful because in the definition of $(\X,\mathcal{M})^{\Prism}$ it allows us to replace all instances of $\Log$ with $\mathcal{A}_{\N^r}$. That is,
\[(\X,\mathcal{M})^\Prism \simeq (\X/\mathcal{A}_{\N^r})^\Prism.\]
As in the lemma statement, we let $(\rA,\I,\M_{\Spf\rA})^a$ be the logarithmic Breuil-Kisin prism constructed in Lemma \ref{lem:standard_log_BK_prism}, so in particular $\X=\Spf \rA/\I$. We will use the diagram
\begin{display}
\X \ar{r}{\mathrm{id}} \ar{d} & \X \ar{d} & \Spf \overline{\B} \ar{l} \ar{d} \\
\Spf\rA \ar{r} & \mathcal{A}_{\N^r} & \Spf \B \ar{l}.
\end{display}
Upon taking relative prismatic stacks on the columns, we obtain maps
\[(\X/\Spf\rA)^\Prism \longrightarrow (\X,\mathcal{M})^\Prism \longleftarrow (\Spf \overline{\B}/\Spf\B)^\Prism.\]
Noting that $(\X/\Spf\rA)^\Prism \simeq \Spf\rA$ and $(\Spf \overline{\B}/\Spf\B)^\Prism \simeq \Spf\B$, equivalently this gives
\[\Spf \rA \longrightarrow (\X,\mathcal{M})^\Prism \longleftarrow \Spf \B.\]
Revisiting the construction in Lemma \ref{lem:log_prism_stack_map}, these are precisely $\rho_{\X,\rA}$ and $\rho_{\X,\B}$.

Even better, using finite limit compatibility of the relative prismatization in this diagram (for $\delta$-pairs by \cite[Lemma 4.20]{AKNRelativeDelta}) we obtain
\[\Spf\rA \times_{(\X,\mathcal{M})^\Prism} \Spf\B \simeq \left(\Spf\overline{\B}/(\Spf\rA\times^{\R}_{\mathcal{A}_{\N^r}}\Spf\B)\right)^\Prism.\]
We will compute the relative prismatic stack on the right, aiming to show it is $\Spf(\rA \coprod_{\ast} \B)$ by simplifying the relative stack until we may identify it with a prismatic envelope. Since $\mathcal{A}_{\N^r}=[\Spf\Z_p\langle \N^r\rangle/\mathbf{G}_m^r]$, we may explicitly compute $\Spf\rA\times^{\R}_{\mathcal{A}_{\N^r}}\Spf\B$ as
\[\Spf\left((\rA\otimes^{\L}_{\Z_p}\B)[v_1^{\pm1},\ldots,v_r^{\pm1}]\dquot (\alpha(m_i)-v_i\T_i)_{i=1}^r\right)^{\wedge}_{(p,\E,d)}.\]
In this presentation, we are using $v_i$ as the coordinates on $\mathbf{G}_m^r$, and the map from $\Spf \overline{\B}$ is given by $v_i \mapsto 1$. We use $\alpha(m_i)$ as in the proof of Lemma \ref{lem:prism_coproduct}. A cursory glance shows that this is \emph{very} similar to $\D^{\mathrm{ex}}$ from Lemma \ref{lem:prism_coproduct}: replacing $\Z_p$ by $\W(k)$ gives exactly $\D^{\mathrm{ex}}$, presented there as
\[\D^{\mathrm{ex}}=\left((\rA\otimes^{\L}_{\W(k)}\B)[v_1^{\pm1},\ldots,v_r^{\pm1}]\dquot (\alpha(m_i)-v_i\T_i)_{i=1}^r\right)^{\wedge}_{(p,\E,d)}.\]
Since the maps $\W(k)\longrightarrow \overline{\B}$ induced by $\rA$ and $\B$ agree and $\W(k)$ is $p$-completely formally \'{e}tale over $\Z_p$, the canonical map
\[\left( \overline{\B}/\D^{\mathrm{ex}} \right)^\Prism \xlongrightarrow{\sim}
\left(\Spf \overline{\B}/(\Spf\rA\times^{\R}_{\mathcal{A}_{\N^r}}\Spf\B)\right)^\Prism\]
is an equivalence by the proof of \cite[Lemma 5.9]{AKNRelativeDelta} (which gives the comparison at the level of stacks).

Now it remains to compute $\left(\overline{\B}/\D^{\mathrm{ex}}\right)^\Prism$. By \cite[Theorem 7.17]{APC2} we then see
\[\Spf\rA\times_{(\X,\mathcal{M})^\Prism}\Spf\B \simeq \Spf \Prism_{\overline{\B}/\D^{\mathrm{ex}}}.\]
This may be computed directly as a prismatic envelope as $\overline{\B}$ is locally cut out by a Koszul regular sequence in $\D^{\mathrm{ex}}/d$, and better the map $\D^{\mathrm{ex}} \to \overline{\B}$ identifies with the map $q$ from Lemma \ref{lem:prism_coproduct}. Thus since the coproduct in Lemma \ref{lem:prism_coproduct} is the prismatic envelope of $\ker(q)/d$ we have shown $\Spf \rA \times_{(\X,\mathcal{M})^\Prism} \Spf \B \simeq \Spf (\rA \coprod_* \B)$.

The assertion with trivial log structures follows by a similar argument.
\end{proof}

\begin{proposition} \label{prop:log_stack_presentation}
Let $(\X,\mathcal{M})=(\Spf\R,\mathcal{M})$ be a framed log affine over $\O_K$, and let
\[\rho_{\X,\rA}: \Spf \rA \longrightarrow (\X,\mathcal{M})^{\Prism}\]
be the fpqc epimorphism of Lemma \ref{lem:log_prism_stack_map} associated with the logarithmic Breuil-Kisin prism of its framing. Then the \v{C}ech nerve of $\Spf \rA$ identifies with $\Spf \rA^{(n)}$ where $\rA^{(n)}$ are the iterated self-coproducts of $(\rA,\I,\M_{\Spf \rA})^a$ in $(\X,\mathcal{M})_\Prism^\op$. In particular
\[(\X,\mathcal{M})^\Prism \simeq \left|\Spf\rA^{(\bullet)}\right|\]
in formal stacks.
\end{proposition}

\begin{proof}
Lemma \ref{lem:mixed_BK_fiber_product}, applied with $\B=\rA$, identifies the first term of the \v{C}ech nerve of $\Spf \rA$ over $(\X,\mathcal{M})^{\Prism}$ with $\Spf(\rA\coprod_*\rA)$. The higher terms are the formal spectra of the iterated coproducts $(\rA^{(1)})^{\widehat{\otimes}_{\rA} n}$, for both the site-theoretic \v{C}ech nerve and the \v{C}ech nerve of $\rho_{\X,\rA}: \Spf \rA \to (\X,\mathcal{M})^{\Prism}$. Using Lemma \ref{lem:log_prism_stack_map}, this is an fpqc epimorphism so the claim follows.
\end{proof}

\begin{proposition} \label{prop:stacky_qcoh}
Let $(\X,\mathcal{M})$ satisfy Assumption \ref{ass:horizontal_nc} over $\O_K$. There is a symmetric monoidal equivalence
\[\iota_{(\X,\mathcal{M})}: \D((\X,\mathcal{M})_\Prism,\O_\Prism) \xlongrightarrow{\sim} \D_{\mathrm{qc}}((\X,\mathcal{M})^\Prism).\]
\end{proposition}

\begin{proof}
We first note that the assertion can be checked strict \'{e}tale locally on a framed log affine, for both the source and target. On the target $\D_{\mathrm{qc}}((\X,\mathcal{M})^\Prism)$ this follows from \'{e}tale descent for $\X^{\Prism}$ in $\X$ from \cite{APC2}, and on the source $\D((\X,\mathcal{M})_\Prism,\O_\Prism)$ it is clear by complete fpqc descent of $\rA \mapsto \widehat{\D}(\rA)$ on derived $(p,I)$-complete rings $\rA$ (see the argument in \cite[Proposition 2.7]{BhattScholzePrismaticFCrystals}).

Now working locally, we may use the logarithmic Breuil-Kisin atlas $(\rA,\I,\M_{\Spf \rA})^a$ from Lemma \ref{lem:standard_log_BK_prism} to obtain
\[\lim_{[n]\in\Delta} \D_{(p,\I^{(n)})}(\rA^{(n)}) \simeq \D((\X,\mathcal{M})_\Prism,\O_\Prism)\]
where $\rA^{(n)}$ is the logarithmic \v{C}ech nerve. The (analogous) proof of Corollary \ref{cor:laur_strict_etale} for ordinary crystals gives strict \'{e}tale descent of the crystal side, which is well-known.

On the other side, we can use Proposition \ref{prop:log_stack_presentation} to identify
\[\D_{\mathrm{qc}}((\X,\mathcal{M})^{\Prism}) \simeq \lim_{[n]\in\Delta}\D_{\mathrm{qc}}(\Spf \rA^{(n)}) \simeq \D((\X,\mathcal{M})_\Prism,\O_\Prism).\]
Indeed, $\rho_{\X,\rA}: \Spf \rA \to (\X,\mathcal{M})^{\Prism}$ is an fpqc epimorphism and we have identified its \v{C}ech nerve with $\Spf \rA^{(n)}$, and thus $\D_{\mathrm{qc}}$ satisfies descent for such maps. See the discussion after \cite[Definition 2.4]{syncrisloc} for more details.
\end{proof}

In the following result we write $\R f_{\Prism,*}$ and $\R f_*$ for the pushforwards on the two categories of Proposition \ref{prop:stacky_qcoh}, respectively.

\begin{proposition} \label{prop:stacky_qcoh_pushforward}
Consider a map $f:(\X,\mathcal{M}) \to (\Y,\mathcal{N})$ where the source and target satisfy Assumption \ref{ass:horizontal_nc} over $\O_K$. Then for every $\mathcal{E} \in \D((\X,\mathcal{M})_\Prism,\O_\Prism)$, there is a natural equivalence
\[\iota_{(\Y,\mathcal{N})} \R f_{\Prism,*}\mathcal{E} \simeq \R f_*\iota_{(\X,\mathcal{M})}\mathcal{E}.\]
\end{proposition}

\begin{proof}
There is a commutative diagram
\begin{display}
\D((\Y,\mathcal{N})_\Prism,\O_\Prism) \ar{r}{\iota_{(\Y,\mathcal{N})}} \ar{d}{f^*} & \D_{\mathrm{qc}}((\Y,\mathcal{N})^\Prism) \ar{d}{(f^\Prism)^*} \\
\D((\X,\mathcal{M})_\Prism,\O_\Prism) \ar{r}{\iota_{(\X,\mathcal{M})}} & \D_{\mathrm{qc}}((\X,\mathcal{M})^\Prism).
\end{display}
By compatibility of the equivalence $\iota$ with pullbacks, upon taking right adjoints the claim follows.
\end{proof}

\begin{remark}
In \cite{KoshikawaYao2025,KoshikawaYaoCorrigendum}, they use the saturated log prismatic site for Laurent $\rF$-crystals. However, by \cite[Lemma 5.9]{Inoue}, the strict subsite is cofinal in the full log prismatic site. When $(\X,\mathcal{M})$ is fs, a strict log prism $(\rA,\I,\M_{\Spf \rA})^a$ in $(\X,\mathcal{M})_{\Prism}^\op$ is saturated; thus $(\X,\mathcal{M})_{\Prism}^{\mathrm{str}} \subset (\X,\mathcal{M})_{\Prism}^{\mathrm{sat}} \subset (\X,\mathcal{M})_{\Prism}^{\mathrm{full}}$ for the strict, saturated, and full log prismatic sites. The first and last have the same crystal categories, so all three do.
\label{rem:crystal_indep}
\end{remark}

We therefore suppress this distinction for $\D((\X,\mathcal{M})_\Prism,\O_\Prism)$ and its perfect and Laurent variants.

\subsection{Solid Laurent crystals}
\label{app:solid_laur}

In this subsection, we define a solid generalization of
\[\Perf((\X,\mathcal{M})_\Prism,\Laur),\]
and show that in nice situations we can still understand it explicitly through \v{C}ech nerves. We first define a useful analytic stack that we've frequently used in the article.

\begin{definition} \label{def:anspf}
Let $(\rA,(d))$ be an orientable prism. Then we write
\[\AnSpf(\rA[1/d]_p^{\wedge}):=\colim_n \Spa((\rA\dquot p^n)[1/d],\rA\dquot p^n),\]
where we use $\Spa$ to denote the analytic stack attached to a Huber pair as in \cite[Example 2.1.5]{Kipp2}.
\end{definition}

\begin{definition} \label{def:laur_category}
Let $(\rA,(d))$ be an orientable prism. Associated to the derived $(p,d)$-complete ring $\rA$ is the analytic adic ring $(\rA,\rA)$. We define
\[\D_{\Box}(\rA[1/d]_p^{\wedge}) := \lim_n \Mod_{(\rA\dquot p^n)[1/d]} \D_{\Box}((\rA,\rA)),\]
where $(\rA\dquot p^n)[1/d]$ is regarded as an algebra object of $\D_{\Box}((\rA,\rA))$. Equivalently, this is $\D_{\Box}$ of $\AnSpf(\rA[1/d]_p^{\wedge})$.

Then
\[\D_{\Box}((\X,\mathcal{M})_\Prism,\Laur) := \lim_{(\rA,(d),\M_{\Spf \rA})^a \in (\X,\mathcal{M})_\Prism^{\mathrm{op}}} \D_{\Box}(\rA[1/d]_p^{\wedge}).\]
\end{definition}

Note there is no issue with working only with principal $d$ as this is true locally, so we obtain the same category of crystals. In this setting we define perfect Laurent crystals as the dualizable objects
\[\Perf((\X,\mathcal{M})_\Prism,\Laur) := \D_{\Box}((\X,\mathcal{M})_\Prism,\Laur)^{\mathrm{dual}},\]
noting one obtains the expected description in Lemma \ref{lem:perf_comparison}.

Our first task is to show that $\D_{\Box}((\X,\mathcal{M})_\Prism,\Laur)$ is computable, as a priori there is no reason for it to be possible to compute this category just using the \v{C}ech nerve of some Breuil-Kisin prism. We will show that this is in fact the case in Theorem \ref{thm:laur_cech}.

This will require us to prove some descendability results for log prismatic envelopes, for which the following criterion will be useful.

\begin{lemma}
Let $\rA\to\B$ be a faithfully flat map of connective animated rings. If $\pi_0\B$ is countably presented as a $\pi_0\rA$-module, then $\rA\to\B$ is descendable of index at most $2$.
\label{lem:countable_flat_descendable}
\end{lemma}

\begin{proof}
See \cite[Corollary 3.33 and the proof of Proposition 3.32]{MathewGaloisStableHomotopy}.
\end{proof}

\begin{lemma} \label{lem:cech_descendable}
Let $(\X,\mathcal{M})$ be a framed log affine over $\O_K$, let $(\rA,\I,\M_{\Spf\rA})^a$ be the logarithmic Breuil-Kisin prism associated with its framing, and let $(\B,(d),\mathcal{N})^a\in (\X,\mathcal{M})_{\Prism}^\op$. Put
\[\rC=\rA\coprod_*\B\]
for the coproduct in $(\X,\mathcal{M})_\Prism^\op$. Then $\Spa(\rC,\rC)\to\Spa(\B,\B)$ is a proper descendable cover.
\end{lemma}

Note that as a degenerate case the same assertions hold when $\X$ is smooth over $\O_K$, when $(\rA,\I)$ is an ordinary Breuil-Kisin atlas.

\begin{proof}
Using Lemma \ref{lem:prism_coproduct}, we can write the coproduct explicitly as
\begin{equation}
\rC \xlongrightarrow{ \sim } \left(\D^{\mathrm{ex}}\{x_1/d,\ldots,x_m/d\}_\delta\right)^{\wedge}_{(p,d)},
\label{eq:prismatic_coproduct_envelope}
\end{equation}
where $(d,x_1,\ldots,x_m) = \ker(\D^{\mathrm{ex}} \to \overline{\B})$. Using \cite[Proposition 5.52, and Remark 5.54]{MaoDerivedCrystalline} we see that $\B \to \rC$ is $(p,d)$-completely faithfully flat. The maps $\B \dquot (p^k,d^k) \to \rC\dquot(p^k,d^k)$ are also faithfully flat and countably presented on $\pi_0$, since we adjoin finitely many variables as a $\delta$-ring, and therefore descendable of index $\le 2$ by Lemma \ref{lem:countable_flat_descendable}.

However for our application, we need to show complete descendability as in \cite[Definition 2.4.1.7]{Kipp} -- this simply means that the canonical map $\rJ:=\fib(\B \to \rC) \to \B$ is $\widehat{\otimes}$-nilpotent in the derived $(p,d)$-complete category (by rigidity of prisms, we can use the same $d$ for the source and target). We will write $\rJ_k=\rJ\dquot(p^k,d^k)$. Derived completeness then yields
\[\RHom_{\B}(\rJ\widehat{\otimes}_{\B} \rJ,\B) \simeq \R\lim_k \RHom_{\B\dquot(p^k,d^k)}(\rJ_k^{\otimes 2},\B\dquot(p^k,d^k)).\]
As we have seen already, the canonical map $\rJ\widehat{\otimes}_{\B}\rJ\to\B$ is nullhomotopic modulo $(p^k,d^k)$. One may then still apply the argument of \cite[Lemma 11.22]{BhattScholzeWittGrassmannian} (note the statement does not directly apply, but the key idea that $\lim^2$ vanishes for an $\N$-indexed limit does) to show that $(\rJ\widehat{\otimes}_{\B}\rJ)^{\otimes 2} \to \B$ is nullhomotopic and therefore that $\rC$ is descendable of index $\le 4$ in the derived $(p,d)$-complete category as desired. Passing to $\D_{\Box}((\B,\B))$ is lax monoidal and thus preserves $\widehat{\otimes}$-nilpotence of $\rJ:=\fib(\B \to \rC) \to \B$, so the first claim follows.

Finally, it remains to prove properness. In our present setting, Example 2.1.13 of \cite{Kipp2} tells us that a map $\Spa(\rC,\rC) \to \Spa(\B,\B)$ is proper if and only if $\rC$ is integral over the subring generated by $\B$ and topologically nilpotent elements of $\rC$. It is therefore enough to prove that
\[\pi_0(\B\dquot(p,d)) \longrightarrow \pi_0(\rC\dquot(p,d))\]
is integral. To see this is really the same condition, note that if a monic polynomial $f$ over $\pi_0\B$ annihilates $c$ modulo $(p,d)$, then $f(t)-f(c)$ is a monic polynomial over $\pi_0\B+(p,d)\pi_0\rC$ annihilating $c$. Now we apply Lemma \ref{lem:prismatic_envelope_integral} to \eqref{eq:prismatic_coproduct_envelope}, using $\D^{\mathrm{ex}}/(d,x_1,\ldots,x_m)\simeq\B/d$.
\end{proof}

The argument for Lemma \ref{lem:prismatic_envelope_integral} arose from a suggestion of Bhargav Bhatt to reduce to \cite[Lemma 12.4]{BhattScholzePrisms}.

\begin{lemma} \label{lem:prismatic_envelope_integral}
Let $\rA$ be an animated $\delta$-ring, let $d\in\pi_0\rA$ be distinguished, and let $\{f_i\}$ be any (arbitrarily large) family of elements of $\pi_0\rA$. Suppose $\rA$ is further derived $(p,d)$ complete. Then the morphism
\[\Spf \left(\rA\{f_i/d\}_\delta\right)^{\wedge}_{(p,d)} \longrightarrow \Spf \rA\]
induced by the structure map is integral in the sense of \cite[2.2.5]{Kipp2}. In other words, on rings it is true in the ordinary sense modulo $(p,d)$.
\end{lemma}

\begin{proof}
As we discuss integrality on formal schemes, we must show 
\[\pi_0(\rA\dquot(p,d)) \longrightarrow \pi_0\left(\left(\rA\{f_i/d\}_\delta\right)^{\wedge}_{(p,d)}\dquot(p,d)\right)\]
is integral. As a first reduction, for just testing integrality it suffices to take a single $f_i$, so we will just work with $f$ from now on. Observe that
\[\pi_0 \rA/(p,d) \longrightarrow \pi_0 \rA_{\perf}/(p,d)\]
is a universal homeomorphism (where we use the perfectoidization in the target). In particular, we can check integrality after base change along this map and assume $\rA=\rA_{\perf}$. Now $\pi_0 \rA_{\perf} = \Ainf(\R)$ for some integral perfectoid ring $\R$, which by Andre's lemma (\cite[Theorem 7.14]{BhattScholzePrisms}) has a perfectoid $p$-completely flat cover $\rS$ where the image $\overline{f} \in \R$ admits compatible $p$-power roots, as well as containing compatible $p$-power roots of unity $\{\zeta_{p^n}\}_{n\ge 0}$. By this second property, the prismatic ideal in $\Ainf(\rS)$ can be taken to be $[p]_q$. Integrality may then be tested after passing to $\rS$, so we check whether
\[\Ainf(\rS) \dquot (p,[p]_q) \longrightarrow \Ainf(\rS)\{f/[p]_q\}_\delta^{\wedge} \dquot (p,[p]_q)\]
is integral. Since the compatible system of $p$-power roots $\overline{f}^\flat \in \rS^\flat$ satisfies $f-[\overline{f}^\flat]\in ([p]_q)$, we can replace the target prismatic envelope with $\Ainf(\rS)\{[\overline{f}^\flat]/[p]_q\}_\delta^{\wedge}$. Following the setup in \cite[\S 12.2]{BhattScholzePrisms}, the claim follows from the description in \cite[Lemma 12.4]{BhattScholzePrisms}.
\end{proof}

\begin{theorem} \label{thm:laur_cech}
Let $(\X,\mathcal{M})$ be a framed log affine over $\O_K$, and let $(\rA,\I,\M_{\Spf\rA})^a$ be the logarithmic Breuil-Kisin prism associated with its framing. If $(\rA^{(m)},\I,\M^{(m)})^a$ is the degree $m$ term of its \v{C}ech nerve in the strict log prismatic site, then
\[\D_{\Box}((\X,\mathcal{M})_\Prism,\Laur) \simeq \lim_{[m]\in\Delta}\D_{\Box}(\rA^{(m)}[1/\I]^{\wedge}_p).\]
\end{theorem}

\begin{proof}
By definition, 
\[\D_{\Box}((\X,\mathcal{M})_\Prism,\Laur) \simeq \lim_{(\B,(d),\M_{\Spf \B})^a \in (\X,\mathcal{M})_\Prism^\op} \D_{\Box}(\B[1/d]_p^{\wedge})\]
where $\D_{\Box}(\B[1/d]_p^{\wedge})$ is interpreted as in Definition \ref{def:laur_category}. Now Lemma \ref{lem:cech_descendable} tells us
\[\D_{\Box}(\B[1/d]^{\wedge}_p) \simeq \lim_{[m]\in\Delta}\D_{\Box}((\rA^{(m)}\coprod_*\B)[1/d]^{\wedge}_p).\]
Now upon commuting limits, we get
\[\D_{\Box}((\X,\mathcal{M})_\Prism,\Laur) \simeq \lim_{[m]\in \Delta} \left( \lim_{(\B,(d),\M_{\Spf \B})^a \in (\X,\mathcal{M})_\Prism^\op} \D_{\Box}((\rA^{(m)}\coprod_*\B)[1/d]^{\wedge}_p) \right).\]
We are now done, as formally we have
\[\lim_{(\B,(d),\M_{\Spf \B})^a \in (\X,\mathcal{M})_\Prism^\op} \D_{\Box}((\rA^{(m)}\coprod_*\B)[1/d]^{\wedge}_p) \simeq \D_{\Box}(\rA^{(m)}[1/\I]_p^{\wedge}),\]
proving the claim. One way to see this last point is that the coproduct functor
\[(\X,\mathcal{M})_\Prism^\op \longrightarrow (\X,\mathcal{M})_{\Prism,/\rA^{(m)}}^\op\]
sending $(\B,(d),\M_{\Spf \B})^a \mapsto (\rA^{(m)}\coprod_*\B)$ regarded as a log prism is left adjoint to the forgetful functor, so by adjunction the limit is
\[\lim_{\rA^{(m)} \to \B} \D_{\Box}(\B[1/d]_p^{\wedge}) \simeq \D_{\Box}(\rA^{(m)}[1/\I]_p^{\wedge}).\]
Thus $\D_{\Box}((\X,\mathcal{M})_\Prism,\Laur) \simeq \lim_{[m]\in\Delta}\D_{\Box}(\rA^{(m)}[1/\I]^{\wedge}_p)$.
\end{proof}

\begin{lemma} \label{lem:perf_comparison}
Let $(\rA,(d))$ be an orientable transversal prism, i.e. $\rA/d$ is $p$-torsionfree. Then we have
\[\D_{\Box}(\rA[1/d]^{\wedge}_p)^{\mathrm{dual}} \simeq \lim_n\Perf((\rA/p^n)[1/d]) \simeq \Perf(\rA[1/d]^{\wedge}_p).\]
In the setting of Theorem \ref{thm:laur_cech}, these give
\[\Perf((\X,\mathcal{M})_\Prism,\Laur) \simeq \lim_{[m]\in\Delta} \Perf(\rA^{(m)}[1/\I]^{\wedge}_p).\]
\end{lemma}

\begin{proof}
This follows by taking dualizable objects in Theorem \ref{thm:laur_cech} and applying \cite[Theorem 14.10]{ScholzeAnalyticGeometry} upon noting that dualizable objects can be detected termwise in a limit diagram.
\end{proof}

\begin{lemma} \label{lem:qdr_laur_cech}
Let $(\Spf\R,\mathcal{M})$ be a framed log affine over $\Spf\W(k)[\zeta_p]$. Let $((\rA^+)^{(m)},[p]_q,\mathcal{M}^{(m)})^a$ be the degree $m$ term of the log $q$-de Rham \v{C}ech nerve in $(\Spf\R,\mathcal{M})_\Prism$, and put $\rA^{(m)}=(\rA^+)^{(m)}[1/[p]_q]^{\wedge}_p$. Then
\[\D_{\Box}((\Spf\R,\mathcal{M})_\Prism,\Laur)\simeq \lim_{[m]\in\Delta}\D_{\Box}(\rA^{(m)}).\]
\end{lemma}

\begin{proof}
In Lemma \ref{lem:standard_log_BK_prism}, the log $q$-de Rham prism is a valid choice of a logarithmic Breuil-Kisin atlas. Theorem \ref{thm:laur_cech} therefore gives the needed \v{C}ech nerve description.
\end{proof}

\begin{corollary} \label{cor:laur_strict_etale}
On log formal schemes satisfying Assumption \ref{ass:horizontal_nc} over $\O_K$, the assignment
\[(\X,\mathcal{M}) \longmapsto \D_{\Box}((\X,\mathcal{M})_\Prism,\Laur)\]
is a stack for the strict \'{e}tale topology. The same holds for $\Perf((\X,\mathcal{M})_\Prism,\Laur)$.
\end{corollary}

\begin{proof}
Consider a strict \'{e}tale cover $\U\to\X$ by framed log affines, where $\X$ itself is also framed. We may give every term $\U^{(n)}$ of its \v{C}ech nerve the framing induced from $\X$, and let $\rA^+_\X$ be the logarithmic Breuil-Kisin prism associated with the framing of $\X$. Because all terms $\U^{\bullet} \to \X$ in the \v{C}ech nerve are \'{e}tale, we obtain unique Breuil-Kisin prisms $\rA^{+}_{\U^{(n)}}$ \'{e}tale over $\rA^+_\X$. We write $\rA^{(m)}_{\U^{(n)}}$ for the (Laurentized) terms in the \v{C}ech nerves of these Breuil-Kisin prisms. Our present situation can be summarized with the following augmented bicosimplicial diagram:
\begin{display}[shift left/.default=+0.8ex]
\D_{\Box}((\X,\mathcal{M})_\Prism,\Laur)
  \ar{r}
  \ar{d}
& \D_{\Box}(\rA_{\X}^{(0)})
  \ar[shift left]{r}
  \ar[shift right]{r}
  \ar{d}
& \D_{\Box}(\rA_{\X}^{(1)})
  \ar{l}
  \ar[shift right, shift right]{r}
  \ar{r}
  \ar[shift left, shift left]{r}
  \ar{d}
& \cdots
  \ar[shift right]{l}
  \ar[shift left]{l} \\
\D_{\Box}((\U,\mathcal{M}|_{\U})_\Prism,\Laur)
  \ar{r}
  \ar[shift left]{d}
  \ar[shift right]{d}
& \D_{\Box}(\rA_{\U}^{(0)})
  \ar[shift left]{r}
  \ar[shift right]{r}
  \ar[shift left]{d}
  \ar[shift right]{d}
& \D_{\Box}(\rA_{\U}^{(1)})
  \ar{l}
  \ar[shift right, shift right]{r}
  \ar{r}
  \ar[shift left, shift left]{r}
  \ar[shift left]{d}
  \ar[shift right]{d}
& \cdots
  \ar[shift right]{l}
  \ar[shift left]{l} \\
\D_{\Box}((\U^{(1)},\mathcal{M}|_{\U^{(1)}})_\Prism,\Laur)
  \ar{r}
  \ar{u}
  \ar[shift right, shift right]{d}
  \ar{d}
  \ar[shift left, shift left]{d}
& \D_{\Box}(\rA_{\U^{(1)}}^{(0)})
  \ar{u}
  \ar[shift left]{r}
  \ar[shift right]{r}
  \ar[shift right, shift right]{d}
  \ar{d}
  \ar[shift left, shift left]{d}
& \D_{\Box}(\rA_{\U^{(1)}}^{(1)})
  \ar{u}
  \ar{l}
  \ar[shift right, shift right]{r}
  \ar{r}
  \ar[shift left, shift left]{r}
  \ar[shift right, shift right]{d}
  \ar{d}
  \ar[shift left, shift left]{d}
& \cdots
  \ar[shift right]{l}
  \ar[shift left]{l} \\
\vdots \ar[shift right]{u} \ar[shift left]{u}
& \vdots \ar[shift right]{u} \ar[shift left]{u}
& \vdots \ar[shift right]{u} \ar[shift left]{u}
& {}
\end{display}
Here we can use Theorem \ref{thm:laur_cech} and descent of $\D_{\Box}(\AnSpf \rA)$ for \'{e}tale covers to deduce each augmented row or column is a limit diagram. For Laurent crystals, as used here, this follows because $\D_{\Box}(\AnSpf \rA) = \lim_k \D_{\Box}(\Spa(\rA^+[1/\I] \dquot p^k, \rA^+ \dquot p^k))$ where the maps are \'{e}tale covers for each $k$ and satisfy descent. Since this is a limit we get the corresponding descent result for $\D_{\Box}(\AnSpf \rA)$.

Then, using that limits commute with limits, we deduce
\[
\begin{aligned}
\lim_{[n]\in\Delta}
\D_{\Box}((\U^{(n)},\mathcal{M}|_{\U^{(n)}})_\Prism,\Laur)
&\simeq
\lim_{[n]\in\Delta}\lim_{[m]\in\Delta}
\D_{\Box}(\rA_{\U^{(n)}}^{(m)})\\
&\simeq
\lim_{[m]\in\Delta}\lim_{[n]\in\Delta}
\D_{\Box}(\rA_{\U^{(n)}}^{(m)})\\
&\simeq
\lim_{[m]\in\Delta}
\D_{\Box}(\rA_{\X}^{(m)})\\
&\simeq
\D_{\Box}((\X,\mathcal{M})_\Prism,\Laur).
\end{aligned}
\]
The claim follows for $\Perf((\X,\mathcal{M})_\Prism,\Laur)$ since dualizability can be tested termwise in a limit of symmetric monoidal categories.
\end{proof}

The canonical $t$-structure on $\D_{\Box}((\X,\mathcal{M})_\Prism,\Laur)$ is defined by taking the connective objects or $\le 0$ part to be the crystals whose evaluations on every prism are connective, and taking the coconnective objects to be the right orthogonal of $\D^{\le -1}_{\Box}((\X,\mathcal{M})_\Prism,\Laur)$.

\begin{corollary} \label{cor:laur_t_struct}
Let $(\X,\mathcal{M})$ satisfy Assumption \ref{ass:horizontal_nc} over $\O_K$. The canonical $t$-structure on $\D_{\Box}((\X,\mathcal{M})_\Prism,\Laur)$ restricts to $\Perf((\X,\mathcal{M})_\Prism,\Laur)$.
\end{corollary}

\begin{proof}
Using a covering by framed log affines via Lemma \ref{lem:strict_etale_local_framing}, we may use their associated logarithmic Breuil-Kisin prisms $(\rA,\I,\M_{\Spf \rA})^a$ from Lemma \ref{lem:standard_log_BK_prism} to produce a covering family (by Theorem \ref{thm:laur_cech}). The argument ends up being slightly technical since it is not clear whether the maps in the Breuil-Kisin \v{C}ech nerve are $t$-exact for $\D_{\Box}(-)$.

First we simply check that there exists a $t$-structure on $\Perf((\X,\mathcal{M})_\Prism,\Laur)$ regarded as ordinary Laurent crystals, defined via orthogonality from the Breuil-Kisin presentation of $\Perf^{\le 0}((\X,\mathcal{M})_\Prism,\Laur)$: we ask for the connective objects to be those which are connective on each \v{C}ech term $\rA^{(n)}[1/\I]_p^{\wedge}$. This follows by the abstract argument of \cite[Footnote 70]{Fgauge}, noting that we actually only need the maps $\rA[1/\I]^{\wedge}_p \to \rA^{(n)}[1/\I]^{\wedge}_p$ to be flat (which they are in the ordinary sense since $\rA[1/\I]^{\wedge}_p$ is Noetherian, using that $p$-complete flatness upgrades to ordinary flatness with a Noetherian source). By Lemma \ref{lem:perf_comparison} this also yields a $t$-structure on the dualizable objects in $\D_{\Box}((\X,\mathcal{M})_\Prism,\Laur)$ as dualizability is detected termwise in a limit.

We must check that any $\mathcal{E}$ in $\Perf^{\le 0}((\X,\mathcal{M})_\Prism,\Laur)$ is also in $\D^{\le 0}_{\Box}((\X,\mathcal{M})_\Prism,\Laur)$, by showing all evaluations on log prisms are connective. By construction of 
\[\Perf^{\le 0}((\X,\mathcal{M})_\Prism,\Laur),\] 
we see the evaluation on $(\rA \coprod_* \B)[1/\I]_p^{\wedge}$ is connective (as it is obtained by pullback from the value on $(\rA,\I,\M_{\Spf \rA})^a$, which by definition is connective). In Lemma \ref{lem:cech_descendable} we argued $(\rA \coprod_* \B)[1/\I]_p^{\wedge}$ is a descendable algebra in $\D_{\Box}(\B[1/\rJ]_p^{\wedge})$. In particular, the base change functor $\B[1/\rJ]_p^{\wedge} \to (\rA \coprod_* \B)[1/\I]_p^{\wedge}$ is conservative: in this situation base change detects connectivity on perfect complexes.\footnote{At this point one might be surprised we use something fancy like descendability rather than just argue by $p$-complete faithful flatness. The issue is that in order  to run this argument, we would want $\B \to \rA \coprod_* \B$ to be faithfully flat (rather than just $(p,\rJ)$-completely faithfully flat) so we can deduce the claim after $p$-complete localization. One can even check $\rho_{\rA}: \Spf \rA \to (\X,\mathcal{M})^{\Prism}$ is fpqc to deduce complete faithful flatness, but without Noetherianity of $\B$ to apply \cite[Lemma 5.15]{BhattCohenMacaulayAIC} it seems unclear how to finish (in fact the claim may very well be false).}

One may also check that $\mathcal{E}\in \Perf^{\ge 0}((\X,\mathcal{M})_\Prism,\Laur)$ is also in $\D^{\ge 0}_{\Box}((\X,\mathcal{M})_\Prism,\Laur)$ for its canonical $t$-structure. Indeed this is defined as the right orthogonal of 
\[\D^{\le -1}_{\Box}((\X,\mathcal{M})_\Prism,\Laur).\]
For $\mathcal{E}' \in \D^{\le -1}_{\Box}((\X,\mathcal{M})_\Prism,\Laur)$, every evaluation on \v{C}ech terms $\rA^{(n)}[1/\I]_p^{\wedge}$ lies in the $\le -1$ part of the solid $t$-structure. Writing $\mathcal{E}_n$ and $\mathcal{E}'_n$ for these evaluations, we get
\[\Map_{\D_{\Box}(\rA^{(n)}[1/\I]_p^{\wedge})}(\mathcal{E}'_n, \mathcal{E}_n) = 0\]
using that $\mathcal{E}_n$ are in the $\ge 0$ part (using that it is perfect here to relate ordinary and solid $t$-structures). But mapping spaces are computed by the limit over the \v{C}ech nerve, so the claim of the corollary follows.
\end{proof}

We denote its heart by $\Coh((\X,\mathcal{M})_\Prism,\Laur)$. Corollary \ref{cor:laur_strict_etale} then gives strict \'{e}tale descent for its heart.

\subsection{Analytifications}

We show in this subsection that the Breuil-Kisin covers give proper effective epimorphisms after analytification, and use this to deduce presentations of $\X^{\Prism,\Box}$ and $(\X,\mathcal{M})^{\Prism,\Box}$ similar to Proposition \ref{prop:log_stack_presentation}.

\begin{proposition} \label{prop:ordinary_BK_analytic_presentation}
Let $\X$ be a smooth $p$-adic formal scheme over $\O_K$, and assume $(\rA,\I)$ is a Breuil-Kisin atlas. The induced map
\[\Spa(\rA,\rA) \longrightarrow \X^{\Prism,\Box}\]
is a representable proper effective epimorphism. If $\rA^{(n)}$ denotes the $n$th iterated self-coproduct in the category of prisms over $\X$, then
\[\X^{\Prism,\Box} \simeq \left|\Spa(\rA^{(\bullet)},\rA^{(\bullet)})\right|\]
in analytic stacks.
\end{proposition}

\begin{proof}
To begin, choosing a quasiregular semiperfectoid cover $\coprod_i \Spf \R_i \to \X$ in the infinite-root/\'{e}tale topology gives an effective epimorphism
\[\coprod_i\Spa(\Prism_{\R_i},\Prism_{\R_i})\longrightarrow\X^{\Prism,\Box}\]
by Construction 3.7.1 and Corollary 3.6.3 of \cite{Kipp2}. It therefore suffices to test the first assertion that $\Spa(\rA,\rA) \longrightarrow \X^{\Prism,\Box}$ is a representable proper effective epimorphism after pullback to each $\Spa(\Prism_{\R_i},\Prism_{\R_i})$. Applying Lemma \ref{lem:mixed_BK_fiber_product} and using finite limit preservation of $(-)^{\Box}$ (Definition \ref{def:box_analytic_realization}) we may therefore identify this pullback with
\[\Spa(\rA\coprod_*\Prism_{\R_i},\rA\coprod_*\Prism_{\R_i})\longrightarrow\Spa(\Prism_{\R_i},\Prism_{\R_i}),\]
which is a proper descendable cover by Lemma \ref{lem:cech_descendable}. Since representability is local on the target, as well as properness by \cite[Remark 2.3.1]{Kipp2}, this proves the first assertion.

Upon using the non-logarithmic case in Lemma \ref{lem:mixed_BK_fiber_product} and finite limit preservation of $(-)^{\Box}$ (as in Definition \ref{def:box_analytic_realization}) we may identify the $n$th \v{C}ech term with $\Spa(\rA^{(n)},\rA^{(n)})$. Since the map from $\Spa(\rA,\rA)$ is an effective epimorphism, we deduce $\X^{\Prism,\Box}\simeq\left|\Spa(\rA^{(\bullet)},\rA^{(\bullet)})\right|$.
\end{proof}

The logarithmic version of this proposition holds, but it is far less obvious because we first have to produce a single effective epimorphism from a log prism to carry out a similar strategy. For a fine monoid $\rP$, let $\mathcal{A}_\rP=[\Spf\Z_p\langle \rP\rangle/\Spf\Z_p\langle \rP^{\mathrm{gp}}\rangle]$ be its Artin cone. This comes equipped with a natural cover $\Spf\Z_p\langle \rP\rangle\to\mathcal{A}_\rP$, and we have an \'{e}tale map $\mathcal{A}_\rP\to\Log$ by \cite[Corollary 5.24]{OlssonLogStacks}.

The purpose of this first lemma is to produce the replacement for the semiperfectoid effective epimorphism.

\begin{lemma} \label{lem:chart_log_prism}
Let $(\X,\mathcal{M})$ be a framed log affine over $\O_K$, let $(\rA,\I,\M_{\Spf\rA})^a$ be the logarithmic Breuil-Kisin prism associated with its framing, and let $\rP=\N^r$ be its framing's chart. There is an orientable bounded fine log prism $(\B_\rP,\I,\mathcal{N}_\rP)^a$ with a strict map
\[(\Spf(\B_\rP/\I),\mathcal{N}_\rP|_{\Spf(\B_\rP/\I)}) \longrightarrow (\X,\mathcal{M})\]
yielding an effective epimorphism $\Spa(\B_\rP,\B_\rP) \longrightarrow (\X,\mathcal{M})^{\Prism,\Box}$. In formal stacks the map
\[\Spf \B_{\rP} \longrightarrow  (\X,\mathcal{M})^{\Prism}\]
is also an effective epimorphism.
\end{lemma}

\begin{proof}
We will first construct the effective epimorphism 
\[\Spf \B_{\rP} \longrightarrow  (\X,\mathcal{M})^{\Prism}\]
of fpqc stacks, and argue that after applying the functor $(-)^{\Box}$ from Definition \ref{def:box_analytic_realization} we obtain an effective epimorphism of analytic stacks.

Since we work with a framed log affine, when constructing $(\X,\mathcal{M})^{\Prism}$ we can realize it as the relative stack $(\X/\mathcal{A}_{\rP})^{\Prism}$. Writing the stack this way reveals that there is an fpqc epimorphism
\[\Spf\rA \times_{\mathcal{A}_{\rP}^\Prism}(\Z_p^\Prism\times\Spf\Z_p\langle \rP\rangle) \longrightarrow (\X,\mathcal{M})^{\Prism}.\]
Indeed, termwise $\Spf \rA \to \X^{\Prism}$ is an fpqc epimorphism, as is $\Z_p^\Prism\times\Spf\Z_p\langle \rP\rangle \to \Z_p^\Prism \times \mathcal{A}_{\rP}$. Thus we must construct an fpqc epimorphism 
\[\Spf \B_{\rP} \to \Spf\rA \times_{\mathcal{A}_{\rP}^\Prism}(\Z_p^\Prism\times\Spf\Z_p\langle \rP\rangle).\]
This follows from rewriting the target stack in a different way. The main issue is that $\mathcal{A}_{\rP}$ is a stack -- if it was a smooth $p$-adic formal scheme, we could try to relate the target to $\Spf$ of some coproduct of prisms, which would give the desired claim.

What fixes this is that $\mathcal{A}_{\rP}$ has affine diagonal. Let us rewrite the target using the Cartesian square

\begin{display}
\Spf\rA \times_{\mathcal{A}_{\rP}^\Prism}(\Z_p^\Prism\times\Spf\Z_p\langle \rP\rangle) \ar{d} \ar{r} & \X^{\Prism} \times_{\X^{\Prism}} \Spf \rA \ar{d} \\
\Z_p\langle \rP \rangle^{\Prism} \times_{\Z_p\langle \rP \rangle^{\Prism}} \Z_p^\Prism\times\Spf\Z_p\langle \rP\rangle \ar{r} & \mathcal{A}_{\rP}^{\Prism} \times_{\Z_p^\Prism} \Z_p^{\Prism}.
\end{display}
The idea is to trivially rewrite each term, so that upon commuting limits with limits we get
\[(\X \times \mathbf{G}_m^r)^{\Prism} \times_{(\X \times \Spf \Z_p\langle \rP \rangle)^{\Prism}} \Spf(\rA \widehat{\otimes} \Z_p\langle \rP \rangle).\]
The key thing we have used here is that $(-)^{\Prism}$ is compatible with fiber products, and $\Spf \Z_p\langle \rP \rangle \times_{\mathcal{A}_{\rP}} \X$ becomes affine. For the proof, in the upper right corner we could have equally well chosen $\Z_p\langle \rP \rangle$ again, and the entire argument still works. We only make this choice to reuse the Breuil-Kisin prism $\rA$ in a moment.

This new presentation is substantially more useful, because there is a Breuil-Kisin prism $\rA\langle u_1^{\pm1},\ldots,u_r^{\pm1}\rangle_{(p,\I)}$ with $\delta(u_i)=0$ giving an fpqc epimorphism
\[\Spf \rA\langle u_1^{\pm1},\ldots,u_r^{\pm1}\rangle_{(p,\I)} \longrightarrow (\X \times \mathbf{G}_m^r)^{\Prism}.\]
Upon base change, we obtain an fpqc epimorphism
\[\Spf \rA\langle u_1^{\pm1},\ldots,u_r^{\pm1}\rangle_{(p,\I)} \times_{(\X \times \Spf \Z_p\langle \rP \rangle)^{\Prism}} \Spf(\rA \widehat{\otimes} \Z_p\langle \rP \rangle) \longrightarrow \Spf\rA\times_{\mathcal{A}_{\rP}^\Prism}(\Z_p^\Prism\times\Spf\Z_p\langle \rP\rangle).\]
The source is now a fiber product of two prisms with their natural maps to $(\X \times \Spf \Z_p\langle \rP \rangle)^{\Prism}$, exactly the (non-logarithmic) situation understood in Lemma \ref{lem:mixed_BK_fiber_product}. Letting the coproduct of these prisms be $\B_{\rP}$ we obtain a composite fpqc epimorphism
\[\Spf \B_{\rP} \longrightarrow \Spf\rA\times_{\mathcal{A}_{\rP}^\Prism}(\Z_p^\Prism\times\Spf\Z_p\langle \rP\rangle) \longrightarrow (\X,\mathcal{M})^{\Prism}.\]
This gives the second claim. Note that $\B_{\rP}$ is also a log prism as it naturally comes with a projection to $\Z_p^\Prism \times \mathcal{A}_{\rP}$ compatible with the map to $\mathcal{A}_{\rP}^\Prism$ using this map.

For the variant on analytic stacks, we apply $(-)^{\Box}$ to both maps. The important property is that it is compatible with finite limits, so in all of the fiber products above we may simply apply it individually to each formal stack. For the first arrow, it is a proper effective epimorphism of analytic stacks by Proposition \ref{prop:ordinary_BK_analytic_presentation}.

For the second arrow, on the first factor we again apply Proposition \ref{prop:ordinary_BK_analytic_presentation}, and for the second factor we only use that
\[\Spf \Z_p\langle \rP\rangle \longrightarrow \mathcal{A}_{\rP}\]
is sent to an effective epimorphism of analytic stacks by $(-)^{\Box}$ which is clear.
\end{proof}

\begin{proposition} \label{prop:log_BK_analytic_presentation}
Let $(\X,\mathcal{M})$ be a framed log affine over $\O_K$, and let $(\rA,\I,\M_{\Spf\rA})^a$ be the logarithmic Breuil-Kisin prism associated with its framing. The induced map
\[\Spa(\rA,\rA)\longrightarrow(\X,\mathcal{M})^{\Prism,\Box}\]
is a representable proper effective epimorphism. If $\rA^{(n)}$ denotes the $n$th iterated self-coproduct in $(\X,\mathcal{M})_\Prism^{\mathrm{op}}$, then
\[(\X,\mathcal{M})^{\Prism,\Box}\simeq\left|\Spa(\rA^{(\bullet)},\rA^{(\bullet)})\right|\]
in analytic stacks.
\end{proposition}

\begin{proof}
The proof is similar to Proposition \ref{prop:ordinary_BK_analytic_presentation}. Letting $\rP=\N^r$ be the framing's chart, by Lemma \ref{lem:chart_log_prism} we get an effective epimorphism
\[\Spa(\B_\rP,\B_\rP)\longrightarrow (\X,\mathcal{M})^{\Prism,\Box}\]
induced by the map $\rho_{\X,\B_{\rP}}$. Then using Lemma \ref{lem:mixed_BK_fiber_product}, we may identify the pullback of $\Spa(\rA,\rA)$ along this map with
\[\Spa(\rA\coprod_*\B_\rP,\rA\coprod_*\B_\rP)\longrightarrow\Spa(\B_\rP,\B_\rP).\]
This is a proper descendable cover by Lemma \ref{lem:cech_descendable}; since representability is local on the target, as is properness by \cite[Remark 2.3.1]{Kipp2}, this proves the first assertion. 

Finally, Lemma \ref{lem:mixed_BK_fiber_product} and finite limit preservation of $(-)^{\Box}$ (as in Definition \ref{def:box_analytic_realization}) identify the $n$th \v{C}ech term with $\Spa(\rA^{(n)},\rA^{(n)})$. Since the map from $\Spa(\rA,\rA)$ is an effective epimorphism, we deduce $(\X,\mathcal{M})^{\Prism,\Box}\simeq\left|\Spa(\rA^{(\bullet)},\rA^{(\bullet)})\right|$.
\end{proof}

\bibliographystyle{amsalpha}
\bibliography{references}

@incollection{Faltings1989,
  author    = {Faltings, Gerd},
  title     = {Crystalline cohomology and {$p$}-adic {Galois}-representations},
  booktitle = {Algebraic Analysis, Geometry, and Number Theory},
  editor    = {Igusa, Jun-Ichi},
  pages     = {25--80},
  publisher = {Johns Hopkins University Press},
  address   = {Baltimore, MD},
  year      = {1989}
}

@incollection{Fontaine1990,
  author    = {Fontaine, Jean-Marc},
  title     = {Repr{\'e}sentations $p$-adiques des corps locaux. {I}},
  booktitle = {The {Grothendieck} Festschrift, Vol. {II}},
  series    = {Progress in Mathematics},
  volume    = {87},
  pages     = {249--309},
  publisher = {Birkh{\"a}user Boston},
  address   = {Boston, MA},
  year      = {1990},
  doi       = {10.1007/978-0-8176-4575-5_6}
}

@incollection{AndreattaIovita2008,
  author    = {Andreatta, Fabrizio and Iovita, Adrian},
  title     = {Global applications of relative {$(\varphi,\Gamma)$}-{Modules} {I}},
  booktitle = {Repr{\'e}sentations $p$-adiques de groupes $p$-adiques {I}: repr{\'e}sentations galoisiennes et $(\varphi,\Gamma)$-modules},
  editor    = {Berger, Laurent and Breuil, Christophe and Colmez, Pierre},
  series    = {Ast{\'e}risque},
  number    = {319},
  pages     = {339--420},
  publisher = {Soci{\'e}t{\'e} math{\'e}matique de France},
  year      = {2008},
  url       = {https://www.numdam.org/item/AST_2008__319__339_0/}
}

@article{MinWangRelativePhiGamma,
  author        = {Min, Yu and Wang, Yupeng},
  title         = {Relative {$(\varphi,\Gamma)$}-modules and prismatic {$F$}-crystals},
  journal       = {Mathematische Zeitschrift},
  volume        = {310},
  number        = {1},
  eid           = {16},
  pages         = {Article No. 16},
  year          = {2025},
  eprint        = {2110.06076},
  archivePrefix = {arXiv},
  primaryClass  = {math.AG},
  doi           = {10.1007/s00209-025-03696-9}
}

@article{AnschutzExtendingTorsors,
  author        = {Ansch{\"u}tz, Johannes},
  title         = {Extending torsors on the punctured {$\mathrm{Spec}(A_{\mathrm{inf}})$}},
  journal       = {Journal f{\"u}r die reine und angewandte Mathematik},
  volume        = {783},
  pages         = {227--268},
  year          = {2022},
  doi           = {10.1515/crelle-2021-0077},
  note          = {Numbering follows arXiv:1804.06356v2},
  eprint        = {1804.06356},
  archivePrefix = {arXiv},
  primaryClass  = {math.NT},
  url           = {https://arxiv.org/abs/1804.06356v2}
}

@article{GrosLeStumQuiros2023,
  author        = {Gros, Michel and Le Stum, Bernard and Quir{\'o}s, Adolfo},
  title         = {Absolute calculus and prismatic crystals on cyclotomic rings},
  journal       = {Algebra \& Number Theory},
  volume        = {20},
  number        = {7},
  pages         = {1301--1384},
  year          = {2026},
  eprint        = {2310.13790},
  archivePrefix = {arXiv},
  primaryClass  = {math.AG},
  doi           = {10.2140/ant.2026.20.1301}
}

@article{GrosLeStumQuirosNegative,
  author        = {Gros, Michel and Le Stum, Bernard and Quir{\'o}s, Adolfo},
  title         = {Twisted differential operators of negative level and prismatic crystals},
  journal       = {Tunisian Journal of Mathematics},
  volume        = {4},
  number        = {1},
  pages         = {19--54},
  year          = {2022},
  eprint        = {2010.04433},
  archivePrefix = {arXiv},
  primaryClass  = {math.AG},
  doi           = {10.2140/tunis.2022.4.19}
}

@article{BhattScholzePrisms,
  author  = {Bhatt, Bhargav and Scholze, Peter},
  title   = {Prisms and prismatic cohomology},
  journal = {Annals of Mathematics},
  volume  = {196},
  number  = {3},
  pages   = {1135--1275},
  year    = {2022},
  doi     = {10.4007/annals.2022.196.3.5}
}

@incollection{bhattframing,
  author    = {Bhatt, Bhargav},
  title     = {Specializing varieties and their cohomology from characteristic {$0$} to characteristic {$p$}},
  booktitle = {Algebraic Geometry: {Salt Lake City} 2015},
  series    = {Proceedings of Symposia in Pure Mathematics},
  volume    = {97.2},
  publisher = {American Mathematical Society},
  year      = {2018},
  pages     = {43--88},
  doi       = {10.1090/pspum/097.2/01699}
}

@article{BhattScholzePrismaticFCrystals,
  author  = {Bhatt, Bhargav and Scholze, Peter},
  title   = {Prismatic $\mathrm{F}$-crystals and crystalline {Galois} representations},
  journal = {Cambridge Journal of Mathematics},
  volume  = {11},
  number  = {2},
  pages   = {507--562},
  year    = {2023},
  doi     = {10.4310/CJM.2023.v11.n2.a3}
}

@unpublished{BhattCohenMacaulayAIC,
  author        = {Bhatt, Bhargav},
  title         = {{Cohen--Macaulayness} of absolute integral closures},
  year          = {2020},
  note          = {Preprint, arXiv:2008.08070},
  eprint        = {2008.08070},
  archivePrefix = {arXiv},
  primaryClass  = {math.AG},
  doi           = {10.48550/arXiv.2008.08070}
}

@article{BhattMathewArc,
  author        = {Bhatt, Bhargav and Mathew, Akhil},
  title         = {The {arc}-topology},
  journal       = {Duke Mathematical Journal},
  volume        = {170},
  number        = {9},
  pages         = {1899--1988},
  year          = {2021},
  eprint        = {1807.04725},
  archivePrefix = {arXiv},
  primaryClass  = {math.AG},
  doi           = {10.1215/00127094-2020-0088}
}

@article{proet,
  author        = {Bhatt, Bhargav and Scholze, Peter},
  title         = {The pro-\'{e}tale topology for schemes},
  journal       = {Ast{\'e}risque},
  number        = {369},
  pages         = {99--201},
  year          = {2015},
  eprint        = {1309.1198},
  archivePrefix = {arXiv},
  primaryClass  = {math.AG}
}

@article{BhattHansenCons,
  author        = {Bhatt, Bhargav and Hansen, David},
  title         = {The six functors for {Zariski}-constructible sheaves in rigid geometry},
  journal       = {Compositio Mathematica},
  volume        = {158},
  number        = {2},
  pages         = {437--482},
  year          = {2022},
  eprint        = {2101.09759},
  archivePrefix = {arXiv},
  primaryClass  = {math.AG},
  doi           = {10.1112/S0010437X22007291}
}

@unpublished{MannPadicSixFunctors,
  author        = {Mann, Lucas},
  title         = {A {$p$}-Adic 6-Functor Formalism in Rigid-Analytic Geometry},
  year          = {2022},
  note          = {Preprint, arXiv:2206.02022},
  eprint        = {2206.02022},
  archivePrefix = {arXiv},
  primaryClass  = {math.AG},
  doi           = {10.48550/arXiv.2206.02022}
}

@unpublished{HeyerMannSixFunctors,
  author        = {Heyer, Claudius and Mann, Lucas},
  title         = {6-{F}unctor {F}ormalisms and {S}mooth {R}epresentations},
  year          = {2024},
  note          = {Preprint, arXiv:2410.13038},
  eprint        = {2410.13038},
  archivePrefix = {arXiv},
  primaryClass  = {math.CT},
  doi           = {10.48550/arXiv.2410.13038}
}

@article{BhattLurieModpRH,
  author        = {Bhatt, Bhargav and Lurie, Jacob},
  title         = {A {Riemann--Hilbert} correspondence in positive characteristic},
  journal       = {Cambridge Journal of Mathematics},
  volume        = {7},
  number        = {1--2},
  pages         = {71--217},
  year          = {2019},
  doi           = {10.4310/CJM.2019.v7.n1.a3},
  eprint        = {1711.04148},
  archivePrefix = {arXiv},
  primaryClass  = {math.AG}
}

@unpublished{APC,
  author        = {Bhatt, Bhargav and Lurie, Jacob},
  title         = {Absolute prismatic cohomology},
  year          = {2022},
  note          = {Preprint, arXiv:2201.06120},
  eprint        = {2201.06120},
  archivePrefix = {arXiv},
  primaryClass  = {math.AG},
  doi           = {10.48550/arXiv.2201.06120}
}

@unpublished{APC2,
  author        = {Bhatt, Bhargav and Lurie, Jacob},
  title         = {The prismatization of {$p$}-adic formal schemes},
  year          = {2022},
  note          = {Preprint, arXiv:2201.06124},
  eprint        = {2201.06124},
  archivePrefix = {arXiv},
  primaryClass  = {math.AG}
}

@unpublished{AKNRelativeDelta,
  author        = {Antieau, Benjamin and Krause, Achim and Nikolaus, Thomas},
  title         = {Prismatic cohomology relative to {$\delta$}-rings},
  year          = {2026},
  note          = {To appear in Annales scientifiques de l'\'{E}cole normale sup\'{e}rieure; arXiv:2310.12770v2},
  eprint        = {2310.12770},
  archivePrefix = {arXiv},
  primaryClass  = {math.AG},
  url           = {https://arxiv.org/abs/2310.12770v2}
}

@unpublished{BhattPadicHodgeNotes,
  author = {Bhatt, Bhargav},
  title  = {Aspects of {$p$}-adic {Hodge} theory},
  year   = {2025},
  note   = {Lecture notes, Princeton University, Fall 2025},
  url    = {https://www.math.ias.edu/~bhatt/teaching/mat517f25/pHT-notes.pdf}
}

@article{BlickleIH,
  author  = {Blickle, Manuel},
  title   = {The intersection homology {$D$}-module in finite characteristic},
  journal = {Mathematische Annalen},
  volume  = {328},
  number  = {3},
  pages   = {425--450},
  year    = {2004},
  doi     = {10.1007/s00208-003-0492-z},
  url     = {https://doi.org/10.1007/s00208-003-0492-z}
}

@article{BlickleGamma,
  author  = {Blickle, Manuel},
  title   = {Minimal {$\gamma$}-sheaves},
  journal = {Algebra \& Number Theory},
  volume  = {2},
  number  = {3},
  pages   = {347--368},
  year    = {2008},
  doi     = {10.2140/ant.2008.2.347},
  url     = {https://doi.org/10.2140/ant.2008.2.347}
}

@book{BoschGuntzerRemmert,
  author    = {Bosch, Siegfried and G{\"u}ntzer, Ulrich and Remmert, Reinhold},
  title     = {Non-Archimedean Analysis: A Systematic Approach to Rigid Analytic Geometry},
  series    = {Grundlehren der mathematischen Wissenschaften},
  volume    = {261},
  publisher = {Springer-Verlag},
  address   = {Berlin},
  year      = {1984}
}

@unpublished{EmertonCoherentRings,
  author    = {Emerton, Matthew},
  title     = {On a class of coherent rings, with applications to the smooth representation theory of {$\mathrm{GL}_2(\mathbf{Q}_p)$} in characteristic {$p$}},
  year      = {2008},
  note      = {Preprint, May 19, 2008},
  url       = {https://math.uchicago.edu/~emerton/pdffiles/frob.pdf}
}

@incollection{EmertonKisin1,
  author    = {Emerton, Matthew and Kisin, Mark},
  title     = {An introduction to the {Riemann--Hilbert} correspondence for unit {$F$}-crystals},
  booktitle = {Geometric Aspects of {Dwork} Theory},
  editor    = {Adolphson, Alan and Baldassarri, Francesco and Berthelot, Pierre and Katz, Nicholas and Loeser, Fran{\c c}ois},
  pages     = {677--700},
  publisher = {Walter de Gruyter},
  address   = {Berlin},
  year      = {2004},
  doi       = {10.1515/9783110198133.2.677}
}

@incollection{GabberTStructures,
  author    = {Gabber, Ofer},
  title     = {Notes on some {$t$}-structures},
  booktitle = {Geometric Aspects of {Dwork} Theory},
  editor    = {Adolphson, Alan and Baldassarri, Francesco and Berthelot, Pierre and Katz, Nicholas and Loeser, Fran{\c c}ois},
  volume    = {II},
  pages     = {711--734},
  publisher = {Walter de Gruyter},
  address   = {Berlin},
  year      = {2004},
  doi       = {10.1515/9783110198133.2.711}
}

@book{EmertonKisin2,
  author    = {Emerton, Matthew and Kisin, Mark},
  title     = {The {Riemann--Hilbert} correspondence for unit {$F$}-crystals},
  series    = {Ast{\'e}risque},
  number    = {293},
  publisher = {Soci{\'e}t{\'e} math{\'e}matique de France},
  year      = {2004},
  pages     = {vi+257},
  isbn      = {2-85629-154-6},
  mrnumber  = {2071510},
  zbl       = {1056.14025},
  doi       = {10.24033/ast.620},
  url       = {https://numdam.org/item/AST_2004__293__1_0/}
}

@unpublished{GuoLiFrobeniusHeight,
  author        = {Guo, Haoyang and Li, Shizhang},
  title         = {Frobenius height of prismatic cohomology with coefficients},
  year          = {2025},
  note          = {Preprint, arXiv:2309.06663v2},
  eprint        = {2309.06663},
  archivePrefix = {arXiv},
  primaryClass  = {math.AG},
  doi           = {10.48550/arXiv.2309.06663}
}

@article{TsujiLogCrystallinePD,
  author  = {Tsuji, Takeshi},
  title   = {Poincar{\'e} duality for logarithmic crystalline cohomology},
  journal = {Compositio Mathematica},
  volume  = {118},
  number  = {1},
  pages   = {11--41},
  year    = {1999},
  doi     = {10.1023/A:1001020809306}
}

@unpublished{syncrisloc,
  author        = {Pentland, Dylan},
  title         = {Syntomification and crystalline local systems},
  year          = {2025},
  note          = {Preprint, arXiv:2510.16961},
  eprint        = {2510.16961},
  archivePrefix = {arXiv},
  primaryClass  = {math.NT},
  doi           = {10.48550/arXiv.2510.16961}
}

@incollection{KatzPAdicProperties,
  author    = {Katz, Nicholas M.},
  title     = {$p$-adic properties of modular schemes and modular forms},
  booktitle = {Modular Functions of One Variable {III}},
  editor    = {Kuijk, Willem and Serre, Jean-Pierre},
  series    = {Lecture Notes in Mathematics},
  volume    = {350},
  pages     = {69--190},
  publisher = {Springer},
  address   = {Berlin},
  year      = {1973},
  doi       = {10.1007/978-3-540-37802-0_3},
  mrnumber  = {0447119}
}

@unpublished{Inoue,
  author        = {Inoue, Kentaro},
  title         = {Log prismatic {$F$}-crystals and realization functors},
  year          = {2025},
  note          = {Preprint, arXiv:2505.01084},
  eprint        = {2505.01084},
  archivePrefix = {arXiv},
  primaryClass  = {math.AG},
  doi           = {10.48550/arXiv.2505.01084},
  url           = {https://arxiv.org/pdf/2505.01084}
}

@article{KoshikawaYao2025,
  author  = {Koshikawa, Teruhisa and Yao, Zijian},
  title   = {Logarithmic prismatic cohomology {II}},
  journal = {Advances in Mathematics},
  volume  = {479},
  pages   = {110446},
  year    = {2025},
  doi     = {10.1016/j.aim.2025.110446},
  eprint  = {2306.00364},
  archivePrefix = {arXiv},
  primaryClass = {math.AG}
}

@misc{KoshikawaYaoCorrigendum,
  author = {Inoue, Kentaro and Koshikawa, Teruhisa and Yao, Zijian},
  title  = {Corrigendum to {``Logarithmic prismatic cohomology {II}''}},
  year   = {2026},
  note   = {Corrigendum, available at \url{https://researchmap.jp/7000017226/misc/53305261}},
  url    = {https://researchmap.jp/7000017226/misc/53305261}
}

@unpublished{KoshikawaLogPrismaticI,
  author        = {Koshikawa, Teruhisa},
  title         = {Logarithmic prismatic cohomology {I}},
  year          = {2022},
  note          = {Preprint, arXiv:2007.14037v3},
  eprint        = {2007.14037},
  archivePrefix = {arXiv},
  primaryClass  = {math.AG}
}

@unpublished{KedlayaLiu2016,
  author        = {Kedlaya, Kiran S. and Liu, Ruochuan},
  title         = {Relative $p$-adic {Hodge} theory, {II}: Imperfect period rings},
  year          = {2016},
  note          = {Preprint, arXiv:1602.06899},
  eprint        = {1602.06899},
  archivePrefix = {arXiv},
  primaryClass  = {math.NT},
  doi           = {10.48550/arXiv.1602.06899}
}

@article{OlssonLogStacks,
  author  = {Olsson, Martin C.},
  title   = {Logarithmic geometry and algebraic stacks},
  journal = {Annales scientifiques de l'\'{E}cole Normale Sup\'{e}rieure},
  series  = {4},
  volume  = {36},
  number  = {5},
  pages   = {747--791},
  year    = {2003},
  doi     = {10.1016/j.ansens.2002.11.001}
}

@misc{OlssonLogPrismaticStack,
  author       = {Olsson, Martin C.},
  title        = {Stacks and Log Prismatic Cohomology},
  howpublished = {Video lecture, Workshop on $p$-adic Arithmetic Geometry,
                  Institute for Advanced Study,
                  \url{https://www.ias.edu/video/stacks-and-log-prismatic-cohomology}},
  year         = {2023},
  note         = {Lecture delivered November 13, 2023}
}

@phdthesis{Kipp,
  author = {Kipp, Niklas},
  title  = {A six-functor formalism for syntomic cohomology of {$p$}-adic formal schemes},
  school = {Universit{\"a}t Regensburg},
  year   = {2026},
  doi    = {10.5283/epub.78438},
  url    = {https://epub.uni-regensburg.de/78438/}
}

@unpublished{Kipp2,
  author        = {Kipp, Niklas},
  title         = {A six-functor formalism for syntomic cohomology},
  year          = {2026},
  note          = {Preprint, arXiv:2608.15744},
  eprint        = {2608.15744},
  archivePrefix = {arXiv},
  doi           = {10.48550/arXiv.2608.15744}
}

@unpublished{HauckStacky,
  author        = {Hauck, Maximilian},
  title         = {A stacky approach to $p$-adic {Hodge} theory},
  year          = {2024},
  note          = {Preprint, arXiv:2409.10557},
  eprint        = {2409.10557},
  archivePrefix = {arXiv},
  primaryClass  = {math.AG},
  doi           = {10.48550/arXiv.2409.10557}
}

@unpublished{HauckLA,
  author        = {Hauck, Maximilian},
  title         = {Locally analytic representations in mixed characteristic via stacks},
  year          = {2026},
  note          = {Preprint, arXiv:2609.28141},
  eprint        = {2609.28141},
  archivePrefix = {arXiv},
  primaryClass  = {math.RT},
  doi           = {10.48550/arXiv.2609.28141}
}

@article{HemoRicharzScholbach2023,
  author        = {Hemo, Tamir and Richarz, Timo and Scholbach, Jakob},
  title         = {Constructible sheaves on schemes},
  journal       = {Advances in Mathematics},
  volume        = {429},
  pages         = {109179},
  year          = {2023},
  doi           = {10.1016/j.aim.2023.109179},
  eprint        = {2305.18131},
  archivePrefix = {arXiv},
  primaryClass  = {math.AG}
}

@article{Porat_locally_analytic,
  author        = {Porat, Gal},
  title         = {Locally analytic vectors and decompletion in mixed characteristic},
  journal       = {Forum of Mathematics, Sigma},
  volume        = {13},
  eid           = {e173},
  pages         = {Paper No. e173, 26 pp.},
  year          = {2025},
  doi           = {10.1017/fms.2025.10121},
  eprint        = {2407.19791},
  archivePrefix = {arXiv},
  primaryClass  = {math.NT}
}

@misc{PoratSolidLocallyAnalytic,
  author        = {Porat, Gal},
  title         = {Solid locally analytic representations in mixed characteristic},
  year          = {2025},
  eprint        = {2510.13673},
  archivePrefix = {arXiv},
  primaryClass  = {math.NT}
}

@article{yamashita,
  author  = {Yamashita, Go},
  title   = {{$p$}-adic {Hodge} theory for open varieties},
  journal = {Comptes Rendus. Math{\'e}matique},
  volume  = {349},
  number  = {21--22},
  pages   = {1127--1130},
  year    = {2011},
  doi     = {10.1016/j.crma.2011.10.016},
  url     = {https://doi.org/10.1016/j.crma.2011.10.016}
}

@article{etCohDiamonds,
  author        = {Scholze, Peter},
  title         = {{\'E}tale cohomology of diamonds},
  journal       = {Ast{\'e}risque},
  volume        = {466},
  pages         = {439--674},
  year          = {2026},
  eprint        = {1709.07343},
  archivePrefix = {arXiv},
  primaryClass  = {math.AG},
  url           = {https://smf.emath.fr/sites/default/files/2026-07/AST-466-web__sample.pdf}
}

@misc{clausenscholze,
  author       = {Clausen, Dustin and Scholze, Peter},
  title        = {Analytic Stacks},
  year         = {2023},
  howpublished = {Lecture course, Institut des Hautes {\'E}tudes Scientifiques and Max Planck Institute for Mathematics},
  url          = {https://people.mpim-bonn.mpg.de/scholze/AnalyticStacks.html}
}

@unpublished{RodriguezCamargoSolidGeometry,
  author = {Rodr{\'i}guez Camargo, Juan Esteban},
  title  = {Notes on Solid Geometry},
  year   = {2024},
  note   = {Seminar notes},
  url    = {https://blogs.cuit.columbia.edu/jr4460/files/2024/09/SeminarSolidGeometrynotes2.pdf}
}

@unpublished{ScholzeAnalyticGeometry,
  author        = {Scholze, Peter},
  title         = {Lectures on Analytic Geometry},
  year          = {2026},
  note          = {Lecture notes, arXiv:2605.03655},
  eprint        = {2605.03655},
  archivePrefix = {arXiv},
  primaryClass  = {math.AG},
  doi           = {10.48550/arXiv.2605.03655}
}

@incollection{DiaoLanLiuZhuLogAdic,
  author        = {Diao, Hansheng and Lan, Kai-Wen and Liu, Ruochuan and Zhu, Xinwen},
  title         = {Logarithmic adic spaces: some foundational results},
  booktitle     = {$p$-adic {Hodge} Theory, Singular Varieties, and Non-Abelian Aspects},
  editor        = {Bhatt, Bhargav and Olsson, Martin},
  series        = {Simons Symposia},
  pages         = {65--182},
  publisher     = {Springer},
  address       = {Cham},
  year          = {2023},
  eprint        = {1912.09836},
  archivePrefix = {arXiv},
  primaryClass  = {math.AG},
  doi           = {10.1007/978-3-031-21550-6_3}
}

@article{MathewModpRH,
  author        = {Mathew, Akhil},
  title         = {The mod-$p$ {Riemann--Hilbert} correspondence and the perfect site},
  journal       = {Tunisian Journal of Mathematics},
  volume        = {5},
  number        = {2},
  pages         = {369--403},
  year          = {2023},
  eprint        = {2205.11657},
  archivePrefix = {arXiv},
  primaryClass  = {math.AG},
  doi           = {10.2140/tunis.2023.5.369}
}

@article{BhattScholzeWittGrassmannian,
  author        = {Bhatt, Bhargav and Scholze, Peter},
  title         = {Projectivity of the {Witt} vector affine {Grassmannian}},
  journal       = {Inventiones Mathematicae},
  volume        = {209},
  number        = {2},
  pages         = {329--423},
  year          = {2017},
  eprint        = {1507.06490},
  archivePrefix = {arXiv},
  primaryClass  = {math.AG},
  doi           = {10.1007/s00222-016-0710-4}
}

@article{MathewGaloisStableHomotopy,
  author  = {Mathew, Akhil},
  title   = {The {Galois} group of a stable homotopy theory},
  journal = {Advances in Mathematics},
  volume  = {291},
  pages   = {403--541},
  year    = {2016},
  doi     = {10.1016/j.aim.2015.12.017}
}

@article{drinfeldmathew,
  author        = {Mathew, Akhil},
  title         = {Faithfully flat descent of almost perfect complexes in rigid geometry},
  journal       = {Journal of Pure and Applied Algebra},
  volume        = {226},
  number        = {5},
  pages         = {106938},
  year          = {2022},
  eprint        = {1912.10968},
  archivePrefix = {arXiv},
  primaryClass  = {math.AG},
  doi           = {10.1016/j.jpaa.2021.106938}
}

@misc{Fgauge,
  author       = {Bhatt, Bhargav},
  title        = {Prismatic $\mathrm{F}$-gauges},
  year         = {2022},
  howpublished = {\url{https://www.math.ias.edu/~bhatt/teaching/mat549f22/lectures.pdf}}
}

@unpublished{MikamiRelativelyDiscrete,
  author        = {Mikami, Yutaro},
  title         = {{$(\varphi,\Gamma)$}-modules over relatively discrete algebras},
  year          = {2026},
  note          = {Preprint, arXiv:2409.14145v3},
  eprint        = {2409.14145},
  archivePrefix = {arXiv},
  primaryClass  = {math.NT},
  doi           = {10.48550/arXiv.2409.14145}
}

@unpublished{MikamiDuality,
  author        = {Mikami, Yutaro},
  title         = {Finiteness and duality of cohomology of {$(\varphi,\Gamma)$}-modules and the 6-functor formalism of locally analytic representations},
  year          = {2025},
  note          = {Preprint, arXiv:2504.01780},
  eprint        = {2504.01780},
  archivePrefix = {arXiv},
  primaryClass  = {math.NT},
  doi           = {10.48550/arXiv.2504.01780}
}

@article{MaoDerivedCrystalline,
  author  = {Mao, Zhouhang},
  title   = {Revisiting derived crystalline cohomology},
  journal = {Bulletin de la Soci\'{e}t\'{e} Math\'{e}matique de France},
  volume  = {152},
  number  = {4},
  pages   = {659--784},
  year    = {2024},
  doi     = {10.24033/bsmf.1524zm},
  eprint  = {2107.02921},
  archivePrefix = {arXiv},
  primaryClass  = {math.AG}
}

@unpublished{BoscoDrinfeld,
  author        = {Bosco, Guido},
  title         = {On the {$p$}-adic pro-\'{e}tale cohomology of {Drinfeld} symmetric spaces},
  year          = {2023},
  note          = {Preprint, arXiv:2110.10683},
  eprint        = {2110.10683},
  archivePrefix = {arXiv},
  primaryClass  = {math.NT},
  doi           = {10.48550/arXiv.2110.10683}
}

@incollection{Schneider,
  author    = {Schneider, Peter and Vign{\'e}ras, Marie-France},
  title     = {A functor from smooth {$\mathfrak{o}$}-torsion representations to {$(\varphi,\Gamma)$}-modules},
  booktitle = {On Certain {$L$}-Functions},
  series    = {Clay Mathematics Proceedings},
  volume    = {13},
  pages     = {525--601},
  publisher = {American Mathematical Society},
  address   = {Providence, RI},
  year      = {2011},
  url       = {https://www.uni-muenster.de/IVV5WS/WebHop/user/pschnei/pap/functor.pdf}
}

@incollection{Brinon,
  author    = {Brinon, Olivier},
  title     = {Filtered {$(\varphi,N)$}-modules and semi-stable representations},
  booktitle = {An excursion into {$p$}-adic {Hodge} theory: from foundations to recent trends},
  series    = {Panoramas et Synth{\`e}ses},
  volume    = {54},
  pages     = {93--129},
  publisher = {Soci{\'e}t{\'e} Math{\'e}matique de France},
  year      = {2019},
  url       = {https://www.math.u-bordeaux.fr/~obrinon/publications/pst.pdf}
}

@misc{StacksProject,
  author       = {{Stacks Project Authors}},
  title        = {Stacks Project},
  year         = {2026},
  howpublished = {\url{https://stacks.math.columbia.edu}}
}

@article{AlonsoJeremiasSaorinTStructures,
  author  = {Alonso Tarr{\'i}o, Leovigildo and Jerem{\'i}as L{\'o}pez, Ana and Saor{\'i}n, Manuel},
  title   = {Compactly generated {$t$}-structures on the derived category of a {Noetherian} ring},
  journal = {Journal of Algebra},
  volume  = {324},
  number  = {3},
  pages   = {313--346},
  year    = {2010},
  doi     = {10.1016/j.jalgebra.2010.04.023},
  url     = {https://arxiv.org/abs/0706.0499}
}

@unpublished{Watanabe,
  author        = {Watanabe, Takumi},
  title         = {A remark on an integral structure of the imperfect coefficient ring of {$(\varphi,\Gamma)$}-modules},
  year          = {2026},
  note          = {Preprint, arXiv:2604.17559; Appendix 2 by Dylan Pentland},
  eprint        = {2604.17559},
  archivePrefix = {arXiv},
  primaryClass  = {math.NT},
  doi           = {10.48550/arXiv.2604.17559}
}

@article{Wigner1973,
  author  = {Wigner, David},
  title   = {Algebraic cohomology of topological groups},
  journal = {Transactions of the American Mathematical Society},
  volume  = {178},
  pages   = {83--93},
  year    = {1973},
  doi     = {10.1090/S0002-9947-1973-0338132-7}
}

@article{lutkebohermert,
  author  = {L{\"u}tkebohmert, Werner},
  title   = {{Riemann}'s existence problem for a {$p$}-adic field},
  journal = {Inventiones Mathematicae},
  volume  = {111},
  pages   = {309--330},
  year    = {1993},
  doi     = {10.1007/BF01231290}
}

@article{KatzNilpMon,
  title={Nilpotent connections and the monodromy theorem: Applications of a result of {Turrittin}},
  author={Katz, Nicholas M},
  journal={Publications math{\'e}matiques de l'IHES},
  volume={39},
  pages={175--232},
  year={1970}
}

@article{GLQcartier,
  title={Cartier transform and prismatic crystals},
  author={Gros, Michel and Le Stum, Bernard and Quir{\'o}s, Adolfo},
  journal={Tunisian Journal of Mathematics},
  volume={5},
  number={3},
  pages={405--456},
  year={2023},
  publisher={Mathematical Sciences Publishers}
}

\end{document}